%% file: GP23.tex
\documentclass[leqno,12pt]{article}
\usepackage{amsmath,amssymb,rotating,a4wide,color}
\usepackage{wasysym}
\usepackage{hyperref}
\usepackage{enumerate}
\usepackage{mathtools}
\usepackage{float}
\usepackage{colortbl}
\usepackage{makecell}
\usepackage[font=footnotesize,labelfont=bf,margin=2cm]{caption}

\usepackage{array}
\newcolumntype{M}{>{\centering\arraybackslash}m{\dimexpr.25\linewidth-2\tabcolsep}}

\usepackage{tikz}
\usetikzlibrary{circuits.logic.US,circuits.logic.IEC,fit}
\usetikzlibrary{cd}
\usetikzlibrary{backgrounds}
\usetikzlibrary{math}
\usepackage{subcaption}
\usepackage{lscape}

\usepackage[all,cmtip]{xy}

\usepackage{verbatim}

\newdir^{ (}{{}*!/-3pt/\dir^{(}}
\newdir_{ (}{{}*!/-3pt/\dir_{(}}

\entrymodifiers={+!!<0pt,\fontdimen22\textfont2>}

\input{GP22-Figures.tex}
\input{GP22-Figures-genus-2.tex}

\def\bigtimes{\mathop{\raise-2pt\hbox{\huge$\times$}}}

\newbox\circbulletbox
\setbox\circbulletbox=\hbox{\,{\large$\circledcirc$}\kern-8.7pt\raise .6pt\hbox{$\bullet$}\,}

\let\le\leqslant
\let\ge\geqslant
\let\leq\leqslant
\let\geq\geqslant

\def\circVbig{\hbox{\text{\it\r{V}}}}
\def\circVscript{\hbox{\scriptsize\text{\it\r{V}}}}
\def\circVscriptscript{\mbox{\tiny\text{\it\r{V}}}}

\def\circVlimits_#1^#2{{\mathchoice%
{\circVbig{}^{\kern2pt #2}_{\kern-2pt #1}}%
{\circVbig{}^{\kern2pt #2}_{\kern-2pt #1}}%
{\scriptstyle\circVscript{}^{\kern1.7pt #2}_{\kern-1pt #1}}%
{\scriptscriptstyle\circVscriptscript{}^{\kern1.5pt #2}_{\kern-1pt #1}}%
}}
\def\circVr_#1{\circVlimits_#1^r}
\def\circVs_#1{\circVlimits_#1^s}

\def\circWbig{\hbox{\text{\it\r{W}}}}
\def\circWscript{\hbox{\scriptsize\text{\it\r{W}}}}
\def\circWscriptscript{\mbox{\tiny\text{\it\r{W}}}}

\def\circWlimits_#1^#2{{\mathchoice%
{\circWbig{}^{\kern2pt #2}_{\kern-2pt #1}}%
{\circWbig{}^{\kern2pt #2}_{\kern-2pt #1}}%
{\scriptstyle\circWscript{}^{\kern1.7pt #2}_{\kern-1pt #1}}%
{\scriptscriptstyle\circWscriptscript{}^{\kern1.5pt #2}_{\kern-1pt #1}}%
}}

\def\OM{\mathchoice
{\rlap{\kern3.2pt$\overline{\phantom{L}}$}M}
{\rlap{\kern3.2pt$\overline{\phantom{L}}$}M}
{\rlap{\kern2.4pt$\scriptstyle\overline{\phantom{L}}$}M}
{\rlap{\kern1.8pt$\scriptscriptstyle\overline{\phantom{L}}$}M}}

\def\mycirc{{\kern1pt\circ\kern2pt}}

\def\Bar#1{\bar{\bar{#1}}}
\let\bbar\hat

\def\Aut{\mathop{\rm Aut}\nolimits}

\def\Pic{\mathop{\rm Pic}\nolimits}

\def\Spec{\mathop{\rm Spec}\nolimits}

\def\deg{\mathop{\rm deg}\nolimits}

\def\PGL{\mathop{\rm PGL}\nolimits}

\def\reg{{\rm reg}}

\let\phi\varphi
\let\theta\vartheta
\let\epsilon\varepsilon
\let\setminus\smallsetminus

\let\emptyset\varnothing

\newcommand{\BP}{{\mathbb{P}}}
\newcommand{\BQ}{{\mathbb{Q}}}
\newcommand{\BR}{{\mathbb{R}}}

\newcommand{\BZ}{{\mathbb{Z}}}

\newcommand{\Fm}{{\mathfrak{m}}}

\newcommand{\Fp}{{\mathfrak{p}}}
\newcommand{\Fq}{{\mathfrak{q}}}

\newcommand{\CC}{{\cal C}}

\newcommand{\CP}{{\cal P}}

\newcommand{\CU}{{\cal U}}

\newcommand{\CX}{{\cal X}}
\newcommand{\CY}{{\cal Y}}
\newcommand{\CZ}{{\cal Z}}

\newcommand{\op}{\operatorname}

\newcommand{\Cst}{C^{\op{st}}}

\newcommand{\wbar}{{\overline{w}}}
\newbox\mybox
\def\arrover#1{\mathrel{
\setbox\mybox=\hbox spread 1.4em
{\hfil$\scriptstyle#1$\hfil}
\vbox{\offinterlineskip\copy\mybox
\hbox to\wd\mybox{\rightarrowfill}}}}

\def\larrover#1{\mathrel{
\setbox\mybox=\hbox spread 1.4em
{\hfil$\scriptstyle#1\vphantom{g}$\hfil}
\vbox{\offinterlineskip\copy\mybox
\hbox to\wd\mybox{\leftarrowfill}}}}

\def\ontoover#1{\mathrel{
\setbox\mybox=\hbox spread 1.4em
{\hfil$\scriptstyle#1\vphantom{g}$\hfil}
\vbox{\offinterlineskip\copy\mybox
\hbox to\wd\mybox{\rightarrowfill\hskip-2.8mm
$\rightarrow$}}}}
\def\leftontoover#1{\mathrel{
\setbox\mybox=\hbox spread 1.4em
{\hfil$\scriptstyle#1\vphantom{g}$\hfil}
\vbox{\offinterlineskip\copy\mybox
\hbox to\wd\mybox{$\leftarrow$\hskip-2.8mm
\leftarrowfill}}}}
\let\longto\longrightarrow
\let\into\hookrightarrow
\let\onto\twoheadrightarrow

\def\isoto{\mathrel{
\setbox\mybox=\hbox spread 0.9em
{\hfil$\scriptstyle\sim$\hfil}
\vbox{\offinterlineskip\copy\mybox
\hbox to\wd\mybox{\rightarrowfill}}}}

\def\Bigskip{\bigskip\bigskip}

\DeclarePairedDelimiter\floor{\lfloor}{\rfloor}

\newtheorem{Thm}{Theorem}[section]
\newtheorem{Prop}[Thm]{Proposition}
\newtheorem{Proposition}[Thm]{Proposition}

\newtheorem{Lem}[Thm]{Lemma}

\newtheorem{Cor}[Thm]{Corollary}

\newtheorem{Def}[Thm]{Definition}

\newtheorem{Rem}[Thm]{Remark}
\newtheorem{Ex}[Thm]{Example}
\newtheorem{Exs}[Thm]{Examples}

\newtheorem{Cons}[Thm]{Construction}

\newtheorem{Alg}[Thm]{Algorithm}
\newtheorem{Sum}[Thm]{Summary}

\numberwithin{Thm}{subsection}
\def\UseTheoremCounterForNextEquation{\setcounter{equation}{\value{Thm}}\addtocounter{Thm}{1}}

\def\qed{{\hskip0pt\unskip\unskip\nobreak\hfil\penalty50
\hskip1em\hbox{}\nobreak\hfil
{$\square$}
\parfillskip=0pt\finalhyphendemerits=0
\par}\medskip}

\newenvironment{Proof}
{\noindent{\bf Proof.}\ }
{\qed}

{\noindent{\bf Proof of #1.}\ }
{\qed}

\begin{document}

\title{\strut
\vskip-80pt
Reduction of Hyperelliptic Curves\\
in Residue Characteristic $2$\\
}
\author{
\begin{minipage}{.3\hsize}
Tim Gehrunger\\[12pt]
\small Department of Mathematics \\
ETH Z\"urich\\
8092 Z\"urich\\
Switzerland \\
tim.gehrunger@math.ethz.ch\\[9pt]
\end{minipage}
\qquad
\begin{minipage}{.3\hsize}
Richard Pink\\[12pt]
\small Department of Mathematics \\
ETH Z\"urich\\
8092 Z\"urich\\
Switzerland \\
pink@math.ethz.ch\\[9pt]
\end{minipage}
}
\date{\today}

\maketitle

\Bigskip

\begin{abstract}
Consider a hyperelliptic curve of genus $g$ over a field $K$ of characteristic zero. After extending $K$ we can view it as a marked curve with its $2g+2$ Weierstrass points. We present an explicit algorithm to compute the stable reduction of this marked curve for a valuation of residue characteristic~$2$ over a finite extension of~$K$. In the cases $g\le2$ we work out relatively simple conditions for the structure of this reduction.
\end{abstract}

{\renewcommand{\thefootnote}{}
}

\newpage
\tableofcontents
\newpage


\section{Introduction}
\label{Intro}

{\bf 1.1 Motivation and strategy:}
Let $K$ be a valued field of characteristic $\not=2$, and let $C$ be a hyperelliptic curve over~$K$, that is, a curve with  the  equation $z^2=f(x)$ for some polynomial~$f$. After extending $K$ if necessary, the curve admits stable reduction. While in principle there is a general algorithm to find a stable model, the goal of this article is to describe this model efficiently.

Our idea is to rigidify the situation using the Weierstrass points of~$C$. These are the ramification points of the canonical double covering $\pi\colon C\onto\bar C$ of a rational curve~$\bar C$. After extending $K$ we suppose that these points are defined over~$K$. We then propose to construct the stable model $\CC$ of $C$ as a marked curve with its Weierstrass points. The stable model of the unmarked curve is easily obtained from this by contracting irreducible components of the special fiber.

Why should this simplify the construction of the stable unmarked model? Because one can use the stable model of the rational curve $\bar C$ marked by the branch points of~$\pi$. These points are the zeros of $f$ and possibly the point $x=\infty$, and the stable model $\bar\CC$ of $\bar C$ as a marked curve can by computed easily. It turns out that the model $\CC$ dominates $\bar\CC$ and that some of its complexity is already present in~$\bar\CC$.  

When the residue characteristic is $\not=2$, we already know from \cite{GehrungerPink2021} and the references therein that $\CC$ is the normalization of $\bar\CC$ in the function field of~$C$. The present article therefore concentrates on the case of residue characteristic~$2$. The complications in that case can all be traced back to the fact that any ramification of a double covering in characteristic~$2$ is wild.

In that case it is not hard to show that $\CC$ is minimal among all semistable models of the marked curve $C$ that dominate~$\bar\CC$ (see Proposition \ref{AbsRelStabMod}). We therefore generalize the situation a little and begin with an arbitrary semistable model $\bar\CC$ of the rational marked curve $\bar C$ and aim to construct the minimal semistable model of the marked curve $C$ that dominates~$\bar\CC$. This has then become a purely local problem over~$\bar\CC$.

Another useful fact is that the hyperelliptic involution $\sigma$ of~$C$, that is, the covering involution of~$\pi$, extends uniquely to an involution $\sigma$ of~$\CC$. It turns out that the quotient scheme $\bbar\CC := \CC/\langle\sigma\rangle$ is a semistable model of $\bar C$ that dominates~$\bar\CC$ (Proposition \ref{BbarModelSemiStab}). The problem thus divides up into the problem of describing $\bbar\CC$ in explicit coordinates and then constructing $\CC$ as the normalization of $\bbar\CC$ in the function field of~$C$.

Over the smooth locus of~$\bar\CC$, all this has been essentially solved by Lehr and Matignon \cite{LehrMatignon2006}. Our main contribution thus lies in describing what to do over a double point of~$\bar\CC$.

\medskip
{\bf 1.2 Overview:}
We now explain the relevant issues in detail. First, to avoid the recurring need for field extensions and the resulting cumbersome changes of notation, we reduce the general case throughout to the case that $K$ is algebraically closed. Let $R$ denote its valuation ring and $k=R/\Fm$ its residue field. Let $v$ denote the valuation on $K$ that is normalized to $v(2)=1$. For any rational number $\alpha$ we choose a suitable fractional power $2^\alpha\in K$ with $v(2^\alpha)=\alpha$. For any Laurent polynomial $f$ over $K$ we let $v(f)$ denote the minimum of the valuations of its coefficients. We let $\bar C_0,\bbar C_0,C_0$ denote the closed fibers of $\bar\CC,\bbar\CC,\CC$, respectively.

As another preparation, for any Laurent polynomial $f\in R[x^{\pm1}]$ we set  
$$w(f)\ :=\ \sup \bigl\{ v(f-h^2) \bigm| h\in R[x^{\pm 1}] \bigr\}\ \in\ \BR\cup\{\infty\},$$
which measures how well $f$ can be approximated by squares. It will turn out that the precise values of $v(f-h^2)$ and $w(f)$ are irrelevant if they exceed~$2$. Accordingly, we call a decomposition of the form $f=h^2+g$ with $g,h\in R[x^{\pm 1}]$ \emph{optimal} if
$$v(g)=w(f)\quad\hbox{or}\quad v(g)>2.$$
We prove that an optimal decomposition of $f$ exists and can be computed effectively (Proposition \ref{OptDecompExists}).

\medskip
The construction of $\CC$ requires finding explicit local coordinates for the normalization of $\bar\CC$ in the function field of~$C$. To explain how suppose first that $x$ is a local coordinate of $\bar\CC$ near a smooth point $\bar p$ of its special fiber. After rescaling $z$ and $f$ by elements of $K^\times$ we may assume that $f$ has integral coefficients and non-zero reduction modulo~$\Fm$. Then an optimal decomposition $f=h^2+g$ with polynomials $g$, $h$ exists (Proposition \ref{OddDecomPol}) and the normalization near $\bar p$ is given by the coordinates $x$ and $t$ with $z=h+2^{\gamma/2}t$ for $\gamma := \min\{2,v(g)\}$ (see Proposition \ref{LehrPropositionCoveringPol} or Lehr \cite[Prop.\,1]{Lehr2001ReductionOP}). Using the resulting equation one can quickly decide where this normalization is smooth and therefore equal to~$\CC$. In particular this happens over any marked point~$\bar p$ (see Proposition \ref{SmoothMarkProp}).

\medskip
At an unmarked smooth point $\bar p$ where the normalization is singular, the theory of Lehr and Matignon \cite{LehrMatignon2006} tells us how to find $\bbar\CC$ by blowing up $\bar\CC$ near $\bar p$ in terms of the zeros of an auxiliary polynomial $S_f$ associated to~$f$. The normalization of this blowup in the function field of~$C$ is then semistable and therefore a local chart of~$\CC$. 
As in the situation of \cite[Thm.\,5.1]{LehrMatignon2006}, the irreducible components of $C_0$ above $\bar p$ are arranged in the form of an oriented tree with the components of genus $>0$ precisely at the ends. In particular, such points do not contribute to any bad reduction of the jacobian of~$C$.

\medskip
Now consider a double point $\bar p$ of~$\bar C_0$. Here another dichotomy occurs: Recall that the number of branch points of $\pi$ is $2g+2$, where $g$ is the genus of~$C$. Consequently there is either an even number of branch points on each side of $\bar p$, or an odd number on each side.
In the second case the normalization of $\bar\CC$ in the function field of~$C$ is already semistable with a unique double point above~$\bar p$ and can be computed explicitly (Proposition \ref{OddPointProp}), giving $\CC$ locally over~$\bar p$.

\medskip
So assume that the number of branch points on each side of $\bar p$ is even. Then locally near $\bar p$ the model $\bar\CC$ is isomorphic to ${\Spec R[x,y]/(xy-2^\alpha)}$ for some $\alpha>0$, which is called the \emph{thickness of~$\bar p$}. Using $y=2^\alpha\kern-1.5pt/x$ we can embed this coordinate ring into the ring of Laurent polynomials $K[x^{\pm 1}]$. This gives us the freedom to rescale $z$ by $K^\times x^\BZ$, which by the equation $z^2=f(x)$ amounts to rescaling $f$ by $K^\times x^{2\BZ}$. Since $\bar p$ is even, we can reduce ourselves to the case that $f$ lies in $R[x,y]/(xy-2^\alpha)$ and has a unit as constant term (Proposition \ref{PropHyperEquationNearP}).

Now observe that the inverse image of $\bar p$ in $\bbar C_0$ may contain irreducible components that lie between the proper transforms of the two irreducible components of $\bar C_0$ that meet at~$\bar p$ and of other irreducible components sticking out from those. 

The former are those that are described by coordinates of the form $x/2^\lambda$ for $0<\lambda<\alpha$. To identify them we study the behavior of optimal decompositions under substitutions of the form $x=2^\lambda u$. Specifically, we consider the function
$$ \wbar\colon\ \BQ\,\cap\, [0, v(a)] \longto\BR,\ \ \lambda\mapsto \wbar(\lambda) := \min \{ 2, w(f(2^\lambda u)) \}.$$
This function can be computed explicitly in terms of optimal decompositions (Propositions \ref{ScalOpt} and \ref{ScalOptDecomp}) and is piecewise linear concave (Proposition \ref{WBarProp}). We prove that the substitution $x=2^\lambda u$ yields an irreducible component of $\bbar C_0$ over $\bar p$ if and only if $\lambda$ is a break point of~$\wbar$ (Proposition \ref{ComponentsTypeBProp}).

Consider the blowup of $\bar\CC$ obtained by adjoining a component with coordinate $x/2^\lambda$ for each break point~$\lambda$. This is a model of $\bar C$ that lies between $\bar\CC$ and $\bbar\CC$, and whose normalization in the function field of $C$ is semistable above all double points. Finding the remaining irreducible components of $\bbar C_0$ above $\bar p$ thus reduces to the earlier problem over a smooth point of~$\bar\CC$.

\medskip
Combining everything this gives an explicit algorithm producing~$\bbar\CC$, from which $\CC$ can be constructed as the normalization in the same way as above. The procedure also provides further details: Propositions \ref{LehrPropositionCovering} and \ref{ComponentsTypeBDetailProp} tell us where the morphism $C_0\onto\bbar C_0$ is inseparable, separable, respectively \'etale and which irreducible components of $\bbar C_0$ decompose in~$C_0$. Propositions \ref{EvenPointNorm2} and \ref{EvenPointNorm<2} determine whether $C_0$ has one or two double points above a double point of $\bbar C_0$. Finally, Proposition \ref{Jacobian} lists some consequences for the reduction behavior of the jacobian of~$C$.


\medskip
{\bf 1.3 Small genus:}
%
Applying these methods to the case of genus $1$, we find that the type of stable reduction only depends on the stable marked reduction of $\bar C$ and the thickness of its double point if it has one. In total, there are $4$ different reduction types in this case.

 In the case of genus $2$, the reduction behavior of $C$ depends on the stable marked reduction of $\bar C$ and the thicknesses of  its  double points as well as on an additional parameter~$\delta$, which is the valuation of a certain expression in the coefficients of~$f$. In total, there are $54$ different reduction types in this case.  For the types of the  unmarked stable reduction and the reduction of the jacobian, our results yield relatively simple conditions  depending only on the thicknesses and $\delta$.

\medskip
{\bf 1.4 Structure of the paper:}
Chapter \ref{GlobalStableModel} contains preparatory material: Section~\ref{App} provides the justification for working over an algebraically closed field with Theorem \ref{NormalNormal}. In Section \ref{SemiStabCurves} we review basic facts about semistable and stable marked curves over~$R$. The following Section \ref{Genus0} concentrates on curves of genus~$0$, with a special emphasis on explicit local coordinates and algorithms for constructing semistable models from others. 
In the final Section \ref{HyperellCurves} we turn to hyperelliptic curves, providing the set-up for the remainder of the article. A summary of all the schemes and morphisms needed in our construction is given in Diagram \ref{AllCPDiagram}. From here on we assume that $R$ has residue characteristic~$2$.

In Chapter \ref{ApproxBySquares} we study optimal approximations of Laurent polynomials. In Section \ref{OptDecomp} we define and characterize them, and in Section \ref{OddDecomp} and Proposition \ref{LaurDecOpt} we construct those that are useful for us. In Section \ref{DoubCov} we show how they arise in computing the normalization of $R[x^{\pm1}]$ in a quadratic extension of $K(x)$.
In Section \ref{ScalingDecomp} we study their behavior under substitutions of the form $x=2^\lambda u$. In Section \ref{SepDecomp} we discuss how to separate positive and negative exponents in optimal decompositions, which can sometimes simplify explicit computations. In Section \ref{LaurDecomp} we discuss the decompositions that are used in the theory of Lehr and Matignon \cite{LehrMatignon2006}.

The local constructions of $\bbar\CC$ and $\CC$ are carried out in Chapter \ref{LocalStableModel}. Sections \ref{SmoothMarkedPoints} through \ref{EvenDoublePoints} describe the situation over each kind of point of $\bar C_0$ in turn: over smooth marked resp.\ unmarked points and over the two types of double points. In Section \ref{Algo} we combine these constructions into explicit algorithms and summarize the resulting properties of irreducible components and double points. We also briefly discuss some consequences for the reduction behavior of the jacobian of~$C$.

In the final Chapter \ref{Examples} we apply these methods to work out the reduction behavior in detail in genus $1$ and~$2$. For genus $2$ we only list the final results, leaving the detailed computations to look up in the associated computer algebra worksheets~\cite{GPComp}. 

\medskip
{\bf 1.5 Relation with other work:}
Semistable reductions of hyperelliptic curves have historically mainly been studied and constructed in residue characteristic $\not =$ 2, see for example the construction of Bosch in \cite{Bosch1980}. A more recent approach is the article \cite{DDMM} by Dokchitser, Dokchitser, Maistret, and Morgan, which describes the special fiber in their notion of cluster pictures. Similarly, in \cite{GehrungerPink2021}, the authors of the present article have given a description of the stable marked reduction.

In a series of articles \cite{Raynaud1970}, \cite{Raynaud1990} and \cite{Raynaud1999}, Raynaud studied the case of mixed characteristic $(0, 2)$ extensively. In particular, his articles provide some properties of the special fiber of semistable models under the additional assumption that the stable marked model of $\bar C$ is smooth. This condition is referred to as the case of equidistant geometry. Even though Raynaud restricted himself to this case, many of his theorems and ideas can be generalized and proved very valuable for understanding the general situation. Lehr and Matignon, building on Raynauds work, fully described the special fiber in the equidistant situation in their article \cite{LehrMatignon2006}. Reading their work was one of our key motivations to write this article.  Arzdorf and Wewers in \cite{ArzdorfWewers} generalize this and construct a semistable model of $C$ using the language of Berkovich analytic spaces. We were made aware of their article only after already completing most of our work. The construction carried out in our article is very similar to theirs, the main differences being that we are concerned with marked models and state everything in the language of schemes.  In a recent preprint \cite{FioreYelton}, Fiore and Yelton define and describe the so-called relatively stable model of hyperelliptic curves using the language of cluster pictures. Their construction is again very similar to that of Arzdorf and Wewers and to the one carried out in this article.  Being made public while we were already writing this article, there is some overlap with our work. In particular, in  \cite{FioreYelton}, the reduction behavior of their relative stable model is studied and explained in some of the cases. The conditions given partly rely on finding a root of a polynomial analogous to our stability polynomial, which we were able to avoid in our work.

Moreover, in  \cite{Liu1993}, Liu gives criteria for the type of the stable reduction of the unmarked curve $C$ in terms of Igusa invariants. This result is first proved in the setting of $\op{char}(k)\neq 2$ and carries over to the wild case by a moduli argument.  

\section{Semistable and stable models of curves}
\label{GlobalStableModel}

\subsection{Reduction to an algebraically closed field}
\label{App}

Throughout this article we let $R_1$ be a complete discrete valuation ring with quotient field~$K_1$. At several places we will need to replace $R_1$ by its integral closure in a finite extension of~$K_1$. Instead of having to say this repeatedly, we find it more convenient to work over an algebraic closure instead. 
So we fix an algebraic closure $K$ of $K_1$ and let $R$ denote the integral closure of $R_1$ in~$K$. Since $R_1$ is complete, the valuation on $K_1$ extends to a unique valuation with values in $\BQ$ on~$K$, whose associated valuation ring is~$R$. As this ring is not noetherian, we have to be careful when dealing with schemes over $\Spec R$.

We will be interested in intermediate fields $K_1\subset K_2\subset K_3\subset K$ that are finite over~$K_1$. For these the integral closure of $R_1$ in $K_2$ is $R_2 := R\cap K_2$ and that in $K_3$ is $R_3 := R\cap K_3$. These are again complete discrete valuation rings and finite extensions of~$R_1$. In particular they are free $R_1$-modules; hence the union $R$ of all $R_2$ is faithfully flat over~$R_1$.

By a scheme over any of these rings we will mean a scheme over the spectrum of this ring. In base extensions we will also drop the symbol $\Spec$. In the rest of this section we will discuss conditions for a normal scheme over $R$ to arise by base extension from a normal scheme over $R_2$ for some~$R_2$. First, by EGA4 \cite[Th.\,8.8.2, Cor.\,8.8.2.5]{EGA4} we have:

\begin{Prop}\label{EGA4882}
\begin{enumerate}
\item[(a)] For any finitely presented scheme $X$ over $R$ there exist $R_2$ as above and a finitely presented scheme $X_2$ over $R_2$ such that $X \cong X_2\!\times_{R_2}\!R$.
\item[(b)] For any $R_2$ as above and any finitely presented schemes $X_2$ and $Y_2$ over~$R_2$, for any morphism 
$\phi\colon X_2\!\times_{R_2}\!R \to Y_2\!\times_{R_2}\!R$
there exist $R_3$ as above such that $\phi$ comes by base extension from a morphism $X_2\!\times_{R_2}\!R_3 \to Y_2\!\times_{R_2}\! R_3$.
\item[(c)] Same as (b) for isomorphisms.
\end{enumerate}
\end{Prop}

\begin{Prop}\label{NormalDown}
In Proposition \ref{EGA4882} (a), if $X$ is integral and normal, then so is~$X_2$.
\end{Prop}

\begin{Proof}
The problem being local on~$X_2$, we may assume that $X_2=\Spec A_2$ for an $R_2$-algebra~$A_2$. Then $X \cong \Spec\,A_2\otimes_{R_2}R$ being integral and normal means that $A_2\otimes_{R_2}R$ is a normal integral domain. 
Next, multiplication by any nonzero element $a\in A_2$ induces nonzero homomorphisms of $A_2$-modules $A_2\onto A_2a \into A_2$. Since $R$ is faithfully flat over~$R_2$, these induce nonzero homomorphisms $A_2\otimes_{R_2}R \onto (A_2\otimes_{R_2}R)(a\otimes1) \into A_2\otimes_{R_2}R$. 
Thus $a\otimes1$ is again nonzero and the natural homomorphism $A_2\to A_2\otimes_{R_2}R$ is injective. In particular $A_2$ is itself integral.
Let $\tilde A_2$ denote its normalization. Then the inclusion $A_2\into \tilde A_2$ induces an integral extension ${A_2\otimes_{R_2}R} \into {\tilde A_2\otimes_{R_2}R}$ of integral domains with the same quotient field. Since $A_2\otimes_{R_2}R$ is already normal by assumption, this integral extension must be trivial. As $R$ is a faithfully flat $R_2$-algebra, it follows that $A_2=\tilde A_2$, and hence $X_2=\Spec \tilde A_2$ is normal.
\end{Proof}

\medskip
In the other direction we want to give conditions for the base change to $R$ of a normal scheme to be normal. We approach this in steps:

\begin{Prop}\label{NormalUp}
Let $X_1$ be a normal integral scheme that is flat of finite type over~$R_1$, whose closed fiber is generically smooth and for which $X_1\!\times_{R_1}\!K$ is integral and normal. Then $X_1\!\times_{R_1}\!R$ is integral and normal.
\end{Prop}

\begin{Proof}
It suffices to show that the scheme $X\times_{R_1} R_2$ is integral and normal for every $R_2$ as above. By assumption and Proposition \ref{NormalDown} this already holds for the generic fiber $X_1\times_{R_1}K_2$. By flatness it follows that $X_1\times_{R_1}R_2$ is integral.

Next, by Serre's criterion \cite[Th.\,5.8.6]{EGA4} a noetherian integral scheme is normal if and only it is (R1) and (S2). By assumption these properties already hold in the generic fiber of $X_1\times_{R_1}R_2$. The remaining points of codimension $1$ are the generic points of the special fiber. At all such points $X_1\times_{R_1}R_2$ is smooth over $R_2$ by assumption; hence it is (R1) there as well. Finally $X_1\times_{R_1}R_2\to X_1$ is flat and finite, so all its fibers are Artin and hence (S2). Since $X_1$ is (S2), it follows from \cite[Cor.\,6.4.2]{EGA4} that $X_1\times_{R_1}R_2$ is (S2) as well. Thus $X_1\times_{R_1}R_2$ is normal, as desired.
\end{Proof}

\begin{Prop}\label{NormalFin}
For any integral scheme $Y_1$ of finite type over~$R_1$, the normalization in any finite extension of the function field of $Y_1$ is finite over~$Y_1$. In particular it is again of finite type over~$R_1$.
\end{Prop}

\begin{Proof}
As $R_1$ is a complete noetherian local ring, it is excellent by \cite[Scholie 7.8.3 (iii)]{EGA4}. Since $Y_1$ is of finite type over~$R_1$, it is itself excellent by [loc.\ cit.\ (ii)]. Thus its normalization is finite over it by [loc.\ cit.\ (vi)], and hence again of finite type over~$R_1$.
\end{Proof}

\begin{Thm}\label{NormalNormal}
Let $Y_1$ be an integral scheme that is flat of finite type over~$R_1$. Let $L_1$ be a finite extension of the function field of~$Y_1$ such that $L_1\otimes_{R_1}K$ is a field. 
Then the normalization of $Y_1$ in $L_1\otimes_{R_1}K$ is finitely presented and arises by base change via ${R_2\into R}$ from the normalization of $Y_2$ in $L_1\otimes_{R_1}K_2$ for some $K_2$ finite over~$K_1$.
\end{Thm}

\begin{Proof}
By flatness the function field of $Y_1$ is an overfield of~$K_1$, and hence so is~$L_1$. Let $\tilde X$ be the normalization of $Y_1\times_{R_1}K_1$ in $L_1\otimes_{R_1}K$. As this is equally the normalization of $Y_1\times_{R_1}K$ in $L_1\otimes_{R_1}K$, which is of finite type over the field~$K$, it follows that $\tilde X$ is finite over $Y_1\times_{R_1}K$ by \cite[Scholie 7.8.3]{EGA4}. Thus $\tilde X$ is of finite type over~$K$ and hence finitely presented over~$R$. By Propositions \ref{EGA4882} (a) and \ref{NormalDown} there therefore exists $K_2$ finite over $K_1$ as above, such that $\tilde X$ arises by base extension from a normal integral scheme over~$K_2$. This means that $\tilde X = X_2\times_{R_2}K$, where $X_2$ is the normalization of $Y_1$ in $L_1\otimes_{R_1}K_2$. By Proposition \ref{NormalFin} this $X_2$ is again of finite type over~$R_2$. Moreover, as $X_2$ is reduced with dense generic fiber, its affine coordinate rings are $R_2$-torsion free; hence $X_2$ is flat over~$R_2$.

Now recall that $R_2$ is excellent by \cite[Scholie 7.8.3 (iii)]{EGA4}. Thus by de Jong \cite[Lemma 2.13]{deJongIHES1996} or Temkin \cite[Thm.\,3.5.5]{Temkin2010}, after replacing $K_2$ by a finite extension and $X_2$ by the corresponding normalization, we can assume that the closed fiber of $X_2\to\Spec R_2$ is generically smooth. Then $X_2$ satisfies the assumptions of Proposition \ref{NormalUp} with $R_2$ in place of~$R_1$, and so $X_2\!\times_{R_2}\!R$ is integral and normal. This is therefore the normalization of $Y_1$ in $L_1\otimes_{R_1}K$ and finitely presented over~$R$.
\end{Proof}

\medskip
In the rest of this article we apply the above results to the case that $Y_1$ is a semistable curve over~$R_1$ with generic fiber $\BP^1_{K_1}$ and $L_1$ is the function field of a hyperelliptic curve~$C_1$. In order to construct a good model of $C_1$ we need to replace $K_1$ at various places by some finite extension and $Y_1$ by a suitable blowup. By Theorem \ref{NormalNormal} we can instead work over the single field $K$ and avoid cumbersome changes of notation.

Let $v$ denote the valuation on $K$, and let $\Fm$ be the maximal ideal and $k:=R/\Fm$ the residue field of~$R$.

\subsection{Semistable curves}
\label{SemiStabCurves}

In this section we review basic known facts about stable marked curves over~$R$. See Knudsen \cite{Knudsen1983} or Liu \cite[\S10.3]{LiuAlgGeo2002} or Temkin \cite{Temkin2010} for the general definition and properties of semistable and stable curves over arbitrary schemes.

\medskip
Let $C$ be a connected smooth proper algebraic curve of genus $g$ over~$K$. By a \emph{model of~$C$} we mean a flat and finitely presented curve $\CC$ over $R$ with generic fiber~$C$. We call such a model \emph{semistable} if the special fiber $C_0$ is smooth except possibly for finitely many ordinary double points. Every double point $p\in C_0$ then possesses an \'etale neighborhood in $\CC$ which is \'etale over $\Spec R[x,y]/(xy-a)$ for some nonzero $a\in \Fm$, such that $p$ corresponds to the point $x=y=0$. Here the valuation $v(a)$ depends only on the local ring of $\CC$ at~$p$, for instance by Liu \cite[\S10.3.2 Cor.\,3.22]{LiuAlgGeo2002}. Following Liu 
\cite[\S10.3.1 Def.\,3.23]{LiuAlgGeo2002} we call $v(a)$ the \emph{thickness of~$p$}. 

\medskip
Any model is an integral separated scheme. Thus for any two models $\CC$ and $\CC'$ over~$R$, the identity morphism on $C$ extends to at most one morphism $\CC\to\CC'$. If this morphism exists, we say that $\CC$ \emph{dominates}~$\CC'$. This defines a partial order on the collection of all models of $C$ up to isomorphism. By blowing up one model one can construct many other models that dominate it. 
Conversely, one can construct the \emph{contraction} of an irreducible component with the following properties:

\begin{Prop}\label{ContractExist}
Assume that $\CC$ is normal with reducible closed fiber~$C_0$, and let $T$ be an irreducible component of~$C_0$. 
\begin{enumerate}
\item[(a)] There exists a normal model $\CC'$ that is dominated by~$\CC$, such that the morphism $\CC\to\CC'$ maps $T$ to a closed point $p'$ and induces an isomorphism $\CC\setminus T \stackrel{\sim}{\to} \CC\setminus\{p'\}$.
\item[(b)] The model $\CC'$ is unique up to unique isomorphism.
\item[(c)] If $\CC$ dominates another model~$\CX$, such that the morphism $\CC\to\CX$ maps $T$ to a closed point, then $\CC'$ dominates~$\CX$.
\end{enumerate}
\end{Prop}

\begin{Proof}
For (a) see \cite[\S6.7 Prop.\,4]{BoschLuetkebohmertRaynaud1990} or the proof of \cite[Prop.\,4.4.6]{Temkin2010}. For (c) see for instance \cite[Prop.\,4.3.2]{Temkin2010}. Finally, (c) implies (b).
\end{Proof}

Any semistable model is normal by \cite[\S10.3.1 Prop.3.15 (c)]{LiuAlgGeo2002}, so Proposition \ref{ContractExist} can be applied to it. We call an irreducible component $T$ \emph{unstable} if it is isomorphic to $\BP^1_k$ and contains at most two double points. 

\begin{Prop}\label{ContractSemiStab}
Suppose that $\CC$ is semistable. 
\begin{enumerate}
\item[(a)] The contraction $\CC'$ is semistable if and only if $T$ is unstable.
\item[(b)] In that case $p'$ is a smooth point if $T$ contains $1$ double point, respectively a double point if it contains $2$ double points.
\end{enumerate}
\end{Prop}

\begin{Proof}
See \cite[\S10.3.2 Lemma\,3.31]{LiuAlgGeo2002} and \cite[Cor.\,B.2]{Temkin2010}.
\end{Proof}

\begin{Prop}\label{RelStabMod}
Consider any model $\CX$ of $C$ over~$R$.
\begin{enumerate}
\item[(a)] Among the semistable models of $C$ that dominate~$\CX$ there exists a \emph{minimal} model~$\CC$, that is, such that every semistable model that dominates $\CX$ also dominates~$\CC$.
\item[(b)] This model is unique up to unique isomorphism.
\item[(c)] The morphism $\CC\onto\CX$ is an isomorphism at all points where $\CX$ is already semistable.
\item[(d)] A semistable model $\CC'$ that dominates $\CX$ is minimal if and only if no fiber of $\CC'\onto\CX$ contains an unstable irreducible component.
\end{enumerate}
\end{Prop}

\begin{Proof}
See Liu \cite[2.3-8]{Liu2006} or Temkin \cite[1.2-5]{Temkin2010}.
An extension of $R$ is rendered unnecessary by the reductions in Section~\ref{App}.
\end{Proof}

\begin{Prop}\label{VarRelStabModExists}
For any models $\CX_0,\ldots,\CX_n$ of $C$ over~$R$, there exists a minimal semistable model $\CC$ that dominates each~$\CX_i$ for all~$i$, and it is unique up to unique isomorphism.
\end{Prop}

\begin{Proof}
Since $\CX_0,\ldots,\CX_n$ are models of~$C$, the diagonal morphism $C \to \CY := \CX_0\times_R \ldots\times_R \CX_n$ is an isomorphism in the generic fiber. Let $\CZ$ denote the normalization of $\CY$ in the function field of~$C$. The fact that the $\CX_i$ are proper over~$R$  then implies that $\CZ$ is proper over~$\CX$. It is therefore a model of $C$ which dominates~$\CX$. By construction, any semistable model of $C$ which dominates $\CX$ and possesses morphisms ${\CC\to\CX_i}$ for all~$i$ must also dominate~$\CZ$. Thus the proposition follows by applying Proposition \ref{RelStabMod} with $\CZ$ in place of~$\CX$.
\end{Proof}

%


\begin{Prop}\label{ContractSeq}
Any morphism of semistable models is the composite of finitely many contractions of unstable irreducible components.
\end{Prop}

\begin{Proof}
Consider a morphism of semistable models $\CC\onto\CX$. If it is not yet an isomorphism, then $\CC$ is not minimal in the sense of Proposition \ref{RelStabMod} (a), so some fiber of $\CC\onto\CX$ contains an unstable irreducible component~$T$. By Proposition \ref{ContractSemiStab} the contraction $\CC'$ of $T$ is then semistable and by Proposition \ref{ContractExist} it dominates~$\CX$. The proposition thus follows by induction on the number of irreducible components of the special fiber of~$\CC$.
\end{Proof}

\begin{Prop}\label{SemistabMorphReg}
Consider a morphism $\pi\colon\CC\onto\CC'$ of semistable models, let $C_0\onto C_0'$ be the induced morphism of closed fibers, and let $Z$ be an irreducible component of~$C_0$ whose image $Z':=\pi(Z)$ is an irreducible component of~$C'_0$. Then $\pi(Z\cap C_0^\reg) \subset Z'\cap C_0^{\prime\reg}$.
\end{Prop}

\begin{Proof}
By Proposition \ref{ContractSeq} and induction it suffices to prove this when $\pi$ is the contraction of an unstable irreducible component $T\not= Z$ of~$C_0$. In that case $\pi$ is a local isomorphism outside $T$. Since $Z\cap C_0^\reg$ is contained in~$Z\setminus T$, its image is therefore contained in $Z'\cap C_0^{\prime\reg}$, as desired.
\end{Proof}

\medskip
Now consider an integer $n\ge0$ and distinct $K$-rational points $P_1,\ldots,P_n \in C(K)$. This turns $C$ into a \emph{smooth semistable marked curve} $(C,P_1,\ldots,P_n)$ over~$K$. If $\CC$ is a semistable model of $C$ such that these points extend to pairwise disjoint sections $\CP_1,\ldots,\CP_n \in \CC(R)$ which avoid all double points of the special fiber, we call $(\CC,\CP_1,\ldots,\CP_n)$ a \emph{semistable model} of $(C,P_1,\ldots,P_n)$ over~$R$. 

\medskip
From now on we assume that $2g+n\ge3$. Then the group of automorphisms of $C$ which preserve the given points is finite, and $(C,P_1,\ldots,P_n)$ is a smooth \emph{stable marked curve}. A \emph{stable model} of $(C,P_1,\ldots,P_n)$ over $R$ is a semistable model such that the group of automorphisms of the closed fiber which preserve the given sections is finite as well. A semistable model is stable if and only if its closed fiber possesses no irreducible component that is isomorphic to $\BP^1_k$ and contains at most two double or marked points. 
The special fiber $(C_0,P_{0,1},\ldots,P_{0,n})$ of a stable model is called \emph{stable reduction of} $(C,P_1,\ldots,P_n)$.

\begin{Prop}\label{AbsStabMod}
\begin{enumerate}
\item[(a)] A stable model $(\CC,\CP_1,\ldots,\CP_n)$ of $(C,P_1,\ldots,P_n)$ exists.
\item[(b)] This model is unique up to unique isomorphism.
\item[(c)] For every semistable model $(\CC',\CP_1',\ldots,\CP_n')$ the model $\CC'$ dominates $\CC$.
\end{enumerate}
\end{Prop}

\begin{Proof}
See Liu \cite[2.19-21]{Liu2006} or Temkin \cite[1.2-5]{Temkin2010} or Cuzub \cite[Th.\,3.4]{Cuzub2018}).
\end{Proof}

\begin{Prop}\label{AbsRelStabMod}
Let $(\CC,\CP_1,\ldots,\CP_n)$ be the stable model of $(C,P_1,\ldots,P_n)$. Let $\CX$ be a model of~$C$, such that $\CC$ dominates~$\CX$ and $P_1,\ldots,P_n$ extend to pairwise disjoint sections of the smooth locus of~$\CX$. Then $\CC$ is the minimal semistable model of $C$ that dominates~$\CX$.
\end{Prop}

\begin{Proof}
Let $\CC'$ be the minimal semistable model of $C$ that dominates $\CX$ from Proposition \ref{RelStabMod}. Then $\CC$ dominates~$\CC'$ by Proposition \ref{RelStabMod} (a).
Conversely, since $\CX$ is already semistable in a neighborhood of the sections extending $P_1,\ldots,P_n$, the morphism $\CC'\to\CX$ is an isomorphism there by \ref{RelStabMod} (c). Thus these points extend to pairwise disjoint sections $\CP_1',\ldots,\CP_n'$ of the smooth locus of~$\CC'$, making $(\CC',\CP_1',\ldots,\CP_n')$ a semistable marked model. Thus $\CC'$ dominates $\CC$ by Proposition \ref{AbsStabMod} (c).
Together this shows that $\CC\cong\CC'$.
\end{Proof}

\begin{Rem}\label{AbsStabModLoc}
\rm In the situation of Proposition \ref{AbsRelStabMod}, the construction of the stable model $(\CC,\CP_1,\ldots,\CP_n)$ becomes a local problem at the points where $\CX$ is not yet semistable, and one can examine these points separately.
\end{Rem}

\subsection{Semistable curves of genus $0$}
\label{Genus0}

In this section we collect a number of special results in genus~$0$, emphasizing facts about explicit coordinates that are hard to find in the existing literature. For this we fix a connected smooth projective algebraic curve $\bar C$ of genus $0$ over~$K$. 

\medskip
Consider a semistable model $\bar\CC$ of $\bar C$ over~$R$ and let $\bar C_0$ denote its closed fiber. Then by assumption $\bar C_0$ is a connected projective curve over $k$ that is smooth except for ordinary double points, and by flatness it has arithmetic genus~$0$. Thus $\bar C_0$ is a union of rational curves isomorphic to~$\BP^1_k$, which meet at the double points and are arranged in the form of a tree (see for instance Cuzub \cite[\S4]{Cuzub2018}).
An irreducible component of $\bar C_0$ is called \emph{stable} if it contains at least $3$ double points; otherwise it is called \emph{unstable}. An irreducible component that contains only one double point is called a \emph{leaf}. Any distinct irreducible components $T$ and $T'$ are connected by a unique shortest path across double points and possibly other irreducible components. We say that these other irreducible components \emph{lie between $T$ and~$T'$}. 

\medskip
Any unstable irreducible component $T$ of $\bar C_0$ can be contracted to a point in another semistable model by Proposition \ref{ContractExist}, and the image of $T$ is a smooth point if $T$ is a leaf, respectively a double point if not. 
By iterating this procedure one can construct many more semistable contractions. For instance, consider any double point $\bar p$ of $\bar C_0$ and let $I$ be the set of irreducible components in one of the two connected components of $\bar C_0\setminus\{\bar p\}$. Then by starting at the leaves in $I$ and iterating one finds a semistable contraction which maps this connected component to a smooth point and is an isomorphism on the complement. 

Iterating this again, for any given irreducible component $T\subset\bar C_0$ one can contract all other irreducible components to smooth points in a semistable model $\Bar\CC$. Then $\Bar\CC$ is a smooth model of $\bar C$ and therefore isomorphic to $\BP^1_R$ (see for instance Liu \cite[Ch.8 Ex. 3.5]{LiuAlgGeo2002}). Since $\bar\CC\onto\Bar\CC$ is an isomorphism over a neighborhood of $T\cap C_0^\reg$, it follows that any smooth point $\bar p\in\bar C_0$ possesses an open neighborhood in $\bar\CC$ that is isomorphic to an open subscheme of $\Spec R[x]$.  

Similarly, let $I$ be the set of irreducible components that do not meet a given double point $\bar p$ of~$\bar C_0$. Then by iterating the above procedure one can find a semistable contraction such that $\bar\CC\onto\Bar\CC$ is an isomorphism over a neighborhood of~$\bar p$ and the closed fiber of $\Bar\CC$ possesses only the two irreducible components adjacent to the image of~$\bar p$. For this there then exist explicit global coordinates $x$ and~$y$, such that
\UseTheoremCounterForNextEquation
\begin{equation}\label{Genus0DoublePointNbhd}
\Bar\CC\ \cong\ \Spec R[\tfrac{1}{x}]\ \cup\ \Spec R[x,y]/(xy-a)\ \cup\ \Spec R[\tfrac{1}{y}]
\end{equation}
for some nonzero $a\in\Fm$, and $\bar p$ corresponds to the point $x=y=0$ in the middle chart
(compare Cuzub \cite[discussion following Def.\,4.7]{Cuzub2018}).
In particular it follows that some open neighborhood of $\bar p$ in $\bar\CC$ is isomorphic to an open subscheme of $\Spec R[x,y]/(xy-a) = \Spec R[x,\tfrac{a}{x}]$.

\medskip
To describe the smooth models of $\bar C$ in terms of coordinates fix a rational function $x$ on $\bar C$ that yields an isomorphism $\bar C\cong\BP^1_K$. Then any other isomorphism $\bar C\cong\BP^1_K$ differs from this by an element of $\Aut_K(\BP^1_K) \cong \PGL_2(K)$. To make this precise abbreviate $A(x) := \smash{\frac{ax+b}{cx+d}}$ for any homothety class $A=\smash{[\binom{a\ b}{c\ d}]} \in \PGL_2(K)$. Then for any $A,B\in\PGL_2(K)$ we have $A(B(x))=(AB)(x)$. Thus for two substitutions $x=A(y)$ and $x=B(z)$ we have $z = (B^{-1}A)(y)$, and this substitution is an automorphism of $\BP^1_R$ if and only if $B^{-1}A\in\PGL_2(R)$. The smooth model with coordinate $y$ therefore depends only on the coset $A\cdot\PGL_2(R)$, and the smooth models up to isomorphism are in bijection with $\PGL_2(K)/\PGL_2(R)$. 

Let $B$ denote the subgroup of upper triangular matrices in $\PGL_2$. Then by the Iwasawa decomposition $\PGL_2(K) = B(K)\cdot\PGL_2(R)$, the inclusion $B\into\PGL_2$ induces a bijection from $B(K)/B(R)$ to $\PGL_2(K)/\PGL_2(R)$. Therefore any smooth model of $\bar C$ can be described by a coordinate $y$ such that $x=ay+b$ for some pair $(a,b)\in K^\times\times K$, which is unique up to 
the action $(a,b)\mapsto (ua,ub+va)$ for all $(u,v)\in R^\times\times R$.

\medskip
Returning to an arbitrary semistable model $\bar\CC$, for every irreducible component of $\bar C_0$ choose a contraction $\bar\CC\onto\bar\CC_i \cong\BP^1_R$ which is an isomorphism generically on this irreducible component. If their number is $n$, the diagonal morphism ${\bar\CC\to\bar\CC_1\times_R\ldots\times_R\bar\CC_n}$ is finite, and since $\bar\CC$ is normal, it follows that $\bar\CC$ is the normalization of $\bar\CC_1\times_R\ldots\times_R\bar\CC_n$ in the function field of~$\bar C$. This shows that $\bar\CC$ is determined by the models~$\bar\CC_i$ and can be constructed explicitly from them.

More generally, we will show how to construct new semistable models from given ones by adjoining irreducible components. We divide such irreducible components into the following types. Consider a semistable model $\bbar\CC$ of $\bar C$ which dominates~$\bar\CC$. Then the morphism of the closed fibers $\bbar C_0 \onto \bar C_0$ maps some irreducible components isomorphically to their images and contracts the others to closed points. 

\begin{Def}\label{Notationabcd}
An irreducible component of $\bbar C_0$ is called
\begin{itemize}
\item of type (a) if it maps isomorphically to an irreducible component of $\bar C_0$;
\item of type (b) if it lies between irreducible components of type (a);
\item of type (c) if it is not of type (a) or (b) and is not a leaf;
\item of type (d) if it is not of type (a) or (b) and is a leaf.
\end{itemize} 
For a sketch of this see Figure~\ref{FigB}.
\end{Def}

\begin{figure}[h] \centering \FigB
\caption{Sketch of the morphism $\bbar C_0 \onto \bar C_0$. Irreducible components of type (a) are drawn in black, those of type (b) in orange, those of type (c) in green, and those of type (d) in blue.}\label{FigB}
\end{figure}


First we look at components above a smooth point $\bar p\in \bar C_0$. For this we identify a neighborhood of $\bar p$ in $\bar\CC$ with an open subscheme of $\Spec R[x]$, such that $\bar p$ corresponds to the point $x=0$.

\begin{Prop}\label{AddCompsTypeCDCoord}
An irreducible component $\bbar T$ of $\bbar C_0$ is a component of type (c) or (d) above $\bar p$ if and only if it is given by a coordinate $y$ with $x=ay+b$ for $a,b\in\Fm$ with $a\neq0$.
\end{Prop}

\begin{Proof}
Choose a coordinate $y$ along $\bbar T$ such that $x=ay+b$ for $a\in K^\times$ and $b\in K$. The valuation at the generic point of $\bbar T$ then satisfies $v(x) = v(ay+b) = \min\{v(a),v(b)\}$. Thus $\bbar T$ maps to $\bar p$ if and only if this number is $>0$. In this case $\bbar T$ must be a component of type (c) or (d), as desired.
\end{Proof}

\begin{Cons}\label{AddCompsTypeCDCons}
\rm Conversely, suppose that for every $i$ in a finite set $I$ we are given a substitution $x=a_ix_i+b_i$ with $a_i,b_i\in\Fm$ and $a_i\neq0$. Let $\bar\CC_i$ denote a smooth model of $\bar C$ with the global coordinate~$x_i$. We will construct a minimal semistable model $\bbar\CC$ of $\bar C$ that dominates $\bar\CC$ as well as all~$\bar\CC_i$.

If $I=\emptyset$ there is nothing to do. Otherwise the number 
$$\alpha\ :=\ \min\bigl(\{v(a_i) \mid i\in I\} \cup \{v(b_i-b_j) \mid i,j\in I,\ i\neq j\}\bigr).$$
is finite and positive. Choose any nonzero $a\in\Fm$ such that $v(a)=\alpha$, and any $b\in R$ such that $v(b_i-b)\ge\alpha$ for all~$i$. Write $x=ay+b$ with a new coordinate $y$ on~$\bar C$. Then the blowup of $\Spec R[x]$ in the ideal $(x-b,a)$ is the union of the affine charts
$$\Spec R[x,\tfrac{a}{x-b}]\ =\ \Spec R[x-b,\tfrac{1}{y}] \quad\hbox{and}\quad \Spec R[y].$$
Gluing this with $\bar\CC\setminus\{\bar p\}$ over a neighborhood of $\bar p$ yields a semistable model $\tilde\CC$ of $\bar C$ that dominates~$\bar\CC$. Its exceptional fiber $E$ is the irreducible component of the closed fiber of $\tilde\CC$ with the coordinate~$y$. For a sketch of this see Figure \ref{FigE} below.
\begin{figure}[h] \centering \FigE
\caption{Sketch of Construction \ref{AddCompsTypeCDCons}.}\label{FigE}
\end{figure}

Now observe that solving the equation $ay+b = x = a_ix_i+b_i$ for $y$ yields the substitution $y = \tfrac{a_i}{a}\, x_i + \tfrac{b_i-b}{a}$, which by construction has coefficients in~$R$. For any $i$ with $v(a_i)=\alpha$ we have $\tfrac{a_i}{a}\in R^\times$; hence $y$ defines the same smooth model as~$x_i$ and we have already constructed the associated irreducible component~$E$. 

For all other $i$ we have $v(a_i)>\alpha$ and hence $\tfrac{a_i}{a}\in\Fm$. We group these indices into finitely many subsets $I_\nu$ according to the residue class of $\smash{\tfrac{b_i-b}{a}}$ modulo~$\Fm$. For each $\nu$ we choose a representative $c_\nu\in R$ of this residue class and consider the substitution $y=z_\nu+c_\nu$  with $\Spec R[y] = \Spec R[z_\nu]$. Then the point $y=c_\nu$ on $E$ is given equivalently by $z_\nu=0$. Also, for all $i\in I_\nu$ the resulting substitutions $z_\nu = \tfrac{a_i}{a}\, x_i + (\tfrac{b_i-b}{a}-c_\nu)$ now have coefficients in~$\Fm$, just as in the original problem. 

Moreover, each $I_\nu$ is now a proper subset of~$I$. Indeed, this is clear if $v(a_i)=\alpha$ for some~$i$, because this index does not lie in~$I_\nu$. Otherwise by construction there exist $i<j$ with $v(b_i-b_j)=\alpha$, so that $\tfrac{b_i-b}{a}$ and $\tfrac{b_j-b}{a}$ are not congruent modulo~$\Fm$. Thus again each $I_\nu$ is a proper subset of~$I$.

By recursion we can therefore assume that for every~$\nu$, we have already constructed a semistable model that dominates $\tilde\CC$, is isomorphic to $\tilde\CC$ outside that point, and contains the desired irreducible components for all $i\in I_\nu$. By gluing these models over $\tilde\CC$ we obtain the desired model $\bbar\CC$. 
\end{Cons}

\begin{Prop}\label{AddCompsTypeCDProp}
The model $\bbar\CC$ constructed in \ref{AddCompsTypeCDCons} is, up to isomorphism, the unique minimal semistable model of $\bar C$ that dominates $\bar\CC$ and whose closed fiber possesses an irreducible component with coordinate $x_i$ for each $i\in I$.
\end{Prop}

\begin{Proof}
By construction $\bbar\CC$ is a semistable model that dominates $\bar\CC$ and whose closed fiber possesses an irreducible component with coordinate $x_i$ for each $i\in I$. 

We prove that $\bbar\CC$ is minimal by induction over $|I|$. In the case $I=\emptyset$ this holds trivially because $\bbar\CC=\bar\CC$. Otherwise let $\tilde\CC$ and $E$ be as in Construction \ref{AddCompsTypeCDCons}. 
If $\bbar\CC$ is not minimal, by Propositions \ref{ContractExist} and \ref{ContractSemiStab} some irreducible component $\bbar T$ of its closed fiber can be contracted, obtaining another semistable model $\check\CC$ with the same properties. Then $\bbar T$ must lie over $\bar p$ and cannot be one of the components with coordinate~$x_i$. Also, by the induction hypothesis $\bbar\CC$ is already minimal among all semistable models that dominate $\tilde\CC$ and possess an irreducible component with coordinate $x_i$ for each $i\in I$. Thus $\check\CC$ cannot dominate~$\tilde\CC$, leaving only the case that $\bbar T$ maps isomorphically to $E\subset\tilde\CC$. Then $E$ does not have a coordinate~$x_i$, and Construction \ref{AddCompsTypeCDCons} shows that $E$ contains at least $3$ double points. But then $\bbar T$ also contains at least $3$ double points, contradicting the semistability of $\check\CC$ in Proposition \ref{ContractSemiStab} (a). We have thus reached a contradiction, proving the minimality of~$\bbar\CC$.

Finally, the uniqueness of $\bbar\CC$ follows by applying Proposition \ref{VarRelStabModExists} to the model $\bar\CC$ and the smooth models associated to the coordinates $x_i$ for all $i\in I$.
\end{Proof}

\begin{Rem}\label{AddCompsTypeCDRem}
\rm Any model that dominates $\bar\CC$ and is isomorphic to $\bar\CC$ outside $\bar p$ can be constructed as in \ref{AddCompsTypeCDCons}, for instance by letting the process run with the coordinates from Proposition \ref{AddCompsTypeCDCoord} for all irreducible components above~$\bar p$. It is also enough apply the process with the components of type (d) only, because the construction automatically adjoins the necessary components of type (c) to ensure semistability.
\end{Rem}

Now we look at components above a double point of~$\bar C_0$. For this we identify a neighborhood of $\bar p$ in $\bar\CC$ with an open subscheme of $\Spec R[x,y]/(xy-a)$ for some nonzero $a\in\Fm$, such that $\bar p$ corresponds to the point $x=y=0$. As before let $\bbar\CC$ be a semistable model of $\bar C$ which dominates~$\bar\CC$ and with closed fiber~$\bbar C_0$.

\begin{Prop}\label{AddCompsTypeBCDCoord}
An irreducible component $\bbar T$ of $\bbar C_0$ is a component 
\begin{itemize}
\item of type (b) over $\bar p$ if and only if it is given by a coordinate $z$ with $x=bz$ for some nonzero $b\in\Fm$ such that $\frac{a}{b}\in\Fm$;
\item of type (c) or (d) over $\bar p$ if and only if it is given by a coordinate $z$ with $x=bz+c$ for nonzero $b,c\in\Fm$ such that $\frac{b}{c},\frac{a}{c}\in\Fm$. 
\end{itemize}
\end{Prop}

\begin{Proof}
Choose a coordinate $z$ on $\bbar T$ such that $x=bz+c$ with $b\in K^\times$ and $c\in K$. The valuation at the generic point of $\bbar T$ then satisfies $v(x) = v(bz+c) = \min\{v(b),v(c)\}$ and hence $v(y) = v(\tfrac{a}{x}) = v(a)-\min\{v(b),v(c)\}$. Thus $\bbar T$ maps to $\bar p$ if and only if both these numbers are $>0$, that is, if $0 < \min\{v(b),v(c)\} < v(a)$. Let us assume this.

Suppose first that $v(b)\le v(c)$. Then we have $\frac{c}{b}\in R$ and can replace $z$ by $z+\frac{c}{b}$, which is also a coordinate for~$\smash{\bbar T}$. Afterwards we have $x=bz$ with $b,\frac{a}{b}\in\Fm$. Let $\bar\CX$ be the blowup of ${\Spec R[x,y]/(xy-a)}$ in the ideal $(x,b)$, which is the union of the affine charts
$$\Spec R[x,\tfrac{b}{x}] \ \cong\ \Spec R[x,w]/(xw-b) \quad\hbox{and}\quad \Spec R[\tfrac{x}{b},\tfrac{a}{x}]\ \cong\ \Spec R[z,y]/(zy-\tfrac{a}{b}).$$
Its exceptional fiber $E$ has the coordinate $z$ and meets the proper transforms of the irreducible components $x=0$ respectively $y=0$ of $\bar C_0$ in a double point each. Thus $E$ is an irreducible component of type (b) above~$\bar p$. Moreover, locally near $\bar p$ the morphism $\bbar\CC\onto\bar\CC$ must factor through~$\bar\CX$ and map $\bbar T$ isomorphically to~$E$. Thus $\bbar T$ is a component of type (b) above~$\bar p$, finishing the first case.

Suppose now that $v(b)>v(c)$. Then $c$ is nonzero with $c,\frac{b}{c},\frac{a}{c}\in\Fm$. Let $\bar\CX$ be the blowup of ${\Spec R[x,y]/(xy-a)}$ in the ideal $(x,c)$, which is the union of the affine charts
$$\Spec R[x,\tfrac{c}{x}] \quad\hbox{and}\quad \Spec R[\tfrac{x}{c},\tfrac{a}{x}].$$
As we have seen above its exceptional fiber $E$ is a component of type (b) above~$\bar p$. Write $x=cw+c$, so that $w$ is a coordinate along $E$ and $w=0$ defines a smooth point $\tilde p$ on it. The equation $cw+c = x = bz+c$ then reduces to $w=\tfrac{b}{c}z$. Thus the blowup of $\bar\CX$ in the ideal $(w,\tfrac{b}{c})$ at $\tilde p$ has another exceptional divisor $E'$ with the coordinate~$w$, which meets the proper transform $\smash{\tilde E}$ of $E$ in a double point that is distinct from the two double points coming from~$\bar\CX$. For a sketch of this see Figure~\ref{FigD}.
\begin{figure}[h] \centering \FigD
\caption{Sketch for the proof of Proposition \ref{AddCompsTypeBCDCoord}.}\label{FigD}
\end{figure}

Gluing this with $\bar\CC\setminus\{\bar p\}$ over a neighborhood of $\bar p$ yields a semistable model $\tilde\CC$ of $\bar C$ that dominates $\bar\CC$ and whose special fiber possesses an irreducible component with coordinate~$z$. Since $\tilde E$ corresponds to a stable irreducible component above~$\bar p$, this $\tilde\CC$ is a minimal semistable model with these properties. By the uniqueness in Proposition \ref{VarRelStabModExists}, it follows that $\bbar\CC$ dominates~$\tilde\CC$. As $\bbar T$ maps isomorphically to~$E'$, it cannot lie between irreducible components of type (b) and is therefore a component of type (c) or (d) above~$\bar p$, finishing the second case.
\end{Proof}

\begin{Cons}\label{AddCompsTypeBCons}
\rm Conversely, suppose that we are given a sequence of nonzero elements $1\,{=}\,b_0,b_1, \allowbreak \ldots,b_r,b_{r+1}\,{=}\,a$ with $r\ge0$ such that $\tfrac{b_i}{b_{i-1}}\in\Fm$ for all $1\le i\le r+1$. Put 
$$\bar\CX_i\ :=\ \Spec R\bigl[\tfrac{x}{b_{i-1}},\tfrac{b_i}{x}\bigr]\ \cong\ \Spec R[x_i,y_i]/(x_iy_i-\tfrac{b_i}{b_{i-1}})$$
and glue these charts together to a scheme $\bar\CX$ over the intersections 
$$\bar\CX_i\cap\bar\CX_{i+1}\ =\ \Spec R\bigl[\tfrac{x}{b_i},\tfrac{b_i}{x}\bigr]$$
for all $1\le i\le r$. Then $\bar\CX$ is a local model of $\bar C$ which dominates $\Spec R[x,\tfrac{a}{x}]$. Its exceptional fiber is consists of $r$ copies of $\BP^1_k$ with coordinates $\tfrac{x}{b_i}$ for all $1\le i\le r$, which are arranged in sequence such that each meets the next and the outer two meet the proper transforms of the respective irreducible components below. Gluing this with $\bar\CC\setminus\{\bar p\}$ over a neighborhood of $\bar p$ yields a semistable model $\bbar\CC$ of $\bar C$ that dominates~$\bar\CC$. By construction this model has $r$ irreducible components of type (b) above~$\bar p$. 
\begin{figure}[h] \centering \FigC
\caption{Sketch of Construction \ref{AddCompsTypeBCons}.}\label{FigC}
\end{figure}
\end{Cons}

\begin{Prop}\label{AddCompsTypeBProp}
The model $\bbar\CC$ constructed in \ref{AddCompsTypeBCons} is, up to isomorphism, the unique minimal semistable model of $\bar C$ that dominates $\bar\CC$ and whose closed fiber possesses an irreducible component with coordinate $\tfrac{x}{b_i}$ for every $1\le i\le r$.
\end{Prop}

\begin{Proof}
Direct consequence of the construction and Proposition \ref{VarRelStabModExists}.
\end{Proof}

\begin{Rem}\label{AddCompsTypeBRem}
\rm Any model that dominates $\bar\CC$ and is isomorphic to $\bar\CC$ outside~$\bar p$ can be constructed by first adjoining all components of type (b) as in \ref{AddCompsTypeBCons} and then all components of type (c) and (d) by the process in \ref{AddCompsTypeCDCons}.
\end{Rem}

\medskip
As a last topic in this section we describe an efficient construction of the stable marked model of a curve of genus~$0$, in which all coordinates are obtained from each other by affine linear coordinate changes. We begin with a local construction.

\begin{Cons}\label{StabMarkedGenus0Local}
\rm Suppose that we are given a chart $\bar\CU=\Spec R[x]$ of some model of $\bar C$ and a nonempty finite set $I\subset R$ whose elements are not all congruent modulo~$\Fm$. These elements represent sections of $\bar\CU$ which may partly but not completely meet in the closed fiber. Partition $I$ into nonempty proper subsets $I_\nu$ according to the residue class modulo~$\Fm$. Then any subset with $|I_\nu|=1$ represents a section that is disjoint from the other sections. For any subset with $|I_\nu|>1$ the value 
$$\alpha_\nu\ :=\ \min\{v(\xi-\xi')\mid \xi,\xi'\in I_\nu\}$$
is finite and $\ge0$. Choose an element $a_\nu\in R\setminus\{0\}$ with $v(a_\nu)=\alpha_\nu$ and an element $b_\nu\in R$ such that $v(\xi-b_\nu)\ge\alpha_\nu$ for all $\xi\in I_\nu$.
In the coordinate $x_\nu=\frac{x-b_\nu}{a_\nu}$ the elements of $I_\nu$ then correspond to the elements of $\smash{I'_\nu := \{\frac{\xi-b_\nu}{a_\nu}\mid \xi\in I_\nu\}}$. By construction this is a nonempty subset of $R$ whose elements are not all congruent modulo~$\Fm$. Consider the blowup of $\bar\CU\cong\Spec R[x]$ in the ideal $(x-b_\nu,a_\nu)$, which in explicit coordinates is given by
$$\Spec R[x_\nu] \cup \Spec R[x,x_\nu^{-1}].$$
Let $\bar\CU' \onto \bar\CU$ be the result of gluing together these local blowups over neighborhoods of the sections $x=b_\nu$ for all $\nu$ with $|I_\nu|>1$. Then the given sections of $\bar\CU$ lift to sections of~$\bar\CU'$, such that those coming from $I_\nu$ and only those land in the chart $\bar\CU'_\nu := \Spec R[x_\nu]$.

We can now repeat the construction with $(\bar\CU'_\nu,I'_\nu)$ in place of $(\bar\CU,I)$, as long as there is a subset with $|I_\nu|>1$. Since $|I'_\nu|<|I|$, this process terminates. Gluing the respective local blowups yields a semistable modification $\tilde\CU \onto \bar\CU$ such that the original sections lift to disjoint sections of the smooth locus of~$\tilde\CU$. Let $\tilde U_0\onto U_0$ denote the respective closed fibers. Then the construction guarantees that every irreducible component of $\tilde U_0$ that maps to a point in $U_0$ contains at least three double or marked points. The assumption also implies that the proper transform of $U_0$ in~$\tilde U_0$, which is still isomorphic to $\Spec k[x]$, contains at least two double or marked points.
\end{Cons}

\begin{Cons}\label{StabMarkedGenus0}
\rm Now suppose that $\bar C$ is marked by $n\ge3$ distinct rational points $\bar P_1,\ldots,\bar P_n$. To start the process from Construction \ref{StabMarkedGenus0Local} we choose a coordinate $x_0$ on~$\bar C$, that is, an isomorphism $\bar C\cong\BP^1_K$, such that $\bar P_1$ corresponds to the point~$\infty$. Then the other points $\bar P_i$ correspond to distinct elements $\xi_i\in K$. Since $n\ge3$, the value 
$$\alpha\ :=\ \min\{v(\xi_i-\xi_j)\mid 2\le i<j\le n\}$$
is finite. Choose an element $a\in K^\times$ with $v(a)=\alpha$ and an element $b\in K$ such that $v(\xi_i-b)\ge\alpha$ for all $2\le i\le n$.
In the coordinate $x=\frac{x_0-b}{a}$ the points $\bar P_2,\ldots,\bar P_n$ then correspond to the elements of $I := \smash{\{\frac{\xi_i-b}{a}\mid 2\le i\le n\}}$. They therefore extend to sections of $\bar\CU:=\Spec R[x]$, while the point $\bar P_1$ extends to the section $x^{-1}=0$ of $\Spec R[x^{-1}]$. Gluing the modification $\tilde\CU \onto \bar\CU$ from Construction \ref{StabMarkedGenus0Local} with $\Spec R[x^{-1}]$ over a neighborhood of the section $x^{-1}=0$ yields a semistable model $\bar\CC$ of~$\bar C$, such that $\bar P_1,\ldots,\bar P_n$ extend to disjoint sections $\bar\CP_1,\ldots,\bar\CP_n$ of the smooth locus of~$\bar\CC$. Moreover, the construction shows that every irreducible component of the closed fiber of $\bar\CC$ contains at least three double or marked points. Thus  $(\bar\CC,\bar\CP_1,\ldots,\bar\CP_n)$ is the stable model of $(\bar C,\bar P_1,\ldots,\bar P_n)$.
\end{Cons}

\begin{Rem}\label{DDMMRef}
\rm A related way of describing semistable curves of genus zero is that of cluster pictures from Dokchitser-Dokchitser-Maistret-Morgan \cite[\S4]{DDMM}.
\end{Rem}

\subsection{Hyperelliptic curves}
\label{HyperellCurves}

Now let $C$ be a \emph{hyperelliptic curve} of genus $g$ over~$K$. Thus $C$ is a connected smooth proper algebraic curve which comes with a double covering $\pi\colon C \onto \bar C$ of a rational curve $\bar C\cong\BP^1_K$. Often the genus $g$ is required to be $\ge2$, but in this article we only assume $g\ge1$.

\medskip
Consider a model $\bar\CC$ of $\bar C$ over~$R$. We say that a model $\CC$ of $C$ over $R$ \emph{dominates} $\bar\CC$ if and only if $\pi$ extends to a morphism $\CC\onto\bar\CC$.
By applying Proposition \ref{RelStabMod} in the case that $\CX$ is the normalization of $\bar\CC$ in the function field of~$C$, there exists a minimal semistable model of $C$ that dominates~$\bar\CC$, and it is unique up to unique isomorphism.

\medskip
Throughout the rest of this article we assume that $K$ has characteristic~$0$. Then the covering $\pi$ is only tamely ramified, and by the Hurwitz formula it is ramified at precisely $2g+2$ closed points, namely, at the Weierstrass points of~$C$. Let $P_1,\ldots,P_{2g+2} \in C(K)$ denote these points and $\bar P_1,\ldots,\bar P_{2g+2} \in \bar C(K)$ their images under~$\pi$. Since $2g+2\ge4$, both $(C,P_1,\ldots,P_{2g+2})$ and $(\bar C,\bar P_1,\ldots,\bar P_{2g+2})$ are stable marked curves. 

\medskip
For the following we fix a semistable model $(\bar\CC,\bar\CP_1,\ldots,\bar\CP_{2g+2})$ of $(\bar C,\bar P_1,\ldots,\bar P_{2g+2})$ over~$R$. We let $\CC$ be the minimal semistable model of $C$ that dominates~$\bar\CC$ and denote the morphism $\CC\onto\bar\CC$ again by~$\pi$. Since $\CC$ is proper over~$R$, each point $P_i$ extends to a unique section $\CP_i$ of~$\CC$.

\begin{Prop}\label{MinMarkModelSemiStab}
$(\CC,\CP_1,\ldots,\CP_{2g+2})$ is a semistable model of $(C,P_1,\ldots,P_{2g+2})$.
\end{Prop}

\begin{Proof}
By construction we have $\pi(\CP_i) = \bar\CP_i$ for all~$i$, and by assumption these sections are pairwise disjoint. Thus the sections $\CP_i$ are pairwise disjoint. 
Also, by assumption the sections $\bar\CP_i$ land in the smooth locus of~$\bar\CC$. Let $\CX$ be the normalization of $\bar\CC$ in the function field of~$C$. Then in Proposition \ref{SmoothMarkProp} below we will show that $\CX$ is smooth over a neighborhood of each~$\bar\CP_i$. Since $\CC$ is the minimal semistable model of $C$ that dominates~$\CX$, Proposition \ref{RelStabMod} (c) implies that $\CC\onto\CX$ is an isomorphism there. Thus each section $\CP_i$ lands in the smooth locus of~$\CC$, and we are done.
\end{Proof}

\medskip
Next let $\sigma$ denote the covering involution of $\pi\colon C\onto\bar C$. By the uniqueness of the minimal semistable model in Proposition \ref{RelStabMod} (b), this extends uniquely to an automorphism of $\CC$ of order~$2$. We denote this extension again by $\sigma$ and consider the quotient $\bbar\CC := \CC/\langle\sigma\rangle$. Since $\CC$ dominates $\bar\CC$, it follows that $\bbar\CC$ dominates~$\bar\CC$. Also $\sigma$ fixes each ramification point~$P_i$ and therefore each section~$\CP_i$. Let $\bbar\CP_i$ denote the section of $\bbar\CC$ that is induced by~$\CP_i$. Let $C_0$ and $\bbar C_0$ denote the closed fibers of $\CC$ and~$\bbar\CC$, respectively.

\begin{Prop}\label{BbarModelSemiStab}
\begin{enumerate}
\item[(a)] $(\bbar\CC,\bbar\CP_1,\ldots,\bbar\CP_{2g+2})$ is a semistable model of $(\bar C,\bar P_1,\ldots,\bar P_{2g+2})$.
\item[(b)] The inverse image of the smooth locus of $\bbar C_0$ is the smooth locus of~$C_0$.
\item[(c)] The inverse image of a double point of thickness $\alpha$ of $\bbar C_0$ is either a double point of thickness $\alpha/2$, or two double points of thickness $\alpha$ that are interchanged by~$\sigma$.
\end{enumerate}
\end{Prop}

\begin{Proof}
The quotient $\bbar\CC$ is semistable by Raynaud \cite[Appendice]{Raynaud1990} and the reduction to a discrete valuation ring in Section \ref{App}. As the morphism $\CC\onto\bbar\CC$ is surjective, the image of a smooth point of the special fiber cannot be a double point. In particular, since the marked sections $\CP_i$ land in the smooth locus of~$\CC$, their images land in the smooth locus of~$\bbar\CC$. By construction they are also disjoint, proving (a).

Next consider a closed point $p\in C_0$ with image $\bar p\in\bar C_0$. Then the inverse image of $\bar p$ is $\{p,\sigma(p)\}$. If $p$ is a smooth point, then so is $\bar p$ by Liu \cite[Prop. 3.48 (a)]{LiuAlgGeo2002}. If $p$ is a double point that is not fixed by~$\sigma$, the covering is \'etale at~$p$, so $\bar p$ is a double point of the same thickness as~$p$. If $p$ is a double point that is fixed by~$\sigma$, by Raynaud \cite[Prop. 2.3.2]{Raynaud1999} its image $\bar p$ is either a double point of twice the thickness as~$p$, or it is a smooth point and there exists a ramification point of the generic fiber which reduces to $p$. As in our situation all ramification points are marked points and reduce to smooth points by semistability, the last case cannot in fact occur, proving (b) and (c).
\end{Proof}

\begin{Prop}\label{BbarModelSemiStabIrrComp}
Let $\bbar T$ be an irreducible component of $\bbar C_0$ and let $T$ be its inverse image in~$C_0$. Then either
\begin{itemize}
\item[(a)] $T$ is isomorphic to $\BP^1_k$ and purely inseparable of degree $2$ over~$\bbar T$, or
\item[(b)] $T$ is irreducible and smooth and separable of degree $2$ over~$\bbar T$, or
\item[(c)] $T$ is isomorphic to $\BP^1_k \sqcup \BP^1_k$, each component mapping isomorphically to~$\bbar T$.
\end{itemize}
In particular the irreducible components of $C_0$ are smooth and have no self-intersections.
\end{Prop}

\begin{Proof}
By Proposition \ref{BbarModelSemiStab} the image of a double point $p$ of $C_0$ is a double point $\bar p$ of~$\bbar C_0$. Since $\bbar C_0$ has genus zero, this point $\bar p$ lies in two distinct irreducible components of $\bbar C_0$. The local surjectivity of $C_0\onto\bbar C_0$ thus implies the same for~$p$. This proves the last sentence of the proposition.

If $T$ is irreducible, this leaves only the possibilities (a) and (b). If $T$ is reducible, each of its irreducible components must map isomorphically to $T\cong\BP^1_k$. Also, these components must be interchanged by the hyperelliptic involution~$\sigma$. In the proof of Proposition \ref{BbarModelSemiStab} we have seen that this rules out that they intersect in a double point, leaving only the possibility (c).
\end{Proof}


\begin{Prop}\label{MinMarkModelStab}
If $(\bar\CC,\bar\CP_1,\ldots,\bar\CP_{2g+2})$ is stable, so is $(\CC,\CP_1,\ldots,\CP_{2g+2})$.
\end{Prop}

\begin{Proof}
Let $\smash{(\CC',\CP_1',\ldots,\CP_{2g+2}')}$ be the stable model of $(C,P_1,\ldots,P_{2g+2})$. By its uniqueness $\sigma$ extends uniquely to an automorphism of $\CC'$ of order~$2$, and the quotient $\bbar\CC{}' := \CC/\langle\sigma\rangle$ is a model of~$\bar C$. Let $\bbar\CP_i'$ denote the section of $\bbar\CC'$ that is induced by~$\CP_i$. The same argument as in the proof of Proposition \ref{BbarModelSemiStab} (a) then shows that  $(\bbar\CC',\bbar\CP_1',\ldots,\bbar\CP_{2g+2}')$ is a semistable model of $(\bar C,\bar P_1,\ldots,\bar P_{2g+2})$. Since $(\bar\CC,\bar\CP_1,\ldots,\bar\CP_{2g+2})$ is stable, the minimality in Proposition \ref{AbsRelStabMod} now implies that $\bbar\CC'$ dominates~$\bar\CC$. Thus $\CC'$ dominates~$\bar\CC$ and therefore also the normalization $\CX$ of $\bar\CC$ in the function field of~$C$. In the proof of Proposition \ref{MinMarkModelSemiStab} we have seen that $P_1,\ldots,P_n$ extend to pairwise disjoint sections of the smooth locus of~$\CX$. By Proposition \ref{AbsRelStabMod} it follows that $\CC'$ is the minimal semistable model of $C$ that dominates~$\CX$. Thus $(\CC',\CP_1',\ldots,\CP_{2g+2}') \cong (\CC,\CP_1,\ldots,\CP_{2g+2})$, and so the latter is stable, as desired.
\end{Proof}

\medskip
Our primary goal in this article is to compute the stable model of $(C,P_1,\ldots,P_{2g+2})$. Proposition \ref{MinMarkModelStab} turns this into a local problem over the stable model of $(\bar C,\bar P_1,\ldots,\bar P_{2g+2})$. After this reduction, the stability condition becomes irrelevant. The same method therefore solves the slightly more general problem of computing the minimal semistable model of $C$ that dominates~$\bar\CC$ for an arbitrary semistable model $(\bar\CC,\bar\CP_1,\ldots,\bar\CP_{2g+2})$ of $(\bar C,\bar P_1,\ldots,\bar P_{2g+2})$. Throughout the following we therefore only work with the semistable marked models introduced above.

\medskip
For reference we collect the schemes and sections we have introduced in the following diagram. Recall that we have natural morphisms $\CC\onto\bbar\CC\onto\bar\CC$ that are compatible with the given sections. We let $(C_0, p_1, \dots, p_{2g+2})$ and $(\bbar C_0,\bbar p_1, \dots,\bbar p_{2g+2})$ and $(\bar C_0,\bar p_1, \dots,\bar p_{2g+2})$ denote the special fibers of $(\CC,\CP_1,\dots,\CP_{2g+2})$ and $(\bbar\CC,\bbar\CP_1,\dots,\bbar\CP_{2g+2})$ and $(\bar\CC,\bar\CP_1,\dots,\bar\CP_{2g+2})$, respectively.
\UseTheoremCounterForNextEquation
\begin{equation}\label{AllCPDiagram}
\vcenter{\xymatrix{
\ C\  \ar@{^{ (}->}[r] \ar@{->>}[d]_-\pi & 
\ \CC\ \ar@{->>}[d]_-{} & 
\ C_0\ \ar@{_{ (}->}[l] \ar@{->>}[d]_-{} &&
\ P_i\ \ar@{^{ (}->}[r] \ar@{|->}[d] & 
\ \CP_i\ \ar@{|->}[d] & 
\ p_i\ \ar@{|->}[d] \ar@{_{ (}->}[l] \\
\ \bbar C\ \ar@{^{ (}->}[r] \ar@{=}[d] & 
\ \bbar\CC\ \ar@{->>}[d] & 
\ \bbar C_0\ \ar@{->>}[d] \ar@{_{ (}->}[l] &&
\ \bbar P_i\ \ar@{^{ (}->}[r] \ar@{=}[d] & 
\ \bbar \CP_i\ \ar@{|->}[d] & 
\ \bbar p_i\ \ar@{|->}[d] \ar@{_{ (}->}[l] \\
\ \bar C\ \ar@{^{ (}->}[r] \ar@{->>}[d] &
\ \bar\CC\ \ar@{->>}[d] & 
\ \bar C_0\ \ar@{_{ (}->}[l] \ar@{->>}[d] &&
\ \bar P_i\ \ar@{^{ (}->}[r] & 
\ \bar\CP_i\ & 
\  \bar p_i\ \ar@{_{ (}->}[l]\\
\ \Spec K\ \ar@{^{ (}->}[r] & 
\ \Spec R\  & 
\ \Spec k\ \ar@{_{ (}->}[l] \\}}
\end{equation}

\medskip
To describe the relation between the special fibers $\bbar C_0$ and~$\bar C_0$, we use the terminology concerning the type of an irreducible component of $\bbar C_0$ from Definition \ref{Notationabcd}. We also divide the double points of $\bar C_0$ into two classes. For this recall that $\bar C_0$ is marked with $2g+2$ distinct points in the smooth locus. As the complement of a double point consists of two connected components, this divides the $2g+2$ marked points into two groups.

\begin{Def}\label{EvenOddDef}
A double point $\bar p$ of $\bar C_0$ is called \emph{even} if each connected component of $\bar C_0 \setminus \{\bar p\}$ contains an even number of the points $\bar p_1,\ldots,\bar p_{2g+2}$. Otherwise, it is called \emph{odd}. 
\end{Def}


When the characteristic of the residue field $k$ of $R$ is not~$2$, in \cite{GehrungerPink2021} we have shown that $\bbar\CC=\bar\CC$ and have given an explicit construction of~$\CC$. In this case the inverse image of an odd double point is a double point of half the thickness, and the inverse image of an even double point consists of two double points of the same thickness. In particular, $C$ has good reduction if and only if $\bar C_0$ is smooth.

\medskip

For the remainder of this article we assume that $k$ has characteristic~$2$. As explained in the introduction, the situation is then much more complicated. 

\medskip
To motivate the content of the next two chapters, let us first consider a smooth point $\bar p\in \bar C_0$. Choose a neighborhood $\bar\CU \subset \bar\CC$ that is isomorphic to an open subscheme of $\Spec R[x]$. To determine the minimal semistable model of $C$ above~$\bar\CU$ we must first compute the normalization of $\bar\CU$ in the function field of~$C$. In the generic fiber this 
 can be described 
by an equation of the form $z^2=f(x)$ for a separable polynomial $f\in K[x]$ of degree $2g+1$ or $2g+2$. After rescaling $f$ and $z$ by $K^\times$ we can assume that $f$ has coefficients in $R$ and is nonzero modulo~$\Fm$. The normalization of $R[x]$ in the function field $K(x,z)$ can then be found by a substitution of the form $z=h+at$ with $h\in R[x]$ and nonzero $a\in R$ and a new variable~$t$. Here $f$ must be approximated in an optimal way by the square~$h^2$. 
Sections \ref{OptDecomp} through \ref{DoubCov}
deal with finding such $h$ and hence computing the normalization. 

Next consider a double point $\bar p\in\bar C_0$ and choose a neighborhood $\bar\CU \subset \bar\CC$ that is isomorphic to an open subscheme of $\Spec R[x,y]/(xy-a)$ for some nonzero $a\in\Fm$. Writing ${R[x,y]/(xy-a)} = {R[x,\tfrac{a}{x}]}$, we will have to find a similar optimal approximation involving Laurent polynomials in~$x$. For simplicity the next chapter therefore deals primarily with Laurent polynomials. 

Where the normalization is not semistable, we will have to construct $\bbar\CC$ as a blowup of~$\bar\CC$ whose normalization in the function field of $C$ is semistable. At a smooth point $\bar p\in\bar C_0$ this problem has been solved by Lehr and Matignon \cite{LehrMatignon2006}. In that article they assume that $\bar C_0$ is smooth everywhere (i.e., that the marked curve $(\bar C,\bar P_1,\ldots,\bar P_{2g+2})$ has good reduction), but their treatment actually applies locally to any smooth point. In particular we obtain an explicit description of all irreducible components of $\bbar C_0$ of type (c) or (d) above~$\bar p$. The polynomial computations for this are done in Section \ref{LaurDecomp}.

Above a double point $\bar p\in\bar C_0$ we may also have irreducible components of type (b). The Laurent polynomial computations required to find these are done in Sections \ref{ScalingDecomp} and \ref{SepDecomp}. After having identified  the irreducible components of type (b),
the remaining irreducible components of type (c) and (d) above $\bar p$ can be found as in \cite{LehrMatignon2006}.

The Laurent polynomial computations for all this are done in Chapter~\ref{ApproxBySquares}. The actual construction of $\bbar\CC$ and $\CC$ is carried out in the respective parts of Chapter~\ref{LocalStableModel}, divided according to the case of a smooth marked or unmarked point, respectively an odd or even double point. The resulting algorithm is presented comprehensively in Section \ref{Algo}.

\section{Approximating Laurent polynomials by squares}
\label{ApproxBySquares}

Recall from Section \ref{App} that we start with a complete discrete valuation ring  $R_1$ with quotient field~$K_1$. We fix an algebraic closure $K$ of $K_1$ and let $R$ denote the integral closure of $R_1$ in~$K$. Since $R_1$ is complete, the valuation on $K_1$ extends to a unique valuation with values in $\BQ$ on~$K$, whose associated valuation ring is~$R$. We let $\Fm$ denote the maximal ideal of~$R$, so that the residue field $k := R/\Fm$ is  algebraically closed. For any $a\in R$ we let $[a]$ denote the residue class in~$k$.

From now on we assume that $K$ has characteristic $0$ and $k$ has characteristic~$2$. We normalize the valuation $v$ on $K$ in such a way that $v(2)=1$. For every integer $n\ge1$ we fix an $n$-th root $2^{1/n}\in K$ in a compatible way, such that for all $n,m\ge1$ we have $(2^{1/mn})^m=2^{1/n}$. For any rational number $\alpha = m/n$ we then set $2^\alpha := (2^{1/n})^m$. This defines a group homomorphism $\BQ \to K^\times$, which by the normalization of $v$ satisfies $v(2^\alpha)=\alpha$.

For any Laurent polynomial $f=\sum_i a_i x^i  \in K[x^{\pm 1}]$ we set
\UseTheoremCounterForNextEquation
\begin{equation}\label{vDef}
v(f)\ :=\ \inf\,\{v(a_i) \,|\, i\in\BZ\}.
\end{equation}
This extends $v$ to a valuation on $K[x^{\pm 1}]$, which by the Gauss lemma satisfies the equation $v(fg)=v(f)+v(g)$ for all $f,g\in K[x^{\pm 1}]$.
The elements $f$ with $v(f)\ge0$ make up the subring $R[x^{\pm1}]$,
and for any such we let $[f]$ denote the residue class in $k[x^{\pm1}]$.

\subsection{Optimal decompositions}
\label{OptDecomp}

For the following sections, we fix a Laurent polynomial $f\in R[x^{\pm 1}]$ with $f \not \equiv 0 \mod\Fm$, in other words with $v(f)=0$. To this we associate the value
\UseTheoremCounterForNextEquation
\begin{equation}\label{wDef}
w(f)\ :=\ \sup \bigl\{ v(f-h^2) \bigm| h\in R[x^{\pm 1}] \bigr\}\ \in\ \BR\cup\{\infty\},
\end{equation}
which measures how well $f$ can be approximated by squares.

\begin{Rem}\label{wfInfty}
\rm This supremum can be $\infty$ without ever being attained. For example, suppose that $f=1+g$ with $v(g)>2$. Simple properties of the binomial coefficients then imply that $\smash{v\bigl(\binom{1/2}{n}g^n\bigr) \ge {n\cdot(v(g)-2)}}$ for all~$n$, so that the binomial series $\smash{\sum_{n\ge0}\binom{1/2}{n} g^n}$ converges coefficientwise to a square root of $f$ in $R[[x]]$. For any $m\ge0$ the partial sum $\smash{h_m := \sum_{n=0}^m\binom{1/2}{n} g^n}$ lies in $R[x]$ and satisfies $v(f-h_m^2) \ge (m+1)\cdot(v(g)-2)$, which goes to $\infty$ for $m\to\infty$. Therefore $w(f)=\infty$, and this supremum is never attained by some $v(f-h^2)$ unless $f$ is already a square.

\medskip
One can avoid this phenomenon by restricting the upper and lower degree of $h$ in (\ref{wDef}). However, for our purposes in Section \ref{DoubCov} and later the precise values of $v(f-h^2)$ and $w(f)$ are irrelevant if they exceed~$2$. We therefore define:
\end{Rem}

\begin{Def}\label{optDef}
A decomposition of the form $f=h^2+g$ with $g,h\in R[x^{\pm 1}]$ is called \emph{optimal} if it satisfies the condition
\UseTheoremCounterForNextEquation
$$v(g)=w(f)\quad\hbox{or}\quad v(g)>2.$$
\end{Def}

Observe that $v(g)>2$ implies that $w(f)>2$, because for any decomposition $f=h^2+g$ we have $w(f)\ge v(g)$ by (\ref{wDef}). Note also that, even if $w(f)\le2$, it is not a priori clear that this supremum is attained and an optimal decomposition exists. But we will prove this in the next section. In the rest of this section we discuss how to recognize an optimal decomposition.

%
%


\begin{Lem}\label{HalfLemma}
For any $h,\tilde h \in R[x^{\pm1}]$ we have
$$\min\{2,v(\tilde h^2-h^2)\}\ =\ 2\cdot\min\{1,v(\tilde h-h)\}.$$
\end{Lem}

\begin{Proof}
We must show that for any $0 \leq\alpha \leq 2$ we have
$$v(\tilde h^2-h^2) \geq \alpha \quad\iff\quad v(\tilde h-h)\geq \alpha/2.$$
For this note first that by assumption we have $v(2h)\geq 1 \geq \alpha/2$. Thus if $v(\tilde h-h) \geq \alpha/2$, it follows that ${v(\tilde h+h)}={v((\tilde h-h)+2h)}\geq \alpha/2$ as well. This implies that ${v(\tilde h^2-h^2)} = {v(\tilde h-h)+v(\tilde h+h)}\geq \alpha$, proving the implication ``$\Leftarrow$''.

Conversely, if ${v(\tilde h^2-h^2)} = {v(\tilde h-h)+(\tilde h+h)}\geq\alpha$, we must have $v(\tilde h-h)\geq \alpha/2$ or $v(\tilde h+h)\geq \alpha/2$. By the same reasoning as above, these inequalities are equivalent, proving the implication ``$\Rightarrow$''.
\end{Proof}

\begin{Proposition}\label{RecogOpt<2}
A decomposition $f=h^2+g$ with $g,h \in R[x^{\pm1}]$ and $v(g)<2$ is optimal if and only if the residue class $[g/2^{v(g)}]$ is not a square in $k[x^{\pm1}]$.
\end{Proposition}

\begin{Proof}
Abbreviate $\alpha := v(g)<2$. If the decomposition is not optimal, there exist $\tilde g,\tilde h \in R[x^{\pm1}]$ with $f=\tilde h^2+\tilde g$ and $v(\tilde g)>\alpha$. Then $v(g-\tilde g)=\alpha<2$, and $h^2+g=f=\tilde h^2+\tilde g$ implies that $\tilde h^2-h^2 = g-\tilde g$. By Lemma \ref{HalfLemma} it follows that $v(\tilde h-h)=\alpha/2$. Write $\tilde h = h+2^{\alpha/2}\ell$ with $\ell\in R[x^{\pm1}]$. Then we have $g-\tilde g = \tilde h^2-h^2 = 2^{1+\alpha/2}h\ell + 2^\alpha\ell^2$ and hence
$$g/2^\alpha-\tilde g/2^\alpha\ =\ 2^{1-\alpha/2}h\ell + \ell^2$$
within $R[x^{\pm1}]$. Here 
$v(\tilde g/2^\alpha) = v(\tilde g)-\alpha >0$ and $v(2^{1-\alpha/2}h\ell) \ge 1-\alpha/2 > 0$ by assumption. Thus the equality implies that $[g/2^\alpha] = [\ell^2]$ is a square in $k[x^{\pm1}]$.

Conversely, if $[g/2^\alpha]$ is a square in $k[x^{\pm1}]$, there exists $\ell\in R[x^{\pm1}]$ with $v(g/2^\alpha-\ell^2)>0$, or equivalently $v(g-2^\alpha\ell^2)>\alpha$. Setting $\tilde h := h+2^{\alpha/2}\ell$ we then deduce that
$$\tilde g\ :=\ f-\tilde h^2\ =\ g - 2^\alpha\ell^2- 2^{1+\alpha/2}h\ell.$$
Since $v(2^{1+\alpha/2}h\ell) \ge 1+\alpha/2 > \alpha$ by assumption, this implies that $v(\tilde g)>\alpha$. Thus $f=\tilde h^2+\tilde g$ is a better decomposition and the decomposition $f=h^2+g$ is not optimal.
\end{Proof}

\begin{Proposition}\label{RecogOpt=2}
A decomposition $f=h^2+g$ with $g,h \in R[x^{\pm1}]$ and $v(g)=2$ is optimal if and only if the equation $[g/4] = t^2+[h]t$ does not have a solution $t\in k[x^{\pm1}]$.
\end{Proposition}

\begin{Proof}
If the decomposition is not optimal, there exist $\tilde g,\tilde h \in R[x^{\pm1}]$ with $f=\tilde h^2+\tilde g$ and $v(\tilde g)>2$. Then $v(g-\tilde g)=2$, and $h^2+g=f=\tilde h^2+\tilde g$ implies that $\tilde h^2-h^2 = g-\tilde g$. By Lemma \ref{HalfLemma} this implies that $v(\tilde h-h)\ge1$. Write $\tilde h = h+2\ell$ with $\ell\in R[x^{\pm1}]$. Then $g-\tilde g = \tilde h^2-h^2 = 4\ell^2+4h\ell$ implies that 
$$g/4-\tilde g/4\ =\ \ell^2+h\ell$$
within $R[x^{\pm1}]$. Here 
$v(\tilde g/4) = v(\tilde g)-2 >0$ by assumption. Thus the equality implies that the equation $[g/4] = t^2+[h]t$ has the solution $t=[\ell]\in k[x^{\pm1}]$.

Conversely, suppose that the equation $[g/4] = t+t^2$ has a solution in $k[x^{\pm1}]$. Then there exists $\ell\in R[x^{\pm1}]$ with $[g/4] = [\ell^2]+[h][\ell]$, or equivalently $v(g-4\ell^2-4h\ell)>2$. Setting $\tilde h := h+2\ell$ we then deduce that
$$\tilde g\ :=\ f-\tilde h^2\ =\ g- 4\ell^2 - 4h\ell$$
satisfies $v(\tilde g)>2$. Thus $f=\tilde h^2+\tilde g$ is a better decomposition and the decomposition $f=h^2+g$ is not optimal.
\end{Proof}

\subsection{Odd decompositions}
\label{OddDecomp}

We call a Laurent polynomial \emph{even} if it possesses only monomials with even exponents, and \emph{odd} if it possesses only monomials with odd exponents. Any $g\in K[x^{\pm1}]$ can be written in a unique way as $g=g_e+g_o$ with $g_e$ even and $g_o$ odd. We call $g_e$ the \emph{even part} and $g_o$ the \emph{odd part of~$g$}.

\begin{Def}\label{OddDecompDef}
A decomposition $f=h^2+g$ with $g,h\in K[x^{\pm 1}]$ is called \emph{odd} if $g$ is odd. 
\end{Def}

\begin{Prop}\label{OddDecompR}
For any odd decomposition both $h$ and $g$ have coefficients in~$R$.
\end{Prop}

\begin{Proof}
By assumption we have $f=\sum_i b_ix^i$ with $b_i\in R$. Write $h=\sum_i c_ix^i$ with $c_i\in K$ and pick an index $i$ with $v(c_i)$ minimal. Then the coefficient of $x^{2i}$ in $g$ is zero, because $g$ is odd. Taking the coefficients of $x^{2i}$ in the equation $f=h^2+g$ thus yields the equation
$$b_{2i}\ =\ c_i^2+2\sum_{j>0} c_{i+j}c_{i-j}.$$ 
By the minimality of $v(c_i)$ and the fact that $v(2)>0$ this implies that $v(b_{2i})=v(c_i^2)$. As $v(b_{2i})\ge0$, it follows that $v(c_i)\ge0$. By minimality again this implies that all coefficients of $h$ lie in~$R$. By the equation $f=h^2+g$ the same then also follows for~$g$. 
\end{Proof}

\medskip
To prove that an odd decomposition exists, we follow Fiore \cite[Prop.\,7.3.9]{Fiore2018}, because his proof is more elegant than our original one and allows for better explicit computation.

\begin{Lem}\label{OddPolynomialDecompositionLemma}
For any even $p\in R[x^{\pm1}]$ there exists $q\in R[x^{\pm 1}]$ with $p(x)=q(x)q(-x)$.   
\end{Lem}
\begin{Proof}
Choose an integer $m$ such that $x^{2m}p(x)$ is a polynomial. Being even, we can write it over the algebraically closed field $K$ in the form
$$x^{2m}p(x)\ =\ c\cdot\prod_{\nu=1}^n (a_\nu-x^2)$$
with $c,a_\nu\in K$. Choose elements $d,b_\nu,i\in K$ with $d^2=c$ and $b_\nu^2=a_\nu$ and $i^2=-1$, and set 
$$q(x)\ :=\ (ix)^{-m} d\prod_{\nu=1}^n (b_\nu-x).$$ 
Then 
$$q(x)q(-x)\ =\ (ix)^{-m}(-ix)^{-m} d^2\prod_{\nu=1}^n (b_\nu-x)(b_\nu+x)
\ =\ x^{-2m} c\prod_{\nu=1}^n (a_\nu-x^2)\ =\ p(x).$$
Thus $q$ solves our problem in $K[x^{\pm1}]$. But since $v(q)=v(q(-x))$, the same formula implies that $2v(q) = v(q)+v(q(-x)) = v(p)\ge0$. Therefore $q\in R[x^{\pm1}]$ and we are done.
\end{Proof}

\begin{Prop} \label{OddDecompExists}
An odd decomposition of $f$ exists. 
\end{Prop}

\begin{Proof}
As $f_e$ is even, by Lemma \ref{OddPolynomialDecompositionLemma} there exists $q\in R[x^{\pm 1}]$ such that $f_e(x)=q(x)q(-x)$.
Writing $q = q_e+q_o$ we observe that $q(-x) = q_e-q_o$. Thus we get $f_e=q_e^2-q_o^2$. 
Next, we choose $i\in R$ with $i^2=-1$ and set $h:=q_e+iq_o\in R[x^{\pm 1}]$. 
Then
$$f-h^2\ =\ f_e+f_o-(q_e+iq_o)^2\ =\ (q_e^2-q_o^2) + f_o - (q_e^2+2iq_eq_o-q_o^2)
\ =\ f_o - 2iq_eq_o.$$
Here we note that, as $q_e$ is even and $q_o$ is odd, their product $q_e q_o$ is odd. Thus $g := f_o - 2iq_eq_o$ is odd, and we have found the odd decomposition $f=h^2+g$, as desired.
\end{Proof}


\begin{Prop} \label{OddDecompOptimal}
For any odd decomposition $f=h^2+g$ we have
$$\min\{2,w(f)\}\ =\ \min\{2,v(g)\}.$$
In particular the decomposition is optimal unless $v(g)=2<w(f)$.
\end{Prop}

\begin{Proof}
By the definition of $w(f)$ we always have $w(f)\ge v(g)$. Thus the equality holds if $v(g)\ge2$, and the decomposition is optimal if $v(g)>2$ by Definition \ref{optDef}. In the case $v(g)<2$ we observe that, since $g$ is odd, the residue class $[g/2^{v(g)}]$ is a nonzero element of $k[x^{\pm1}]$ that possesses only monomials with odd exponents. It is therefore not a square in $k[x^{\pm1}]$, and so the decomposition is optimal by Proposition \ref{RecogOpt<2}. Thus $w(f)=v(g)$, and again the equality follows. Together this also shows that the decomposition is optimal unless $v(g)=2<w(f)$.
\end{Proof}


\begin{Prop}\label{OptDecompExists}
An optimal decomposition of $f$ exists and can be computed effectively.
\end{Prop}

\begin{Proof}
The proof of Proposition \ref{OddDecompExists} gives an effective construction of an odd decomposition $f=h^2+g$. By Proposition \ref{OddDecompOptimal} this is optimal unless $v(g)=2$. In that case we can use Proposition \ref{RecogOpt=2} to effectively decide whether the decomposition is optimal, and if not, its proof yields an effective procedure to produce a better decomposition $f=\tilde h^2+\tilde g$ with $v(\tilde g)>2$. That decomposition is then optimal by Definition \ref{optDef}.
\end{Proof}


\begin{Prop}\label{OddDecomRange}
Consider any odd decomposition $f=h^2+g$. If $f$ possesses only monomials with exponents in an interval $[d_2,d_1]$, then $h$ possesses only monomials with exponents in the interval $[d_2/2,d_1/2]$.
\end{Prop}

\begin{Proof}
Write $f = \sum_ia_ix^i$ and $h = \sum_ib_ix^i$ with $a_i,b_i\in R$. Let $i$ be maximal such that $b_i\neq0$. Then the coefficient of $x^{2i}$ in $h^2$ is $b_i^2$, and the coefficient in $g$ vanishes because $g$ is odd. Thus we have $a_{2i}=b_i^2\neq0$. Therefore $2i\le d_1$ and hence $i\le d_1/2$. The analogous argument shows that the minimal index $j$ with $b_j\neq0$ satisfies $j\ge d_2/2$.
\end{Proof}


\begin{Prop}\label{OddDecomPol}
If $f\in R[x]$, an optimal decomposition $f=h^2+g$ with $h,g\in R[x]$ exists and can be computed effectively.
\end{Prop}

\begin{Proof}
Same as for Proposition \ref{OptDecompExists}, taking into account Proposition \ref{OddDecomRange} for $d_2=0$.
\end{Proof}


\begin{Prop} \label{OddDecompOptimal2}
If $[f]\in k$, then any odd decomposition $f=h^2+g$ is optimal.
\end{Prop}

\begin{Proof}
By Proposition \ref{OddDecompOptimal} the decomposition is optimal unless $v(g)=2$. In that case the residue class $[g/4]\in k[x^{\pm1}]$ is nonzero and odd. Moreover, since $v(g)>0$, the equation $f=h^2+g$ and the assumption $[f]\in k$ implies that $[h]\in k$ as well. Together this implies that the equation $[g/4] = t^2+[h]t$ cannot have a solution $t\in k[x^{\pm1}]$. By Proposition \ref{RecogOpt=2} the decomposition is therefore optimal in this case as well.
\end{Proof}


\begin{Ex}\label{OptDecompExample}
\rm Consider $f = 1+ax+bx^2$ with $a,b\in R$. Choose $c\in R$ with $c^2=b$ and set $h := 1+cx$. Then $g := f-h^2 = (a-2c)x$ is odd, yielding an odd decomposition $f=h^2+g$ with $v(g) = v(a-2c)$. By Proposition \ref{OddDecompOptimal} this is optimal unless $v(a-2c)=2$. In that case abbreviate $d := \frac{a-2c}{4} \in R^\times$. Then the equation $[g/4] = t^2+[h]t$ from Proposition  \ref{RecogOpt=2} boils down to the equation 
$$[d]x\ =\ t^2+[1+cx]t\ =\ t\cdot(t+1+[c]x)$$
with $[d]\neq0$ in~$k$. For any solution $t\in k[x^{\pm1}]$, unique factorization in the ring $k[x^{\pm1}]$ shows that both $t$ and $t+1+[c]x$ must be pure monomials in~$x$. As $k$ has characteristic~$2$, this is only possible if $[d]=[c]$ with $t=1$ or $t=[c]x$. Note that $[d]=[c]$ is equivalent to $v(a-6c) = v(4d-4c)>2$. Proposition \ref{RecogOpt=2} thus shows that the given odd decomposition is optimal except if $v(a-6c)>2=v(a-2c)$.

In that case, setting $\tilde h := 1-cx$ yields another odd decomposition $f=\tilde h^2+\tilde g$ with $\tilde g = (a+2c)x$ and $v(\tilde g) = v(a+2c) = v((a-6c)+8c) > 2$, which is therefore optimal. Alternatively, setting $\hat h := 1+\frac{a}{2}x$ yields a truncated power series decomposition $f=\hat h^2+\hat g$ as in Section \ref{LaurDecomp} with $\hat g = (c^2-\frac{a^2}{4})x^2$ and $v(\hat g) = v(c^2-\frac{a^2}{4}) = v\bigl((a-2c)(a+2c)/4\bigr) > 2$, which is again optimal.
\end{Ex}

\subsection{Double covers of the affine line}
\label{DoubCov}

Keeping $f$ as above, in this section we assume in addition that $f$ is neither a square nor divisible by the square of a non-unit in $K[x^{\pm 1}]$. This means that all zeros of $f$ in $K^\times$ are simple and there is at least one. It implies that the equation $z^2=f$ defines a quadratic field extension $K(x,z)$ of $K(x)$. We want to compute the normalization of $R[x^{\pm1}]$ in $K(x,z)$.

For this we fix an optimal decomposition $f=h^2+g$. We set $\gamma := \min\{2,w(f)\}$ and substitute $z=h+2^{\gamma/2}t$ with a new variable~$t$. Computing $h^2+2^{1+\gamma/2}ht+2^\gamma t^2 = z^2 = f = h^2+g$ and dividing by $2^\gamma$ then yields the equation
\UseTheoremCounterForNextEquation
\begin{equation}\label{tEq}
2^{1-\gamma/2}ht+t^2\ =\ g/2^\gamma.
\end{equation}
Here the term $2^{1-\gamma/2}h$ lies in $R[x^{\pm1}]$ because $\gamma\le2$, and the term $g/2^\gamma$ lies in $R[x^{\pm1}]$ because $v(g)\ge\gamma$ by Definition \ref{optDef}. The equation thus has coefficients in~$R$, showing that $t$ is integral over $R[x^{\pm1}]$. The following proposition is a variant of Lehr \cite[Prop.\,1]{Lehr2001ReductionOP}.

\begin{Prop}\label{LehrPropositionCovering}
The normalization of $A := R[x^{\pm1}]$ in $K(x,z)$ is flat over $R$ and isomorphic to
$$B\ :=\ R[x^{\pm1}][t]\!\bigm/\!\bigl(2^{1-\gamma/2}ht+t^2-g/2^\gamma\bigr).$$
\begin{enumerate}[(a)]
\item In the case $w(f)<2$ the equation (\ref{tEq}) modulo $\Fm$ has the form $t^2 = [g/2^\gamma]$. The curve $\Spec B/\Fm B$ is irreducible and smooth over $k$ outside finitely many points where $[\frac{dg}{dx}/2^\gamma]=0$, and the double covering $\Spec B/\Fm B \to \Spec A/\Fm A$ is purely inseparable. 
\item In the case $w(f)=2$ the equation (\ref{tEq}) modulo $\Fm$ has the form  $[h]t+t^2 = [g/4]$. The curve $\Spec B/\Fm B$ is irreducible and smooth over $k$ outside finitely many points where $[f]=0$, and the double covering $\Spec B/\Fm B \to \Spec A/\Fm A$ is separable. 
\item In the case $w(f)>2$ the equation (\ref{tEq}) modulo $\Fm$ has the form  $([h]+t)\cdot t = 0$. The curve $\Spec B/\Fm B$ is the union of two distinct rational curves and smooth over $k$ outside finitely many points where $[f]=0$, and each irreducible component maps isomorphically to $\Spec A/\Fm A$.
\end{enumerate}
\end{Prop}

\begin{Proof}
We first prove (a) through (c).

In the case $w(f)<2$ we have $\gamma = w(f) < 2$ and hence $[2^{1-\gamma/2}h]=0$. The equation (\ref{tEq}) modulo $\Fm$ therefore has the form $t^2 = [g/2^\gamma]$. By optimality we also have $\gamma=v(g)$, and by Proposition \ref{RecogOpt<2} the residue class $[g/2^\gamma]$ is not a square in $k[x^{\pm1}]$. Thus its derivative $[\frac{dg}{dx}/2^\gamma]$ is a nonzero element of $k[x^{\pm1}]$ with at most finitely zeros. Away from these, the curve $\Spec B/\Fm B$ is smooth over~$k$, and the morphism $\Spec B/\Fm B \to \Spec A/\Fm A$ is totally inseparable, proving (a).

In the case $w(f)\ge2$ we have $\gamma = 2$ and hence $[2^{1-\gamma/2}h]=[h]$. The equation (\ref{tEq}) modulo $\Fm$ therefore has the form $[h]t+t^2 = [g/4]$. By optimality we also have $v(g)\ge2$ and hence $[f]=[h^2+g]=[h]^2$. Since $[f]\neq0$, this has at most finitely many zeros. Away from these, the equation shows that the morphism $\Spec B/\Fm B \to \Spec A/\Fm A$ is \'etale. In particular the curve $\Spec B/\Fm B$ is smooth over $k$ outside the zeros of~$[f]$.

In the case $w(f)=2$ we know in addition from Proposition \ref{RecogOpt=2} that the equation modulo $\Fm$ is irreducible, proving (b).

In the case $w(f)>2$ by optimality we have $v(g)>2$ and hence $[g/4]=0$. The equation modulo $\Fm$ thus has the form $([h]+t)\cdot t = [h]t+t^2 = 0$. It follows that $\Spec B/\Fm B$ is the union of two distinct rational curves, each mapping isomorphically to $\Spec A/\Fm A$, proving (c).

To show the first statement consider any subfield $K_2\subset K$ with $R_2 := R\cap K_2$ such that $f_2\in R_2[x^{\pm1}]$. Then $B = B_2\otimes_{R_2}R$ with the noetherian ring
$$B_2\ :=\ R_2[x^{\pm1}][t]\!\bigm/\!\bigl(2^{1-\gamma/2}ht+t^2-g/2^\gamma\bigr).$$
As this is free of rank $2$ over $R_2[x^{\pm1}]$, it is flat over~$R_2$; hence $B$ is flat over~$R$. Also, since $f$ has only simple zeros in~$K$, the spectrum of
$$B_2\otimes_{R_2}{K_2}\ \cong\ K_2[x^{\pm1}][z]/(z^2-f)$$ 
is smooth over~$K_2$. Moreover, in each case we have seen that the closed fiber is generically smooth. Thus $\Spec B_2 \to \Spec R_2$ is smooth outside codimension~$2$, and so $B_2$ is regular in codimension~$1$. 
On the other hand, since $\Spec B_2 \to \Spec R_2[x^{\pm1}]$ is finite, all its fibers are Artin and hence (S2). As $\Spec R_2[x^{\pm1}]$ is (S2), it follows from EGA4 \cite[Cor.\,6.4.2]{EGA4} that $\Spec B_2$ is (S2) as well. By Serre's criterion \cite[Th.\,5.8.6]{EGA4} it is therefore normal; hence $B_2$ is normal. Finally, since $B$ is the union of the rings $B_2$ as $K_2$ varies, it is normal as well, and we are done.
\end{Proof}

\begin{Prop}\label{LehrPropositionCoveringPol}
Assume in addition that $f,h,g\in R[x]$. Then the normalization of $A := R[x]$ in $K(x,z)$ is flat over $R$ and isomorphic to
$$B\ :=\ R[x][t]\!\bigm/\!\bigl(2^{1-\gamma/2}ht+t^2-g/2^\gamma\bigr),$$
and the analogues of Proposition \ref{LehrPropositionCovering} (a) through (c) hold accordingly.
\end{Prop}

\begin{Proof}
Same as for Proposition \ref{LehrPropositionCovering}.
\end{Proof}

\subsection{Behavior under scaling}
\label{ScalingDecomp}

Let $f$ be as in Section \ref{OptDecomp}. Writing $f=\sum_ia_ix^i$ with $a_i\in R$, we now assume that the constant coefficient $a_0$ is a unit and we are given a positive number $\alpha\in\BQ$ satisfying 
\UseTheoremCounterForNextEquation
\begin{equation}\label{AlphaDef}
v(a_i)\ge\alpha|i|\ \hbox{for all}\ i<0. 
\end{equation}
Since $a_0$ is a unit, (\ref{AlphaDef}) means that all negative slopes of the Newton polygon of $f$ are $\le-\alpha$. The condition implies that for any $\lambda\in\BQ\cap[0,\alpha]$, the Laurent polynomial
$$f(2^\lambda u)\ =\ \sum_i a_i2^{\lambda i}u^i$$
in the variable $u$ again has coefficients in~$R$. Fix an odd decomposition $f=h^2+g$.

\begin{Proposition}\label{ScalOpt}
For every $\lambda\in\BQ\,\cap\, ]\kern1pt0, \alpha\kern1pt[$ the decomposition $f(2^\lambda u)=h(2^\lambda u)^2+g(2^\lambda u)$ is optimal.
\end{Proposition}

\begin{Proof}
The assumption $0<\lambda<\alpha$ implies that all coefficients except the constant coefficient of $f(2^\lambda u)$ lie in~$\Fm$. As the decomposition $f(2^\lambda u)=h(2^\lambda u)^2+g(2^\lambda u)$ is odd, it is therefore optimal by Proposition \ref{OddDecompOptimal2}.
\end{Proof}

\medskip
Now we consider the function
\UseTheoremCounterForNextEquation
\begin{equation}\label{WBarDef}
\wbar\colon\ \BQ\,\cap\, [0, \alpha] \longto\BR,\ \ \lambda\mapsto 
\wbar(\lambda) := \min \{ 2, w(f(2^\lambda u)) \}.
\end{equation}
Recall that a \emph{break point} of a continuous piecewise linear function is a point were the slope changes.

\begin{Proposition}\label{WBarProp}
\begin{enumerate}
\item[(a)] The function $\wbar$ is continuous and piecewise linear concave.
\item[(b)] If $\wbar$ has a horizontal segment, its value there is $2$ and $w(f(2^\lambda u))>2$ over its interior.
\item[(c)] Consider a point $\lambda\in\BQ\,\cap\, ]\kern1pt0, \alpha\kern1pt[$ and set $\gamma := v(g(2^\lambda u))$. Then $\lambda$ is a break point of~$\wbar$ if and only if either $\gamma=2$, or $\gamma<2$ and $[g(2^\lambda u)/2^\gamma]$ is not a monomial.
\end{enumerate}
\end{Proposition}

\begin{Proof}
By Definition \ref{optDef} and Proposition \ref{ScalOpt} we have 
$$\wbar(\lambda)\ =\ \min\bigl\{2,\, v(g(2^\lambda u))\bigr\},$$
which shows that the image of $\wbar$ is contained in $\BQ$. 
Write $g=\sum_i c_ix^i$ with $c_i\in R$. Then for any $\lambda$ we have $g(2^\lambda u) = \sum_i c_i2^{\lambda i}u^i$ and hence 
\UseTheoremCounterForNextEquation
\begin{equation}\label{WBarProp2}
v(g(2^\lambda u))\ =\ \min\bigl\{ v(c_i)+\lambda i \bigm| i\in\BZ \bigr\}.
\end{equation}
If $g\not=0$, this is the minimum of a non-empty finite collection of affine-linear functions. In particular this proves (a). 

Next, as $g$ is odd, the coefficient $c_0$ vanishes; hence all slopes in (\ref{WBarProp2}) are non-zero. Thus if $\wbar$ has a horizontal segment, its value there is $2$ and $v(g(2^\lambda u))>2$ over its interior. Since we always have $w(f(2^\lambda u)) \ge v(g(2^\lambda u))$, this proves (b).

Now consider $\lambda\in\BQ\,\cap\, ]\kern1pt0, \alpha\kern1pt[$ and set $\gamma := v(g(2^\lambda u))$. If $\gamma>2$, then $\lambda$ lies in the interior of a segment where $\wbar$ has the constant value~$2$; hence it is not a break point. If $\gamma=2$, the fact that all slopes of the function $v(g(2^\lambda u))$ are nonzero implies that the value decreases strictly on at least one side of~$\lambda$. On the other side the value either also decreases or keeps the constant value~$2$, and in both cases $\lambda$ is a break point of~$\wbar$. If $\gamma<2$, then $\lambda$ is a break point if and only if the minimum in (\ref{WBarProp2}) is attained for at least two distinct indices. This means precisely that $[g(2^\lambda u)/2^\gamma]$ is not a monomial. Together this proves (c).
\end{Proof}

\begin{Rem}\label{WBarBPNP}
\rm In the situation of Proposition \ref{WBarProp}, a direct computation based on (\ref{WBarProp2}) shows that $[g(2^\lambda u)/2^\gamma]$ fails to be a monomial if and only if $-\lambda$ is a slope of the Newton polygon of~$g(x)$.
\end{Rem}

%


\begin{Ex}\label{WExample}
\rm Here is a polynomial with an odd decomposition:
$$f\ =\ x^4 + 3 x^3 + 3 x^2 + 4 x + 1+8x^{-1}
\ =\ (x^2+x+1)^2+(x^3+2x+8x^{-1}).$$
The number $\alpha=3$ satisfies the condition \eqref{AlphaDef}, and the function $\wbar$ can be read off from Proposition \ref{ScalOpt}, yielding
$$\wbar(\lambda)\ =\ 
{\scriptstyle\left\{
{\textstyle\begin{array}{cl}
3 \lambda &\text{ if }\ 0\leq \lambda\leq \frac{1}{2}, \\[3pt]
\lambda+1 &\text{ if }\ \frac{1}{2}\leq \lambda\leq 1, \\[3pt]
3-\lambda &\text{ if }\ 1\leq \lambda\leq 3,
\end{array}}\right.}$$
as shown in Figure~\ref{FigExampleW}.
\begin{figure}[h] \centering \FigExampleW
\caption{Plot of $\wbar$ as in Example \ref{WExample}}	
\label{FigExampleW}
\end{figure}
\end{Ex}

\subsection{Separating positive and negative exponents}
\label{SepDecomp}

We keep $f$ as in Section \ref{ScalingDecomp}, so that $f\in R^\times +R[x]x+ \Fm[x^{-1}]x^{-1}$. We can sometimes simplify explicit computations with optimal decompositions by separating the positive and negative parts of~$f$.

\begin{Prop} \label{HenselforLaurent}
There exist $f_1\in R^\times +R[x]x$ and $f_2\in R^\times +\Fm[x^{-1}]x^{-1}$ with $f_1 f_2=f$. Moreover, such a factorization is unique up to multiplication by units. 
\end{Prop}

\begin{Proof}
Choose $n\ge0$ such that $\tilde f := x^nf$ is a polynomial. Then by assumption its residue class satisfies $[\tilde f] = x^n\cdot [f]$ with coprime factors $x^n,[f]\in k[x]$. By Hensel's Lemma there therefore exist polynomials $f_0, f_1\in R[x]$ with $f_0f_1=\tilde f$ and $[f_1]=[f]$ as well as $[f_0]=x^n$ and $\deg(f_0)=n$. The equality $[f_1]=[f]$ implies that $f_1$ lies in $R^\times +R[x]x$, and the conditions on $f_0$ imply that $f_2 := x^{-n}f_0$ lies in $R^\times +\Fm[x^{-1}]x^{-1}$. This yields the desired factorization $f_1f_2=f$.

Conversely, for any such factorization, setting $f_0 := x^nf_2$ yields a factorization $f_0f_1=\tilde f$ as in Hensel's lemma. As that is unique up to multiplication by units, the same follows for the factorization $f_1f_2=f$.
\end{Proof}

\medskip
In the following we fix a factorization $f=f_1f_2$ as in Proposition \ref{HenselforLaurent}. 

\begin{Lem}\label{wfwf1wf2Lem1}
For any decompositions $f_1=h_1^2+g_1$ and $f_2=h_2^2+g_2$, setting $h:=h_1h_2$ yields a decomposition $f=h^2+g$ with $v(g)\ge\min\{v(g_1), v(g_2)\}$.
\end{Lem}

\begin{Proof}
By construction we have $g = f_1f_2-h_1^2h_2^2 = f_1g_2+g_1h_2^2$. Since $v(f_1)$, $v(h_2)\ge0$, this implies that $v(g)\ge\min\{v(g_1), v(g_2)\}$, as desired.
\end{Proof}

\begin{Lem}\label{wfwf1wf2Lem2}
Given any decomposition $f=h^2+g$ there exist decompositions $f_1=h_1^2+g_1$ and $f_2=h_2^2+g_2$ with $v(g)\le\min\{v(g_1), v(g_2)\}$.
\end{Lem}

\begin{Proof}
This is trivial if $v(g)$ is zero, so let us assume that $v(g)>0$. Then $f$ is congruent to $h^2$ modulo~$\Fm$; hence $h$ also lies in $R^\times +R[x]x+ \Fm[x^{-1}]x^{-1}$. Using Proposition \ref{HenselforLaurent} we write $h=\tilde{h}_1 \tilde{h}_2 $ with  $\tilde h_1\in R^\times +R[x]x$ and $\tilde h_2\in R^\times +\Fm[x^{-1}]x^{-1}$. Next we choose units $c_\nu\in R^\times $ such that  $c_\nu^2\tilde h_\nu^2$ has the same constant coefficient as $f_\nu$ for $\nu=1,2$. Setting $h_\nu:=c_\nu \tilde{h}_\nu$ we then have $g_1 := f_1-h_1^2\in R[x]x$ and $g_2 := f_2-h_2^2\in \Fm[x^{-1}]x^{-1}$. With $c:=\left( c_1 c_2\right )^{-1}\in R^\times$ we also have $h=c h_1 h_2$.

We will show that $v(g_1),v(g_2)\ge v(g)$. First we observe that $\tilde g\ :=\ g_2+(1-c^2)+g_1$ has constant coefficient $1-c^2$ and its parts with positive, respectively negative exponents are $g_1$ and~$g_2$. We will compare $\tilde g$ with the Laurent polynomial 
\UseTheoremCounterForNextEquation
\begin{equation}\label{wfwf1wf2Lem2A}
g\ =\ f-h^2\ =\ f_1f_2-c^2h_1^2h_2^2
\ =\  f_1\cdot g_2+h_1^2 h_2^2\cdot (1-c^2)+h_2^2\cdot g_1.
\end{equation}
For this we write $g=\sum_ia_ix^i$ and $\tilde g=\sum_i\tilde a_ix^i$ as well as $f_1 = \sum_i b_{1,i}x^i$ and $h_1^2h_2^2 = \sum_ib_{2,i}x^i$ and $h_2^2 = \sum_i b_{3,i}x^i$ with all coefficients in~$R$. For any integers $i$ and $j$ we put
$$M_{i,j}\ := \ \Biggl\{\begin{array}{ll}
b_{1,i-j} &\text { if } j<0, \\
b_{2,i-j} &\text{ if }j=0, \\
b_{3,i-j} &\text{ if } j>0.
\end{array}$$ 
Then the equation (\ref{wfwf1wf2Lem2A}) means that $a_i = \sum_{j} M_{i,j} \cdot \tilde a_j$ for all $i$. Choose an integer $d\ge0$ such that $a_i=\tilde a_i=0$ for all $i$ not in the interval $[-d,d]$. Then the equation means that the vector $(a_i)_i$ is obtained by multiplying the vector $(\tilde a_i)_i$ with the matrix $M:= (M_{i,j})_{i,j}$, where the indices run from $-d$ to~$d$.

By definition this matrix has coefficients in~$R$. To determine its reduction modulo~$\Fm$ we note that by construction we have $f_1,h_1\in R^\times +R[x]x$ and $h_2\in R^\times +\Fm[x^{-1}]x^{-1}$ and hence $[f_1], [h_1]\in k^\times + k[x]x$ and $[h_2]\in k^\times$. Therefore $[f_1]$, $[h_1^2h_2^2]$, $[h_2^2]$ all lie in $k^\times + k[x]x$. Equivalently this means that $[b_{1,i-j}] = [b_{2,i-j}] = [b_{3,i-j}] = 0$ for all $i<j$ and that $[b_{1,0}]$, $[b_{2,0}]$, $[b_{3,0}]$ are all nonzero. Thus $M$ modulo~$\Fm$ is a lower triangular matrix with nonzero coefficients on the diagonal. This implies that $\det(M)\in R$ is nonzero modulo~$\Fm$. It is therefore a unit, and so the matrix $M$ is invertible over~$R$.

Finally, recall that $\min\{v(a_i)\,| -d\le i\le d\} = v(g)$. As $M$ is invertible over~$R$, it follows that $v(\tilde g) = \min\{v(\tilde a_i)\,| -d\le i\le d\} = v(g)$ as well. Since the positive and negative parts of $\tilde g$ are just $g_1$ and $g_2$, it follows that $v(g_1),v(g_2)\ge v(g)$, finishing the proof.
\end{Proof}

\begin{Prop}\label{WisTheMinimumOfThePolynomialWs}
\begin{enumerate}
\item[(a)] We have  $w(f)=\min\{w(f_1), w(f_2)\}$.
\item[(b)] For any optimal decompositions $f_1=h_1^2+g_1$ and $f_2=h_2^2+g_2$, setting $h:=h_1h_2$ yields an optimal decomposition $f=h^2+g$.
\end{enumerate}
\end{Prop} 

\begin{Proof}
Part (a) follows directly from Lemmas \ref{wfwf1wf2Lem1} and \ref{wfwf1wf2Lem2} and the definition (\ref{wDef}). In the situation of (b), Lemma \ref{wfwf1wf2Lem1} says that $v(g)\ge\min\{v(g_1), v(g_2)\}$. If this minimum is $>2$, the decomposition $f=h^2+g$ is directly optimal by (\ref{wDef}). Otherwise by optimality the smaller of the values $v(g_\nu)$ satisfies $w(f_\nu) = v(g_\nu) \le v(g_{3-\nu}) \le w(f_{3-\nu})$, which by (a) implies that $v(g) \ge \min\{v(g_1), v(g_2)\} = \min\{w(f_1), w(f_2)\} = w(f)$. Thus again $f=h^2+g$ is optimal by (\ref{wDef}).
\end{Proof}

\medskip
Now observe that the Newton polygon of $f$ is the concatenation of the Newton polygons of $f_1$ and~$f_2$. Since $f$ satisfies the condition (\ref{AlphaDef}), it follows that $f_1$ and $f_2$ do so as well. 

\begin{Prop}\label{ScalOptDecomp}
Take any odd decompositions $f_1=h_1^2+g_1$ and $f_2=h_2^2+g_2$ and consider the decomposition $f=h^2+g$ with $h:=h_1h_2$. Then for every $\lambda\in\BQ\,\cap\, ]\kern1pt0, \alpha\kern1pt[$ the decomposition $f(2^\lambda u)=h(2^\lambda u)^2+g(2^\lambda u)$ is optimal, and Proposition \ref{WBarProp} (c) holds for this~$g$ as well.
\end{Prop}

\begin{Proof}
For $\nu=1,2$, applying Proposition \ref{ScalOpt} to $f_\nu$ in place of $f$ shows that the decomposition $f_\nu(2^\lambda u)=h_\nu(2^\lambda u)^2+g_\nu(2^\lambda u)$ is optimal for every $\lambda\in\BQ\,\cap\, ]\kern1pt0, \alpha\kern1pt[$. By Proposition \ref{WisTheMinimumOfThePolynomialWs} (b) it follows that the decomposition $f(2^\lambda u)=h(2^\lambda u)^2+g(2^\lambda u)$ is optimal.
In particular we have 
$$\wbar(f)\ =\ \min\bigl\{2,\, v(g(2^\lambda u))\bigr\}.$$
In the case $\gamma\neq2$ the statement of Proposition \ref{WBarProp} (c) follows exactly as in the proof given there. In the case $\gamma=2$ consider an odd decomposition $f=\tilde h^2+\tilde g^2$. Then the fact that both decompositions $f(2^\lambda u)=h(2^\lambda u)^2+g(2^\lambda u)$ and $f(2^\lambda u)=\tilde h(2^\lambda u)^2+\tilde g(2^\lambda u)$ are optimal implies that $v(\tilde g(2^\lambda u))=2$ as well. Thus $\lambda$ is a break point of $\wbar$ by Proposition \ref{WBarProp} (c), and we are done.
\end{Proof}

\subsection{Truncated power series decompositions}
\label{LaurDecomp}

There is another construction that sometimes yields optimal decompositions, which has been used by Lehr and Matignon \cite{LehrMatignon2006}. To explain this we first consider an arbitrary nonzero polynomial $f\in K[x]$ of degree~$d$. We are interested in square roots of $f$ modulo terms of degree $>d/2$.

\begin{Def}\label{LaurDecDef}
A \emph{truncated power series decomposition} of $f$ is a decomposition of the form $f=h^2+g$ with $h,g\in K[x]$, such that 
$h$ possesses only monomials of degrees $\le d/2$ and
$g$ possesses only monomials of degrees $>d/2$.
\end{Def}

\begin{Prop}\label{LaurDecExists}
If $f$ has nonzero constant term, a truncated power series decomposition of $f$ exists and is unique up to $h\mapsto \pm h$. 
\end{Prop}

\begin{Proof}
Let $a\in K^\times$ be the constant coefficient of~$f$. Then $a^{-1}f$ has constant coefficient~$1$; hence it has a unique formal square root in $K[[x]]$ that also has constant coefficient~$1$. The desired $h$ must be a square root of~$a$ times the truncation of this power series modulo terms of degree $>{d/2}$, and this determines~$g$.
\end{Proof}

\begin{Prop}\label{LaurDecOpt}
Assume that $f$  lies in $R[x]$ and has unit constant term. If $w(f)\ge2$, then any truncated power series decomposition of $f$ has coefficients in $R$ and is optimal.
\end{Prop}

\begin{Proof}
We claim that under the given assumptions, for every integer $1\le i\le \floor{d/2}+1$ there exists an optimal decomposition $f=h^2+g$, 
such that $h$ possesses only monomials of degrees $\le d/2$ and $g$ possesses only monomials of degrees $\ge i$.
In the case $i=\floor{d/2}+1$ this is a truncated power series decomposition, and by uniqueness up to sign it follows that every truncated power series decomposition is optimal.

To prove the claim for $i=1$ we begin with an odd decomposition $f=h^2+g$. Here Proposition \ref{OddDecomRange} shows that $h$ is a polynomial of degree $\le{d/2}$, so we are done if the decomposition is optimal. Otherwise by Proposition \ref{OddDecompOptimal} we have $v(g)=2<w(f)$ and by Proposition \ref{RecogOpt=2} the equation $[g/4] = t^2+[h]t$ has a solution $t\in k[x^{\pm1}]$. Since $[g/4]$ lies in $k[x]$ and has degree $\le d$, this solution $t$ must also lie in $k[x]$ and have degree $\le{d/2}$. Moreover, since $g$ is odd, the polynomial $[g/4]$ is divisible by~$x$, and so after possibly replacing $t$ by $t+[h]$ we can assume that $t$ is also divisible by~$x$. Lifting $t$ coefficientwise to $R$ yields a polynomial $\ell\in R[x]$ of degree $\le{d/2}$ and without constant term, which satisfies $v(g-4\ell^2-4h\ell)>2$. Setting $\tilde h := h+2\ell$ and $\tilde g := f-\tilde h^2$ we obtain a decomposition such that $\tilde h$ has degree $\le{d/2}$ and $\tilde g$ has no constant term and satisfies $v(\tilde g)>2$. This has all the desired properties for $i=1$.

Now take a decomposition $f=h^2+g$ satisfying the claim for some integer $1\le i\le \floor{d/2}$. Since $f$ has unit constant term and $v(g)\ge2$, the constant term $b$ of $h$ is also a unit. Let $c$ be the coefficient of $x^i$ in~$g$ and set $\tilde h := h + \tfrac{c}{2b}x^i$. Since $v(c)\ge v(g)\ge2$, this is another polynomial in $R[x]$ of degree $\le\floor{d/2}$. By construction the polynomial
$$\begin{array}{rl}
\tilde g\ :=\ f-\tilde h^2 & =\ h^2+g - h^2 - \tfrac{c}{b}hx^i - ( \tfrac{c}{2b})^2x^{2i} \\[3pt]
&=\ (g-cx^i) - (h-b)\tfrac{c}{b}x^i - ( \tfrac{c}{2b})^2x^{2i}
\end{array}$$
possesses only monomials of degrees $>i$. Moreover, the fact that $v(c)\ge v(g)\ge2$ and $v(b)=0$ implies that $v((\tfrac{c}{2b})^2) =2v(c)-2 \ge v(c) \ge v(g)$ as well. Thus we have $v(\tilde g)\ge v(g)$, and the decomposition $f=\tilde h^2+\tilde g^2$ satisfies the claim for $i+1$ in place of~$i$.
By induction on $i$ the claim thus follows for all~$i$, and we are done.
\end{Proof}

\medskip
Lehr and Matignon \cite[\S3]{LehrMatignon2006} have combined the truncated power series decompositions for all polynomials obtained from $f$ by linear substitutions into a single object that they call a \emph{special decomposition}. Let $K(x_0)$ denote the field of rational functions in a new variable~$x_0$ and consider another new variable~$y$.

\begin{Prop}\label{ExpCons}
There exist unique polynomials in $y$ of the form
$$H_f(x_0,y)\ =\ \!\! \sum_{0\le k\le \frac{d}{2}}\!\! H_{f,k}(x_0)\, y^k
\qquad\hbox{and}\qquad
G_f(x_0,y)\ =\ \!\! \sum_{\frac{d}{2}< k\le d}\!\! G_{f,k}(x_0)\, y^k$$ 
with coefficients $H_{f,k}(x_0)$, $G_{f,k}(x_0) \in K(x_0)$, such that $H_{f,0}(x_0)=1$ and 
$$\frac{f(x_0+y)}{f(x_0)}\ =\ H_f(x_0,y)^2 + G_f(x_0,y).$$
These satisfy $f(x_0)^k H_{f,k}(x_0) \in K[x_0]$ and $f(x_0)^k G_{f,k}(x_0) \in K[x_0]$ for all relevant~$k$.
\end{Prop}

\begin{Proof}
The Taylor expansion of $f(x_0+y)$ yields the formula
\UseTheoremCounterForNextEquation
\begin{equation}\label{ExpCons1}
\frac{f(x_0+y)}{f(x_0)}\ =\ 1+\sum_{n=1}^d \frac{f^{(n)}(x_0)}{f(x_0)}\cdot\frac{y^n}{n!}.\end{equation}
Using the general binomial series we write its formal square root in $K(x_0)[[y]]$ in the form
\UseTheoremCounterForNextEquation
\begin{equation}\label{ExpCons2}
\sqrt{\frac{f(x_0+y)}{f(x_0)}}\ \ =\ \ 
\sum_{m=0}^\infty \binom{\frac{1}{2}}{m}\cdot\left[\sum_{n=1}^d \frac{f^{(n)}(x_0)}{f(x_0)}\cdot\frac{y^n}{n!}\right]^m
\ =\ \ \sum_{k=0}^\infty H_{f,k}(x_0)\, y^k
\end{equation}
with $H_{f,k}(x_0)\in K(x_0)$ and $H_{f,0}(x_0)=1$. Its truncation $H_f(x_0,y) := \smash{\sum_{0\le k\le d/2} H_{f,k}(x_0) y^k}$ is the unique square root with constant term $1$ modulo terms of exponent $>d/2$. Since $H_f(x_0)^2$ has degree $\le d=\deg(f)$ in~$y$, the difference $G_f(x_0,y) := f(x_0+y)/f(x_0)-H_f(x_0,y)^2$ has the desired form. This proves the existence and uniqueness of the decomposition. From (\ref{ExpCons2}) we can also see that $f(x_0)^k H_{f,k}(x_0) \in K[x_0]$ for all~$k$. This together with (\ref{ExpCons1}) implies that $f(x_0)^k G_{f,k}(x_0) \in K[x_0]$ for all $k$ as well.
\end{Proof}


\begin{Def}\label{StabPolDef}
The \emph{stability polynomial} associated to $f$ is 
$$S_f(x_0)\ :=\ f(x_0)^{2^m} G_{f,2^m}(x_0)\ \in\ K[x_0],$$
where $m$ is the unique integer such that $\frac{d}{2} < 2^m\le d$.
\end{Def}

Lehr and Matignon \cite[Def\,3.4]{LehrMatignon2006} call this---up to a constant factor---the \emph{monodromy polynomial} of~$f$, because they are interested in determining the Galois group of the field extension over which the stable model of $C$ is defined. We prefer a name that refers more directly to the stable model, for whose construction the zeros of the stability polynomial 
play a similar role as a level structure does.

\begin{Prop}\label{StabPolLinSub}
For any $a,c\in K^\times$ and $b\in K$ we have
$$\begin{array}{rl}
G_{cf(ax+b)}(x_0,y) &=\ G_f(ax_0+b,ay) \qquad\hbox{and} \\[3pt]
S_{cf(ax+b)}(x_0) &=\ (ac)^{2^m}\cdot S_f(ax_0+b).
\end{array}$$
\end{Prop}

\begin{Proof}
Replace $(x_0,y)$ by $(ax_0+b,ay)$ in the defining formula $\frac{f(x_0+y)}{f(x_0)} = H_f(x_0,y)^2 + G_f(x_0,y)$ and compute.
\end{Proof}

\medskip

Now assume that $f\in K[x]$ is neither a square nor divisible by the square of a non-unit in $K[x]$. Further, assume that $w(f)<2$. As in Section \ref{DoubCov} we are interested in the quadratic field extension $K(x,z)$ of $K(x)$ that is defined by the equation $z^2=f$. We also assume that $f\in R[x] \setminus \Fm[x]$ with $f \not \equiv 0 \mod\Fm$. Consider any $\xi_0\in R$ with $f(\xi_0)\in R^\times$, and take a new variable~$y$.

\begin{Lem}\label{StabPolXi0Lem}
There exists a unique $\epsilon\in\BQ^{>0}$ such that $v(G_f(\xi_0,2^\epsilon y))=2$.
\end{Lem}

\begin{Proof}
By assumption $f(\xi_0+y)/f(\xi_0) = H_f(\xi_0,y)^2 + G_f(\xi_0,y)$ is not a square in $K(y)$; hence there exists $d/2<k\le d$ such that $G_{f,k}(\xi_0)\neq0$. Thus 
$$\BQ\longto\BQ,\ \epsilon \longmapsto\ v(G_f(\xi_0,2^\epsilon y))\ =\ 
\min \bigl\{ v(G_{f,k}(\xi_0)) +\epsilon k \bigm| \tfrac{d}{2}< k\le d \bigr\}$$
is a piecewise linear function with finitely many strictly positive slopes, which therefore takes the value $2$ at a unique $\epsilon\in\BQ$.
Assume that $\epsilon \leq 0$, which is equivalent to $v(G_f(\xi_0,y))\geq 2$. Note that by assumption, we have $f(\xi_0)\in R^\times$. Hence  $\frac{f(\xi_0+y)}{f(\xi_0)} \in R[y]$ and thus $H_f(\xi_0,y)$ is in $R[y]$ as well because its square can be written as the sum of two elements of $R[y]$. Let $c$ be a square root of $f(\xi_0)$. Then 
$$f(\xi_0+y)=(cH_f(\xi_0,y))^2+c^2G_f(\xi_0,y)$$ 
is a decomposition with $v(c^2G_f(\xi_0,y))\geq 2$. This contradicts $w(f)<2$.
\end{Proof}

\begin{Thm}\label{LMThm}
For any $\epsilon\in\BQ^{>0}$, set $\bar\CX:=\BP^1_R$ with the coordinate $y$ after the substitution $x=\xi_0+2^\epsilon y$, and let $\CX$ denote the normalization of $\bar\CX$ in $K(x,z)$. 
Then the following are equivalent:
\begin{enumerate}
\item[(a)] The closed fiber of $\CX$ has geometric genus $>0$.
\item[(b)] We have $v(G_f(\xi_0,2^\epsilon y))=2$ and there exists $\xi_1\in R$ with $v(\xi_1-\xi_0)\ge\epsilon$ and $S_f(\xi_1)=0$, and the substitution $x=\xi_1+2^\epsilon y$ yields the same models $\bar\CX$ and~$\CX$.
\end{enumerate}
\end{Thm}

\begin{Proof}
Lehr and Matignon proved this in \cite[Thm\,5.1]{LehrMatignon2006} under special assumptions on~$f$. But their proof carries over to our situation with no essential changes. Specifically, the first assertion of \cite[Lemma 3.3 (iii)]{LehrMatignon2006} is not needed in the proof, so the proof works regardless of the parity of~$d$.
\end{Proof}


\begin{Prop}\label{StabPolNonZero}
Under the above assumptions we have $S_f\not=0$.
\end{Prop}

\begin{Proof}
By Proposition \ref{StabPolLinSub} it suffices to prove this after an arbitrary linear substitution $x\rightsquigarrow bx+c$. Since $d>0$, we may thus assume that $f(0)=0$ and that all other zeros of $f$ have valuation $\ge0$. After multiplying $f$ by a constant we may then assume that $f\in R[x]$ and that the reduction $[f]$ has a simple zero at~$0$.

Now consider an arbitrary $\xi_0\in R$ with $f(\xi_0)\in R^\times$, and let $\epsilon\in\BQ$ be as in Lemma \ref{StabPolXi0Lem}. Then we have
$$\frac{f(\xi_0+2^\epsilon y)}{f(x_0)}\ =\ H_f(\xi_0,2^\epsilon y)^2 + G_f(\xi_0,2^\epsilon y)$$
in $K[y]$, where the last term on the right hand side lies in $4R[y]$. Suppose that $\epsilon\le0$. Then applying the reverse substitution $y=2^{-\epsilon}(x-\xi_0) \in R[x]$ we deduce that 
$$\frac{f(x)}{f(x_0)}\ =\ H_f(\xi_0,x-\xi_0)^2 + G_f(\xi_0,x-\xi_0)$$
in $K[x]$, where the last term on the right hand side lies in $4R[x]$. As the left hand side lies in $R[x]$ by assumption, it follows 
that $f(x)$ is a square modulo $4R[x]$. As this contradicts the assumptions on~$f$, we deduce that $\epsilon>0$. 

Now set $\bar\CX_{\xi_0}:=\BP^1_R$ with the coordinate $y$ after the substitution $x=\xi_0+2^\epsilon y$, and let $\CX_{\xi_0}$ denote the normalization of $\bar\CX_{\xi_0}$ in $K(x,z)$. If $S_f=0$, the condition in Theorem \ref{LMThm} is always satisfied with $\xi_1=\xi_0$. Thus for any choice of $\xi_0$ as above, the closed fiber of $\CX_{\xi_0}$ has geometric genus $>0$. 

Now observe that two pairs $(\xi_0,\epsilon)$ and $(\xi'_0,\epsilon')$ yield the same model $\CX_{\xi_0}$ if and only if $v(\xi'_0-\xi_0)>\epsilon$, and in that case we have $\epsilon=\epsilon'$. Since all $\epsilon>0$, we can thus choose arbitrarily many inequivalent pairs $(\xi_0,\epsilon)$ which yield a curve of genus $>0$. For any finite number of them, consider the minimal semistable model $\CC$ of $C$ that dominates each $\CX_{\xi_0}$ from Proposition \ref{VarRelStabModExists}. Then the closed fiber $C_0$ of $\CC$ must contain an irreducible component mapping isomorphically to the closed fiber of any one of the $\CX_{\xi_0}$ and to a point in all others. But this is not possible for arbitrarily many $\xi_0$, because the arithmetic genus of $C_0$ is the genus of~$C$, which is fixed. We have obtained a contradiction, proving that $S_f\not=0$.
\end{Proof}


\begin{Rem}\label{StabPolDegRem}
\rm If $d$ is odd, as in Lehr-Matignon \cite[Lemma 3.3 (iii)]{LehrMatignon2006} one can show that $S_f$ has degree $2^m(d-1)$.
If $d$ is even, one can show that $S_f$ has degree $<2^m(d-1)$.
\end{Rem}

\begin{Ex}\label{StabPolEx}
\rm In the case $d=1$ we have $H_f=1$ and hence $G_f=\frac{f'(x_0)}{f(x_0)}\cdot y$. Since $2^m=1$ we get $S_f=f'(x_0)$ in this case. Concretely for $f=a+bx$ we obtain $S_f=b$. 
In the case $2\le d\le3$ we have $H_f=1+\frac{f'(x_0)}{f(x_0)}\cdot \frac{y}{2}$ and hence 
$$G_f\ =\ \frac{f''(x_0)}{f(x_0)}\cdot \frac{y^2}{2} \allowbreak + \frac{f'''(x_0)}{f(x_0)}\cdot \frac{y^3}{6} - \Bigl(\frac{f'(x_0)}{f(x_0)}\cdot \frac{y}{2}\Bigr)^2.$$ 
Since $2^m=2$ we get $S_f = {(2f''(x_0)f(x_0)} - {f'(x_0)^2})/4$. 
Concretely for $f=a+bx+cx^2$ we obtain $S_f=-(b^2-4ac)/4$, and for $f=a+bx+cx^2+dx^3$ we obtain 
$$S_f\ =\ \frac{3d^2x_0^4 + 4cdx_0^3 + 6bdx_0^2 + 12adx_0 + (4ac-b^2)}{4}.$$
Note that for $1\le d\le2$ this shows that $S_f\in K^\times$, so that $S_f$ has no zero. In this case the function field $K(x,z)$ with $z^2=f$ corresponds to a rational curve, whose reduction cannot have an irreducible component of genus~$>0$.
For more details in the case $d=3$ see Section \ref{Genus1}.
\end{Ex}

\section{Local constructions}
\label{LocalStableModel}

We return to the situation and notation of Section \ref{HyperellCurves}. We fix a closed point $\bar p$ of~$\bar C_0$ and identify a neighborhood $\bar\CU \subset \bar\CC$ with an open subscheme of $\Spec R[x]$ or $\Spec R[x,y]/(xy-a)$ 
for some nonzero $a\in\Fm$. To determine the minimal semistable model of $C$ above~$\bar\CU$ we must first compute the normalization of $\bar\CU$ in the function field of~$C$. In the generic fiber this is given by an equation of the form $z^2=f(x)$ for a separable polynomial $f\in K[x]$ of degree $2g+1$ or $2g+2$. After rescaling $f$ and $z$ by $K^\times$ we can assume that $f$ has coefficients in $R$ and is nonzero modulo~$\Fm$. Recall that the zeros of $f$ in the generic fiber $\bar C\cong\BP^1_K$ are the marked points $\bar P_i$ that are different from $\infty$.


\subsection{Smooth marked points}
\label{SmoothMarkedPoints}

In this section we assume that $\bar p$ is the reduction of a section~$\bar\CP_i$. By semistability it is thus a smooth point of~$\bar C_0$. The following result holds for arbitrary residue characteristic.

\begin{Prop}\label{SmoothMarkProp}
There exists a neighborhood $\bar\CU \subset \bar\CC$ of~$\bar p$, whose normalization in the function field of $C$ is smooth over~$R$ and equal to $\CC$ over~$\bar\CU$.
\end{Prop}

\begin{Proof}
Choose a neighborhood $\bar\CU$ and an embedding $\bar\CU\into\Spec R[x]$ such that $\bar\CP_i$ is given by $x=0$. As the sections $\bar\CP_j$ are all disjoint, we then have $f(x) = xg(x)$ for a polynomial $g\in R[x]$ with $g(0)\in R^\times$. Consider the ring $B := R[x,z]/(z^2-f(x))$. This is flat and integral over $A := R[x]$ and thus contained in the normalization of $A$ in the function field of~$C$. The inverse image of the section $\bar\CP_i$ under $\Spec B\onto\Spec A$ is given by $x=z=0$. Along this section we have
$$d(z^2-f(x))\ =\ 2z\,dz-f'(x)\,dx\ =\ 2z\,dz-xg'(x)\,dx-g(x)\,dx
\ \equiv\ -g(0)dx.$$
As $g(0)$ is a unit, the jacobian criterion implies that $\Spec B$ is smooth over $R$ along this section. In particular, an open neighborhood of this section is normal and hence coincides with the normalization of $\Spec A$ in the function field of~$C$. This proves the first statement of the proposition, and the second follows from Proposition \ref{RelStabMod} (c).
\end{Proof}

\begin{Rem}\label{SmoothMarkRem}
\rm The description of $f$ in the above proof shows that  the residue class  $[f]\in k[x]$ is not a square; so by Proposition \ref{RecogOpt<2} we have $w(f)=0$. We therefore have the case (a) of Proposition \ref{LehrPropositionCoveringPol}; hence the irreducible component of $C_0$ that dominates the closed fiber $\bar U_0$ of $\bar\CU$ is inseparable over $\bar U_0$ and has genus zero.
\end{Rem}

\subsection{Smooth unmarked points}
\label{SmoothUnmarkedPoints}

In this section we assume that $\bar p$ is a smooth point of~$\bar C_0$ that does not lie in any of the sections~$\bar\CP_i$. We identify a neighborhood $\bar\CU \subset \bar\CC$ with an open subscheme of $\Spec R[x]$ and suppose that the function field of $C$ is $K(x,z)$ with $z^2=f(x)$ for a separable polynomial $f\in R[x]$ with $v(f)=0$. Then by assumption $\bar p$ does not lie in the zero locus of~$f$.

With Proposition \ref{LehrPropositionCoveringPol} we can compute the normalization $\CU$ of $\bar\CU$ in $K(x,z)$. In particular we can determine where $\CU$ is smooth, and by Proposition \ref{RelStabMod} (c) it is equal to $\CC$ there. 
So assume that $\CU$ is singular above~$\bar p$. Proposition \ref{LehrPropositionCoveringPol} then shows that $w(f)<2$ and the closed fiber of $\CU$ is inseparable over the closed fiber of~$\bar\CU$, and that $\bbar C_0$ possesses at least one irreducible component of type (d) above~$\bar p$ in the terminology of Definition \ref{Notationabcd}.

\begin{Prop}\label{RaynaudProp}
Let $\bbar T$ be an irreducible component of $\bbar C_0$ above~$\bar p$. Then its inverse image $T$ in $C_0$ is irreducible and 
\begin{itemize}
\item either $\bbar T$ is of type (d) and $T$ has genus $>0$ and $2$-rank $0$ and is separable over~$\bbar T$,
\item or $\bbar T$ is of type (c) and $T$ has genus $0$ and is purely inseparable over~$\bbar T$.
\end{itemize}
\end{Prop}

\begin{Proof}
By Raynaud \cite[Th.\,2]{Raynaud1990}, the dual graph of the inverse image of $\bar p$ in $C_0$ is a tree. 
The hyperelliptic involution therefore stabilizes each irreducible component of $C_0$ above~$\bar p$, and so the quotient morphism $C_0\onto\bbar C_0$ induces a bijection on irreducible components above~$\bar p$. In particular, the inverse image $T\subset C_0$ of any irreducible component $\bbar T \subset \bbar C_0$ above~$\bar p$ is irreducible. Moreover, by Raynaud \cite[Lemme 3.1.2]{Raynaud1999} 
the morphism $T\onto\bbar T$ is separable if and only if $T$ corresponds to a leaf of the tree, that is, if and only if $\bbar T$ is an irreducible component of type (d). In that case $T$ has genus $>0$ by \cite[Prop.\,2 (iii)]{Raynaud1990}, and otherwise $T\onto\bbar T$ being purely inseparable implies that $T$ has genus~$0$.
Finally, by the last statement of \cite[Th.\,2]{Raynaud1990} and its proof on page 187 of \cite{Raynaud1990} any such irreducible component has $2$-rank~$0$.
\end{Proof}

\medskip
As in Lehr and Matignon \cite[Thm\,5.1]{LehrMatignon2006} we can now conclude:

\begin{Thm}\label{StapPolDTypeD}
Let $S_f$ be the stability polynomial associated to $f$ by Definition \ref{StabPolDef}. Consider any $\xi_0\in R$ such that $\bar p\in\bar\CU$ is given by $x=\xi_0$, and assume that $S_f(\xi_0)=0$.
  Then the substitution $x=\xi_0+2^\epsilon y$ yields an irreducible component of $\bbar C_0$ of type (d) over $\bar p$ with coordinate~$y$. Conversely, every irreducible component of $\bbar C_0$ of type (d) over $\bar p$ arises in this way.
\end{Thm}

\begin{Proof}
Take any $\xi_0$ and $\epsilon$ with the stated properties. With $\bar\CX$ and $\CX$ as in Theorem \ref{LMThm}, the closed fiber $X_0$ of $\CX$ then has geometric genus $>0$. Let $\CC'$ be the minimal semistable model of $C$ from Proposition \ref{VarRelStabModExists} which dominates both $\CC$ and~$\CX$, and let $T'$ be the irreducible component of the closed fiber of $\CC'$ that surjects to~$X_0$. Write $\CC'\onto\CC$ as the composite of finitely many contractions of unstable irreducible components as in Proposition \ref{ContractSeq}. Then since $T'$ has geometric genus $>0$, its image can never be contracted in this sequence; hence it surjects to an irreducible component of genus $>0$ of~$C_0$. By construction this irreducible component lies above~$\bar p$; hence by Proposition \ref{RaynaudProp} it must surject to an irreducible component of $\bbar C_0$ of type (d). This proves the first statement of the proposition. 

The last statement follows by combining Proposition \ref{RaynaudProp} and the implication (a)$\Rightarrow$(b) of Theorem \ref{LMThm}.
\end{Proof}


\begin{Cons}\label{SUPCDdown}
\rm To construct $\bbar\CC$ above $\bar p$ one first makes a list of all zeros $\xi_i\in R$ of $S_f$ such that $\bar p\in\bar\CU$ is given by $x=\xi_0$. For each of these one computes the associated $\epsilon_i:=\epsilon\in\BQ$ from Lemma \ref{StabPolXi0Lem}.
For efficiency, for any $i$ one may want to remove all $\xi_j$ with $j\not=i$ for which $v(\xi_j-\xi_i)\ge\epsilon_i$, because they yield the same irreducible component of type (d) as~$\xi_i$. Then one applies Construction \ref{AddCompsTypeCDCons} to the substitutions $x=\xi_i+2^{\epsilon_i} y$. By Theorem \ref{StapPolDTypeD} this yields $\bbar\CC_0$ over a neighborhood of~$\bar p$.
\end{Cons}

\begin{Rem}\label{SUPCDup}
\rm To construct $\CC$ above $\bar p$ one first constructs $\bbar\CC$. Over the smooth locus of $\bbar\CC$ one can find explicit coordinates of $\CC$ by Proposition \ref{LehrPropositionCoveringPol}. Over any double point of $\bbar\CC$ above $\bar p$ one finds explicit coordinates of $\CC$ by Proposition \ref{EvenPointNorm<2} below. In fact, such a double point is always even, because by construction all marked points lie on one side of it, namely in $C_0\setminus\{\bar p\}$. Also, Proposition \ref{EvenPointNorm<2} applies in this situation by Proposition \ref{ComponentsTypeBProp}, because we already know that there are no irreducible components of type (b) above that double point.
\end{Rem}

\subsection{Double points}
\label{DoublePoints}

Now we assume that $\bar p$ is a double point of~$\bar C_0$ and identify a neighborhood $\bar\CU \subset \bar\CC$ with an open subscheme of $\Spec R[x,y]/(xy-a)$ for some nonzero $a\in\Fm$, such that $\bar p$ is given by $x=y=0$. After multiplying $a$ by a unit we may assume that $a=2^\alpha$ for some $\alpha>0$, where $\alpha$ is the thickness of~$\bar p$.
Identifying $y$ with $\frac{2^\alpha}{x}$, the ring in question is isomorphic to the subring $R[x,\tfrac{2^\alpha}{x}]$ of the ring of Laurent polynomials $K[x^{\pm1}]$.

Recall from Definition \ref{EvenOddDef} that the double point $\bar p$ is called \emph{even} if each connected component of $\bar C_0 \setminus \{\bar p\}$ contains an even number of the points $\bar p_1,\ldots,\bar p_{2g+2}$. Otherwise, it is called \emph{odd}. 

\begin{Prop}\label{PropHyperEquationNearP}
There exist polynomials $f_1\in R[x]$ and $f_2\in R[y]$ with constant terms~$1$ such that an equation for $C$ is given by
$$z^2\ =\ \biggl\{\!\!\begin{array}{cl}
f_1(x)f_2(\tfrac{2^\alpha}{x}) &\hbox{if $\bar p$ is even,}\\[3pt]
xf_1(x)f_2(\tfrac{2^\alpha}{x}) &\hbox{if $\bar p$ is odd.}
\end{array}$$
\end{Prop}

\begin{Proof}
For any $i$ let $\xi_i\in K\cup\{\infty\}$ be the $x$-coordinate of the point~$\bar P_i$. After reordering we may assume that $v(\xi_i)>0$ if and only if $1\le i\le r$ and that $\xi_i=\infty$ at most for $i=2g+2$. 

Observe that a point $\xi\in K$ reduces to the point $x=y=0$ in $\Spec R[x,y]/(xy-2^\alpha)$ if and only if both $\xi$ and $\smash{\tfrac{2^\alpha}{\xi}}$ lie in~$\Fm$, or equivalently if $0<v(\xi)<\alpha$. By the semistability assumption on~$\bar\CC$ this does not happen for any of the points~$\xi_i$. Thus for all $1\le i\le r$ we have  $v(\xi_i)\ge\alpha$. These points therefore reduce to a point on $\bar C_0$ on one side of the double point~$\bar p$, and all others, including~$\infty$, to a point on the other side. Thus the number of marked points $\bar p_i$ on one side of $\bar p$ is~$r$, and so $\bar p$ is odd if and only if $r$ is odd.

Now set 
\vskip-15pt
$$f_2(y)\ :=\ \prod_{i=1}^r \bigl(1-\tfrac{\xi_i}{2^{\alpha}} \cdot y\bigr)
\quad\hbox{and}\quad
f_1(x)\:=\ \prod_{i=r+1}^{2g+2} \bigl(1-\tfrac{1}{\xi_i}\cdot x\bigr).$$
By construction both are polynomials with coefficients in $R$ and constant coefficient~$1$. Moreover 
$$x^rf_1(x)f_2(\tfrac{2^\alpha}{x})\ =\ 
\prod_{i=1}^r \bigl(x-\xi_i\bigr) \cdot \prod_{i=r+1}^{2g+2} \bigl(1-\tfrac{1}{\xi_i}\cdot x\bigr)$$
is a polynomial with a simple zero at each $\xi_i\neq\infty$ and no other zeros. Thus the hyperelliptic curve $C$ can be described by the equation $z^2=x^rf_1(x)f_2(\tfrac{2^\alpha}{x})$. Finally write $r=2s+e$ with $s\in\BZ$ and $e\in\{0,1\}$, so that $e=0$ if $\bar p$ is even and $e=1$ if $\bar p$ is odd. Using the substitution $z=x^sz$ we can then rewrite the equation for $C$ in the form $z^2=x^ef_1(x)f_2(\tfrac{2^\alpha}{x})$, as desired.
\end{Proof}

\subsection{Odd double points}
\label{OddDoublePoints}

In this section we assume that $\bar p$ is an odd double point. 

\begin{Prop}\label{OddPointProp}
There exists a neighborhood $\bar\CU \subset \bar\CC$ of~$\bar p$, whose normalization $\CU$ in the function field of $C$ is equal to $\CC$ over~$\bar\CU$ and possesses a unique double point $p$ over $\bar p$ of thickness~$\tfrac{\alpha}{2}$. Moreover, both irreducible components of $C_0$ at $p$ are purely inseparable over the respective irreducible components of~$\bar C_0$.
\end{Prop}

\begin{Proof}
By Proposition \ref{PropHyperEquationNearP} there exist polynomials $f_1\in R[x]$ and $f_2\in R[y]$ with constant terms~$1$ such that the function field of $C$ is $K(x,z)$ with $z^2 = xf_1(x)f_2(\tfrac{2^\alpha}{x})$. Choose $R_2$ finite over $R_1$ as in Section \ref{App} such that $2^{\alpha/2}$ and all coefficients of $f_1$ and $f_2$ lie in~$R_2$. Setting $w:= \smash{\frac{2^{\alpha/2}}{x}z}$, we then have $zw=2^{\alpha/2}f_1(x)f_2(\tfrac{2^\alpha}{x})$ and $w^2 = y f_1(x)f_2(\tfrac{2^\alpha}{x})$.
In particular this shows that $z$ and $w$ are integral over the ring 
$$A\ :=\ R_2\bigl[x,\tfrac{2^\alpha}{x}\bigr]\ \cong\  R_2[x,y]/(xy-2^\alpha).$$

\begin{Lem}\label{OddPointLem1}
The $R_2$-subalgebra $B\subset K(x,z)$ generated by $x,y,z,w$ has the presentation
$$B\ =\ R_2\bigl[x,y,z,w\bigr]
\!\Bigm/\! \biggl(\!\begin{array}{lll}
xy-2^\alpha, & xw-2^{\alpha/2}z, & yz-2^{\alpha/2}w, \\[3pt]
z^2-xf_1 f_2, & zw-2^{\alpha/2}f_1 f_2, & w^2-yf_1f_2
\end{array}\!\biggr).$$
\end{Lem}

\begin{Proof}
Consider the $A$-submodule $M\subset K(x,z)$ that is generated by $z$ and~$w$. By the definition of $w$ this is equal to $I\frac{z}{x}$ for the ideal $I:=(x,2^{\alpha/2})$ of~$A$. A simple computation shows that all relations between the generators of this ideal are linear combinations of the relations $x\cdot2^{\alpha/2} = 2^{\alpha/2}\cdot x$ and $y\cdot x = 2^{\alpha/2}\cdot2^{\alpha/2}$. Thus all $A$-linear relations between the generators of $M$ are linear combinations of the relations  $x w= 2^{\alpha/2} z$ and $y z = 2^{\alpha/2} w$. Also, as the hyperelliptic involution $\sigma$ acts by $-1$ on~$M$, the sum $A+M$ is direct. Moreover, the above relations for $z^2$ and $zw$ and $w^2$ show that $A+M$ is already the subring in question. Together this yields the stated presentation.
\end{Proof}

\medskip
Next let $r$ be a uniformizer of the complete discrete valuation ring~$R_2$. Then our double point $\bar p$ corresponds to the maximal ideal $\Fp := (r,x,y) \subset A$, and the only maximal ideal above~$\Fp$ is $\Fq := (r,x,y,z,w) \subset B$. 

\begin{Lem}\label{OddPointLem2}
The completion $\hat B$ of $B$ at $\Fq$ is isomorphic to $R_2[[z,u]]/(zu-2^{\alpha/2})$.
\end{Lem}

\begin{Proof}
The polynomials $f_1,f_2\in R_2[x,y]$ have constant term~$1$, so they represent units in~$\hat B$. The equations $z^2=xf_1 f_2$ and $w^2=yf_1f_2$ thus imply that $x$ and $y$ lie in the square of the maximal ideal $\hat\Fq$ of~$\hat B$. Since $u:= w/f_1f_2$ is congruent to $w$ modulo~$\hat\Fq$, it follows that $\hat\Fq$ is already generated by $r,z,u$. As $\hat B$ is a noetherian complete local ring, this shows that the natural homomorphism $R_2[[z,u]]\to \hat B$ is surjective. In particular, there exists a power series $g\in R_2[[z,u]]$ with unit constant term such that $f_1f_2=g$ in~$\hat B$. Thus we have $x=z^2g^{-1}$ and $y=w^2g^{-1} = u^2g$ in $\hat B$ and can eliminate the variables~$x$ and~$y$. Moreover the relation $zw=2^{\alpha/2}g$ is equivalent to $zu=2^{\alpha/2}$ and is quickly seen to imply all the remaining relations in Lemma \ref{OddPointLem1}. Thus we have the desired isomorphism.
\end{Proof}

\medskip
Lemma \ref{OddPointLem2} implies that $\hat B$ is normal. Since $\hat B$ is faithfully flat over the localization~$B_\Fq$, from Liu \cite[\S4.1.2 Exerc.\,1.16]{LiuAlgGeo2002} it follows that $B_\Fq$ is normal as well. On the other hand observe that by construction $B$ is finite over~$A$ and therefore contained in the normalization $\tilde B\subset K(x,z)$ of~$A$. By Proposition \ref{NormalFin} this normalization is finite over~$A$ and hence over~$B$, and so by Liu \cite[\S4.1.2 Prop.\,1.29]{LiuAlgGeo2002} the normal locus is open in $\Spec B$. Together this implies that $B[\frac{1}{b}] = \tilde B[\frac{1}{b}]$ for some $b\in B\setminus\Fq$. 

Since $b$ is a unit in~$\hat B$, Lemma \ref{OddPointLem2} implies that $\Spec B$ is semistable with exactly one ordinary double point over $\Fp$ of thickness~$\alpha/2$. After base change from $\Spec R_2$ to $\Spec R$ we conclude that 
the normalization $\CU$ is semistable with a unique double point of thickness $\frac{\alpha}{2}$ near~$\bar p$. By Proposition \ref{RelStabMod} (c) the morphism $\CC\onto\CX$ is an isomorphism there.

\medskip
Finally, the presentation in Lemma \ref{OddPointLem1} implies that $B[\frac{1}{x}] \cong R_2[x^{\pm1},z]/(z^2-xf_1f_2)$. As the equation $z^2=xf_1f_2)$ is inseparable modulo $\Fm$, it follows that the irreducible component of $C_0$ that meets $p$ and on which $x\not=0$ is inseparable over the corresponding irreducible component of~$\bar C_0$. By symmetry the same follows for the other irreducible component of $C_0$ that meets~$p$. This finishes the proof.
\end{Proof}

\subsection{Even double points}
\label{EvenDoublePoints}

In this section we assume that $\bar p$ is an odd double point. We fix polynomials $f_1\in R[x]$ and $f_2\in R[\tfrac{2^\alpha}{x}]$ with constant terms~$1$ as in Proposition \ref{PropHyperEquationNearP}, such that the function field of $C$ is $K(x,z)$ with $z^2 = f_1f_2$. Then $f:= f_1f_2$ and $\alpha$ satisfy the condition (\ref{AlphaDef}) from Section \ref{ScalingDecomp}.

By Proposition \ref{AddCompsTypeBCDCoord} the irreducible components of $\bbar C_0$ of type (b) above $\bar p$ are given by coordinates $z=\frac{x}{b}$ for nonzero $b\in\Fm$ such that $\frac{2^\alpha}{b}\in\Fm$. After rescaling $z$ by a unit we can describe them with $b=2^\lambda$ for $\lambda\in\BQ\,\cap\, ]\kern1pt0, \alpha\kern1pt[$. Our first job is to decide which values of $\lambda$ occur. For this consider the function $\wbar$ from (\ref{WBarDef}).

\begin{Prop}\label{ComponentsTypeBProp} 
The substitution $x=2^\lambda u$ yields an irreducible component of $\bbar C_0$ of type (b) over $\bar p$ if and only if $\lambda$ is a break point of~$\wbar$. 
\end{Prop}

\begin{Proof}
The proof requires some preparation. Let $\bar\CC'$ be the semistable model of $\bar C$ obtained from $\bar\CC$ by adjoining an irreducible component with coordinate $u=x/2^\lambda$ as in Construction \ref{AddCompsTypeBCons}. 
Let $\CX'$ be the normalization of $\bar\CC'$ in the function field of~$C$, and let  $\CC'$ be the minimal semistable model of $C$ which dominates~$\CX'$ according to Proposition \ref{RelStabMod}. 
Consider the open chart $\bar\CU := \Spec R[u^{\pm1}]$ of $\bar\CC'$ and let $\CU$ denote its inverse image in~$\CX'$. Let $\bar U_0= \Spec k[u^{\pm1}]$
and $U_0$ denote their respective closed fibers.
Fix an odd decomposition $f=h^2+g$ and set $\gamma := v(g(2^\lambda u))$.

\begin{Lem}\label{AdaptedCoveringLemma}  
\begin{enumerate}[(a)]
\item If $\gamma<2$, then $U_0$ is irreducible and purely inseparable over~$\bar U_0$, and it is singular if $\lambda$ is a break point of~$\wbar$, respectively isomorphic to $\BP^1_k\setminus\{0,\infty\}$ if not.
\item If $\gamma=2$, then $U_0$ is irreducible and smooth and \'etale over~$\bar U_0$, and it either has genus $>0$ or is isomorphic to $\BP^1_k$ minus at least three points.
\item If $\gamma>2$, then $U_0$ is the disjoint union of two copies of $\BP^1_k\setminus\{0,\infty\}$, each mapping isomorphically to~$\bar U_0$.
\end{enumerate}
\end{Lem}

\begin{Proof}
By construction $\CU$ is the normalization of $\Spec R[u^{\pm1}]$ in the function field $K(x,z)$.  By Proposition \ref{ScalOpt} the decomposition $f(2^\lambda u)=h(2^\lambda u)^2+g(2^\lambda u)$ is optimal, and so $\wbar(\lambda) = \min\{2,\gamma\}$. We can therefore compute $U_0$ using Proposition \ref{LehrPropositionCovering}. Observe that since $f_1$ and $f_2$ have constant coefficients~$1$, the assumption $0<\lambda<\alpha$ implies that $[f_1(2^\lambda u)] = [f_2(2^\lambda u)] = 1$ and hence $[f(2^\lambda u)] =1$. Moreover, since $g$ has constant coefficient~$0$, it follows that $[g(2^\lambda u)] = 0$ and therefore $[h(2^\lambda u)] = 1$.

\medskip
In the case $\gamma<2$ Proposition \ref{LehrPropositionCovering} (a) implies that
\UseTheoremCounterForNextEquation
\begin{equation}\label{ACLa}
U_0\ \cong\ \Spec k[u^{\pm 1},t]/(t^2-\ell)
\end{equation}
with the nonzero odd Laurent polynomial $\ell := [g(2^\lambda u)/2^\gamma] \in k[u^{\pm1}]$. In particular it is purely inseparable over $\bar U_0= \Spec k[u^{\pm1}]$ and therefore irreducible. Next we know from Proposition \ref{WBarProp} (c) that $\ell$ is a monomial if and only if $\lambda$ is not a break point of~$\wbar$. In that case we can write $\ell=cu^{2r+1}$ with some $c\in k^\times$ and deduce that $U_0 \cong \Spec k[(tu^{-r})^{\pm1}] \cong \BP^1_k\setminus\{0,\infty\}$. Otherwise observe that since $\ell$ is odd and $k$ has characteristic~$2$, we have $\ell=um^2$ for some Laurent polynomial $m\in k[u^{\pm1}]$ that is not a monomial. Thus $\ell$ has a multiple zero at some point in~$k^\times$ and $U_0$ is singular there. This proves (a).

\medskip
In the case $\gamma=2$, Proposition \ref{LehrPropositionCovering} (b) and the fact that $[h(2^\lambda u)] = 1$ imply that
\UseTheoremCounterForNextEquation
\begin{equation}\label{ACLb}
U_0\ \cong\ \Spec k[u^{\pm 1},t]/(t^2+t-\ell)
\end{equation}
with the irreducible odd Laurent polynomial $\ell := [g(2^\lambda u)/4] \in k[u^{\pm1}]$. Thus $U_0$ is irreducible and \'etale over $\Spec k[u^{\pm1}]$ and therefore smooth. Consider the associated covering of smooth projective curves $\pi\colon X_0 \onto \BP^1_k$. 
If $\ell$ lies in $k[u]$, then $\pi$ is unramified over the point $u=0$ and $X_0$ has two points over it. Similarly, if $\ell$ lies in $k[u^{-1}]$, then $\pi$ is unramified over the point $u=\infty$ and $X_0$ has two points over it. In both cases $\pi^{-1}(\{0,\infty\})$ consists of at least three points.
If $\ell$ lies neither in $k[u]$ nor in $k[u^{-1}]$, then $\pi$ is ramified over both $0$ and~$\infty$. As the ramification is wild, the ramification divisor then has multiplicity $\ge2$ at both points. By the Hurwitz formula $X_0$ then has genus $>0$, finishing the proof of (b).

\medskip
In the case $\gamma>2$, Proposition \ref{LehrPropositionCovering} (c) and the fact that $[h(2^\lambda u)] = 1$ imply that
\UseTheoremCounterForNextEquation
\begin{equation}\label{ACLc}
U_0\ \cong\ \Spec k[u^{\pm 1},t]/(t(t+1)),
\end{equation}
proving (c).
\end{Proof}

\medskip
Next let $\bar C'_0$ denote the closed fiber of~$\bar\CC'$, let $E\subset\bar C'_0$ be the exceptional fiber of $\bar\CC'\onto\bar\CC$, and observe that $E$ consists of $\bar U_0$ and two double points of~$\bar C'_0$. Let $C'_0$ denote the closed fiber of~$\CC'$, and consider an irreducible component $Z$ of $C'_0$ which surjects to~$E$.

\begin{Lem}\label{MorphRegLem}
The morphism $\CC'\onto\CX'$ sends $Z\cap C_0^{\prime\reg}$ isomorphically to an irreducible component of~$U_0^\reg$.
\end{Lem}

\begin{Proof}
By construction the hyperelliptic involution $\sigma$ extends to $\CX'$ and thus, by the uniqueness of the minimal semistable model in Proposition \ref{RelStabMod} (b), also to~$\CC'$. By Proposition \ref{BbarModelSemiStab} (a) the quotient $\bbar\CC'$ is a semistable model of $\bar C$, which by construction dominates~$\bar\CC'$. Let $\bbar Z$ denote the image of $Z$ in the closed fiber $\bbar C'_0$ of~$\bbar\CC'$. Then the image of $Z\cap C_0^{\prime\reg}$ in $\bbar C'_0$ is contained in $\bbar Z\cap \bbar C_0^{\prime\reg}$ by Proposition \ref{BbarModelSemiStab} (b). Moreover, applying Proposition \ref{SemistabMorphReg} to $\bbar\CC'\onto\bar\CC'$ shows that the image of $\bbar Z\cap \bbar C_0^{\prime\reg}$ in $\bar C_0'$ is contained in $E\cap\bar C_0^{\prime\reg}=\bar U_0$. By the definition of $\CX'$ the image of $Z\cap C_0^{\prime\reg}$ under the morphism $\CC'\onto\CX'$ is therefore contained in~$U_0$.

Next suppose that $U_0$ possesses a singular point~$q$. By Lemma \ref{AdaptedCoveringLemma} we then have $\gamma<2$ and $U_0\onto\bar U_0\cong \Spec k[u^{\pm1}]$ is purely inseparable.
Thus $q$ is not an ordinary double point, hence $\CX'$ is not semistable there, and so $\CC'\onto\CX'$ is not an isomorphism over~$q$. Any point of $Z$ above $q$ must then be a double point of $C_0'$ where $Z$ meets the exceptional fiber. The image of $Z\cap C_0^{\prime\reg}$ under the morphism $\CC'\onto\CX'$ is therefore contained in~$U_0^\reg$.

On the other hand we know from Proposition \ref{RelStabMod} (c) that the morphism $\CC'\onto\CX'$ is an isomorphism at all points where $\CX'$ is already smooth. In particular it is an isomorphism over~$U_0^\reg$. Thus $Z\cap C_0^{\prime\reg}$ is the inverse image of $U_0^\reg$ in $Z$ and maps isomorphically to an irreducible component of~$U_0^\reg$, as desired.
\end{Proof}

\medskip
Next recall that an irreducible component of $C'_0$ is called unstable if it is isomorphic to $\BP^1_k$ and contains at most two double points. 

\begin{Lem}\label{UnstabLem}
The above $Z$ is unstable if and only if $\lambda$ is not a break point of~$\wbar$. 
\end{Lem}

\begin{Proof}
Suppose that $Z$ is unstable. Then Lemma \ref{MorphRegLem} implies that some irreducible component of $ U_0^\reg$ is isomorphic to $\BP^1_k$ minus at most two closed points. But since the closure of this irreducible component in $\CX'$ surjects to~$E$, it already contains at least two distinct points above the two points $0,\infty\in E\setminus\bar U_0$. Thus some irreducible component of $\bar U_0$ must be isomorphic to $\BP^1_k$ minus exactly two closed points. By Lemma \ref{AdaptedCoveringLemma} and Proposition \ref{WBarProp} (c) this happens only if $\gamma$ is not a break point of~$\wbar$, as desired.

Conversely suppose that $\lambda$ is not a break point of~$\wbar$. Then by Proposition \ref{WBarProp} (c) we have $\gamma\not=2$, and combining Lemmas \ref{AdaptedCoveringLemma} and \ref{MorphRegLem} shows that $Z\cap C_0^{\prime\reg} \cong \BP^1_k\setminus\{0,\infty\}$. On the other hand, since $Z$ surjects to~$E$, it already contains at least two distinct points above the two points $0,\infty\in E\setminus\bar U_0$, and by Lemma \ref{MorphRegLem} these must be double points of~$C_0'$. Thus $Z$ is a semistable rational curve without self-intersection and therefore isomorphic to~$\BP^1_k$. 
Together this implies that $Z$ is unstable, as desired.
\end{Proof}

\medskip
Now we can prove Proposition \ref{ComponentsTypeBProp}. For this recall that in Section \ref{HyperellCurves} we had defined $\CX$ as the normalization of $\bar\CC$ in the function field of~$C$ and then $\CC$ as the minimal semistable model of $C$ that dominates~$\CX$. Since $\bar\CC'$ dominates~$\bar\CC$, it follows that $\CX'$ dominates~$\CX$, and so the minimality of $\CC$ implies that $\CC'$ dominates~$\CC$. Also, the minimality of $\CC'$ implies that the morphism $\CC'\onto\CC$ is an isomorphism if and only if $\CC$ already dominates~$\CX'$. By construction this is so if and only if $\bbar\CC$ dominates~$\bar\CC'$, that is, if and only if the value $\lambda$ occurs for an irreducible component of $\bbar C_0$ of type (b) over~$\bar p$.

Assume now that $\lambda$ is not a break point of~$\wbar$. Then any irreducible component $Z$ of $C'_0$ which surjects to~$E$ is unstable by Lemma \ref{UnstabLem}. By Proposition \ref{ContractSemiStab} (a) it can thus be contracted to a point in another semistable model $\CC''$ of~$C$. Since $Z$ maps to the closed point $\bar p\in\bar\CC$, it also maps to a closed point of~$\CX$; hence $\CC''$ dominates~$\CX$ by Proposition \ref{ContractExist} (c). By the minimality of $\CC$ it follows that $\CC''$ dominates~$\CC$, and so $\CC'\onto\CC$ is not an isomorphism. Thus $\lambda$ does not occur for an irreducible component of~$\bbar C_0$, proving one direction of the desired equivalence.

Conversely assume that $\lambda$ does not occur for an irreducible component of~$\bbar C_0$. Then $\CC'\onto\CC$ is not an isomorphism. By Proposition \ref{ContractSeq} there then exists an unstable irreducible component $Z$ of $C_0'$ that maps to a closed point in~$\CC$. Now recall that by construction $\bar\CC'\onto\bar\CC$ is an isomorphism over $\bar\CC\setminus\{\bar p\}$. Thus $\CX'\onto\CX$ and hence $\CC'\onto\CC$ are isomorphisms over $\bar\CC\setminus\{\bar p\}$. Therefore $Z$ must lie over the closed point $\bar p \in \bar\CC$. If $Z$ were to map to a closed point of~$\bar\CC'$, it would also map to a closed point of~$\CX'$, and the contraction of $Z$ would be another semistable model of $C$ that dominates~$\CX'$ by Proposition \ref{ContractExist} (c), contradicting the minimality of~$\CC'$. This leaves only the possibility that $Z$ surjects to~$E$. Since $Z$ is unstable, Lemma \ref{UnstabLem} then implies that $\lambda$ is not a break point of~$\wbar$, proving the other direction of the desired equivalence.
\end{Proof} 


\begin{Prop}\label{ComponentsTypeBDetailProp}
Let $\lambda$ be a break point of~$\wbar$, let $\bbar T$ be the associated irreducible component of type (b) of $\bbar C_0$, and let $T$ be an irreducible component of $C_0$ over~$\bbar T$.
\begin{enumerate}[(a)]
\item If $\wbar(\lambda)<2$, then $T\onto\bbar T$ is inseparable of degree $2$ 
and there is at least one irreducible component of type (c) or (d) of $\bbar C_0$ that meets~$\bbar T$.
\item If $\wbar(\lambda)=2$, then $T\onto\bbar T$ is separable of degree $2$ and there is no irreducible component of type (c) or (d) of $\bbar C_0$ that meets~$\bbar T$.
\end{enumerate}
\end{Prop}

\begin{Proof}
Keeping the notation from the proof of Proposition \ref{ComponentsTypeBProp}, the cases correspond to the first two cases of Lemma \ref{AdaptedCoveringLemma}, in whose proof we have seen that $\wbar(\lambda)=\gamma$. Since $\lambda$ occurs, the model $\CC$ dominates the model~$\CX'$, and by Proposition \ref{RelStabMod} (c) the morphism $\CC\onto\CX'$ is an isomorphism over the regular locus~$U_0^\reg$. Since $U_0^\reg$ is irreducible in both cases of Lemma \ref{AdaptedCoveringLemma}, it follows that $T$ contains $U_0^\reg$ as an open subscheme and meets no other irreducible component of $C_0$ there. 

In the case $\gamma=2$ the assertions of Lemma \ref{AdaptedCoveringLemma} (b) now finish the proof of (b). 
In the case $\gamma<2$ we know from Lemma \ref{AdaptedCoveringLemma} (a) that $T\onto\bbar T$ is purely inseparable and that $U_0$ possesses at least one singular point. At this point $U_0$ is not semistable; hence $\CC\onto\CX'$ cannot be an isomorphism there, which means that $C_0$ must possess another irreducible component over it. This corresponds to an irreducible component of type (c) or (d) of $\bbar C_0$ that meets~$\bbar T$, finishing the proof of (a).
\end{Proof}


\begin{Rem}\label{EvenDoublePointsCD}
\rm Once we have identified all irreducible components of $\bbar C_0$ of type (b) above~$\bar p$, we can find those of type (c) and (d) as follows. Let $0<\lambda_1<\ldots<\lambda_r<\alpha$ be the break points of~$\wbar$, and construct a model $\tilde\CC$ of $\bar C$ from $\bar\CC$ by applying Construction \ref{AddCompsTypeBCons} to the elements $1,2^{\lambda_1},\ldots,2^{\lambda_r},2^\alpha$. Then $\bbar\CC$ dominates~$\tilde\CC$, and the morphism $\bbar\CC\onto\tilde\CC$ is an isomorphism at all double points of $\tilde\CC$ above~$\bar p$. All the remaining irreducible components above~$\bar p$ therefore lie over smooth points of the closed fiber of~$\tilde\CC$. We can find these by applying the method of Section \ref{SmoothUnmarkedPoints} with $\tilde\CC$ in place of~$\CC$.
\end{Rem}

In the rest of this section we discuss how to compute the minimal semistable model $\CC$ of $C$ over a neighborhood of $\bar p$ under the assumption that $\bbar\CC\onto\bar\CC$ is an isomorphism. This can be applied in particular after replacing $\bar\CC$ by~$\bbar\CC$.

Since every irreducible component of type (c) or (d) above a double point is connected to an irreducible component of type (b), the assumption means that there are no irreducible components of type (b). By Proposition \ref{ComponentsTypeBProp} this means that $\wbar$ has no break points. By the definition (\ref{WBarDef}) of $\wbar$ this leads to one of two cases: Either $\wbar$ is constant with value~$2$, or $\wbar$ is linear with a nonzero slope and has values $\le2$.

\begin{Prop}\label{EvenPointNorm2}
If $\wbar$ is constant with value~$2$, then $\CC\onto\bar\CC$ is \'etale over~$\bar p$ and $\CC$ has two double points of the same thickness over~$\bar p$.
\end{Prop}

\begin{Proof}
Fix a decomposition $f=h^2+g$ which is odd or arises from odd decompositions of $f_1$ and $f_2$ as in Proposition \ref{ScalOptDecomp}. Then for every $\lambda\in\BQ\,\cap\, ]\kern1pt0, \alpha\kern1pt[$ the decomposition $f(2^\lambda u)=h(2^\lambda u)^2+g(2^\lambda u)$ is optimal and by Proposition \ref{WBarProp} (b) we have $v(g(2^\lambda u))>2$. As this is a continuous function of~$\lambda$, it follows that $v(g)\ge2$ and $v(g(2^\alpha u))\ge2$. Writing $g=\sum_ic_ix^i$ this means that $v(c_i)\ge2$ for all $i\ge0$, respectively $v(c_i)+\alpha i\ge2$ for all $i<0$. Therefore $g/4$ lies in the ring $A := R[x,\tfrac{2^\alpha}{x}]$. Moreover, the inequality $v(g(2^{\alpha/2} u))>2$ implies that the constant coefficient of $g/4$ lies in~$\Fm$; hence $g/4$ lies in the maximal ideal $\Fp := (\Fm,x,\tfrac{2^\alpha}{x})$ of~$A$ that corresponds to the point $\bar p\in\bar C_0$. On the other hand the fact that $f_1,f_2\in A$ have constant terms~$1$ implies that $f= f_1f_2$ is congruent to $1$ modulo~$\Fp$. Since $g\in\Fp$, the equation $f=h^2+g$ thus implies that $h\equiv1$ modulo $\Fp$ as well.

As in (\ref{tEq}) the substitution $z=h+2t$ with a new variable~$t$ now yields the equation $ht+t^2=g/4$, showing that $t$ is integral over~$A$. Moreover, as $h$ is a unit at~$\Fp$, the equation is \'etale there. Setting
\UseTheoremCounterForNextEquation
\begin{equation}\label{EvenPointNorm2B}
B\ :=\ A[t]/(ht+t^2-g/4)\ \cong\ R[x,y,t]/(xy-2^\alpha,\ ht+t^2-g/4),
\end{equation}
it follows that $B[a^{-1}]$ is \'etale over $A$ for some $a\in A\setminus\Fp$ and therefore the normalization of $A[a^{-1}]$ in $K(x,z)$. In particular $\Spec B[a^{-1}]$ is semistable, and so it is a local chart of the minimal semistable model $\CC$ over~$\bar\CC$. Thus $\CC\onto\bar\CC$ is \'etale of degree $2$ over~$\bar p$ and $\CC$ has precisely two double points of the same thickness over~$\bar p$.
\end{Proof}


\begin{Prop}\label{EvenPointNorm<2}
If $\wbar$ is linear with a nonzero slope, then $\CC$ has a single double point over~$\bar p$ of half the thickness.
\end{Prop}

\begin{Proof}
Suppose first that $\wbar$ is increasing, so that $w(f) = \wbar(0) < 2$. Let $\bar T$ be the irreducible component of $\bar C_0$ whose intersection with the chart $\Spec R[x,y]/(xy-a)$ is given by $y=0$. Then Proposition \ref{LehrPropositionCovering} (a) implies that $C_0$ possesses a unique irreducible component $T$ that surjects to $\bar T$ and that $T\onto\bar T$ is inseparable of degree~$2$.
Thus $T$ has a unique point above~$\bar p$, and this must be a double point of $C_0$ of half the thickness as $\bar p$ by Proposition \ref{BbarModelSemiStab} (c). If $\wbar$ is decreasing, the same argument applies with the equation $x=0$.
\end{Proof}

\medskip
We give a second proof which at the same time produces local equations for $\CC$.

\medskip
\begin{Proof}
Fix a decomposition $f=h^2+g$ which is odd or arises from odd decompositions of $f_1$ and $f_2$ as in Proposition \ref{ScalOptDecomp}. Choose $R_2$ finite over $R_1$ as in Section \ref{App} such that $2^{\alpha/2}$ and all coefficients of $f,g,h$ lie in~$R_2$. As before we identify $A := R_2[x,\tfrac{2^\alpha}{x}]$ with the ring $R_2[x,y]/(xy-2^\alpha)$ for $y=\tfrac{2^\alpha}{x}$. 

Suppose first that $\wbar$ is increasing, so that $\gamma := \wbar(0) < \delta := \wbar(\alpha) \le 2$. Writing $g=\sum_ic_ix^i$ with $c_i\in R_2$, by the proof of Proposition \ref{WBarProp} we then have 
$$\wbar(\lambda)\ =\  v(g(2^\lambda u))\ =\ \min\bigl\{ v(c_i)+\lambda i \bigm| i\in\BZ \bigr\}$$
for all $\lambda\in \BQ\,\cap\, [0, \alpha]$. Since $\wbar$ is linear with positive slope, the minimum must be attained for a fixed index $i>0$. Moreover, as this conclusion holds for some odd decomposition, this index $i$ is necessarily odd (even if we perform the actual computation with the decomposition from \ref{ScalOptDecomp}).

Now $\gamma = v(g)$ shows that $v(c_i)=\gamma$ and $v(c_j)\ge\gamma$ for all $j>0$, and $v(g(2^\alpha u)) \ge\gamma$ shows that $v(c_j)+\alpha j\ge\gamma$ for all $j<0$. Thus we have $g=2^\gamma x^i\ell$ for some $\ell\in A$ whose constant term $c_i/2^\gamma$ is a unit. The assumptions also imply that $\gamma+\alpha i = v(g(2^\alpha u)) = \delta \le 2$.

Since $i$ is odd, we can substitute $z=h+2^{\gamma/2}x^{(i-1)/2}u$ with a new variable~$u$. Plugging this into the equation $z^2=f=h^2+g=h^2+2^\gamma x^i\ell$ and simplifying yields the equation 
\UseTheoremCounterForNextEquation
\begin{equation}\label{EPN21}
u^2+2^{(2-\gamma)/2}x^{-(i-1)/2}hu\ =\ x\ell.
\end{equation}
Since $\gamma+\alpha i = \delta$, we can rewrite $2^{(2-\gamma)/2}x^{-(i-1)/2}h = 2^{\alpha/2}m$ with
$$m\ :=\ 2^{(2-\delta)/2}y^{(i-1)/2}h\ \in A.$$
The equation (\ref{EPN21}) thus simplifies to $u^2+2^{\alpha/2}mu = x\ell$.
Setting $v := \smash{\frac{2^{\alpha/2}}{x}u} = \smash{\frac{y}{2^{\alpha/2}}u}$, we then have $uv+ymu=2^{\alpha/2}\ell$ and $v^2+ymv=y\ell$. In particular this shows that $u$ and $v$ are integral over~$A$.

\begin{Lem}\label{EvenPointNorm<2Lem1}
The $R_2$-subalgebra $B\subset K(x,z)$ generated by $x,y,u,v$ has the presentation
$$B\ =\ R_2\bigl[x,y,u,v\bigr]
\!\Bigm/\! \biggl(\!\begin{array}{lll}
xy-2^\alpha, & xv-2^{\alpha/2}u, & yu-2^{\alpha/2}v, \\[3pt]
u^2+2^{\alpha/2}mu-x\ell, & uv+ymu-2^{\alpha/2}\ell, & v^2+ymv-y\ell
\end{array}\!\biggr).$$
\end{Lem}

\begin{Proof}
Same as for Lemma \ref{OddPointLem1}: The relations $xv=2^{\alpha/2}u$ and $yu=2^{\alpha/2}v$ give a presentation for the $A$-submodule $M := Au+Av\subset K(x,z)$, the sum $A+M$ is direct, and $A+M$ is already a subring with the last three relations.
\end{Proof}

\medskip
Next let $r$ be a uniformizer of the complete discrete valuation ring~$R_2$. Then our double point $\bar p$ corresponds to the maximal ideal $\Fp := (r,x,y) \subset A$, and the only maximal ideal above~$\Fp$ is $\Fq := (r,x,y,u,v) \subset B$. 

\begin{Lem}\label{EvenPointNorm<2Lem2}
The completion $\hat B$ of $B$ at $\Fq$ is isomorphic to $R_2[[u,w]]/(uw-2^{\alpha/2})$.
\end{Lem}

\begin{Proof}
Same as for Lemma \ref{OddPointLem2}: By construction $\ell$ is a unit at~$\Fp$, so we can use the equations $u^2+2^{\alpha/2}mu=x\ell$ and $v^2+ymv=y\ell$ to eliminate the variables $x$ and~$y$ and replace the variable $v$ by $w := (v+ym)\ell^{-1}$, showing that $\hat B \cong R_2[[u,w]]/(uw-2^{\alpha/2})$.
\end{Proof}

\medskip
By the same argument as at the end of the proof of Proposition \ref{OddPointProp} it follows that $B[\frac{1}{b}]$ is normal for some $b\in B\setminus\Fq$. Finally, as $b$ is a unit in~$\hat B$, Lemma \ref{OddPointLem2} implies that $\Spec B$ is semistable with exactly one ordinary double point over $\Fp$ of thickness~$\alpha/2$. After base change from $\Spec R_2$ to $\Spec R$ the proposition follows in the case that $\wbar$ is increasing.

\medskip
If $\wbar$ is decreasing, we can apply the same arguments to $g(\tfrac{2^\alpha}{x})$ in place of~$g$. With $\gamma := \wbar(\alpha)$ we then get $\delta := \wbar(0) =\gamma+\alpha i\le 2$ for an odd integer $i>0$ and can write $g=2^\gamma y^i\ell$ for some $\ell\in A$ whose constant term is a unit. Substituting $z=h+2^{\gamma/2}y^{(i-1)/2}u$ and setting $v := \smash{\frac{2^{\alpha/2}}{y}u} = \smash{\frac{x}{2^{\alpha/2}}u}$ and $m := 2^{(2-\delta)/2}x^{(i-1)/2}h$, we then get the same equations as in Lemma \ref{EvenPointNorm<2Lem1} except that $x$ and $y$ are interchanged. The rest of the proof is exactly the same.
\end{Proof}

\subsection{Algorithms and summaries}
\label{Algo}

Now we summarize our results and combine them into explicit algorithms. Recall that we start with a semistable model $(\bar\CC,\bar\CP_1,\ldots,\bar\CP_{2g+2})$ of $\bar C\cong\BP^1_K$, where the marked points in the generic fiber are the branch points of $C\onto\bar C$. This may or may not be the stable marked model $\bar\CC$ from Construction \ref{StabMarkedGenus0}.

\begin{Alg}\label{AlgbbarC}
\rm To compute $\bbar\CC$ we can proceed as follows: 

\medskip\emph{Step 1:} We begin with a simplification that will allow us to compute the stability polynomial only once, provided that we stick to linear substitutions later on. Choose a coordinate $\tilde x$ on $\bar C$ such that one of the marked points is given by~$\tilde x=\infty$. Let $\tilde f\in K[\tilde x]$ be the polynomial of degree $2g+1$ whose zeros are the remaining branch points. Make a list of all zeros $\tilde\xi_i\in K$ of the stability polynomial $S_{\tilde f}$ from Definition \ref{StabPolDef}. For each $i$ compute the unique $\tilde\alpha_i\in\BQ$ with $\smash{v(G_{\tilde f}(\tilde\xi_i,2^{\tilde\alpha_i} y))=2}$ from Lemma \ref{StabPolXi0Lem}.

\medskip\emph{Step 2:} For every even double point $\bar p\in\bar C_0$ identify a neighborhood with an open subscheme of $\Spec R[x,y]/(xy-2^\alpha)$ for some $\alpha\in\BQ^{>0}$. This is possible with a linear substitution $\tilde x=ax+b$ for some $a\in K^\times$ and $b\in K$. Then find $f_1$ and $f_2$ as in Proposition \ref{PropHyperEquationNearP} and compute the break points of the function $\wbar$ from (\ref{WBarDef}) associated to $f := f_1f_2$, using either an odd decomposition of $f$ and Proposition \ref{WBarProp} (c), or odd decompositions of $f_1$ and $f_2$ and Proposition \ref{ScalOptDecomp}. By Proposition \ref{ComponentsTypeBProp} this determines all irreducible components of $\bbar C_0$ of type (b) over $\bar p$. Adjoin them to $\bar\CC$ as in Remark \ref{EvenDoublePointsCD}. This yields a model $\tilde\CC$ of $\bar C$ that lies between $\bbar\CC$ and $\bar\CC$ and contains precisely the irreducible components of $\bar C_0$ and those of $\bbar C_0$ of type (b).

\medskip\emph{Step 3:} Cover the smooth locus of $\tilde\CC\setminus\bigcup_{i=1}^{2g+2}\bar\CP_i$ with charts $\bar\CU$ isomorphic to open subschemes of $\Spec R[x]$. Since one of the marked points is given by~$\tilde x=\infty$, this is possible with a linear substitution $\tilde x=ax+b$ for some $a\in K^\times$ and $b\in K$. Choose $c\in K^\times$ such that $f(x) := c\tilde f(ax+b)$ satisfies $v(f)=0$. Then Proposition \ref{StabPolLinSub} implies that the zeros of the stability polynomial $S_f$ are the points $\xi_i := (\tilde\xi_i-b)/a\in K^\times$ for all $i$ as in Step 1 and that the numbers $\alpha_i := \tilde\alpha_i-v(a)\in\BQ$ satisfy $v(G_f(\xi_i,2^{\alpha_i} y))=2$. From these $\xi_i$ select those with $\xi_i\in R$ and $\alpha_i>0$, such that $\xi_i$ defines a point in the closed fiber of~$\bar\CU$. By Theorem \ref{StapPolDTypeD} the substitutions $x=\xi_i+2^{\alpha_i} y$ then yield precisely all irreducible components of $\bbar C_0$ of type (d) over a point of~$\bar\CU$. Adjoin them and the necessary irreducible components of type (c) to $\tilde\CC$ as in Constructions \ref{SUPCDdown} and Construction \ref{AddCompsTypeCDCons}, obtaining~$\bbar\CC$.
\end{Alg}

\begin{Sum}\label{SummIrrComp}
Irreducible components:
\rm Let $\bbar T$ be an irreducible component of $\bbar C_0$ and let $T$ denote its inverse image in~$C_0$. By Proposition \ref{BbarModelSemiStabIrrComp} there are three possibilities for~$T$. With Proposition \ref{LehrPropositionCovering} or \ref{LehrPropositionCoveringPol} we can decide which one occurs. More specifically we have the following cases:
\begin{itemize}
\item[(a)]  Either  $T$ is isomorphic to $\BP^1_k$ and purely inseparable of degree $2$ over~$\bbar T$. This can happen when $\bbar T$ is of type (a), in particular when $\bbar T$ contains a marked point by Remark \ref{SmoothMarkRem} or an odd double point by Proposition \ref{OddPointProp}. It also happens whenever $\bbar T$ is of type (b) and $\wbar(\lambda)<2$ by Proposition \ref{ComponentsTypeBDetailProp} (a), and whenever $\bbar T$ is of type (c) by Proposition \ref{RaynaudProp}.
\item[(b)]  Or  $T$ is irreducible, smooth and separable of degree $2$ over~$\bbar T$. This can happen when $\bbar T$ is of type (a). It also happens whenever $\bbar T$ is of type (b) and $\wbar(\lambda)=2$ by Proposition \ref{ComponentsTypeBDetailProp} (b). Finally it happens whenever $\bbar T$ is of type (d) by Proposition \ref{RaynaudProp}, and in that case $T$ has $2$-rank~$0$.
\item[(c)] 
 Or  $T$ is isomorphic to $\BP^1_k \sqcup \BP^1_k$, each component mapping isomorphically to~$\bbar T$. This can only happen when $\bbar T$ is of type (a).
\end{itemize}
\end{Sum}

\begin{Sum}\label{SummDoubPts}
Double points:
\rm Let $\bbar p$ be a double point of thickness $\alpha$ of~$\bbar C_0$. By Proposition \ref{BbarModelSemiStab} (c) there are two possibilities for $C_0$ above~$\bbar p$:
\begin{itemize}
\item[(a)] There is a unique double point of thickness $\alpha/2$ above~$\bbar p$. This happens when $\bbar p$ is an odd double point by Proposition \ref{OddPointProp}, or when $\bbar p$ is an even double point and $\wbar$ is linear with a nonzero slope by Proposition \ref{EvenPointNorm<2}.
\item[(b)] There are two double points of thickness $\alpha$ above~$\bbar p$ that are interchanged by~$\sigma$. This happens when $\bbar p$ is an even double point and  $\wbar$ is constant with value~$2$ by Proposition \ref{EvenPointNorm2}.
\end{itemize}
\end{Sum}

\begin{Alg}\label{AlgC}
\rm To compute explicit coordinates for $\CC$ one first computes~$\bbar\CC$. Then locally over the smooth locus of $\bbar\CC$ one finds coordinates for $\CC$ using Proposition \ref{LehrPropositionCovering} or \ref{LehrPropositionCoveringPol}. Near a marked section one can directly use the proof of Proposition \ref{SmoothMarkProp}.
Local coordinates over an odd double point are produced by the proof of Proposition \ref{OddPointProp}, and over an even double point by the proof of Proposition \ref{EvenPointNorm2}, respectively the second proof of Proposition \ref{EvenPointNorm<2}.
\end{Alg}

\begin{Rem}\label{DualGraph}
\rm From the above information on irreducible components and double points one can determine the dual graph of~$C_0$.
In particular one can determine whether the dual graph is a tree and hence whether the jacobian of $C$ has good reduction. Specifically, let $n_2$ be the number of irreducible components of $\bbar C_0$ whose inverse image in $C_0$ consists of two irreducible components, and let $m_2$ be the number of double points of $\bbar C_0$ whose inverse image in $C_0$ consists of two double points.
\end{Rem}

\begin{Prop}\label{BettiNumber}
The first Betti number $h^1$ of the dual graph of $C_0$ is $m_2-n_2$.
\end{Prop}

\begin{Proof}
Let $n$ be the number of irreducible components of $\bbar C_0$ and $m$ the number of double points of~$\bbar C_0$. Then the dual graph of $\bbar C_0$ has $n$ vertices and $m$ edges. As this dual graph is a tree, it follows that $m=n-1$. By Summaries \ref{SummIrrComp} and \ref{SummDoubPts} the number of irreducible components of $C_0$ is $n+n_2$ and the number of double points of $C_0$ is $m+m_2$. Thus the dual graph of $C_0$ has $n+n_2$ vertices and $m+m_2$ edges. But this graph is also connected, because $C_0$ is connected. Thus its $h^1$ is $(m+m_2)-(n+n_2)+1 = m_2-n_2$.
\end{Proof}

\begin{Prop}\label{Jacobian}
\begin{itemize}
\item[(a)] The stable reduction of the jacobian of $C$ is an extension of an abelian variety of dimension $g-(m_2-n_2)$ with a torus of dimension $m_2-n_2$.
\item[(b)] The jacobian of $C$ has good reduction if and only if $m_2=n_2$.
\end{itemize}
\end{Prop}

\begin{Proof}
\def\Pic{{\underline{\rm Pic}}}
By Raynaud \cite[Th.\,8.2.1]{Raynaud1970} the relative $\Pic^0_{\CC/R}$ is representable and separated, the condition $(N)^*$ from \cite[D\'ef.\,6.1.4]{Raynaud1970} being satisfied in our case. Moreover $\Pic^0_{\CC/R}$ is smooth by  \cite[Cor.\,2.3.2]{Raynaud1970}. By Proposition \ref{BettiNumber} and the computation in Bosch-L\"utkebohmert-Raynaud  \cite[\S9.2 Example~8]{BoschLuetkebohmertRaynaud1990} its closed fiber is an extension of an abelian variety of dimension $g-(m_2-n_2)$ with a torus of dimension $m_2-n_2$. Thus $\Pic^0_{\CC/R}$ is the stable model of the jacobian of $C$ and its reduction has the stated properties.
\end{Proof}


\section{Examples}
\label{Examples}

\subsection{Genus 1} 
\label{Genus1}

In this section we present the results for  genus~$1$,  leaving most computations to the reader.  In this case  $C\onto\bar C$ has $4$ branch points, and we identify $(\bar C,\bar P_1,\ldots,\bar P_4)$ with $(\BP^1_K,0,1,\infty,a)$ for some $a\in K\setminus\{0,1\}$. After  reordering the branch points and  applying a M\"obius transformation we may without loss of generality assume that $a \in R \setminus \left( 1+\Fm\right)$. Set $\alpha := v(a) \in \BQ^{\ge0}$ and let $(\bar\CC,\bar\CP_1,\ldots,\bar\CP_4)$ be the stable model of $(\BP^1_K,0,1,\infty,a)$ over~$R$. 
The elliptic curve $C$ is  then  defined by the Legendre equation
\UseTheoremCounterForNextEquation
\begin{equation}\label{LegEq}
z^2 =  f(x) :=  x(x-1)(x-a),
\end{equation}
where the ramification points $P_1,\ldots, P_4$ are the $2$-division points. Let $\CC$ be the minimal semistable model of $C$ that dominates~$\bar\CC$. Then the points $P_i$ extend to sections $\CP_i$ of $\CC$ and $(\CC,\CP_1,\ldots,\CP_4)$ is the stable model of $(C,P_1,\ldots,P_4)$ over~$R$ by Proposition \ref{MinMarkModelStab}.
Computing the stability polynomial for $f$ as in Example \ref{StabPolEx} yields
$$S_f(x_0)\ =\ \tfrac{1}{4} \cdot \bigl( 3x_0^4 - 4(a+1)x_0^3 + 6a x_0^2 - a^2\bigr).$$
It turns out that the zeros of $S_f$ are precisely the $x$-coordinates of the nonzero $3$-division points of~$C$, if the point with $x=\infty$ is taken as the identity element.

\medskip

Suppose first that $\alpha=0$. Then $\bar C_0$ is irreducible. Let $\bbar T$ be the irreducible component of $\bbar C_0$ that maps isomorphically to~$\bar C_0$. As it contains marked points, its inverse image $T$ in $C_0$ is an irreducible component that is purely inseparable over it. As $\bar C_0$ does not contain any double point and the total arithmetic genus of $C_0$ is~$1$, there must be exactly one irreducible component $\bbar T'$ of type (d) in~$\bbar C_0$. The inverse image $T'$ of $\bbar T'$ in $C_0$ has genus~$1$  and meets $T$  in a unique one double point.  


\medskip
Now suppose that $\alpha>0$. Then $\bar \CC$ has a unique double point~$\bar p$,   which has thickness $\alpha$. Let $\bbar T_1$ and $\bbar T_2$ be the two irreducible components of $\smash{\bbar C_0}$ that map onto the two irreducible components of~$\bar C_0$. 
As they contains marked points, their inverse images $T_1$ and $T_2$ in $C_0$ are irreducible components that are purely inseparable over them.
For a suitable choice of an equation for $C$ near the double point, computations yield $\wbar(\lambda) = \min\{2,\lambda,\alpha-\lambda\}$. 

\medskip
In the case $\alpha\le4$ the graph of $\wbar$ consists of two oblique line segments
with exactly one breakpoint of value $\alpha/2$ at the midpoint $\lambda=\alpha/2$. Hence $\bbar C_0$ has exactly one component of type (b) over~$\bar p$, which meets $\bbar T_1$ and $\bbar T_2$ in one double point each. Let $T$ denote the inverse image of this component in~$C_0$. Then $T$ meets $T_1$ and $T_2$ in one double point each. 

\medskip
In the case $0<\alpha<4$ we have $\wbar(\alpha/2)=\alpha/2<2$. By Proposition \ref{ComponentsTypeBDetailProp} (a) it follows that $T$ is of genus $0$ and inseparable over~$\bbar T$ and that $\bbar C_0$ possesses at least one irreducible component $\bbar T'$ of type (d) above~$\bar p$. As in the case $\alpha=0$ above this component is unique and meets~$\bbar T$, and its inverse image $T$ in $C_0$ must have genus~$1$.


\medskip
In the case $\alpha=4$ we have $\wbar(\alpha/2)=\alpha/2=2$. By Proposition \ref{ComponentsTypeBDetailProp} (b) it follows that $T$  is irreducible and separable of degree $2$ over~$\bbar T$ and that $\bbar C_0$ has no component of type (c) or (d) that meets~$\bbar T$. As the arithmetic genus of $C_0$ is $1$, this $T$ must then have genus~$1$.


\medskip
In the remaining case $\alpha>4$ the graph of $\wbar$ consists of three line segments
with 
two breakpoints of value $2$ at $\lambda\in\{2,\alpha-2\}$. This means that $\bbar C_0$ possesses exactly two irreducible components $\bbar T_1'$ and $\bbar T_2'$ of type (b) over~$\bar p$, such that each $\bbar T_i'$ meets $\bbar T_i$ at a double point 
and $\bbar T_1'$ meets $\bbar T_2'$ at a double point. 
 Let $T_1'$ and $T_2'$ be the irreducible components of $C_0$ over $\bbar T_1'$ and~$\bbar T_2'$. Then  Proposition \ref{EvenPointNorm2} shows that  $T_1'$ and $T_2'$  meet in two double points.  Since each $T'_i$ meets $T_i$ in a single double point, this implies that the dual graph of $C_0$ contains a loop.  As the total arithmetic genus of $C_0$ is~$1$, it follows that all irreducible components of $C_0$ are rational,  that there are no other irreducible components, and  that the dual graph has no other loops. 

\begin{figure}[h] \centering \FigA
	\caption{Stable reduction in genus $1$.}\label{FigA}
\end{figure}

\medskip
 In all cases the  above results show that  the combinatorial structure of $C_0$  depends only on~$\alpha$. In particular $C_0$ possesses an irreducible component of genus $1$ if and only if ${\alpha\le4}$. Thus $C$ has good unmarked reduction if and only if $\alpha\le 4$. Moreover, computations show that the reduction is supersingular for $\alpha<4$ and ordinary for $\alpha=4$. This  corresponds  to the well-known fact that $C$ has the $j$-invariant ${2^8(a^2-a+1)^3/(a^2-a)^2}$, whose valuation is $\ge0$ if and only if $\alpha\le4$ and $=0$ if and only if $\alpha=4$. In \cite{Yelton2021}, Yelton constructed a semistable model of $C$ using elementary methods and also found the threshold $\alpha=4$.

\medskip
 The results are depicted in Figure  \ref{FigA}, using the colors from Figure~\ref{FigB}. 
An irreducible component of genus $g' > 0$ is labeled by $g=g'$, while those of genus $0$ remain unlabeled.
\subsection{Genus 2}
\label{Genus2}

In this section we describe the results for $g=2$, leaving the detailed computations to look up in \cite{GPComp}. Here $C\onto\bar C$ has $6$ branch points and we start with the associated stable marked model $(\bar\CC,\bar\CP_1,\ldots,\bar\CP_6)$ of~$\bar C$. The seven possibilities for the combinatorial structure of its closed fiber $(\bar C_0,\bar p_1,\ldots,\bar p_6)$ are shown in Figure~\ref{FigTwoAll}. Filled circles signify even double points and empty circles odd double points. The arrows indicate the ways that one type can degenerate into another. 

\begin{figure}[h] \centering 
\FigTwoAll
\caption{The possibilities for $(\bar C_0,\bar p_1,\ldots,\bar p_6)$ in genus~$2$.}\label{FigTwoAll}
\end{figure}

We observe that $\bar C_0$ has at most three even double points and at most one odd double point, and that the combinatorial structure is invariant under symmetries interchanging the former transitively. Let $\alpha\ge\beta\ge\gamma\ge0$ denote the thicknesses of the even double points and $\epsilon\ge0$ that of the odd double point, where we interpret $0$ as the thickness of a double point that does not exist. 
Moreover, 
each even double point is connected to a unique leaf component with exactly two marked points. After possibly interchanging the marked points we can assume that $\bar\CP_1$ meets this component for the double point of thickness $\alpha$ if $\alpha>0$. We also assume that $\bar\CP_2$ has maximal distance from~$\bar\CP_1$, that is, that the sum of the thicknesses of the double points between them is maximal. Then $\bar\CP_2$ must meet the leaf component containing the double point of thickness $\beta$ if $\beta>0$. We identify $\bar C$  with $\BP^1_K$ in such a way that $\bar P_1$ is identified with~$0$ and $\bar P_2$ with~$\infty$. Then $C$ is defined by an equation of the form 
\UseTheoremCounterForNextEquation
\begin{equation}\label{Genus2Equation}
z^2\ =\ f(x)\ =\ ax+bx^2+cx^3+dx^4+ex^5
\end{equation}
with $f\in K[x]$ separable of degree~$5$. Write $\operatorname{disc}(f)$ for the discriminant of $f$. 
Rescaling $x$ and $z$ by factors in $K^\times$\!, we can now arrange to have $v(f)=0$ and
$$\left\{\kern-4pt\begin{array}{rl}
	v(a) \kern-4pt&=\, \alpha+2\epsilon \\[3pt]
	v(e) \kern-4pt&=\, \beta \\[3pt]
	v(\operatorname{disc}(f)) \kern-4pt&=\, 2\alpha + 2\gamma + 6 \epsilon \\[3pt]
	v(b) \kern-4pt&\geq \, \epsilon \\[3pt]
	v(b) \kern-4pt&= \, \epsilon\ \ \hbox{if $\alpha>0$}
\end{array}\kern-4pt\right\}
\text{ or equivalently }
\left\{\kern-3pt\begin{array}{rl}
	\alpha \kern-4pt&=\, v(a)-2\epsilon \\[3pt]
	\beta \kern-4pt&=\, v(e) \\[3pt]
	\gamma \kern-4pt&=\, \frac{1}{2}v(\operatorname{disc}(f))-\alpha-3\epsilon \\[3pt]
	\epsilon \kern-4pt&= \, \min\{v(b),\frac{1}{2}v(\operatorname{disc}(f))-v(a)\}
\end{array}\kern-4pt\right\}
$$
With the equation in this form we choose a square root $\sqrt{bd}$ of $bd$ and  set $\delta:=v\bigl(c-2\sqrt{bd}\kern2pt\bigr)$.
It turns out that the combinatorial structure of  $(C_0,p_1,\ldots,p_6)$  depends only on the values of $\alpha, \beta, \gamma, \delta$ and $\epsilon$, which are all $\geq 0$. 
Computations show that we always have $\delta \geq \min \{2, \gamma\}$, with equality if $\gamma<\min \{\beta, 2\}$. 


\medskip

All in all, the seven cases 
from Figure \ref{FigTwoAll} divide into 54 subcases for the combinatorial structure of $(C_0,p_1,\ldots,p_6)$.
In each subcase we draw the irreducible components using  the same colors as in Figure~\ref{FigB}. We also label any irreducible component of genus $g' > 0$ by $g=g'$, while all irreducible components of genus $0$ remain unlabeled.

\medskip
 The computations were done using Algorithm \ref{AlgbbarC} with the exception of step 3, which became too impractical due to the degree of the stability polynomial. Instead we used indirect arguments and the good reduction part of the stable reduction criteria for the unmarked curve $C$ from Liu \cite[Th.\,1]{Liu1993} in terms of Igusa invariants. In an upcoming article of the first author \cite{Gehrunger2024}, the use of Liu's criterion will be replaced by elementary arguments, making the proof of the classification independent of Liu's criterion. 


\newpage
\medskip
{\bf Case (A):} This is the case of ``equidistant geometry'' of Lehr-Matignon, that is, where $\bar C_0$ is smooth. Hence we have $\alpha=\beta=\gamma=\epsilon=0$, and so the combinatorial structure of $C_0$ depends only on~$\delta$. The 3 possible subcases are sketched in Figure \ref{FigTwoCaseA}.

\begin{figure}[H] \centering 
\FigTwoCaseA
\caption{The possibilities for $(C_0,p_1,\ldots,p_6)$ in the case (A).}
\label{FigTwoCaseA}
\end{figure}

\newpage

\medskip
{\bf Case (B):} Here $\bar C_0$ has exactly one even double point of thickness $\alpha>0$, and we have $\beta=\gamma=\epsilon=0$. The combinatorial structure of $C_0$ depends only on $\alpha$ and $\delta$. The 11 possible subcases are sketched in Figure \ref{FigTwoCaseB}.
\begin{figure}[h!] \centering 
\FigTwoCaseB
\caption{The possibilities for $(C_0,p_1,\ldots,p_6)$ in the case (B).}\label{FigTwoCaseB}
\end{figure}

\newpage
\medskip
{\bf Case (C):} Here $\bar C_0$ has two even double points of respective thicknesses $\alpha\ge\beta>0$, and we have $\gamma=\delta=\epsilon=0$. 
It turns out that the combinatorial structure $C_0$ over each double point is the same as in the case of genus~$1$ and depends only on the thickness of  that  double point. The 6 possible subcases are sketched in Figure \ref{FigTwoCaseC}.

\begin{figure}[h] \centering 
\FigTwoCaseC
\caption{The possibilities  of $(C_0,p_1,\ldots,p_6)$ in the case (C).}\label{FigTwoCaseC}
\end{figure}

{\bf Case (D):} 
Here $\bar C_0$ has one irreducible component without marked points in the middle, which is connected by double points of thicknesses $\alpha\ge\beta\ge\gamma>0$ to leaf components with two marked points each. We also have $\epsilon=0$,
and the combinatorial structure of $C_0$ depends only on $\alpha, \beta, \gamma$ and  $\delta$. The 24 possible subcases are sketched in Figures \ref{FigTwoCaseDone} and \ref{FigTwoCaseDtwo}.

\newpage

\begin{figure}[H] \centering 
\FigTwoCaseDone
\caption{The first 12 possibilities for $(C_0,p_1,\ldots,p_6)$ in the case (D).}\label{FigTwoCaseDone}
\end{figure}

\newpage

\begin{figure}[H] \centering 
\FigTwoCaseDtwo
\caption{The remaining 12 possibilities for $(C_0,p_1,\ldots,p_6)$ in the case (D).}\label{FigTwoCaseDtwo}
\end{figure}

%

\medskip
{\bf Cases (E--G):} Here $\bar C_0$ has an odd double point~$\bar p$ of thickness $\epsilon>0$, and we always have $\gamma=\delta=0$. The combinatorial structure of $C_0$ depends only on the values $\alpha\ge\beta\ge0$. It turns out that the situation on each side of~$\bar p$ is the same as for the reduction of a curve of genus~$1$, and that the two sides are independent of each other.
\medskip
{\bf Case (E):} Here we have $\alpha=\beta=0$, and there is a single subcase only, which is sketched in Figure~\ref{FigTwoCaseE}.

\begin{figure}[H] \centering 
\FigTwoCaseE
\caption{The single possibility for $(C_0,p_1,\ldots,p_6)$ in the case (E).}\label{FigTwoCaseE}
\end{figure}

\medskip
{\bf Case (F):} Here we have $\alpha>\beta=0$. The 3 possible subcases are sketched in Figure~\ref{FigTwoCaseF}.

\begin{figure}[H] \centering 
\FigTwoCaseF
\caption{The possibilities for $(C_0,p_1,\ldots,p_6)$ in the case (F).}\label{FigTwoCaseF}
\end{figure}

\medskip
{\bf Case (G):} Here we have $\alpha\ge\beta>0$. The 6 possible subcases are sketched in Figure~\ref{FigTwoCaseG}.

\begin{figure}[H] \centering 
\FigTwoCaseG
\caption{The possibilities for $(C_0,p_1,\ldots,p_6)$ in the case (G).}\label{FigTwoCaseG}
\end{figure}

{}From the above results, one can also determine the closed fiber $\Cst_0$ of the stable reduction of the unmarked curve~$C$. The list of possible cases and their names are taken from Liu \cite[Th.\,1]{Lehr2001ReductionOP}. There the case distinctions were given in terms of Igusa invariants, which are complicated polynomials in the coefficients of~$f$. Our results yield relatively simple conditions in terms of the numbers $\alpha,\beta,\gamma,\delta$ alone. The seven cases for $\Cst_0$ are sketched in Figure \ref{FigTwoCaseUnmarked}. 

\begin{figure}[H] \centering 
\FigTwoUnmarked
\caption{The possibilities for the stable reduction of the unmarked curve.}\label{FigTwoCaseUnmarked}
\end{figure}


 From these results one can also deduce  the reduction behavior of the jacobian:

\begin{Cor}\label{JacRed}
The reduction of $J(C)$  
\begin{enumerate}[(a)]
\item is good  if and only if $\alpha+\delta \leq 4$ or $3\alpha+\gamma \leq 8$; 
\item has toric rank $1$ if and only if $3\alpha+\gamma>8$ and $\delta+\alpha>4$ and $\beta+\gamma \leq 4$;
\item has toric rank $2$ if and only if $\beta+\gamma>4$. 
\end{enumerate}
\end{Cor}

\begin{Exs}\label{Genus2Examples}
\rm Finally, we provide some examples over $K=\bar \BQ_2$. In each case
the values of $\alpha,\ldots,\epsilon$ can be determined easily from the equation, if necessary after a linear substitution.
The first curve was first considered  by Diophantus as problem 17 in book VI of his {\it{Arithmetica}} \cite{Arithmetica}. 
The last three curves are precisely the hyperelliptic curves of genus $2$ with many automorphisms (compare \cite{MuellerPink2022}).
 \begin{enumerate}[(a)]
 	\item For the curve $z^2=x^6+x^2+1$ we have case (D) with $\alpha=\beta=\gamma=\delta=1$. Hence the stable marked reduction is of type (D22).
 	\item For the curve $z^2=2x^5+2x^3+1$ we have case (A) with $\delta =2/5$. Hence the stable marked reduction is of type (A2). 
\item For the curve $z^2 = x^4 + 3 x^3 + 3 x^2 + 4 x + 1+8x^{-1}$ we have case (B) with $\alpha=3$ and $\delta=1$. Hence the stable  marked reduction is of type (B5). The computation in Example \ref{WExample} already showed that there are two components of type (b) above the double point of $\bar C_0$ of thickness~$3$.
\item For the curve $z^2=16x^5 + x^4 + 4x^3 - x^2 + 8x$ we have case (D) with $\alpha=4$ and $\beta=3$ and $\gamma=1$. Hence the stable  marked reduction is of type (D9).
\item For the curve $z^2=x^5 + x^4 - 4x^3 - 10x^2 + 12x$ we have
case (G) with $\alpha=\beta=1$. Hence the stable marked reduction is of type (G1).
\item For the curve $z^2=x^5-1$ we have case (A) with $\delta=1$. Hence the stable marked reduction is of type (A3). 
\item For the curve $z^2=x^5-x$ we have case (D) with $\alpha=\beta=\gamma=\delta=1/2$. Hence the stable  marked reduction is of type (D22).
\item For the curve $z^2=x^6-1$ we have case (D) with $\alpha=\beta=\gamma=\delta=1$. Hence the stable marked reduction is of type (D22).
\end{enumerate}
\end{Exs}


\end{document}

%% file: GP22-Figures.tex
\def\FigA{
\def\geqone{g\kern1pt{=}\kern1pt1}
\def\arraystretch{0.75}
\begin{tabular}{ | c | M | M | M |}
\hline
\rule{0pt}{2.5ex}    $\alpha$    & $\bar C_0$ & $\bbar{C}_0$  & $C_0$ \\
\hline

$\alpha=0$ & 

\begin{tikzpicture}[baseline=0, scale=0.7]
\draw [thick,black] plot  coordinates {(0,0) (4,1) };
\begin{scope}[shift={(1/8,1/32)}]
\begin{scope}[shift={(3/4,3/16)}]
 \draw [thick,black] plot  coordinates {(-1/32,1/8) (1/32,-1/8) };
  \draw (-2/32,1/8) node [above] {\hbox{\footnotesize$0$}};
\end{scope}   
\begin{scope}[shift={(6/4,6/16)}]
\draw [thick,black] plot  coordinates {(-1/32,1/8) (1/32,-1/8) };
  \draw (-2/32,1/8) node [above] {\hbox{\footnotesize$a$}};
\end{scope}
\begin{scope}[shift={(9/4,9/16)}]
\draw [thick,black] plot  coordinates {(-1/32,1/8) (1/32,-1/8) };
  \draw (-2/32,1/8) node [above] {\hbox{\footnotesize$1$}};
\end{scope}
\begin{scope}[shift={(12/4,12/16)}]
\draw [thick,black] plot  coordinates {(-1/32,1/8) (1/32,-1/8) };
  \draw (-2/32,1/8) node [above] {\hbox{\footnotesize$\infty$}};
\end{scope}
\end{scope}
\end{tikzpicture} &

\begin{center}\begin{tikzpicture}[baseline=0, scale=0.7]
\draw [thick,black] plot  coordinates {(0,0) (4,1) };
\draw [thick,blue] plot  coordinates {(0,24/12) (8/12,-8/12) };
\fill [black] (8/17,2/17) circle (0.08);
\begin{scope}[shift={(2/8,2/32)}]
\begin{scope}[shift={(3/4,3/16)}]
\draw [thick,black] plot  coordinates {(-1/32,1/8) (1/32,-1/8) };
\draw (-2/32,1/8) node [above] {\hbox{\footnotesize$0$}};
\end{scope}   
\begin{scope}[shift={(6/4,6/16)}]
\draw [thick,black] plot  coordinates {(-1/32,1/8) (1/32,-1/8) };
\draw (-2/32,1/8) node [above] {\hbox{\footnotesize$a$}};
\end{scope}
\begin{scope}[shift={(9/4,9/16)}]
\draw [thick,black] plot  coordinates {(-1/32,1/8) (1/32,-1/8) };
\draw (-2/32,1/8) node [above] {\hbox{\footnotesize$1$}};
\end{scope}
\begin{scope}[shift={(12/4,12/16)}]
\draw [thick,black] plot  coordinates {(-1/32,1/8) (1/32,-1/8) };
\draw (-2/32,1/8) node [above] {\hbox{\footnotesize$\infty$}};
\end{scope}
\end{scope}
\end{tikzpicture}\end{center}&

\begin{center}\begin{tikzpicture}[baseline=0, scale=0.7]
\draw [thick,black] plot  coordinates {(0,0) (4,1) };
\draw [thick,blue] plot  coordinates {(0,24/12) (8/12,-8/12) };
\draw (1/12,20/12) node [right] {\hbox{\footnotesize\,$\geqone$\;}};
\fill [black] (8/17,2/17) circle (0.08);
\begin{scope}[shift={(2/8,2/32)}]
\begin{scope}[shift={(3/4,3/16)}]
\draw [thick,black] plot  coordinates {(-1/32,1/8) (1/32,-1/8) };
\draw (-2/32,1/8) node [above] {\hbox{\footnotesize$0$}};
\end{scope}   
\begin{scope}[shift={(6/4,6/16)}]
\draw [thick,black] plot  coordinates {(-1/32,1/8) (1/32,-1/8) };
\draw (-2/32,1/8) node [above] {\hbox{\footnotesize$a$}};
\end{scope}
\begin{scope}[shift={(9/4,9/16)}]
\draw [thick,black] plot  coordinates {(-1/32,1/8) (1/32,-1/8) };
\draw (-2/32,1/8) node [above] {\hbox{\footnotesize$1$}};
\end{scope}
\begin{scope}[shift={(12/4,12/16)}]
\draw [thick,black] plot  coordinates {(-1/32,1/8) (1/32,-1/8) };
\draw (-2/32,1/8) node [above] {\hbox{\footnotesize$\infty$}};
\end{scope}
\end{scope}
\end{tikzpicture}\end{center} \\

\hline
$0 < \alpha <4$ & 

\begin{center}\begin{tikzpicture}[baseline=0, scale=0.7]		
\draw [thick,black] plot  coordinates {(-2,-2) (1/2,1/2 ) };
\draw [thick,black] plot  coordinates {(2,-2) (-1/2,1/2 ) };
\draw [thick,black] plot  coordinates {(3/2-1/8   +1/8,-3/2+1/8   +1/8) (3/2-1/8   -1/8,-3/2+1/8   -1/8) };
\draw [thick,black] plot  coordinates {(1/2 +1/8    +1/8,-1/2 -1/8  +1/8) (1/2 +1/8  -1/8,-1/2-1/8   -1/8) };
\draw [thick,black] plot  coordinates {(-1/2 -1/8    -1/8,-1/2 -1/8  +1/8) (-1/2 -1/8  +1/8,-1/2-1/8   -1/8) };
\draw [thick,black] plot  coordinates {(-3/2+1/8   -1/8,-3/2+1/8   +1/8) (-3/2+1/8   +1/8,-3/2+1/8   -1/8) };
\fill [black] (0,0) circle (0.08);
\draw (-1/2 -1/8 ,-1/4-1/16 ) node [left] {\hbox{\footnotesize$0$}};
\draw (1/2 +1/8 ,-1/4-1/16 ) node [right] {\hbox{\footnotesize$1$}};
\draw (-3/2 + 1/8 ,-5/4+1/8 +1/16 ) node [left] {\hbox{\footnotesize$a$}};
\draw (+3/2 - 1/8 ,-5/4+1/8 +1/16 ) node [right] {\hbox{\footnotesize$\infty$}};
\end{tikzpicture}\end{center} &

\begin{center}\begin{tikzpicture}[baseline=0, scale=0.7]		
\draw [thick,black] plot  coordinates {(-15/8,-14/8) (-3/4,1/2 ) };
\draw [thick,black] plot  coordinates {(15/8,-14/8) (3/4,1/2 ) };
\draw [thick,orange] plot coordinates {(-7/4,0) (7/4,0) };
\draw [thick,blue] plot coordinates {(0,7/4-1/4) (0,-1) };
\fill [black] (0,0) circle (0.08);
\fill [black] (-1,0) circle (0.08);
\fill [black] (1,0) circle (0.08);
\begin{scope}[shift={(-1-2/6,-2/3)}]
\draw [thick,black] plot  coordinates {(-1/8,1/16) (1/8,-1/16) };
\draw (-0.05,0.2) node [left] {\hbox{\footnotesize$0$}};
\end{scope}
\begin{scope}[shift={(-1-4/6,-4/3)}]
\draw [thick,black] plot  coordinates {(-1/8,1/16) (1/8,-1/16) };
\draw (-0.05,0.2) node [left] {\hbox{\footnotesize$a$}};
\end{scope}
\begin{scope}[shift={(1+2/6,-2/3)}]
\draw [thick,black] plot  coordinates {(1/8,1/16) (-1/8,-1/16) };
\draw (0.05,0.15) node [right] {\hbox{\footnotesize$1$}};
\end{scope}
\begin{scope}[shift={(1+4/6,-4/3)}]
\draw [thick,black] plot  coordinates {(1/8,1/16) (-1/8,-1/16) };
\draw (0.05,0.15) node [right] {\hbox{\footnotesize$\infty$}};
\end{scope}
\end{tikzpicture}\end{center} &

\begin{center}\begin{tikzpicture}[baseline=0, scale=0.7]		
\draw [thick,black] plot  coordinates {(-15/8,-14/8) (-3/4,1/2 ) };
\draw [thick,black] plot  coordinates {(15/8,-14/8) (3/4,1/2 ) };
\draw [thick,orange] plot coordinates {(-7/4,0) (7/4,0) };
\draw [thick,blue] plot coordinates {(0,7/4-1/4) (0,-1) };
\fill [black] (0,0) circle (0.08);
\fill [black] (-1,0) circle (0.08);
\fill [black] (1,0) circle (0.08);
\draw (0,6/4-3/8) node [right] {\hbox{\footnotesize$\geqone$}};   
\begin{scope}[shift={(-1-2/6,-2/3)}]
\draw [thick,black] plot  coordinates {(-1/8,1/16) (1/8,-1/16) };
\draw (-0.05,0.2) node [left] {\hbox{\footnotesize$0$}};
\end{scope}
\begin{scope}[shift={(-1-4/6,-4/3)}]
\draw [thick,black] plot  coordinates {(-1/8,1/16) (1/8,-1/16) };
\draw (-0.05,0.2) node [left] {\hbox{\footnotesize$a$}};
\end{scope}
\begin{scope}[shift={(1+2/6,-2/3)}]
\draw [thick,black] plot  coordinates {(1/8,1/16) (-1/8,-1/16) };
\draw (0.05,0.15) node [right] {\hbox{\footnotesize$1$}};
\end{scope}
\begin{scope}[shift={(1+4/6,-4/3)}]
\draw [thick,black] plot  coordinates {(1/8,1/16) (-1/8,-1/16) };
\draw (0.05,0.15) node [right] {\hbox{\footnotesize$\infty$}};
\end{scope}
\end{tikzpicture}\end{center} \\

\hline
$\alpha =4$ & 

\begin{center}\begin{tikzpicture}[baseline=0, scale=0.7]
\draw [thick,black] plot  coordinates {(-2,-2) (1/2,1/2 ) };
\draw [thick,black] plot  coordinates {(2,-2) (-1/2,1/2 ) };
\draw [thick,black] plot  coordinates {(3/2-1/8   +1/8,-3/2+1/8   +1/8) (3/2-1/8   -1/8,-3/2+1/8   -1/8) };
\draw [thick,black] plot  coordinates {(1/2 +1/8    +1/8,-1/2 -1/8  +1/8) (1/2 +1/8  -1/8,-1/2-1/8   -1/8) };
\draw [thick,black] plot  coordinates {(-1/2 -1/8    -1/8,-1/2 -1/8  +1/8) (-1/2 -1/8  +1/8,-1/2-1/8   -1/8) };
\draw [thick,black] plot  coordinates {(-3/2+1/8   -1/8,-3/2+1/8   +1/8) (-3/2+1/8   +1/8,-3/2+1/8   -1/8) };
\fill [black] (0,0) circle (0.08);
\draw (-1/2 -1/8 ,-1/4-1/16 ) node [left] {\hbox{\footnotesize$0$}};
\draw (1/2 +1/8 ,-1/4-1/16 ) node [right] {\hbox{\footnotesize$1$}};
\draw (-3/2 + 1/8 ,-5/4+1/8 +1/16 ) node [left] {\hbox{\footnotesize$a$}};
\draw (+3/2 - 1/8 ,-5/4+1/8 +1/16 ) node [right] {\hbox{\footnotesize$\infty$}};
\end{tikzpicture}\end{center} &

\begin{center}\begin{tikzpicture}[baseline=0, scale=0.7]		
\draw [thick,black] plot  coordinates {(-15/8,-14/8) (-3/4,1/2 ) };
\draw [thick,black] plot  coordinates {(15/8,-14/8) (3/4,1/2 ) };
\draw [thick,orange] plot coordinates {(-7/4,0) (7/4,0) };
\fill [black] (-1,0) circle (0.08);
\fill [black] (1,0) circle (0.08);
\begin{scope}[shift={(-1-2/6,-2/3)}]
\draw [thick,black] plot  coordinates {(-1/8,1/16) (1/8,-1/16) };
\draw (-0.05,0.2) node [left] {\hbox{\footnotesize$0$}};
\end{scope}
\begin{scope}[shift={(-1-4/6,-4/3)}]
\draw [thick,black] plot  coordinates {(-1/8,1/16) (1/8,-1/16) };
\draw (-0.05,0.2) node [left] {\hbox{\footnotesize$a$}};
\end{scope}
\begin{scope}[shift={(1+2/6,-2/3)}]
\draw [thick,black] plot  coordinates {(1/8,1/16) (-1/8,-1/16) };
\draw (0.05,0.15) node [right] {\hbox{\footnotesize$1$}};
\end{scope}
\begin{scope}[shift={(1+4/6,-4/3)}]
\draw [thick,black] plot  coordinates {(1/8,1/16) (-1/8,-1/16) };
\draw (0.05,0.15) node [right] {\hbox{\footnotesize$\infty$}};
\end{scope}
\end{tikzpicture}\end{center} &

\begin{center}\begin{tikzpicture}[baseline=0, scale=0.7]		
\draw [thick,black] plot  coordinates {(-15/8,-14/8) (-3/4,1/2 ) };
\draw [thick,black] plot  coordinates {(15/8,-14/8) (3/4,1/2 ) };
\draw [thick,orange] plot coordinates {(-7/4,0) (7/4,0) };
\fill [black] (-1,0) circle (0.08);
\fill [black] (1,0) circle (0.08);
\draw (0,0) node [above] {\hbox{\footnotesize$\geqone$}};   
\begin{scope}[shift={(-1-2/6,-2/3)}]
\draw [thick,black] plot  coordinates {(-1/8,1/16) (1/8,-1/16) };
\draw (-0.05,0.2) node [left] {\hbox{\footnotesize$0$}};
\end{scope}
\begin{scope}[shift={(-1-4/6,-4/3)}]
\draw [thick,black] plot  coordinates {(-1/8,1/16) (1/8,-1/16) };
\draw (-0.05,0.2) node [left] {\hbox{\footnotesize$a$}};
\end{scope}
\begin{scope}[shift={(1+2/6,-2/3)}]
\draw [thick,black] plot  coordinates {(1/8,1/16) (-1/8,-1/16) };
\draw (0.05,0.15) node [right] {\hbox{\footnotesize$1$}};
\end{scope}
\begin{scope}[shift={(1+4/6,-4/3)}]
\draw [thick,black] plot  coordinates {(1/8,1/16) (-1/8,-1/16) };
\draw (0.05,0.15) node [right] {\hbox{\footnotesize$\infty$}};
\end{scope}
\end{tikzpicture}\end{center} \\

\hline 
$\alpha >4$ & 

\begin{center}\begin{tikzpicture}[baseline=0, scale=0.7]
\draw [thick,black] plot  coordinates {(-2,-2) (1/2,1/2 ) };
\draw [thick,black] plot  coordinates {(2,-2) (-1/2,1/2 ) };
\draw [thick,black] plot  coordinates {(3/2-1/8   +1/8,-3/2+1/8   +1/8) (3/2-1/8   -1/8,-3/2+1/8   -1/8) };
\draw [thick,black] plot  coordinates {(1/2 +1/8    +1/8,-1/2 -1/8  +1/8) (1/2 +1/8  -1/8,-1/2-1/8   -1/8) };
\draw [thick,black] plot  coordinates {(-1/2 -1/8    -1/8,-1/2 -1/8  +1/8) (-1/2 -1/8  +1/8,-1/2-1/8   -1/8) };
\draw [thick,black] plot  coordinates {(-3/2+1/8   -1/8,-3/2+1/8   +1/8) (-3/2+1/8   +1/8,-3/2+1/8   -1/8) };
\fill [black] (0,0) circle (0.08);
\draw (-1/2 -1/8 ,-1/4-1/16 ) node [left] {\hbox{\footnotesize$0$}};
\draw (1/2 +1/8 ,-1/4-1/16 ) node [right] {\hbox{\footnotesize$1$}};
\draw (-3/2 + 1/8 ,-5/4+1/8 +1/16 ) node [left] {\hbox{\footnotesize$a$}};
\draw (+3/2 - 1/8 ,-5/4+1/8 +1/16 ) node [right] {\hbox{\footnotesize$\infty$}};
\end{tikzpicture}\end{center} &

\begin{center}\begin{tikzpicture}[baseline=0, scale=0.7]		
\draw [thick,black] plot  coordinates {(-15/8,-14/8) (-3/4,1/2 ) };
\draw [thick,black] plot  coordinates {(15/8,-14/8) (3/4,1/2 ) };
\draw [thick,orange] plot coordinates {(-18/12,4/12) (6/12,-12/12) };
\draw [thick,orange] plot coordinates {(18/12,4/12) (-6/12,-12/12) };
\fill [black] (-1,0) circle (0.08);
\fill [black] (1,0) circle (0.08);
\fill [black] (0,-2/3) circle (0.08);
\begin{scope}[shift={(-1-2/6,-2/3)}]
\draw [thick,black] plot  coordinates {(-1/8,1/16) (1/8,-1/16) };
\draw (-0.05,0.2) node [left] {\hbox{\footnotesize$0$}};
\end{scope}
\begin{scope}[shift={(-1-4/6,-4/3)}]
\draw [thick,black] plot  coordinates {(-1/8,1/16) (1/8,-1/16) };
\draw (-0.05,0.2) node [left] {\hbox{\footnotesize$a$}};
\end{scope}
\begin{scope}[shift={(1+2/6,-2/3)}]
\draw [thick,black] plot  coordinates {(1/8,1/16) (-1/8,-1/16) };
\draw (0.05,0.15) node [right] {\hbox{\footnotesize$1$}};
\end{scope}
\begin{scope}[shift={(1+4/6,-4/3)}]
\draw [thick,black] plot  coordinates {(1/8,1/16) (-1/8,-1/16) };
\draw (0.05,0.15) node [right] {\hbox{\footnotesize$\infty$}};
\end{scope}
\end{tikzpicture}\end{center} &

\begin{center}\begin{tikzpicture}[baseline=0, scale=0.7]		
\draw [thick,black] plot  coordinates {(-15/8,-14/8) (-6/8,4/8 ) };
\draw [thick,black] plot  coordinates {(15/8,-14/8) (3/4,1/2 ) };
\draw [thick,orange] plot [smooth, tension=1] coordinates {(-3/2,2/8) (1/4,-4/8) (-1/3,-10/8) };
\draw [thick,orange] plot [smooth, tension=1] coordinates {(3/2,2/8) (-1/4,-4/8) (1/3,-10/8) };
\fill [black] (0.095/2-1,0.095) circle (0.08);
\fill [black] (-0.095/2+1,0.095) circle (0.08);
\fill [black] (0,-0.30) circle (0.08);
\fill [black] (0,-1.085) circle (0.08);
\begin{scope}[shift={(-1-2/6,-2/3)}]
\draw [thick,black] plot  coordinates {(-1/8,1/16) (1/8,-1/16) };
\draw (-0.05,0.2) node [left] {\hbox{\footnotesize$0$}};
\end{scope}
\begin{scope}[shift={(-1-4/6,-4/3)}]
\draw [thick,black] plot  coordinates {(-1/8,1/16) (1/8,-1/16) };
\draw (-0.05,0.2) node [left] {\hbox{\footnotesize$a$}};
\end{scope}
\begin{scope}[shift={(1+2/6,-2/3)}]
\draw [thick,black] plot  coordinates {(1/8,1/16) (-1/8,-1/16) };
\draw (0.05,0.15) node [right] {\hbox{\footnotesize$1$}};
\end{scope}
\begin{scope}[shift={(1+4/6,-4/3)}]
\draw [thick,black] plot  coordinates {(1/8,1/16) (-1/8,-1/16) };
\draw (0.05,0.15) node [right] {\hbox{\footnotesize$\infty$}};
\end{scope}
\end{tikzpicture}\end{center} \\

\hline
\end{tabular}}

\def\FigB{
\begin{tikzpicture}[
outer frame sep=0ex,%
/tikz/inner frame xsep=2ex,
/tikz/inner frame ysep=2ex,
show background top,%
show background bottom,%
show background left,%
show background right]

\draw[thick] plot (0,0) -- (5/4,1+1/4);

\begin{scope}[xscale=-1, shift={(-6,0)}]
\draw[thick] plot (0,0) -- (5/4,1+1/4);
\end{scope}

\begin{scope}[xscale=-1]
\draw[thick,orange] plot (-2,0) -- (-3/4,5/4);
\end{scope}
\fill [black]   (1,1) circle (0.05);%

\begin{scope}[xscale=1]
\draw[thick,orange] plot (6/4,0) -- (11/4,5/4);
\fill [black]   (7/4,1/4) circle (0.05);%

\end{scope}

\begin{scope}[xscale=1, shift={(2+3/4,3/4)}]
\draw[thick] (1,1/3-1/12) node [scale=1] {$\dots\dots\dots$};
\end{scope}

\begin{scope}[xscale=1, shift={(4,0)}]

\draw[thick,orange] plot (0,0) -- (5/4,1+1/4);
\fill [black]   (1,1) circle (0.05);%
\end{scope}

\begin{scope}[xscale=1]
\draw[thick,green] plot (10/4,2/4) -- (5/4,12/4);
\fill [black]   (14/6,5/6) circle (0.05);%

\begin{scope}[xscale=1]
\draw[thick,green] plot (6/4,7/4) -- (12/4,10/4);
\fill [black]   (36/20,38/20) circle (0.05);%

\begin{scope}[xscale=1]
\draw[thick,blue] plot (19/8,12/4) -- (12/4,7/4);
\fill [black]   (54/20,47/20) circle (0.05);%
\end{scope}
\end{scope}

\begin{scope}[xscale=1]
\draw[thick,blue] plot (1/4,8/4) -- (7/4,11/4);
\fill [black]   (6/5+1/4,13/5) circle (0.05);%
\end{scope}
\end{scope}

\begin{scope}[shift={(6+3/4,1)}]
\draw[very thin,-stealth](0,0) -- (6/4,0);
\end{scope}

\begin{scope}[shift={(9,0)}]

\draw[thick] plot (0,0) -- (5/4,1+1/4);

\begin{scope}[xscale=-1, shift={(-2,0)}]
\draw[thick] plot (0,0) -- (5/4,1+1/4);
\end{scope}
\fill [black]   (1,1) circle (0.05);%
\end{scope}
\end{tikzpicture}

}

\def\FigC{
\begin{tikzpicture}[
outer frame sep=0ex,%
/tikz/inner frame xsep=2ex,
/tikz/inner frame ysep=3ex,
show background top,%
show background bottom,%
show background left,%
show background right]

\draw[thick] plot (0,0) -- (5/4,5/4);
\begin{scope}[xscale=-1, shift={(-12/2,0)}]
\draw[thick] plot (0,0) -- (5/4,5/4);
\end{scope}

\begin{scope}[xscale=1]
\draw[thick,orange] plot (2,0) -- (3/4,5/4);
\end{scope}
\fill [black]   (1,1) circle (0.05);%

\begin{scope}[xscale=1]
\draw[thick,orange] plot (3/2,0) -- (11/4,5/4);
\fill [black]   (7/4,1/4) circle (0.05);%
\end{scope}

\begin{scope}[xscale=1, shift={(2,2/4)}]
\draw (17/12,1/8) node [scale=1] {$\dots\dots$};
\end{scope}

\begin{scope}[xscale=1, shift={(4,0)}]
\draw[thick,orange] plot (0,0) -- (5/4,1+1/4);
\fill [black]   (1,1) circle (0.05);%
\end{scope}

\begin{scope}[shift={(6+1/4,3/4)}]
\draw [very thin,-stealth](0,0) -- (6/4,0);
\end{scope}

\begin{scope}[shift={(8,0)}]

\draw[thick] plot (0,0) -- (5/4,1+1/4);
\fill [black]   (1,1) circle (0.05);%
\begin{scope}[xscale=-1, shift={(-2,0)}]
\draw[thick] plot (0,0) -- (5/4,1+1/4);

\end{scope}
\end{scope}
\end{tikzpicture}

}

\def\FigD{
\begin{tikzpicture}[
outer frame sep=0ex,%
/tikz/inner frame xsep=2ex,
/tikz/inner frame ysep=2ex,
show background top,%
show background bottom,%
show background left,%
show background right]

\draw[thick] plot (0,0) -- (5/4,1+1/4);
\fill [black]   (1,1) circle (0.05);%
\begin{scope}[xscale=-1, shift={(-4,0)}]
\draw[thick] plot (0,0) -- (5/4,1+1/4);

\end{scope}
\begin{scope}[xscale=1, shift={(1,1)}]
\draw[thick,orange] plot (-1/2,0) -- (2+1/2,0);
\fill [black]   (2,0) circle (0.05);%
\draw (-1/2,1/3) node {\hbox{\footnotesize $\tilde E$}};
\end{scope}
\begin{scope}[xscale=1, shift={(2,1)}]
\draw[thick,blue] plot (0,-1/2) -- (0,3/2);
\fill [black]   (0,0) circle (0.05);%
\draw (1/3,5/4) node {\hbox{\footnotesize $E'$}};

\end{scope}

\begin{scope}[shift={(4+3/4,1)}]
\draw [very thin,-stealth](0,0) -- (6/4,0);
\end{scope}

\begin{scope}[shift={(7,0)}]

\draw[thick] plot (0,0) -- (5/4,1+1/4);
\fill [black]   (1,1) circle (0.05);%
\begin{scope}[xscale=-1, shift={(-2,0)}]
\draw[thick] plot (0,0) -- (5/4,1+1/4);

\end{scope}
\end{scope}
\end{tikzpicture}

}

\def\FigE{
\begin{tikzpicture}[
outer frame sep=0ex,%
/tikz/inner frame xsep=2ex,
/tikz/inner frame ysep=2ex,
show background top,%
show background bottom,%
show background left,%
show background right]

\begin{scope}[xscale=1, shift={(2,1)}]
\draw[thick]  plot (-1/2,0) -- (3+1/2,0);

\end{scope}
\begin{scope}[xscale=1, shift={(3+1/2,1)}]
\draw[thick,green] plot (0,-1/2) -- (0,2);
\fill [black]   (0,0) circle (0.05);%
\begin{scope}[xscale=1, shift={(0,3/2)}]
\draw[thick,blue] plot (-1/2,0) -- (3/2,0);
\fill [black]   (0,0) circle (0.05);%
\draw (-3/4,0) node {\hbox{\footnotesize $c_1$}};
\end{scope}
\begin{scope}[xscale=1, shift={(0,2/2)}]
\draw[thick,blue] plot (-1/2,0) -- (3/2,0);
\fill [black]   (0,0) circle (0.05);%
\draw (-3/4,0) node {\hbox{\footnotesize $c_2$}};
\end{scope}

\end{scope}

\begin{scope}[shift={(6+1/4,3/2)}]
\draw [very thin,-stealth](0,0) -- (6/4,0);
\end{scope}

\begin{scope}[shift={(17/2,1)}]

\draw[thick] plot (0,0) -- (4,0);
\fill [black]   (2,0) circle (0.05);%

\end{scope}
\end{tikzpicture}

}

\def\FigExampleW{
\begin{tikzpicture}[
outer frame sep=0ex,%
/tikz/inner frame xsep=2ex,
/tikz/inner frame ysep=2ex,
show background top,%
show background bottom,%
show background left,%
show background right]
\tikzmath{\x1 = 0.75;}

\draw[->] (0, 0) -- (6.4, 0*\x1) node[right] {$\beta$};
\draw[->] (0, 0) -- (0, 4.6*\x1) node[above] {$\wbar$};

\draw[red] (0, 0) -- (2/2, 6/2*\x1);
\draw[red] (2/2, 6/2*\x1) -- (2,4*\x1) ;
\draw[red] (2,4*\x1) -- (6,0) ;

\draw[dashed] (2/2, 6/2*\x1) -- (0, 6/2*\x1) node[left] {\hbox{\footnotesize $\frac{3}{2}$}};
\draw[dashed] (2, 4*\x1) -- (0, 4*\x1) node[left] {\hbox{\footnotesize $2$}};

\draw[dashed] (2/2, 6/2*\x1) -- (2/2, 0*\x1) node[below] {\hbox{\footnotesize $\frac{1}{2}$}};
\draw[dashed] (2, 4*\x1) -- (2, 0) node[below] {\hbox{\footnotesize $1$}};

\draw (0,-0.225) node[left] {\hbox{\footnotesize $0$}};
\draw (6,0) node[below] {\hbox{\footnotesize $3$}};

\end{tikzpicture}
}



%% file: GP22-Figures-genus-2.tex
\def\geqone{g\kern1pt{=}\kern1pt1}
\def\geqtwo{g\kern1pt{=}\kern1pt2}

\def\FigTwoAllA{\begin{tikzpicture}[scale=0.7]
\draw [thick,black] plot  coordinates {(-4/8,-1/8) (24/8,6/8) };
\begin{scope}[shift={(0,0)}]
\draw [thick,black] plot  coordinates {(-1/32,1/8) (1/32,-1/8) };
\end{scope}
\begin{scope}[shift={(1/2,1/8)}]
\draw [thick,black] plot  coordinates {(-1/32,1/8) (1/32,-1/8) };
\end{scope}
\begin{scope}[shift={(2/2,2/8)}]
\draw [thick,black] plot  coordinates {(-1/32,1/8) (1/32,-1/8) };
\end{scope}
\begin{scope}[shift={(3/2,3/8)}]
\draw [thick,black] plot  coordinates {(-1/32,1/8) (1/32,-1/8) };
\end{scope}
\begin{scope}[shift={(4/2,4/8)}]
\draw [thick,black] plot  coordinates {(-1/32,1/8) (1/32,-1/8) };
\end{scope}
\begin{scope}[shift={(5/2,5/8)}]
\draw [thick,black] plot  coordinates {(-1/32,1/8) (1/32,-1/8) };
\end{scope}
\end{tikzpicture}}

\def\FigTwoAllB{\begin{tikzpicture}[scale=0.7]
\draw [thick,black] plot  coordinates {(-48/24,-32/24) (3/8,2/8) };
\draw [thick,black] plot  coordinates {(6/8,-9/8) (-4/24,6/24) };
\fill [black] (0,0) circle (0.08);
\begin{scope}[shift={(-3/8,-2/8)}]
\draw [thick,black] plot  coordinates {(-2/32,3/32) (2/32,-3/32) };
\end{scope}
\begin{scope}[shift={(-6/8,-4/8)}]
\draw [thick,black] plot  coordinates {(-2/32,3/32) (2/32,-3/32) };
\end{scope}
\begin{scope}[shift={(-9/8,-6/8)}]
\draw [thick,black] plot  coordinates {(-2/32,3/32) (2/32,-3/32) };
\end{scope}
\begin{scope}[shift={(-12/8,-8/8)}]
\draw [thick,black] plot  coordinates {(-2/32,3/32) (2/32,-3/32) };
\end{scope}
\begin{scope}[shift={(2/8,-3/8)}]
\draw [thick,black] plot  coordinates {(3/32,2/32) (-3/32,-2/32) };
\end{scope}
\begin{scope}[shift={(4/8,-6/8)}]
\draw [thick,black] plot  coordinates {(3/32,2/32) (-3/32,-2/32) };
\end{scope}
\end{tikzpicture}}

\def\FigTwoAllE{\begin{tikzpicture}[scale=0.7]
\draw (7/4,0/4) node {};
\draw [thick,black] plot  coordinates {(-11/8,-11/8) (1/4,1/4) };
\draw [thick,black] plot  coordinates {(11/8,-11/8) (-1/4,1/4) };
\begin{scope}[shift={(0,0)}]
\fill [white] (0,0) circle (0.1);
\draw [thick] (0,0) circle (0.1);
\end{scope}
\begin{scope}[shift={(1/3,-1/3)}]
\draw [thick,black] plot  coordinates {(-1/12,-1/12) (1/12,1/12) };
\end{scope}
\begin{scope}[shift={(2/3,-2/3)}]
\draw [thick,black] plot  coordinates {(-1/12,-1/12) (1/12,1/12) };
\end{scope}
\begin{scope}[shift={(3/3,-3/3)}]
\draw [thick,black] plot  coordinates {(-1/12,-1/12) (1/12,1/12) };
\end{scope}
\begin{scope}[shift={(-1/3,-1/3)}]
\draw [thick,black] plot  coordinates {(-1/12,1/12) (1/12,-1/12) };
\end{scope}
\begin{scope}[shift={(-2/3,-2/3)}]
\draw [thick,black] plot  coordinates {(-1/12,1/12) (1/12,-1/12) };
\end{scope}
\begin{scope}[shift={(-3/3,-3/3)}]
\draw [thick,black] plot  coordinates {(-1/12,1/12) (1/12,-1/12) };
\end{scope}
\end{tikzpicture}}

\def\FigTwoAllC{\begin{tikzpicture}[scale=0.7]
\draw [thick,black] plot  coordinates {(-10/8,-10/8) (1/4,1/4) };
\draw [thick,black] plot  coordinates {(10/8,-10/8) (-1/4,1/4) };
\draw [thick,black] plot  coordinates {(6/8,-10/8) (18/8,2/8) };
\fill [black] (0,0) circle (0.08);
\fill [black] (4/4,-4/4) circle (0.08);
\begin{scope}[shift={(-1/2,-1/2)}]
\draw [thick,black] plot  coordinates {(-1/12,1/12) (1/12,-1/12) };
\end{scope}
\begin{scope}[shift={(-5/6,-5/6)}]
\draw [thick,black] plot  coordinates {(-1/12,1/12) (1/12,-1/12) };
\end{scope}
\begin{scope}[shift={(1/3,-1/3)}]
\draw [thick,black] plot  coordinates {(-1/12,-1/12) (1/12,1/12) };
\end{scope}
\begin{scope}[shift={(2/3,-2/3)}]
\draw [thick,black] plot  coordinates {(-1/12,-1/12) (1/12,1/12) };
\end{scope}
\begin{scope}[shift={(2-1/6,-1/6)}]
\draw [thick,black] plot  coordinates {(-1/12,1/12) (1/12,-1/12) };
\end{scope}
\begin{scope}[shift={(2-3/6,-3/6)}]
\draw [thick,black] plot  coordinates {(-1/12,1/12) (1/12,-1/12) };
\end{scope}
\end{tikzpicture}}

\def\FigTwoAllF{\begin{tikzpicture}[scale=0.7]
\draw [thick,black] plot  coordinates {(-10/8,-10/8) (1/4,1/4) };
\draw [thick,black] plot  coordinates {(10/8,-10/8) (-1/4,1/4) };
\draw [thick,black] plot  coordinates {(6/8,-10/8) (18/8,2/8) };
\begin{scope}[shift={(0,0)}]
\fill [white] (0,0) circle (0.1);
\draw [thick] (0,0) circle (0.1);
\end{scope}
\fill [black] (4/4,-4/4) circle (0.08);
\begin{scope}[shift={(-1/3,-1/3)}]
\draw [thick,black] plot  coordinates {(-1/12,1/12) (1/12,-1/12) };
\end{scope}
\begin{scope}[shift={(-2/3,-2/3)}]
\draw [thick,black] plot  coordinates {(-1/12,1/12) (1/12,-1/12) };
\end{scope}
\begin{scope}[shift={(-3/3,-3/3)}]
\draw [thick,black] plot  coordinates {(-1/12,1/12) (1/12,-1/12) };
\end{scope}
\begin{scope}[shift={(1/2,-1/2)}]
\draw [thick,black] plot  coordinates {(-1/12,-1/12) (1/12,1/12) };
\end{scope}
\begin{scope}[shift={(2-1/6,-1/6)}]
\draw [thick,black] plot  coordinates {(-1/12,1/12) (1/12,-1/12) };
\end{scope}
\begin{scope}[shift={(2-3/6,-3/6)}]
\draw [thick,black] plot  coordinates {(-1/12,1/12) (1/12,-1/12) };
\end{scope}
\end{tikzpicture}}

\def\FigTwoAllD{\begin{tikzpicture}[scale=0.7]
\draw [thick,black] plot  coordinates {(4/8,-1/8) (-20/8,5/8) };
\begin{scope}[shift={(0,0)}]
\draw [thick,black] plot  coordinates {(2/24,8/24) (-5/16,-20/16) };
\fill [black] (0,0) circle (0.08);
\begin{scope}[shift={(-3/32,-12/32)}]
\draw [thick,black] plot  coordinates {(-1/8,1/32) (1/8,-1/32) };
\end{scope}
\begin{scope}[shift={(-6/32,-24/32)}]
\draw [thick,black] plot  coordinates {(-1/8,1/32) (1/8,-1/32) };
\end{scope}
\end{scope}
\begin{scope}[shift={(-4/4,1/4)}]
\draw [thick,black] plot  coordinates {(2/24,8/24) (-5/16,-20/16) };
\fill [black] (0,0) circle (0.08);
\begin{scope}[shift={(-3/32,-12/32)}]
\draw [thick,black] plot  coordinates {(-1/8,1/32) (1/8,-1/32) };
\end{scope}
\begin{scope}[shift={(-6/32,-24/32)}]
\draw [thick,black] plot  coordinates {(-1/8,1/32) (1/8,-1/32) };
\end{scope}
\end{scope}
\begin{scope}[shift={(-8/4,2/4)}]
\draw [thick,black] plot  coordinates {(2/24,8/24) (-5/16,-20/16) };
\fill [black] (0,0) circle (0.08);
\begin{scope}[shift={(-3/32,-12/32)}]
\draw [thick,black] plot  coordinates {(-1/8,1/32) (1/8,-1/32) };
\end{scope}
\begin{scope}[shift={(-6/32,-24/32)}]
\draw [thick,black] plot  coordinates {(-1/8,1/32) (1/8,-1/32) };
\end{scope}
\end{scope}
\end{tikzpicture}}

\def\FigTwoAllG{\begin{tikzpicture}[scale=0.7]
\draw [thick,black] plot  coordinates {(-10/8,-10/8) (1/4,1/4) };
\draw [thick,black] plot  coordinates {(8/8,-8/8) (-1/4,1/4) };
\draw [thick,black] plot  coordinates {(4/8,-8/8) (14/8,2/8) };
\draw [thick,black] plot  coordinates {(22/8,-10/8) (5/4,1/4) };
\fill [black] (0,0) circle (0.08);
\begin{scope}[shift={(3/4,-3/4)}]
\fill [white] (0,0) circle (0.1);
\draw [thick] (0,0) circle (0.1);
\end{scope}
\fill [black] (6/4,0) circle (0.08);
\begin{scope}[shift={(-1/2,-1/2)}]
\draw [thick,black] plot  coordinates {(-1/12,1/12) (1/12,-1/12) };
\end{scope}
\begin{scope}[shift={(-5/6,-5/6)}]
\draw [thick,black] plot  coordinates {(-1/12,1/12) (1/12,-1/12) };
\end{scope}
\begin{scope}[shift={(3/8,-3/8)}]
\draw [thick,black] plot  coordinates {(-1/12,-1/12) (1/12,1/12) };
\end{scope}
\begin{scope}[shift={(3/2-3/8,-3/8)}]
\draw [thick,black] plot  coordinates {(-1/12,1/12) (1/12,-1/12) };
\end{scope}
\begin{scope}[shift={(2,-1/2)}]
\draw [thick,black] plot  coordinates {(-1/12,-1/12) (1/12,1/12) };
\end{scope}
\begin{scope}[shift={(2+2/6,-5/6)}]
\draw [thick,black] plot  coordinates {(-1/12,-1/12) (1/12,1/12) };
\end{scope}
\end{tikzpicture}}

\def\includefig#1#2{
\vcenter{\hbox{\makebox(110pt,55pt){
\fbox{\makebox(0pt,45pt)[lt]{\makebox(20pt,20pt){(\small #1)}}
\makebox(90pt,45pt){#2}}}}}}

\def\FigTwoAll{
$\xymatrix@C+0pt@R+0pt{
\includefig{A}{\FigTwoAllA} \ar[d] \ar[dr] \\
\includefig{B}{\FigTwoAllB} \ar[d] \ar[dr] & \includefig{E}{\FigTwoAllE} \ar[d] \\
\includefig{C}{\FigTwoAllC}  \ar[d] \ar[dr] & \includefig{F}{\FigTwoAllF} \ar[d] \\
\includefig{D}{\hskip15pt \FigTwoAllD} & \includefig{G}{\FigTwoAllG} }$}

\def\FigTwoCaseAa{\begin{tikzpicture}[baseline=25pt]
\draw (7/4,0/4) node {};
\draw [thick,black] plot  coordinates {(-15/8,-15/8) (3/4,3/4) };
\draw [thick,blue] plot  coordinates {(12/8,-10/8) (-1/8,3/8) };
\fill [black] (1/8,1/8) circle (0.06);
\draw (11/8-2/3,-11/8+2/4) node [below] {\hbox{\footnotesize$\geqone$}};
\begin{scope}[shift={(2/4,2/4)}]
\draw [thick,blue] plot  coordinates {(11/8,-11/8) (-1/4,1/4) };
\fill [black] (0,0) circle (0.06);
\draw (11/8,-11/8+5/12) node [above] {\hbox{\footnotesize$\geqone$}};
\end{scope}
\begin{scope}[shift={(-1/4,-1/4)}]
\draw [thick,black] plot  coordinates {(-1/12,1/12) (1/12,-1/12) };
\end{scope}
\begin{scope}[shift={(-2/4,-2/4)}]
\draw [thick,black] plot  coordinates {(-1/12,1/12) (1/12,-1/12) };
\end{scope}
\begin{scope}[shift={(-3/4,-3/4)}]
\draw [thick,black] plot  coordinates {(-1/12,1/12) (1/12,-1/12) };
\end{scope}
\begin{scope}[shift={(-4/4,-4/4)}]
\draw [thick,black] plot  coordinates {(-1/12,1/12) (1/12,-1/12) };
\end{scope}
\begin{scope}[shift={(-5/4,-5/4)}]
\draw [thick,black] plot  coordinates {(-1/12,1/12) (1/12,-1/12) };
\end{scope}
\begin{scope}[shift={(-6/4,-6/4)}]
\draw [thick,black] plot  coordinates {(-1/12,1/12) (1/12,-1/12) };
\end{scope}
\end{tikzpicture}}

\def\FigTwoCaseAb{\begin{tikzpicture}[baseline=23pt]
\draw (7/4,0/4) node {};
\draw [thick,black] plot  coordinates {(-14/8,-14/8) (4/8,4/8) };
\draw [thick,green] plot  coordinates {(15/8,-11/8) (0,2/4) };
\fill [black] (1/4,1/4) circle (0.06);
\begin{scope}[shift={(3/4,-3/4)}]
\draw [thick,blue] plot  coordinates {(-5/8,-5/8) (2/4,2/4) };
\fill [black] (1/4,1/4) circle (0.06);
\draw (-5/8-0.20,-5/8) node [below] {\hbox{\footnotesize$\geqone$}};
\end{scope}
\begin{scope}[shift={(5/4,-5/4)}]
\draw [thick,blue] plot  coordinates {(-5/8,-5/8) (2/4,2/4) };
\fill [black] (1/4,1/4) circle (0.06);
\draw (-5/8+0.05,-5/8) node [right] {\hbox{\footnotesize$\geqone$}};
\end{scope}
\begin{scope}[shift={(1/8,1/8)}]
\begin{scope}[shift={(-1/4,-1/4)}]
\draw [thick,black] plot  coordinates {(-1/12,1/12) (1/12,-1/12) };
\end{scope}
\begin{scope}[shift={(-2/4,-2/4)}]
\draw [thick,black] plot  coordinates {(-1/12,1/12) (1/12,-1/12) };
\end{scope}
\begin{scope}[shift={(-3/4,-3/4)}]
\draw [thick,black] plot  coordinates {(-1/12,1/12) (1/12,-1/12) };
\end{scope}
\begin{scope}[shift={(-4/4,-4/4)}]
\draw [thick,black] plot  coordinates {(-1/12,1/12) (1/12,-1/12) };
\end{scope}
\begin{scope}[shift={(-5/4,-5/4)}]
\draw [thick,black] plot  coordinates {(-1/12,1/12) (1/12,-1/12) };
\end{scope}
\begin{scope}[shift={(-6/4,-6/4)}]
\draw [thick,black] plot  coordinates {(-1/12,1/12) (1/12,-1/12) };
\end{scope}
\end{scope}
\end{tikzpicture}}

\def\FigTwoCaseAc{\begin{tikzpicture}[baseline=20pt]
\draw (7/4,0/4) node {};
\draw [thick,black] plot  coordinates {(-14/8,-14/8) (4/8,4/8) };
\begin{scope}[shift={(1/4,1/4)}]
\draw [thick,blue] plot  coordinates {(11/8,-11/8) (-1/4,1/4) };
\fill [black] (0,0) circle (0.06);
\draw (11/8-2/3,-11/8+1/3) node [below] {\hbox{\footnotesize$\geqtwo$}};
\end{scope}
\begin{scope}[shift={(1/8,1/8)}]
\begin{scope}[shift={(-1/4,-1/4)}]
\draw [thick,black] plot  coordinates {(-1/12,1/12) (1/12,-1/12) };
\end{scope}
\begin{scope}[shift={(-2/4,-2/4)}]
\draw [thick,black] plot  coordinates {(-1/12,1/12) (1/12,-1/12) };
\end{scope}
\begin{scope}[shift={(-3/4,-3/4)}]
\draw [thick,black] plot  coordinates {(-1/12,1/12) (1/12,-1/12) };
\end{scope}
\begin{scope}[shift={(-4/4,-4/4)}]
\draw [thick,black] plot  coordinates {(-1/12,1/12) (1/12,-1/12) };
\end{scope}
\begin{scope}[shift={(-5/4,-5/4)}]
\draw [thick,black] plot  coordinates {(-1/12,1/12) (1/12,-1/12) };
\end{scope}
\begin{scope}[shift={(-6/4,-6/4)}]
\draw [thick,black] plot  coordinates {(-1/12,1/12) (1/12,-1/12) };
\end{scope}
\end{scope}
\end{tikzpicture}}

\def\includefigA#1{
\makebox(130pt,80pt){#1}}

\def\FigTwoCaseA{
\begin{tabular}{|c|c|c|}
	\hline
{\large\strut}	(A1) & (A2) & (A3) \\
\hline
{\Large\strut} $\delta=0$ & $0<\delta<4/5$ & $\delta\ge4/5$ \\
\hline
\includefigA{\FigTwoCaseAa} & 
\includefigA{\FigTwoCaseAb} & 
\includefigA{\FigTwoCaseAc} \\
\hline
\end{tabular}}

\def\FigTwoCaseBi{\begin{tikzpicture}[baseline=13pt]
\draw [thick,black] plot  coordinates {(-16/12,-16/8) (1/6,1/4) };
\begin{scope}[shift={(-2/9,-3/9)}]
\draw [thick,black] plot  coordinates {(-3/32,2/32) (3/32,-2/32) };
\end{scope}
\begin{scope}[shift={(-4/9,-6/9)}]
\draw [thick,black] plot  coordinates {(-3/32,2/32) (3/32,-2/32) };
\end{scope}
\begin{scope}[shift={(-6/9,-9/9)}]
\draw [thick,black] plot  coordinates {(-3/32,2/32) (3/32,-2/32) };
\end{scope}
\begin{scope}[shift={(-8/9,-12/9)}]
\draw [thick,black] plot  coordinates {(-3/32,2/32) (3/32,-2/32) };
\end{scope}
\draw [thick,orange] plot  coordinates {(-1/2,0) (4/2,0) };
\fill [black] (0,0) circle (0.06);
\draw [thick,blue] plot  coordinates {(3/4,1/3) (3/4,-6/4) };
\draw (3/4,-3/2) node [right] {\hbox{\footnotesize$\geqone$}};
\fill [black] (3/4,0) circle (0.06);
\begin{scope}[shift={(3/2,0)}]
\draw [thick,black] plot  coordinates {(8/12,-8/8) (-1/6,1/4) };
\fill [black] (0,0) circle (0.06);
\begin{scope}[shift={(2/9,-3/9)}]
\draw [thick,black] plot  coordinates {(3/32,2/32) (-3/32,-2/32) };
\end{scope}
\begin{scope}[shift={(4/9,-6/9)}]
\draw [thick,black] plot  coordinates {(3/32,2/32) (-3/32,-2/32) };
\end{scope}
\end{scope}
\draw [thick,blue] plot  coordinates {(-10/9-5/6,-15/9+5/4) (-10/9+1/6,-15/9-1/4) };
\fill [black] (-10/9,-15/9) circle (0.06);
\draw (-10/9-5/6,-15/9+5/4+1/6) node [right] {\hbox{\footnotesize$\geqone$}};
\end{tikzpicture}}

\def\FigTwoCaseBii{\begin{tikzpicture}[baseline=13pt]
\draw [thick,black] plot  coordinates {(-16/12,-16/8) (1/6,1/4) };
\begin{scope}[shift={(-2/9,-3/9)}]
\draw [thick,black] plot  coordinates {(-3/32,2/32) (3/32,-2/32) };
\end{scope}
\begin{scope}[shift={(-4/9,-6/9)}]
\draw [thick,black] plot  coordinates {(-3/32,2/32) (3/32,-2/32) };
\end{scope}
\begin{scope}[shift={(-6/9,-9/9)}]
\draw [thick,black] plot  coordinates {(-3/32,2/32) (3/32,-2/32) };
\end{scope}
\begin{scope}[shift={(-8/9,-12/9)}]
\draw [thick,black] plot  coordinates {(-3/32,2/32) (3/32,-2/32) };
\end{scope}
\draw [thick,orange] plot  coordinates {(-1/2,0) (4/2,0) };
\draw (3/4,0) node [below] {\hbox{\footnotesize$\geqone$}};
\fill [black] (0,0) circle (0.06);
\begin{scope}[shift={(3/2,0)}]
\draw [thick,black] plot  coordinates {(8/12,-8/8) (-1/6,1/4) };
\fill [black] (0,0) circle (0.06);
\begin{scope}[shift={(2/9,-3/9)}]
\draw [thick,black] plot  coordinates {(3/32,2/32) (-3/32,-2/32) };
\end{scope}
\begin{scope}[shift={(4/9,-6/9)}]
\draw [thick,black] plot  coordinates {(3/32,2/32) (-3/32,-2/32) };
\end{scope}
\end{scope}
\draw [thick,blue] plot  coordinates {(-10/9-5/6,-15/9+5/4) (-10/9+1/6,-15/9-1/4) };
\fill [black] (-10/9,-15/9) circle (0.06);
\draw (-10/9-5/6,-15/9+5/4+1/6) node [right] {\hbox{\footnotesize$\geqone$}};
\end{tikzpicture}}

\def\FigTwoCaseBiii{\begin{tikzpicture}[baseline=13pt]
\draw [thick,black] plot  coordinates {(-16/12,-16/8) (1/6,1/4) };
\begin{scope}[shift={(-2/9,-3/9)}]
\draw [thick,black] plot  coordinates {(-3/32,2/32) (3/32,-2/32) };
\end{scope}
\begin{scope}[shift={(-4/9,-6/9)}]
\draw [thick,black] plot  coordinates {(-3/32,2/32) (3/32,-2/32) };
\end{scope}
\begin{scope}[shift={(-6/9,-9/9)}]
\draw [thick,black] plot  coordinates {(-3/32,2/32) (3/32,-2/32) };
\end{scope}
\begin{scope}[shift={(-8/9,-12/9)}]
\draw [thick,black] plot  coordinates {(-3/32,2/32) (3/32,-2/32) };
\end{scope}
\begin{scope}[shift={(3/4,0)}]
\draw [thick,orange] plot [smooth, tension=1] coordinates {(-3/3,2/12) (1/4,-6/12) (-1/3,-14/12) };
\draw [thick,orange] plot [smooth, tension=1] coordinates {(3/3,2/12) (-1/4,-6/12) (1/3,-14/12) };
\fill [black] (0,-0.26) circle (0.06);
\fill [black] (0,-0.99) circle (0.06);
\fill [black] (-3/4+0.14/3,0.14/2) circle (0.06);
\fill [black] (3/4-0.14/3,0.14/2) circle (0.06);
\end{scope}
\begin{scope}[shift={(3/2,0)}]
\draw [thick,black] plot  coordinates {(8/12,-8/8) (-1/6,1/4) };
\begin{scope}[shift={(2/9,-3/9)}]
\draw [thick,black] plot  coordinates {(3/32,2/32) (-3/32,-2/32) };
\end{scope}
\begin{scope}[shift={(4/9,-6/9)}]
\draw [thick,black] plot  coordinates {(3/32,2/32) (-3/32,-2/32) };
\end{scope}
\end{scope}
\draw [thick,blue] plot  coordinates {(-10/9-5/6,-15/9+5/4) (-10/9+1/6,-15/9-1/4) };
\fill [black] (-10/9,-15/9) circle (0.06);
\draw (-10/9-5/6,-15/9+5/4+1/6) node [right] {\hbox{\footnotesize$\geqone$}};
\end{tikzpicture}}

\def\FigTwoCaseBvii{\begin{tikzpicture}[baseline=13pt]
\draw [thick,black] plot  coordinates {(-14/12,-14/8) (1/6,1/4) };
\begin{scope}[shift={(-2/9,-3/9)}]
\draw [thick,black] plot  coordinates {(-3/32,2/32) (3/32,-2/32) };
\end{scope}
\begin{scope}[shift={(-4/9,-6/9)}]
\draw [thick,black] plot  coordinates {(-3/32,2/32) (3/32,-2/32) };
\end{scope}
\begin{scope}[shift={(-6/9,-9/9)}]
\draw [thick,black] plot  coordinates {(-3/32,2/32) (3/32,-2/32) };
\end{scope}
\begin{scope}[shift={(-8/9,-12/9)}]
\draw [thick,black] plot  coordinates {(-3/32,2/32) (3/32,-2/32) };
\end{scope}
\begin{scope}[shift={(3/4,0)}]
\draw [thick,orange] plot [smooth, tension=1] coordinates {(-3/3,2/12) (1,-8/12) (-1/2,-20/12) };
\draw [thick,orange] plot [smooth, tension=.8] coordinates {(1/2+3/3,2/12) (5/12,-8/12) (11/12,-16/12) };
\fill [black] (0.69,-0.36) circle (0.06);
\fill [black] (0.66,-1.16) circle (0.06);
\draw (-1/2,-20/12) node [above] {\hbox{\footnotesize$\geqone$}};
\end{scope}
\fill [black] (0.195/3,0.195/2) circle (0.06);
\begin{scope}[shift={(2,0)}]
\draw [thick,black] plot  coordinates {(8/12,-8/8) (-1/6,1/4) };
\fill [black] (-0.035/3,0.035/2) circle (0.06);
\begin{scope}[shift={(2/9,-3/9)}]
\draw [thick,black] plot  coordinates {(3/32,2/32) (-3/32,-2/32) };
\end{scope}
\begin{scope}[shift={(4/9,-6/9)}]
\draw [thick,black] plot  coordinates {(3/32,2/32) (-3/32,-2/32) };
\end{scope}
\end{scope}
\end{tikzpicture}}

\def\FigTwoCaseBvi{\begin{tikzpicture}[baseline=13pt]
\draw [thick,black] plot  coordinates {(-1-14/12,-14/8) (-1+1/6,1/4) };
\begin{scope}[shift={(-3/36,-4/36)}]
\begin{scope}[shift={(-1-2/9,-3/9)}]
\draw [thick,black] plot  coordinates {(-3/32,2/32) (3/32,-2/32) };
\end{scope}
\begin{scope}[shift={(-1-4/9,-6/9)}]
\draw [thick,black] plot  coordinates {(-3/32,2/32) (3/32,-2/32) };
\end{scope}
\begin{scope}[shift={(-1-6/9,-9/9)}]
\draw [thick,black] plot  coordinates {(-3/32,2/32) (3/32,-2/32) };
\end{scope}
\begin{scope}[shift={(-1-8/9,-12/9)}]
\draw [thick,black] plot  coordinates {(-3/32,2/32) (3/32,-2/32) };
\end{scope}
\end{scope}
\draw [thick,orange] plot [smooth, tension=1] coordinates {(4/12,1/4) (-6/12,-1/4) (-16/12,1/4) };
\begin{scope}[shift={(3/4,0)}]
\draw [thick,orange] plot [smooth, tension=1] coordinates {(-3/3,2/12) (1/4,-6/12) (-1/3,-14/12) };
\draw [thick,orange] plot [smooth, tension=1] coordinates {(3/3,2/12) (-1/4,-6/12) (1/3,-14/12) };
\fill [black] (0,-0.26) circle (0.06);
\fill [black] (0,-0.99) circle (0.06);
\fill [black] (3/4-0.14/3,0.14/2) circle (0.06);
\end{scope}
\fill [black] (-5/6-0.6/3,1/4-0.6/2) circle (0.06);
\fill [black] (0.14,0.036) circle (0.06);
\begin{scope}[shift={(3/2,0)}]
\draw [thick,black] plot  coordinates {(8/12,-8/8) (-1/6,1/4) };
\begin{scope}[shift={(2/9,-3/9)}]
\draw [thick,black] plot  coordinates {(3/32,2/32) (-3/32,-2/32) };
\end{scope}
\begin{scope}[shift={(4/9,-6/9)}]
\draw [thick,black] plot  coordinates {(3/32,2/32) (-3/32,-2/32) };
\end{scope}
\end{scope}
\draw [thick,blue] plot  coordinates {(-6/12,-15/9+7/4) (-6/12,-15/9+1/4) };
\fill [black] (-6/12,-.25) circle (0.06);
\draw (-6/6,-10/6) node [right] {\hbox{\footnotesize$\geqone$}};
\end{tikzpicture}}

\def\FigTwoCaseBiv{\begin{tikzpicture}[baseline=13pt]
\draw [thick,black] plot  coordinates {(-14/12,-14/8) (1/6,1/4) };
\begin{scope}[shift={(-2/9,-3/9)}]
\draw [thick,black] plot  coordinates {(-3/32,2/32) (3/32,-2/32) };
\end{scope}
\begin{scope}[shift={(-4/9,-6/9)}]
\draw [thick,black] plot  coordinates {(-3/32,2/32) (3/32,-2/32) };
\end{scope}
\begin{scope}[shift={(-6/9,-9/9)}]
\draw [thick,black] plot  coordinates {(-3/32,2/32) (3/32,-2/32) };
\end{scope}
\begin{scope}[shift={(-8/9,-12/9)}]
\draw [thick,black] plot  coordinates {(-3/32,2/32) (3/32,-2/32) };
\end{scope}
\draw [thick,orange] plot  coordinates {(-1/4,1/6) (6/4,-6/6) };
\draw [thick,orange] plot  coordinates {(5/2+1/4,1/6) (5/2-6/4,-6/6) };
\fill [black] (0,0) circle (0.06);
\fill [black] (5/4,-5/6) circle (0.06);
\begin{scope}[shift={(3/5,-2/5)}]
\draw [thick,blue] plot coordinates {(1/12,1/3) (-7/24,-7/6) };
\fill [black] (0,0) circle (0.06);
\draw (-6.5/24,-6.5/6) node [below] {\hbox{\footnotesize$\geqone$}};
\end{scope}
\begin{scope}[shift={(5/2-3/5,-2/5)}]
\draw [thick,blue] plot coordinates {(-1/12,1/3) (7/24,-7/6) };
\fill [black] (0,0) circle (0.06);
\draw (6.5/24,-6.5/6) node [below] {\hbox{\footnotesize$\geqone$}};
\end{scope}
\begin{scope}[shift={(5/2,0)}]
\draw [thick,black] plot  coordinates {(8/12,-8/8) (-1/6,1/4) };
\fill [black] (0,0) circle (0.06);
\begin{scope}[shift={(2/9,-3/9)}]
\draw [thick,black] plot  coordinates {(3/32,2/32) (-3/32,-2/32) };
\end{scope}
\begin{scope}[shift={(4/9,-6/9)}]
\draw [thick,black] plot  coordinates {(3/32,2/32) (-3/32,-2/32) };
\end{scope}
\end{scope}
\end{tikzpicture}}

\def\FigTwoCaseBv{\begin{tikzpicture}[baseline=13pt]
\draw [thick,black] plot  coordinates {(-14/12,-14/8) (1/6,1/4) };
\begin{scope}[shift={(-2/9,-3/9)}]
\draw [thick,black] plot  coordinates {(-3/32,2/32) (3/32,-2/32) };
\end{scope}
\begin{scope}[shift={(-4/9,-6/9)}]
\draw [thick,black] plot  coordinates {(-3/32,2/32) (3/32,-2/32) };
\end{scope}
\begin{scope}[shift={(-6/9,-9/9)}]
\draw [thick,black] plot  coordinates {(-3/32,2/32) (3/32,-2/32) };
\end{scope}
\begin{scope}[shift={(-8/9,-12/9)}]
\draw [thick,black] plot  coordinates {(-3/32,2/32) (3/32,-2/32) };
\end{scope}
\draw [thick,orange] plot  coordinates {(-1/4,1/6) (6/4,-6/6) };
\draw [thick,orange] plot  coordinates {(5/2+1/4,1/6) (5/2-6/4,-6/6) };
\fill [black] (0,0) circle (0.06);
\fill [black] (5/4,-5/6) circle (0.06);
\begin{scope}[shift={(3/5,-2/5)}]
\draw [thick,blue] plot coordinates {(1/12,1/3) (-7/24,-7/6) };
\fill [black] (0,0) circle (0.06);
\draw (-6.5/24,-6.5/6) node [below] {\hbox{\footnotesize$\geqone$}};
\end{scope}
\draw (1+4/4+1/12,-1+4/6) node [below] {\hbox{\footnotesize$\geqone$}};
\begin{scope}[shift={(5/2,0)}]
\draw [thick,black] plot  coordinates {(8/12,-8/8) (-1/6,1/4) };
\fill [black] (0,0) circle (0.06);
\begin{scope}[shift={(2/9,-3/9)}]
\draw [thick,black] plot  coordinates {(3/32,2/32) (-3/32,-2/32) };
\end{scope}
\begin{scope}[shift={(4/9,-6/9)}]
\draw [thick,black] plot  coordinates {(3/32,2/32) (-3/32,-2/32) };
\end{scope}
\end{scope}
\end{tikzpicture}}

\def\FigTwoCaseBviii{\begin{tikzpicture}[baseline=15pt]
\draw [thick,black] plot  coordinates {(-14/12,-14/8) (1/6,1/4) };
\begin{scope}[shift={(-2/9,-3/9)}]
\draw [thick,black] plot  coordinates {(-3/32,2/32) (3/32,-2/32) };
\end{scope}
\begin{scope}[shift={(-4/9,-6/9)}]
\draw [thick,black] plot  coordinates {(-3/32,2/32) (3/32,-2/32) };
\end{scope}
\begin{scope}[shift={(-6/9,-9/9)}]
\draw [thick,black] plot  coordinates {(-3/32,2/32) (3/32,-2/32) };
\end{scope}
\begin{scope}[shift={(-8/9,-12/9)}]
\draw [thick,black] plot  coordinates {(-3/32,2/32) (3/32,-2/32) };
\end{scope}
\draw [thick,orange] plot  coordinates {(-1/2,0) (5/2,0) };
\fill [black] (0,0) circle (0.06);
\begin{scope}[shift={(4/2,0)}]
\draw [thick,black] plot  coordinates {(8/12,-8/8) (-1/6,1/4) };
\fill [black] (0,0) circle (0.06);
\begin{scope}[shift={(2/9,-3/9)}]
\draw [thick,black] plot  coordinates {(3/32,2/32) (-3/32,-2/32) };
\end{scope}
\begin{scope}[shift={(4/9,-6/9)}]
\draw [thick,black] plot  coordinates {(3/32,2/32) (-3/32,-2/32) };
\end{scope}
\end{scope}
\begin{scope}[shift={(2/3,0)}]
\draw [thick,blue] plot  coordinates {(0,1/3) (0,-6/4) };
\fill [black] (0,0) circle (0.06);
\draw (0,-6/4) node [left] {\hbox{\footnotesize$\geqone$\kern-1pt}};
\end{scope}
\begin{scope}[shift={(4/3,0)}]
\draw [thick,blue] plot  coordinates {(0,1/3) (0,-6/4) };
\fill [black] (0,0) circle (0.06);
\draw (0,-6/4) node [right] {\hbox{\footnotesize$\geqone$}};
\end{scope}
\end{tikzpicture}}

\def\FigTwoCaseBix{\begin{tikzpicture}[baseline=15pt]
\draw [thick,black] plot  coordinates {(-14/12,-14/8) (1/6,1/4) };
\begin{scope}[shift={(-2/9,-3/9)}]
\draw [thick,black] plot  coordinates {(-3/32,2/32) (3/32,-2/32) };
\end{scope}
\begin{scope}[shift={(-4/9,-6/9)}]
\draw [thick,black] plot  coordinates {(-3/32,2/32) (3/32,-2/32) };
\end{scope}
\begin{scope}[shift={(-6/9,-9/9)}]
\draw [thick,black] plot  coordinates {(-3/32,2/32) (3/32,-2/32) };
\end{scope}
\begin{scope}[shift={(-8/9,-12/9)}]
\draw [thick,black] plot  coordinates {(-3/32,2/32) (3/32,-2/32) };
\end{scope}
\draw [thick,orange] plot  coordinates {(-1/2,0) (5/2,0) };
\fill [black] (0,0) circle (0.06);
\begin{scope}[shift={(4/2,0)}]
\draw [thick,black] plot  coordinates {(8/12,-8/8) (-1/6,1/4) };
\fill [black] (0,0) circle (0.06);
\begin{scope}[shift={(2/9,-3/9)}]
\draw [thick,black] plot  coordinates {(3/32,2/32) (-3/32,-2/32) };
\end{scope}
\begin{scope}[shift={(4/9,-6/9)}]
\draw [thick,black] plot  coordinates {(3/32,2/32) (-3/32,-2/32) };
\end{scope}
\end{scope}
\begin{scope}[shift={(1/2,0)}]
\draw [thick,green] plot  coordinates {(0,1/3) (0,-7/4) };
\fill [black] (0,0) circle (0.06);
\end{scope}
\begin{scope}[shift={(1/2,-2/3)}]
\draw [thick,blue] plot  coordinates {(-1/3,0) (4/3,0)};
\fill [black] (0,0) circle (0.06);
\draw (2/2,0) node [below] {\hbox{\footnotesize$\geqone$\kern-1pt}};
\end{scope}
\begin{scope}[shift={(1/2,-4/3)}]
\draw [thick,blue] plot  coordinates {(-1/3,0) (4/3,0)};
\fill [black] (0,0) circle (0.06);
\draw (2/2,0) node [below] {\hbox{\footnotesize$\geqone$\kern-1pt}};
\end{scope}
\end{tikzpicture}}

\def\FigTwoCaseBx{\begin{tikzpicture}[baseline=15pt]
\draw [thick,black] plot  coordinates {(-14/12,-14/8) (1/6,1/4) };
\begin{scope}[shift={(-2/9,-3/9)}]
\draw [thick,black] plot  coordinates {(-3/32,2/32) (3/32,-2/32) };
\end{scope}
\begin{scope}[shift={(-4/9,-6/9)}]
\draw [thick,black] plot  coordinates {(-3/32,2/32) (3/32,-2/32) };
\end{scope}
\begin{scope}[shift={(-6/9,-9/9)}]
\draw [thick,black] plot  coordinates {(-3/32,2/32) (3/32,-2/32) };
\end{scope}
\begin{scope}[shift={(-8/9,-12/9)}]
\draw [thick,black] plot  coordinates {(-3/32,2/32) (3/32,-2/32) };
\end{scope}
\draw [thick,orange] plot  coordinates {(-1/2,0) (5/2,0) };
\fill [black] (0,0) circle (0.06);
\begin{scope}[shift={(4/2,0)}]
\draw [thick,black] plot  coordinates {(8/12,-8/8) (-1/6,1/4) };
\fill [black] (0,0) circle (0.06);
\begin{scope}[shift={(2/9,-3/9)}]
\draw [thick,black] plot  coordinates {(3/32,2/32) (-3/32,-2/32) };
\end{scope}
\begin{scope}[shift={(4/9,-6/9)}]
\draw [thick,black] plot  coordinates {(3/32,2/32) (-3/32,-2/32) };
\end{scope}
\end{scope}
\begin{scope}[shift={(1,0)}]
\draw [thick,blue] plot  coordinates {(0,1/3) (0,-6/4) };
\fill [black] (0,0) circle (0.06);
\draw (0,-6/4) node [right] {\hbox{\footnotesize$\geqtwo$}};
\end{scope}
\end{tikzpicture}}

\def\FigTwoCaseBxi{\begin{tikzpicture}[baseline=15pt]
\draw [thick,black] plot  coordinates {(-14/12,-14/8) (1/6,1/4) };
\begin{scope}[shift={(-2/9,-3/9)}]
\draw [thick,black] plot  coordinates {(-3/32,2/32) (3/32,-2/32) };
\end{scope}
\begin{scope}[shift={(-4/9,-6/9)}]
\draw [thick,black] plot  coordinates {(-3/32,2/32) (3/32,-2/32) };
\end{scope}
\begin{scope}[shift={(-6/9,-9/9)}]
\draw [thick,black] plot  coordinates {(-3/32,2/32) (3/32,-2/32) };
\end{scope}
\begin{scope}[shift={(-8/9,-12/9)}]
\draw [thick,black] plot  coordinates {(-3/32,2/32) (3/32,-2/32) };
\end{scope}
\draw [thick,orange] plot  coordinates {(-1/2,0) (5/2,0) };
\fill [black] (0,0) circle (0.06);
\begin{scope}[shift={(4/2,0)}]
\draw [thick,black] plot  coordinates {(8/12,-8/8) (-1/6,1/4) };
\fill [black] (0,0) circle (0.06);
\begin{scope}[shift={(2/9,-3/9)}]
\draw [thick,black] plot  coordinates {(3/32,2/32) (-3/32,-2/32) };
\end{scope}
\begin{scope}[shift={(4/9,-6/9)}]
\draw [thick,black] plot  coordinates {(3/32,2/32) (-3/32,-2/32) };
\end{scope}
\end{scope}
\draw (1,0) node [below] {\hbox{\footnotesize$\geqtwo$}};
\end{tikzpicture}}

\def\includefigB#1{
\makebox(130pt,75pt){#1}}

\def\FigTwoCaseB{
\begin{tabular}{|c|c|c|}
\hline
{\large\strut} (B1) & (B2) & (B3) \\
\hline
${\Large\strut}
\delta=0\ \land \ \alpha<4$ 
& $\delta=0\ \land \ \alpha=4$ 
& $\delta=0\ \land \ \alpha>4$ \\
\hline
\includefigB{\FigTwoCaseBi} & 
\includefigB{\FigTwoCaseBii} & 
\includefigB{\FigTwoCaseBiii} \\
\hline\hline
{\large\strut} (B4) & (B5) & (B6) \\
\hline
${\Large\strut}
0<2\delta<\alpha \ \land \ \alpha+\delta<4$ 
& $0<2\delta<\alpha\ \land \ \alpha+\delta=4$ 
& $0<\delta<\frac{4}{3}\ \land \ \alpha+\delta>4$ \\
\hline
\includefigB{\FigTwoCaseBiv} &
\includefigB{\FigTwoCaseBv} & 
\includefigB{\FigTwoCaseBvi} \\
\hline\hline
{\large\strut} (B7)& (B8)& (B9) \\ 
\hline
{\Large\strut} 
  $\delta\ge\frac{4}{3}\ \land \ \alpha>\frac{8}{3}$ 
&  $\alpha=2\delta<\frac{8}{3}$
& $\alpha<2\delta<\frac{8+2\alpha}{5}$ 
\\
\hline
\includefigB{\FigTwoCaseBvii} &
\includefigB{\FigTwoCaseBviii} &
\includefigB{\FigTwoCaseBix} 
 \\ 
\hline 
 \noalign{\global\arrayrulewidth=2.0pt}
\arrayrulecolor{white}\hline
 \noalign{\global\arrayrulewidth=0.4pt}
 \arrayrulecolor{black} \cline{1-2}
{\large\strut} (B10) & (B11) & \multicolumn{1}{c}{}  \\
\cline{1-2}
{\Large\strut}
$\alpha<\frac{8+2\alpha}{5}\leq 2 \delta$ 
& $\alpha=\frac{8}{3} \leq 2 \delta$ 
& \multicolumn{1}{c}{}  \\
\cline{1-2}
\includefigB{\FigTwoCaseBx} &
\includefigB{\FigTwoCaseBxi} & 
\multicolumn{1}{c}{} \\
\cline{1-2}

\end{tabular}}

\def\FigTwoCaseCi{\begin{tikzpicture}[baseline=13pt]
		\begin{scope}[shift={(3/4+1/2,1/2)}]
			\draw [thick, orange] plot coordinates {(-3/4+0.14/3-3/8,0.14/2) (3/4-0.14/3+3/8,0.14/2) };
			\draw (0,-7/8) node [below] {\hbox{\footnotesize$\geqone$}};
			\draw [thick, blue] plot coordinates {(0,0.14/2+3/8) (0,0.14/2-8/8) } ;
			\fill [black] (-3/4+0.14/3,0.14/2) circle (0.06);
			\fill [black] (3/4-0.14/3,0.14/2) circle (0.06);
			\fill [black] (0,0.14/2) circle (0.06);
			
		\end{scope}
		\begin{scope}[shift={(-2/2,-1)}, yscale=-1]
			\draw [thick, orange] plot coordinates {(-3/4+0.14/3-3/8,0.14/2) (3/4-0.14/3+3/8,0.14/2) };
			\draw (0,-7/8) node [above] {\hbox{\footnotesize$\geqone$}};
			\draw [thick, blue] plot coordinates {(0,0.14/2+3/8) (0,0.14/2-8/8) } ;
			\fill [black] (-3/4+0.14/3,0.14/2) circle (0.06);
			\fill [black] (3/4-0.14/3,0.14/2) circle (0.06);
				\fill [black] (0,0.14/2) circle (0.06);
		\end{scope}
		
		\draw [thick, black, shorten >= -0.25cm, shorten <=-0.25cm] (0.14/3+1/2,0.14/2+1/2)--(-0.14/3 -1/4,-0.14/2-1);	
		\begin{scope} [shift={(-1/8+1/4+1/4*253/492,-1/4+1/4 )}]
			\draw [thick,black] plot  coordinates {(492/253*0.061,-1*0.061) (-492/253*0.061,1*0.061) };
	\end{scope} 
	
		\begin{scope} [shift={(-1/8+1/4-1/4*253/492,-1/4-1/4 )}]
			\draw [thick,black] plot  coordinates {(492/253*0.061,-1*0.061) (-492/253*0.061,1*0.061) };
	\end{scope}

		\begin{scope}[shift={(3/2+1/2,1/2)}]
			\begin{scope}[shift={(-1/48,0)}]
				\draw [thick,black] plot  coordinates {(4/12,-8/8) (-1/12,1/4) };
				\begin{scope}[shift={(1/9,-3/9)}]
					\draw [thick,black] plot  coordinates {(8/64,3/64) (-8/64,-3/64) };
				\end{scope}		
				\begin{scope}[shift={(2/9,-6/9)}]
					\draw [thick,black] plot  coordinates {(8/64,3/64) (-8/64,-3/64) };
				\end{scope}
			\end{scope}
		\end{scope}
		\begin{scope}[shift={(-3/2-1/4,-1)}, yscale=-1, xscale=-1]
			\begin{scope}[shift={(-1/48,0)}]
				\draw [thick,black] plot  coordinates {(4/12,-8/8) (-1/12,1/4) };
				\begin{scope}[shift={(1/9,-3/9)}]
					\draw [thick,black] plot  coordinates {(8/64,3/64) (-8/64,-3/64) };
				\end{scope}		
				\begin{scope}[shift={(2/9,-6/9)}]
					\draw [thick,black] plot  coordinates {(8/64,3/64) (-8/64,-3/64) };
				\end{scope}
			\end{scope}
		\end{scope}
\end{tikzpicture}}

\def\FigTwoCaseCii{\begin{tikzpicture}[baseline=13pt]
		\begin{scope}[shift={(3/4+1/2,0+1/2)}]
			\draw [thick, orange] plot coordinates {(-3/4+0.14/3-3/8,0.14/2) (3/4-0.14/3+3/8,0.14/2) };
			\draw (0,0) node [below] {\hbox{\footnotesize$\geqone$}};
			\fill [black] (-3/4+0.14/3,0.14/2) circle (0.06);
			\fill [black] (3/4-0.14/3,0.14/2) circle (0.06);
		\end{scope}
		\begin{scope}[shift={(-2/2,-1)}, yscale=-1]
			\draw [thick, orange] plot coordinates {(-3/4+0.14/3-3/8,0.14/2) (3/4-0.14/3+3/8,0.14/2) };
			\draw (0,-7/8) node [above] {\hbox{\footnotesize$\geqone$}};
			\draw [thick, blue] plot coordinates {(0,0.14/2+3/8) (0,0.14/2-8/8) } ;
			\fill [black] (-3/4+0.14/3,0.14/2) circle (0.06);
			\fill [black] (3/4-0.14/3,0.14/2) circle (0.06);
			\fill [black] (0,0.14/2) circle (0.06);
		\end{scope}
	\draw [thick, black, shorten >= -0.25cm, shorten <=-0.25cm] (0.14/3+1/2,0.14/2+1/2)--(-0.14/3 -1/4,-0.14/2-1);	
	\begin{scope} [shift={(-1/8+1/4+1/4*253/492,-1/4+1/4 )}]
		\draw [thick,black] plot  coordinates {(492/253*0.061,-1*0.061) (-492/253*0.061,1*0.061) };
	\end{scope} 
	
	\begin{scope} [shift={(-1/8+1/4-1/4*253/492,-1/4-1/4 )}]
		\draw [thick,black] plot  coordinates {(492/253*0.061,-1*0.061) (-492/253*0.061,1*0.061) };
	\end{scope}

		\begin{scope}[shift={(3/2+1/2,1/2)}]
			\begin{scope}[shift={(-1/48,0)}]
				\draw [thick,black] plot  coordinates {(4/12,-8/8) (-1/12,1/4) };
				\begin{scope}[shift={(1/9,-3/9)}]
					\draw [thick,black] plot  coordinates {(8/64,3/64) (-8/64,-3/64) };
				\end{scope}		
				\begin{scope}[shift={(2/9,-6/9)}]
					\draw [thick,black] plot  coordinates {(8/64,3/64) (-8/64,-3/64) };
				\end{scope}
			\end{scope}
		\end{scope}
		\begin{scope}[shift={(-3/2-1/4,-1)}, yscale=-1, xscale=-1]
			\begin{scope}[shift={(-1/48,0)}]
				\draw [thick,black] plot  coordinates {(4/12,-8/8) (-1/12,1/4) };
				\begin{scope}[shift={(1/9,-3/9)}]
					\draw [thick,black] plot  coordinates {(8/64,3/64) (-8/64,-3/64) };
				\end{scope}		
				\begin{scope}[shift={(2/9,-6/9)}]
					\draw [thick,black] plot  coordinates {(8/64,3/64) (-8/64,-3/64) };
				\end{scope}
			\end{scope}
		\end{scope}
\end{tikzpicture}}

\def\FigTwoCaseCiii{\begin{tikzpicture}[baseline=13pt]
		\begin{scope}[shift={(3/4+1/2,0+1/2)}]
			\draw [thick, orange] plot coordinates {(-3/4+0.14/3-3/8,0.14/2) (3/4-0.14/3+3/8,0.14/2) };
			\draw (0,0) node [below] {\hbox{\footnotesize$\geqone$}};
			\fill [black] (-3/4+0.14/3,0.14/2) circle (0.06);
			\fill [black] (3/4-0.14/3,0.14/2) circle (0.06);
		\end{scope}
		\begin{scope}[shift={(-2/2,-1)}, yscale=-1]
			\draw [thick, orange] plot coordinates {(-3/4+0.14/3-3/8,0.14/2) (3/4-0.14/3+3/8,0.14/2) };
			\draw (0,0) node [above] {\hbox{\footnotesize$\geqone$}};
			\fill [black] (-3/4+0.14/3,0.14/2) circle (0.06);
			\fill [black] (3/4-0.14/3,0.14/2) circle (0.06);
		\end{scope}
	\draw [thick, black, shorten >= -0.25cm, shorten <=-0.25cm] (0.14/3+1/2,0.14/2+1/2)--(-0.14/3 -1/4,-0.14/2-1);	
	\begin{scope} [shift={(-1/8+1/4+1/4*253/492,-1/4+1/4 )}]
		\draw [thick,black] plot  coordinates {(492/253*0.061,-1*0.061) (-492/253*0.061,1*0.061) };
	\end{scope} 
	
	\begin{scope} [shift={(-1/8+1/4-1/4*253/492,-1/4-1/4 )}]
		\draw [thick,black] plot  coordinates {(492/253*0.061,-1*0.061) (-492/253*0.061,1*0.061) };
	\end{scope} 
		
		\begin{scope}[shift={(3/2+1/2,0+1/2)}]
			\begin{scope}[shift={(-1/48,0)}]
				\draw [thick,black] plot  coordinates {(4/12,-8/8) (-1/12,1/4) };
				\begin{scope}[shift={(1/9,-3/9)}]
					\draw [thick,black] plot  coordinates {(8/64,3/64) (-8/64,-3/64) };
				\end{scope}		
				\begin{scope}[shift={(2/9,-6/9)}]
					\draw [thick,black] plot  coordinates {(8/64,3/64) (-8/64,-3/64) };
				\end{scope}
			\end{scope}
		\end{scope}
		\begin{scope}[shift={(-3/2-1/4,-1)}, yscale=-1, xscale=-1]
			\begin{scope}[shift={(-1/48,0)}]
				\draw [thick,black] plot  coordinates {(4/12,-8/8) (-1/12,1/4) };
				\begin{scope}[shift={(1/9,-3/9)}]
					\draw [thick,black] plot  coordinates {(8/64,3/64) (-8/64,-3/64) };
				\end{scope}		
				\begin{scope}[shift={(2/9,-6/9)}]
					\draw [thick,black] plot  coordinates {(8/64,3/64) (-8/64,-3/64) };
				\end{scope}
			\end{scope}
		\end{scope}
\end{tikzpicture}}

\def\FigTwoCaseCiv{\begin{tikzpicture}[baseline=13pt]
		\begin{scope}[shift={(3/4+1/2,0+1/2)}]
			\draw [thick,orange] plot [smooth, tension=1] coordinates {(-3/3,2/12) (1/4,-6/12) (-1/3,-14/12) };
			\draw [thick,orange] plot [smooth, tension=1] coordinates {(3/3,2/12) (-1/4,-6/12) (1/3,-14/12) };
			\fill [black] (0,-0.26) circle (0.06);
			\fill [black] (0,-0.99) circle (0.06);
			\fill [black] (-3/4+0.14/3,0.14/2) circle (0.06);
			\fill [black] (3/4-0.14/3,0.14/2) circle (0.06);
		\end{scope}
		\begin{scope}[shift={(-2/2,-1)}, yscale=-1]
			\draw [thick, orange] plot coordinates {(-3/4+0.14/3-3/8,0.14/2) (3/4-0.14/3+3/8,0.14/2) };
			\draw (0,-7/8) node [above] {\hbox{\footnotesize$\geqone$}};
			\draw [thick, blue] plot coordinates {(0,0.14/2+3/8) (0,0.14/2-8/8) } ;
			\fill [black] (-3/4+0.14/3,0.14/2) circle (0.06);
			\fill [black] (3/4-0.14/3,0.14/2) circle (0.06);
			\fill [black] (0,0.14/2) circle (0.06);
		\end{scope}
	\draw [thick, black, shorten >= -0.25cm, shorten <=-0.25cm] (0.14/3+1/2,0.14/2+1/2)--(-0.14/3 -1/4,-0.14/2-1);	
	\begin{scope} [shift={(-1/8+1/4+1/4*253/492,-1/4+1/4 )}]
		\draw [thick,black] plot  coordinates {(492/253*0.061,-1*0.061) (-492/253*0.061,1*0.061) };
	\end{scope} 
	
	\begin{scope} [shift={(-1/8+1/4-1/4*253/492,-1/4-1/4 )}]
		\draw [thick,black] plot  coordinates {(492/253*0.061,-1*0.061) (-492/253*0.061,1*0.061) };
	\end{scope} 
		\begin{scope}[shift={(3/2+1/2,0+1/2)}]
			\begin{scope}[shift={(-1/48,0)}]
				\draw [thick,black] plot  coordinates {(4/12,-8/8) (-1/12,1/4) };
				\begin{scope}[shift={(1/9,-3/9)}]
					\draw [thick,black] plot  coordinates {(8/64,3/64) (-8/64,-3/64) };
				\end{scope}		
				\begin{scope}[shift={(2/9,-6/9)}]
					\draw [thick,black] plot  coordinates {(8/64,3/64) (-8/64,-3/64) };
				\end{scope}
			\end{scope}
		\end{scope}
		\begin{scope}[shift={(-3/2-1/4,-1)}, yscale=-1, xscale=-1]
			\begin{scope}[shift={(-1/48,0)}]
				\draw [thick,black] plot  coordinates {(4/12,-8/8) (-1/12,1/4) };
				\begin{scope}[shift={(1/9,-3/9)}]
					\draw [thick,black] plot  coordinates {(8/64,3/64) (-8/64,-3/64) };
				\end{scope}		
				\begin{scope}[shift={(2/9,-6/9)}]
					\draw [thick,black] plot  coordinates {(8/64,3/64) (-8/64,-3/64) };
				\end{scope}
			\end{scope}
		\end{scope}
\end{tikzpicture}}

\def\FigTwoCaseCv{\begin{tikzpicture}[baseline=13pt]
		\begin{scope}[shift={(3/4+1/2,0+1/2)}]
			\draw [thick,orange] plot [smooth, tension=1] coordinates {(-3/3,2/12) (1/4,-6/12) (-1/3,-14/12) };
			\draw [thick,orange] plot [smooth, tension=1] coordinates {(3/3,2/12) (-1/4,-6/12) (1/3,-14/12) };
			\fill [black] (0,-0.26) circle (0.06);
			\fill [black] (0,-0.99) circle (0.06);
			\fill [black] (-3/4+0.14/3,0.14/2) circle (0.06);
			\fill [black] (3/4-0.14/3,0.14/2) circle (0.06);
		\end{scope}
		\begin{scope}[shift={(-2/2,-1)}, yscale=-1]
			\draw [thick, orange] plot coordinates {(-3/4+0.14/3-3/8,0.14/2) (3/4-0.14/3+3/8,0.14/2) };
			\draw (0,0) node [above] {\hbox{\footnotesize$\geqone$}};
			\fill [black] (-3/4+0.14/3,0.14/2) circle (0.06);
			\fill [black] (3/4-0.14/3,0.14/2) circle (0.06);
		\end{scope}
		\draw [thick, black, shorten >= -0.25cm, shorten <=-0.25cm] (0.14/3+1/2,0.14/2+1/2)--(-0.14/3 -1/4,-0.14/2-1);	
		\begin{scope} [shift={(-1/8+1/4+1/4*253/492,-1/4+1/4 )}]
			\draw [thick,black] plot  coordinates {(492/253*0.061,-1*0.061) (-492/253*0.061,1*0.061) };
		\end{scope} 
		
		\begin{scope} [shift={(-1/8+1/4-1/4*253/492,-1/4-1/4 )}]
			\draw [thick,black] plot  coordinates {(492/253*0.061,-1*0.061) (-492/253*0.061,1*0.061) };
		\end{scope} 
	
		\begin{scope}[shift={(3/2+1/2,0+1/2)}]
			\begin{scope}[shift={(-1/48,0)}]
				\draw [thick,black] plot  coordinates {(4/12,-8/8) (-1/12,1/4) };
				\begin{scope}[shift={(1/9,-3/9)}]
					\draw [thick,black] plot  coordinates {(8/64,3/64) (-8/64,-3/64) };
				\end{scope}		
				\begin{scope}[shift={(2/9,-6/9)}]
					\draw [thick,black] plot  coordinates {(8/64,3/64) (-8/64,-3/64) };
				\end{scope}
			\end{scope}
		\end{scope}
		\begin{scope}[shift={(-3/2-1/4,-1)}, yscale=-1, xscale=-1]
			\begin{scope}[shift={(-1/48,0)}]
				\draw [thick,black] plot  coordinates {(4/12,-8/8) (-1/12,1/4) };
				\begin{scope}[shift={(1/9,-3/9)}]
					\draw [thick,black] plot  coordinates {(8/64,3/64) (-8/64,-3/64) };
				\end{scope}		
				\begin{scope}[shift={(2/9,-6/9)}]
					\draw [thick,black] plot  coordinates {(8/64,3/64) (-8/64,-3/64) };
				\end{scope}
			\end{scope}
		\end{scope}
\end{tikzpicture}}

\def\FigTwoCaseCvi{\begin{tikzpicture}[baseline=13pt]
		\begin{scope}[shift={(3/4+1/2,0+1/2)}]
			\draw [thick,orange] plot [smooth, tension=1] coordinates {(-3/3,2/12) (1/4,-6/12) (-1/3,-14/12) };
			\draw [thick,orange] plot [smooth, tension=1] coordinates {(3/3,2/12) (-1/4,-6/12) (1/3,-14/12) };
			\fill [black] (0,-0.26) circle (0.06);
			\fill [black] (0,-0.99) circle (0.06);
			\fill [black] (-3/4+0.14/3,0.14/2) circle (0.06);
			\fill [black] (3/4-0.14/3,0.14/2) circle (0.06);
		\end{scope}
		\begin{scope}[shift={(-2/2,-1)}, yscale=-1]
			\draw [thick,orange] plot [smooth, tension=1] coordinates {(-3/3,2/12) (1/4,-6/12) (-1/3,-14/12) };
			\draw [thick,orange] plot [smooth, tension=1] coordinates {(3/3,2/12) (-1/4,-6/12) (1/3,-14/12) };
			\fill [black] (0,-0.26) circle (0.06);
			\fill [black] (0,-0.99) circle (0.06);
			\fill [black] (-3/4+0.14/3,0.14/2) circle (0.06);
			\fill [black] (3/4-0.14/3,0.14/2) circle (0.06);
		\end{scope}
	\draw [thick, black, shorten >= -0.25cm, shorten <=-0.25cm] (0.14/3+1/2,0.14/2+1/2)--(-0.14/3 -1/4,-0.14/2-1);	
	\begin{scope} [shift={(-1/8+1/4+1/4*253/492,-1/4+1/4 )}]
		\draw [thick,black] plot  coordinates {(492/253*0.061,-1*0.061) (-492/253*0.061,1*0.061) };
	\end{scope} 
	
	\begin{scope} [shift={(-1/8+1/4-1/4*253/492,-1/4-1/4 )}]
		\draw [thick,black] plot  coordinates {(492/253*0.061,-1*0.061) (-492/253*0.061,1*0.061) };
	\end{scope} 
		
		\begin{scope}[shift={(3/2+1/2,0+1/2)}]
			\begin{scope}[shift={(-1/48,0)}]
				\draw [thick,black] plot  coordinates {(4/12,-8/8) (-1/12,1/4) };
				\begin{scope}[shift={(1/9,-3/9)}]
					\draw [thick,black] plot  coordinates {(8/64,3/64) (-8/64,-3/64) };
				\end{scope}		
				\begin{scope}[shift={(2/9,-6/9)}]
					\draw [thick,black] plot  coordinates {(8/64,3/64) (-8/64,-3/64) };
				\end{scope}
			\end{scope}
		\end{scope}
		\begin{scope}[shift={(-3/2-1/4,-1)}, yscale=-1, xscale=-1]
			\begin{scope}[shift={(-1/48,0)}]
				\draw [thick,black] plot  coordinates {(4/12,-8/8) (-1/12,1/4) };
				\begin{scope}[shift={(1/9,-3/9)}]
					\draw [thick,black] plot  coordinates {(8/64,3/64) (-8/64,-3/64) };
				\end{scope}		
				\begin{scope}[shift={(2/9,-6/9)}]
					\draw [thick,black] plot  coordinates {(8/64,3/64) (-8/64,-3/64) };
				\end{scope}
			\end{scope}
		\end{scope}
\end{tikzpicture}}

\def\includefigC#1{
\makebox(130pt,70pt){#1}}

\def\FigTwoCaseC{
\begin{tabular}{|c|c|c|}
	\hline
{\large\strut}	(C1) & (C2) & (C3) \\
\hline
{\Large\strut} $4>\alpha\ge\beta$ & $4=\alpha>\beta$ & $\alpha=4=\beta$ \\
\hline
\includefigC{\FigTwoCaseCi} & 
\includefigC{\FigTwoCaseCii} & 
\includefigC{\FigTwoCaseCiii} \\
\hline\hline
{\large\strut}	(C4) & (C5) & (C6) \\
\hline
{\Large\strut} $\alpha>4>\beta$ & $\alpha>4=\beta$ & $\alpha\ge\beta>4$ \\
\hline
\includefigC{\FigTwoCaseCiv} & 
\includefigC{\FigTwoCaseCv} & 
\includefigC{\FigTwoCaseCvi} \\
\hline
\end{tabular}}

\def\FigTwoCaseDi{\begin{tikzpicture}[baseline=13pt]
\draw [thick,black] plot  coordinates {(-1/2,0) (4/2,0) };
\begin{scope}[shift={(3/2,0)}]
	\draw [thick,black] plot  coordinates {(-4/6,-4/4) (1/6,1/4) };
	\fill [black] (0,0) circle (0.06);
	\begin{scope}[shift={(-2/9,-3/9)}]
		\draw [thick,black] plot  coordinates {(-3/32,2/32) (3/32,-2/32) };
	\end{scope}
	\begin{scope}[shift={(-4/9,-6/9)}]
		\draw [thick,black] plot  coordinates {(-3/32,2/32) (3/32,-2/32) };
	\end{scope}
\end{scope}
\begin{scope}[shift={(3/4,0)}]
	\draw [thick,black] plot  coordinates {(-4/6,-4/4) (1/6,1/4) };
	\fill [black] (0,0) circle (0.06);
	\begin{scope}[shift={(-2/9,-3/9)}]
		\draw [thick,black] plot  coordinates {(-3/32,2/32) (3/32,-2/32) };
	\end{scope}
	\begin{scope}[shift={(-4/9,-6/9)}]
		\draw [thick,black] plot  coordinates {(-3/32,2/32) (3/32,-2/32) };
	\end{scope}
\end{scope}
\begin{scope}[shift={(0,0)}]
	\draw [thick,black] plot  coordinates {(-4/6,-4/4) (1/6,1/4) };
	\fill [black] (0,0) circle (0.06);
	\begin{scope}[shift={(-2/9,-3/9)}]
		\draw [thick,black] plot  coordinates {(-3/32,2/32) (3/32,-2/32) };
	\end{scope}
	\begin{scope}[shift={(-4/9,-6/9)}]
		\draw [thick,black] plot  coordinates {(-3/32,2/32) (3/32,-2/32) };
	\end{scope}
\end{scope}
\draw (-2/4,0) node [left] {\hbox{\footnotesize$\geqtwo$}};
\end{tikzpicture}}

\def\FigTwoCaseDii{\begin{tikzpicture}[baseline=13pt]
		\draw [thick,black] plot  coordinates {(-1/2,0) (3/2+1/8+1/100,0) };
		\draw [thick, black] (3/2+1/8,0) arc (270:450:0.25);
		\draw [thick, black] plot coordinates {(3/2+1/8+1/100,0.5) (2/2+1/8,0.5)};
		\begin{scope}[shift={(3/2,0)}]
			\draw [thick,orange] plot  coordinates {(0,6/8) (0,0) };
			\draw [thick,orange] plot [smooth, tension=1] coordinates {(0,0) (1/16,-4/8) (3/16+1/32,-7/8)  };
			\fill [black] (0,1/2) circle (0.06);
			\fill [black] (0,0) circle (0.06);
		\end{scope}
		\begin{scope}[shift={(3/2+1/16,-1/2)}]
			\draw [thick,black] plot  coordinates {(-4/6,-4/4) (1/6,1/4) };
			\fill [black] (0,0) circle (0.06);
			\begin{scope}[shift={(-2/9,-3/9)}]
				\draw [thick,black] plot  coordinates {(-3/32,2/32) (3/32,-2/32) };
			\end{scope}
			\begin{scope}[shift={(-4/9,-6/9)}]
				\draw [thick,black] plot  coordinates {(-3/32,2/32) (3/32,-2/32) };
			\end{scope}
		\end{scope}
		\begin{scope}[shift={(3/4,0)}]
			\draw [thick,black] plot  coordinates {(-4/6,-4/4) (1/6,1/4) };
			\fill [black] (0,0) circle (0.06);
			\begin{scope}[shift={(-2/9,-3/9)}]
				\draw [thick,black] plot  coordinates {(-3/32,2/32) (3/32,-2/32) };
			\end{scope}
			\begin{scope}[shift={(-4/9,-6/9)}]
				\draw [thick,black] plot  coordinates {(-3/32,2/32) (3/32,-2/32) };
			\end{scope}
		\end{scope}
		\begin{scope}[shift={(0,0)}]
			\draw [thick,black] plot  coordinates {(-4/6,-4/4) (1/6,1/4) };
			\fill [black] (0,0) circle (0.06);
			\begin{scope}[shift={(-2/9,-3/9)}]
				\draw [thick,black] plot  coordinates {(-3/32,2/32) (3/32,-2/32) };
			\end{scope}
			\begin{scope}[shift={(-4/9,-6/9)}]
				\draw [thick,black] plot  coordinates {(-3/32,2/32) (3/32,-2/32) };
			\end{scope}
		\end{scope}
		\draw (-2/4,0) node [left] {\hbox{\footnotesize$\geqone$}};
\end{tikzpicture}}

\def\FigTwoCaseDiii{\begin{tikzpicture}[baseline=13pt]
			\draw [thick,black] plot  coordinates {(-1/2,0) (3/2+1/8+1/100,0) };
		\draw [thick, black] (3/2+1/8,0) arc (270:450:0.25);
		\draw [thick, black] plot coordinates {(3/2+1/8+1/100,0.5) (-1/2,0.5)};
		\begin{scope}[shift={(3/2,0)}]
	\draw [thick,orange] plot  coordinates {(0,6/8) (0,0) };
\draw [thick,orange] plot [smooth, tension=1] coordinates {(0,0) (1/16,-4/8) (3/16+1/32,-7/8)  };
			\fill [black] (0,1/2) circle (0.06);
			\fill [black] (0,0) circle (0.06);
		\end{scope}
	\begin{scope}[shift={(3/2-6/4+3/4,0)}]
	\draw [thick,orange] plot  coordinates {(0,6/8) (0,0) };
\draw [thick,orange] plot [smooth, tension=1] coordinates {(0,0) (1/16,-4/8) (3/16+1/32,-7/8)  };
		\fill [black] (0,1/2) circle (0.06);
		\fill [black] (0,0) circle (0.06);
	\end{scope}
		\begin{scope}[shift={(3/2+1/16,-1/2)}]
			\draw [thick,black] plot  coordinates {(-4/6,-4/4) (1/6,1/4) };
			\fill [black] (0,0) circle (0.06);
			\begin{scope}[shift={(-2/9,-3/9)}]
				\draw [thick,black] plot  coordinates {(-3/32,2/32) (3/32,-2/32) };
			\end{scope}
			\begin{scope}[shift={(-4/9,-6/9)}]
				\draw [thick,black] plot  coordinates {(-3/32,2/32) (3/32,-2/32) };
			\end{scope}
		\end{scope}
	\begin{scope}[shift={(3/4+1/16,-1/2)}]
		\draw [thick,black] plot  coordinates {(-4/6,-4/4) (1/6,1/4) };
		\fill [black] (0,0) circle (0.06);
		\begin{scope}[shift={(-2/9,-3/9)}]
			\draw [thick,black] plot  coordinates {(-3/32,2/32) (3/32,-2/32) };
		\end{scope}
		\begin{scope}[shift={(-4/9,-6/9)}]
			\draw [thick,black] plot  coordinates {(-3/32,2/32) (3/32,-2/32) };
		\end{scope}
	\end{scope}
		\begin{scope}[shift={(0+0/16,0)}]
			\draw [thick,black] plot  coordinates {(-4/6,-4/4) (1/6,1/4) };
			\fill [black] (0,0) circle (0.06);
			\begin{scope}[shift={(-2/9,-3/9)}]
				\draw [thick,black] plot  coordinates {(-3/32,2/32) (3/32,-2/32) };
			\end{scope}
			\begin{scope}[shift={(-4/9,-6/9)}]
				\draw [thick,black] plot  coordinates {(-3/32,2/32) (3/32,-2/32) };
			\end{scope}
		\end{scope}
\end{tikzpicture}}

\def\FigTwoCaseDiv{\begin{tikzpicture}[baseline=13pt]
		\draw [thick,black] plot  coordinates {(-1/2,0) (3/2+3/8,0) };
	\draw [thick, black] plot coordinates {(3/2+3/8,0.5) (-1/2,0.5)};
	\begin{scope}[shift={(3/2,0)}]
\draw [thick,orange] plot  coordinates {(0,6/8) (0,0) };
\draw [thick,orange] plot [smooth, tension=1] coordinates {(0,0) (1/16,-4/8) (3/16+1/32,-7/8)  };
		\fill [black] (0,1/2) circle (0.06);
		\fill [black] (0,0) circle (0.06);
	\end{scope}
	\begin{scope}[shift={(3/2-3/4,0)}]
\draw [thick,orange] plot  coordinates {(0,6/8) (0,0) };
\draw [thick,orange] plot [smooth, tension=1] coordinates {(0,0) (1/16,-4/8) (3/16+1/32,-7/8)  };
		\fill [black] (0,1/2) circle (0.06);
		\fill [black] (0,0) circle (0.06);
	\end{scope}
	\begin{scope}[shift={(3/2+1/16,-1/2)}]
		\draw [thick,black] plot  coordinates {(-4/6,-4/4) (1/6,1/4) };
		\fill [black] (0,0) circle (0.06);
		\begin{scope}[shift={(-2/9,-3/9)}]
			\draw [thick,black] plot  coordinates {(-3/32,2/32) (3/32,-2/32) };
		\end{scope}
		\begin{scope}[shift={(-4/9,-6/9)}]
			\draw [thick,black] plot  coordinates {(-3/32,2/32) (3/32,-2/32) };
		\end{scope}
	\end{scope}
	\begin{scope}[shift={(3/4+1/16,-1/2)}]
		\draw [thick,black] plot  coordinates {(-4/6,-4/4) (1/6,1/4) };
		\fill [black] (0,0) circle (0.06);
		\begin{scope}[shift={(-2/9,-3/9)}]
			\draw [thick,black] plot  coordinates {(-3/32,2/32) (3/32,-2/32) };
		\end{scope}
		\begin{scope}[shift={(-4/9,-6/9)}]
			\draw [thick,black] plot  coordinates {(-3/32,2/32) (3/32,-2/32) };
		\end{scope}
	\end{scope}
\begin{scope}[shift={(0,0)}]
\draw [thick,orange] plot  coordinates {(0,6/8) (0,0) };
\draw [thick,orange] plot [smooth, tension=1] coordinates {(0,0) (1/16,-4/8) (3/16+1/32,-7/8)  };
	\fill [black] (0,1/2) circle (0.06);
	\fill [black] (0,0) circle (0.06);
\end{scope}
	\begin{scope}[shift={(0+1/16,-1/2)}]
		\draw [thick,black] plot  coordinates {(-4/6,-4/4) (1/6,1/4) };
		\fill [black] (0,0) circle (0.06);
		\begin{scope}[shift={(-2/9,-3/9)}]
			\draw [thick,black] plot  coordinates {(-3/32,2/32) (3/32,-2/32) };
		\end{scope}
		\begin{scope}[shift={(-4/9,-6/9)}]
			\draw [thick,black] plot  coordinates {(-3/32,2/32) (3/32,-2/32) };
		\end{scope}
	\end{scope}
\end{tikzpicture}}

\def\FigTwoCaseDx{\begin{tikzpicture}[baseline=13pt]
		\begin{scope}[shift={(3/4+1/2,0+1/2)}]
			\draw [thick,orange] plot [smooth, tension=1] coordinates {(-3/3,2/12) (1/4,-6/12) (-1/3,-14/12) };
			\draw [thick,orange] plot [smooth, tension=1] coordinates {(3/3,2/12) (-1/4,-6/12) (1/3,-14/12) };
			\fill [black] (0,-0.26) circle (0.06);
			\fill [black] (0,-0.99) circle (0.06);
			\fill [black] (-3/4+0.14/3,0.14/2) circle (0.06);
			\fill [black] (3/4-0.14/3,0.14/2) circle (0.06);
		\end{scope}
		\begin{scope}[shift={(-2/2,-1)}, yscale=-1]
		\draw [thick,orange] plot [smooth, tension=1] coordinates {(-3/3,2/12) (1/4,-6/12) (-1/3,-14/12) };
		\draw [thick,orange] plot [smooth, tension=1] coordinates {(3/3,2/12) (-1/4,-6/12) (1/3,-14/12) };
		\fill [black] (0,-0.26) circle (0.06);
		\fill [black] (0,-0.99) circle (0.06);
		\fill [black] (-3/4+0.14/3,0.14/2) circle (0.06);
		\fill [black] (3/4-0.14/3,0.14/2) circle (0.06);
	\end{scope}
	\draw [thick, black, shorten >= -0.25cm, shorten <=-0.25cm] (0.14/3+1/2,0.14/2+1/2)--(-0.14/3 -1/4,-0.14/2-1);
	
	\begin{scope} [shift={(-1/8+1/4,-1/4 )}]
		\fill [black] (0,0) circle (0.06);
	\draw [thick,black] plot  coordinates {(4/12,-8/8) (-1/12,1/4) };
	\begin{scope}[shift={(1/9,-3/9)}]
		\draw [thick,black] plot  coordinates {(8/64,3/64) (-8/64,-3/64) };
	\end{scope}		
	\begin{scope}[shift={(2/9,-6/9)}]
		\draw [thick,black] plot  coordinates {(8/64,3/64) (-8/64,-3/64) };
	\end{scope}
	\end{scope}

		\begin{scope}[shift={(3/2+1/2,0+1/2)}]
			\begin{scope}[shift={(-1/48,0)}]
				\draw [thick,black] plot  coordinates {(4/12,-8/8) (-1/12,1/4) };
				\begin{scope}[shift={(1/9,-3/9)}]
						\draw [thick,black] plot  coordinates {(8/64,3/64) (-8/64,-3/64) };
				\end{scope}		
				\begin{scope}[shift={(2/9,-6/9)}]
					\draw [thick,black] plot  coordinates {(8/64,3/64) (-8/64,-3/64) };
				\end{scope}
			\end{scope}
		\end{scope}
		\begin{scope}[shift={(-3/2-1/4,-1)}, yscale=-1, xscale=-1]
			\begin{scope}[shift={(-1/48,0)}]
			\draw [thick,black] plot  coordinates {(4/12,-8/8) (-1/12,1/4) };
			\begin{scope}[shift={(1/9,-3/9)}]
				\draw [thick,black] plot  coordinates {(8/64,3/64) (-8/64,-3/64) };
			\end{scope}		
			\begin{scope}[shift={(2/9,-6/9)}]
				\draw [thick,black] plot  coordinates {(8/64,3/64) (-8/64,-3/64) };
			\end{scope}
		\end{scope}
	\end{scope}
	\end{tikzpicture}}

\def\FigTwoCaseDix{\begin{tikzpicture}[baseline=13pt]
		\begin{scope}[shift={(3/4+1/2,0+1/2)}]
			\draw [thick,orange] plot [smooth, tension=1] coordinates {(-3/3,2/12) (1/4,-6/12) (-1/3,-14/12) };
			\draw [thick,orange] plot [smooth, tension=1] coordinates {(3/3,2/12) (-1/4,-6/12) (1/3,-14/12) };
			\fill [black] (0,-0.26) circle (0.06);
			\fill [black] (0,-0.99) circle (0.06);
			\fill [black] (-3/4+0.14/3,0.14/2) circle (0.06);
			\fill [black] (3/4-0.14/3,0.14/2) circle (0.06);
		\end{scope}
		\begin{scope}[shift={(-2/2,-1)}, yscale=-1]
			\draw [thick, orange] plot coordinates {(-3/4+0.14/3-3/8,0.14/2) (3/4-0.14/3+3/8,0.14/2) };
			\draw (0,0) node [above] {\hbox{\footnotesize$\geqone$}};
			\fill [black] (-3/4+0.14/3,0.14/2) circle (0.06);
			\fill [black] (3/4-0.14/3,0.14/2) circle (0.06);
		\end{scope}
	\draw [thick, black, shorten >= -0.25cm, shorten <=-0.25cm] (0.14/3+1/2,0.14/2+1/2)--(-0.14/3 -1/4,-0.14/2-1);
	
		\begin{scope} [shift={(-1/8+1/4,-1/4 )}]
		\fill [black] (0,0) circle (0.06);
		\draw [thick,black] plot  coordinates {(4/12,-8/8) (-1/12,1/4) };
		\begin{scope}[shift={(1/9,-3/9)}]
			\draw [thick,black] plot  coordinates {(8/64,3/64) (-8/64,-3/64) };
		\end{scope}		
		\begin{scope}[shift={(2/9,-6/9)}]
			\draw [thick,black] plot  coordinates {(8/64,3/64) (-8/64,-3/64) };
		\end{scope}
	\end{scope}
		
		\begin{scope}[shift={(3/2+1/2,0+1/2)}]
					\begin{scope}[shift={(-1/48,0)}]
			\draw [thick,black] plot  coordinates {(4/12,-8/8) (-1/12,1/4) };
			\begin{scope}[shift={(1/9,-3/9)}]
				\draw [thick,black] plot  coordinates {(8/64,3/64) (-8/64,-3/64) };
			\end{scope}		
			\begin{scope}[shift={(2/9,-6/9)}]
				\draw [thick,black] plot  coordinates {(8/64,3/64) (-8/64,-3/64) };
			\end{scope}
		\end{scope}
		\end{scope}
		\begin{scope}[shift={(-3/2-1/4,-1)}, yscale=-1, xscale=-1]
						\begin{scope}[shift={(-1/48,0)}]
				\draw [thick,black] plot  coordinates {(4/12,-8/8) (-1/12,1/4) };
				\begin{scope}[shift={(1/9,-3/9)}]
					\draw [thick,black] plot  coordinates {(8/64,3/64) (-8/64,-3/64) };
				\end{scope}		
				\begin{scope}[shift={(2/9,-6/9)}]
					\draw [thick,black] plot  coordinates {(8/64,3/64) (-8/64,-3/64) };
				\end{scope}
			\end{scope}
		\end{scope}
\end{tikzpicture}}

\def\FigTwoCaseDviii{\begin{tikzpicture}[baseline=13pt]
		\begin{scope}[shift={(3/4+1/2,0+1/2)}]
			\draw [thick,orange] plot [smooth, tension=1] coordinates {(-3/3,2/12) (1/4,-6/12) (-1/3,-14/12) };
			\draw [thick,orange] plot [smooth, tension=1] coordinates {(3/3,2/12) (-1/4,-6/12) (1/3,-14/12) };
			\fill [black] (0,-0.26) circle (0.06);
			\fill [black] (0,-0.99) circle (0.06);
			\fill [black] (-3/4+0.14/3,0.14/2) circle (0.06);
			\fill [black] (3/4-0.14/3,0.14/2) circle (0.06);
		\end{scope}
		\begin{scope}[shift={(-2/2,-1)}, yscale=-1]
			\draw [thick, orange] plot coordinates {(-3/4+0.14/3-3/8,0.14/2) (3/4-0.14/3+3/8,0.14/2) };
			\draw (0,-7/8) node [above] {\hbox{\footnotesize$\geqone$}};
			\draw [thick, blue] plot coordinates {(0,0.14/2+3/8) (0,0.14/2-8/8) } ;
			\fill [black] (-3/4+0.14/3,0.14/2) circle (0.06);
			\fill [black] (3/4-0.14/3,0.14/2) circle (0.06);
			\fill [black] (0,0.14/2) circle (0.06);
		\end{scope}
		\draw [thick, black, shorten >= -0.25cm, shorten <=-0.25cm] (0.14/3+1/2,0.14/2+1/2)--(-0.14/3 -1/4,-0.14/2-1);
			\begin{scope} [shift={(-1/8+1/4,-1/4 )}]
			\fill [black] (0,0) circle (0.06);
			\draw [thick,black] plot  coordinates {(4/12,-8/8) (-1/12,1/4) };
			\begin{scope}[shift={(1/9,-3/9)}]
				\draw [thick,black] plot  coordinates {(8/64,3/64) (-8/64,-3/64) };
			\end{scope}		
			\begin{scope}[shift={(2/9,-6/9)}]
				\draw [thick,black] plot  coordinates {(8/64,3/64) (-8/64,-3/64) };
			\end{scope}
		\end{scope}
		\begin{scope}[shift={(3/2+1/2,0+1/2)}]
	\begin{scope}[shift={(-1/48,0)}]
	\draw [thick,black] plot  coordinates {(4/12,-8/8) (-1/12,1/4) };
	\begin{scope}[shift={(1/9,-3/9)}]
		\draw [thick,black] plot  coordinates {(8/64,3/64) (-8/64,-3/64) };
	\end{scope}		
	\begin{scope}[shift={(2/9,-6/9)}]
		\draw [thick,black] plot  coordinates {(8/64,3/64) (-8/64,-3/64) };
	\end{scope}
\end{scope}
		\end{scope}
		\begin{scope}[shift={(-3/2-1/4,-1)}, yscale=-1, xscale=-1]
				\begin{scope}[shift={(-1/48,0)}]
				\draw [thick,black] plot  coordinates {(4/12,-8/8) (-1/12,1/4) };
				\begin{scope}[shift={(1/9,-3/9)}]
					\draw [thick,black] plot  coordinates {(8/64,3/64) (-8/64,-3/64) };
				\end{scope}		
				\begin{scope}[shift={(2/9,-6/9)}]
					\draw [thick,black] plot  coordinates {(8/64,3/64) (-8/64,-3/64) };
				\end{scope}
			\end{scope}
		\end{scope}
\end{tikzpicture}}

\def\FigTwoCaseDvii{\begin{tikzpicture}[baseline=13pt]
		\begin{scope}[shift={(3/4+1/2,0+1/2)}]
			\draw [thick, orange] plot coordinates {(-3/4+0.14/3-3/8,0.14/2) (3/4-0.14/3+3/8,0.14/2) };
			\draw (0,0) node [below] {\hbox{\footnotesize$\geqone$}};
			\fill [black] (-3/4+0.14/3,0.14/2) circle (0.06);
			\fill [black] (3/4-0.14/3,0.14/2) circle (0.06);
		\end{scope}
		\begin{scope}[shift={(-2/2,-1)}, yscale=-1]
			\draw [thick, orange] plot coordinates {(-3/4+0.14/3-3/8,0.14/2) (3/4-0.14/3+3/8,0.14/2) };
			\draw (0,0) node [above] {\hbox{\footnotesize$\geqone$}};
			\fill [black] (-3/4+0.14/3,0.14/2) circle (0.06);
			\fill [black] (3/4-0.14/3,0.14/2) circle (0.06);
		\end{scope}
		\draw [thick, black, shorten >= -0.25cm, shorten <=-0.25cm] (0.14/3+1/2,0.14/2+1/2)--(-0.14/3 -1/4,-0.14/2-1);
		
			\begin{scope} [shift={(-1/8+1/4,-1/4 )}]
			\fill [black] (0,0) circle (0.06);
			\draw [thick,black] plot  coordinates {(4/12,-8/8) (-1/12,1/4) };
			\begin{scope}[shift={(1/9,-3/9)}]
				\draw [thick,black] plot  coordinates {(8/64,3/64) (-8/64,-3/64) };
			\end{scope}		
			\begin{scope}[shift={(2/9,-6/9)}]
				\draw [thick,black] plot  coordinates {(8/64,3/64) (-8/64,-3/64) };
			\end{scope}
		\end{scope}
		
		\begin{scope}[shift={(3/2+1/2,0+1/2)}]
					\begin{scope}[shift={(-1/48,0)}]
			\draw [thick,black] plot  coordinates {(4/12,-8/8) (-1/12,1/4) };
			\begin{scope}[shift={(1/9,-3/9)}]
				\draw [thick,black] plot  coordinates {(8/64,3/64) (-8/64,-3/64) };
			\end{scope}		
			\begin{scope}[shift={(2/9,-6/9)}]
				\draw [thick,black] plot  coordinates {(8/64,3/64) (-8/64,-3/64) };
			\end{scope}
		\end{scope}
		\end{scope}
		\begin{scope}[shift={(-3/2-1/4,-1)}, yscale=-1, xscale=-1]
			\begin{scope}[shift={(-1/48,0)}]
	\draw [thick,black] plot  coordinates {(4/12,-8/8) (-1/12,1/4) };
	\begin{scope}[shift={(1/9,-3/9)}]
		\draw [thick,black] plot  coordinates {(8/64,3/64) (-8/64,-3/64) };
	\end{scope}		
	\begin{scope}[shift={(2/9,-6/9)}]
		\draw [thick,black] plot  coordinates {(8/64,3/64) (-8/64,-3/64) };
	\end{scope}
\end{scope}
		\end{scope}
\end{tikzpicture}}

\def\FigTwoCaseDvi{\begin{tikzpicture}[baseline=13pt]
		\begin{scope}[shift={(3/4+1/2,0+1/2)}]
	\draw [thick, orange] plot coordinates {(-3/4+0.14/3-3/8,0.14/2) (3/4-0.14/3+3/8,0.14/2) };
	\draw (0,0) node [below] {\hbox{\footnotesize$\geqone$}};
	\fill [black] (-3/4+0.14/3,0.14/2) circle (0.06);
	\fill [black] (3/4-0.14/3,0.14/2) circle (0.06);
		\end{scope}
	\begin{scope}[shift={(-2/2,-1)}, yscale=-1]
				\draw [thick, orange] plot coordinates {(-3/4+0.14/3-3/8,0.14/2) (3/4-0.14/3+3/8,0.14/2) };
		\draw (0,-7/8) node [above] {\hbox{\footnotesize$\geqone$}};
		\draw [thick, blue] plot coordinates {(0,0.14/2+3/8) (0,0.14/2-8/8) } ;
		\fill [black] (-3/4+0.14/3,0.14/2) circle (0.06);
		\fill [black] (3/4-0.14/3,0.14/2) circle (0.06);
		\fill [black] (0,0.14/2) circle (0.06);
	\end{scope}
\draw [thick, black, shorten >= -0.25cm, shorten <=-0.25cm] (0.14/3+1/2,0.14/2+1/2)--(-0.14/3 -1/4,-0.14/2-1);

	\begin{scope} [shift={(-1/8+1/4,-1/4 )}]
	\fill [black] (0,0) circle (0.06);
	\draw [thick,black] plot  coordinates {(4/12,-8/8) (-1/12,1/4) };
	\begin{scope}[shift={(1/9,-3/9)}]
		\draw [thick,black] plot  coordinates {(8/64,3/64) (-8/64,-3/64) };
	\end{scope}		
	\begin{scope}[shift={(2/9,-6/9)}]
		\draw [thick,black] plot  coordinates {(8/64,3/64) (-8/64,-3/64) };
	\end{scope}
\end{scope}

		\begin{scope}[shift={(3/2+1/2,1/2)}]
						\begin{scope}[shift={(-1/48,0)}]
				\draw [thick,black] plot  coordinates {(4/12,-8/8) (-1/12,1/4) };
				\begin{scope}[shift={(1/9,-3/9)}]
					\draw [thick,black] plot  coordinates {(8/64,3/64) (-8/64,-3/64) };
				\end{scope}		
				\begin{scope}[shift={(2/9,-6/9)}]
					\draw [thick,black] plot  coordinates {(8/64,3/64) (-8/64,-3/64) };
				\end{scope}
			\end{scope}
			\end{scope}
		\begin{scope}[shift={(-3/2-1/4,-1)}, yscale=-1, xscale=-1]
						\begin{scope}[shift={(-1/48,0)}]
				\draw [thick,black] plot  coordinates {(4/12,-8/8) (-1/12,1/4) };
				\begin{scope}[shift={(1/9,-3/9)}]
					\draw [thick,black] plot  coordinates {(8/64,3/64) (-8/64,-3/64) };
				\end{scope}		
				\begin{scope}[shift={(2/9,-6/9)}]
					\draw [thick,black] plot  coordinates {(8/64,3/64) (-8/64,-3/64) };
				\end{scope}
			\end{scope}
		\end{scope}
\end{tikzpicture}}

\def\FigTwoCaseDv{\begin{tikzpicture}[baseline=13pt]
		\begin{scope}[shift={(3/4+1/2,1/2)}]
			\draw [thick, orange] plot coordinates {(-3/4+0.14/3-3/8,0.14/2) (3/4-0.14/3+3/8,0.14/2) };
			\draw (0,-7/8) node [below] {\hbox{\footnotesize$\geqone$}};
			\draw [thick, blue] plot coordinates {(0,0.14/2+3/8) (0,0.14/2-8/8) } ;
			\fill [black] (-3/4+0.14/3,0.14/2) circle (0.06);
			\fill [black] (3/4-0.14/3,0.14/2) circle (0.06);
			\fill [black] (0,0.14/2) circle (0.06);
			
		\end{scope}
		\begin{scope}[shift={(-2/2,-1)}, yscale=-1]
				\draw [thick, orange] plot coordinates {(-3/4+0.14/3-3/8,0.14/2) (3/4-0.14/3+3/8,0.14/2) };
			\draw (0,-7/8) node [above] {\hbox{\footnotesize$\geqone$}};
			\draw [thick, blue] plot coordinates {(0,0.14/2+3/8) (0,0.14/2-8/8) } ;
			\fill [black] (-3/4+0.14/3,0.14/2) circle (0.06);
			\fill [black] (3/4-0.14/3,0.14/2) circle (0.06);
				\fill [black] (0,0.14/2) circle (0.06);
		\end{scope}

\draw [thick, black, shorten >= -0.25cm, shorten <=-0.25cm] (0.14/3+1/2,0.14/2+1/2)--(-0.14/3 -1/4,-0.14/2-1);

	\begin{scope} [shift={(-1/8+1/4,-1/4 )}]
	\fill [black] (0,0) circle (0.06);
	\draw [thick,black] plot  coordinates {(4/12,-8/8) (-1/12,1/4) };
	\begin{scope}[shift={(1/9,-3/9)}]
		\draw [thick,black] plot  coordinates {(8/64,3/64) (-8/64,-3/64) };
	\end{scope}		
	\begin{scope}[shift={(2/9,-6/9)}]
		\draw [thick,black] plot  coordinates {(8/64,3/64) (-8/64,-3/64) };
	\end{scope}
\end{scope}

		\begin{scope}[shift={(3/2+1/2,1/2)}]
					\begin{scope}[shift={(-1/48,0)}]
			\draw [thick,black] plot  coordinates {(4/12,-8/8) (-1/12,1/4) };
			\begin{scope}[shift={(1/9,-3/9)}]
				\draw [thick,black] plot  coordinates {(8/64,3/64) (-8/64,-3/64) };
			\end{scope}		
			\begin{scope}[shift={(2/9,-6/9)}]
				\draw [thick,black] plot  coordinates {(8/64,3/64) (-8/64,-3/64) };
			\end{scope}
		\end{scope}
		\end{scope}
		\begin{scope}[shift={(-3/2-1/4,-1)}, yscale=-1, xscale=-1]
						\begin{scope}[shift={(-1/48,0)}]
				\draw [thick,black] plot  coordinates {(4/12,-8/8) (-1/12,1/4) };
				\begin{scope}[shift={(1/9,-3/9)}]
					\draw [thick,black] plot  coordinates {(8/64,3/64) (-8/64,-3/64) };
				\end{scope}		
				\begin{scope}[shift={(2/9,-6/9)}]
					\draw [thick,black] plot  coordinates {(8/64,3/64) (-8/64,-3/64) };
				\end{scope}
			\end{scope}
		\end{scope}
\end{tikzpicture}}

\def\FigTwoCaseDxi{\begin{tikzpicture}[baseline=13pt]
		\draw [thick,black] plot  coordinates {(-8/6,-8/4) (1/6,1/4) };
			\begin{scope} [shift={(-10/27,-15/27 )}]
				\draw [thick,black] plot  coordinates {(8/12,-8/8) (-1/6,1/4) };
			\fill [black] (0,0) circle (0.06);
			\begin{scope}[shift={(2/9,-3/9)}]
				\draw [thick,black] plot  coordinates {(3/32,2/32) (-3/32,-2/32) };
			\end{scope}
			\begin{scope}[shift={(4/9,-6/9)}]
				\draw [thick,black] plot  coordinates {(3/32,2/32) (-3/32,-2/32) };
			\end{scope}
		\end{scope}
			\begin{scope} [shift={(-20/27,-30/27 )}]
				\draw [thick,black] plot  coordinates {(8/12,-8/8) (-1/6,1/4) };
			\fill [black] (0,0) circle (0.06);
			\begin{scope}[shift={(2/9,-3/9)}]
				\draw [thick,black] plot  coordinates {(3/32,2/32) (-3/32,-2/32) };
			\end{scope}
			\begin{scope}[shift={(4/9,-6/9)}]
				\draw [thick,black] plot  coordinates {(3/32,2/32) (-3/32,-2/32) };
			\end{scope}
		\end{scope}
		\draw [thick,orange] plot  coordinates {(-1/2,0) (4/2,0) };
		\fill [black] (0,0) circle (0.06);
		\draw [thick,blue] plot  coordinates {(3/4,1/3) (3/4,-6/4) };
		\draw (3/4,-3/2) node [right] {\hbox{\footnotesize$\geqone$}};
		\fill [black] (3/4,0) circle (0.06);
		\begin{scope}[shift={(3/2,0)}]
			\draw [thick,black] plot  coordinates {(8/12,-8/8) (-1/6,1/4) };
			\fill [black] (0,0) circle (0.06);
			\begin{scope}[shift={(2/9,-3/9)}]
				\draw [thick,black] plot  coordinates {(3/32,2/32) (-3/32,-2/32) };
			\end{scope}
			\begin{scope}[shift={(4/9,-6/9)}]
				\draw [thick,black] plot  coordinates {(3/32,2/32) (-3/32,-2/32) };
			\end{scope}
		\end{scope}
		\draw [thick,blue] plot  coordinates {(-10/9-5/6,-15/9+5/4) (-10/9+1/6,-15/9-1/4) };
		\fill [black] (-10/9,-15/9) circle (0.06);
		\draw (-10/9-5/6,-15/9+5/4+1/6) node [] {\hbox{\footnotesize$\geqone$}};
\end{tikzpicture}}

\def\FigTwoCaseDxii{\begin{tikzpicture}[baseline=13pt]
		\draw [thick,black] plot  coordinates {(-16/12,-16/8) (1/6,1/4) };
		\begin{scope} [shift={(-10/27,-15/27 )}]
			\draw [thick,black] plot  coordinates {(8/12,-8/8) (-1/6,1/4) };
			\fill [black] (0,0) circle (0.06);
			\begin{scope}[shift={(2/9,-3/9)}]
				\draw [thick,black] plot  coordinates {(3/32,2/32) (-3/32,-2/32) };
			\end{scope}
			\begin{scope}[shift={(4/9,-6/9)}]
				\draw [thick,black] plot  coordinates {(3/32,2/32) (-3/32,-2/32) };
			\end{scope}
		\end{scope}
		\begin{scope} [shift={(-20/27,-30/27 )}]
			\draw [thick,black] plot  coordinates {(8/12,-8/8) (-1/6,1/4) };
			\fill [black] (0,0) circle (0.06);
			\begin{scope}[shift={(2/9,-3/9)}]
				\draw [thick,black] plot  coordinates {(3/32,2/32) (-3/32,-2/32) };
			\end{scope}
			\begin{scope}[shift={(4/9,-6/9)}]
				\draw [thick,black] plot  coordinates {(3/32,2/32) (-3/32,-2/32) };
			\end{scope}
		\end{scope}
		\draw [thick,orange] plot  coordinates {(-1/2,0) (4/2,0) };
		\draw (3/4,0) node [below] {\hbox{\footnotesize$\geqone$}};
		\fill [black] (0,0) circle (0.06);
		\begin{scope}[shift={(3/2,0)}]
			\draw [thick,black] plot  coordinates {(8/12,-8/8) (-1/6,1/4) };
			\fill [black] (0,0) circle (0.06);
			\begin{scope}[shift={(2/9,-3/9)}]
				\draw [thick,black] plot  coordinates {(3/32,2/32) (-3/32,-2/32) };
			\end{scope}
			\begin{scope}[shift={(4/9,-6/9)}]
				\draw [thick,black] plot  coordinates {(3/32,2/32) (-3/32,-2/32) };
			\end{scope}
		\end{scope}
		\draw [thick,blue] plot  coordinates {(-10/9-5/6,-15/9+5/4) (-10/9+1/6,-15/9-1/4) };
		\fill [black] (-10/9,-15/9) circle (0.06);
		\draw (-10/9-5/6,-15/9+5/4+1/6) node [] {\hbox{\footnotesize$\geqone$}};
\end{tikzpicture}}

\def\FigTwoCaseDxiii{\begin{tikzpicture}[baseline=13pt]
		\draw [thick,black] plot  coordinates {(-16/12,-16/8) (1/6,1/4) };
\begin{scope} [shift={(-10/27,-15/27 )}]
	\draw [thick,black] plot  coordinates {(8/12,-8/8) (-1/6,1/4) };
	\fill [black] (0,0) circle (0.06);
	\begin{scope}[shift={(2/9,-3/9)}]
		\draw [thick,black] plot  coordinates {(3/32,2/32) (-3/32,-2/32) };
	\end{scope}
	\begin{scope}[shift={(4/9,-6/9)}]
		\draw [thick,black] plot  coordinates {(3/32,2/32) (-3/32,-2/32) };
	\end{scope}
\end{scope}
\begin{scope} [shift={(-20/27,-30/27 )}]
	\draw [thick,black] plot  coordinates {(8/12,-8/8) (-1/6,1/4) };
	\fill [black] (0,0) circle (0.06);
	\begin{scope}[shift={(2/9,-3/9)}]
		\draw [thick,black] plot  coordinates {(3/32,2/32) (-3/32,-2/32) };
	\end{scope}
	\begin{scope}[shift={(4/9,-6/9)}]
		\draw [thick,black] plot  coordinates {(3/32,2/32) (-3/32,-2/32) };
	\end{scope}
\end{scope}
		\begin{scope}[shift={(3/4,0)}]
			\draw [thick,orange] plot [smooth, tension=1] coordinates {(-3/3,2/12) (1/4,-6/12) (-1/3,-14/12) };
			\draw [thick,orange] plot [smooth, tension=1] coordinates {(3/3,2/12) (-1/4,-6/12) (1/3,-14/12) };
			\fill [black] (0,-0.26) circle (0.06);
			\fill [black] (0,-0.99) circle (0.06);
			\fill [black] (-3/4+0.14/3,0.14/2) circle (0.06);
			\fill [black] (3/4-0.14/3,0.14/2) circle (0.06);
		\end{scope}
		\begin{scope}[shift={(3/2,0)}]
			\draw [thick,black] plot  coordinates {(8/12,-8/8) (-1/6,1/4) };
			\begin{scope}[shift={(2/9,-3/9)}]
				\draw [thick,black] plot  coordinates {(3/32,2/32) (-3/32,-2/32) };
			\end{scope}
			\begin{scope}[shift={(4/9,-6/9)}]
				\draw [thick,black] plot  coordinates {(3/32,2/32) (-3/32,-2/32) };
			\end{scope}
		\end{scope}
		\draw [thick,blue] plot  coordinates {(-10/9-5/6,-15/9+5/4) (-10/9+1/6,-15/9-1/4) };
		\fill [black] (-10/9,-15/9) circle (0.06);
		\draw (-10/9-5/6,-15/9+5/4+1/6) node [] {\hbox{\footnotesize$\geqone$}};
\end{tikzpicture}}

\def\FigTwoCaseDxiv{\begin{tikzpicture}[baseline=13pt]
		\draw [thick,black] plot  coordinates {(-14/12,-14/8) (1/6,1/4) };
		\begin{scope} [shift={(-4/9,-6/9 )}, xscale=-1, yscale=-1]
			\draw [thick,black] plot  coordinates {(8/12,-8/8) (-1/6,1/4) };
			\fill [black] (0,0) circle (0.06);
			\begin{scope}[shift={(2/9,-3/9)}]
				\draw [thick,black] plot  coordinates {(3/32,2/32) (-3/32,-2/32) };
			\end{scope}
			\begin{scope}[shift={(4/9,-6/9)}]
				\draw [thick,black] plot  coordinates {(3/32,2/32) (-3/32,-2/32) };
			\end{scope}
		\end{scope}
		\begin{scope} [shift={(-7/9,-10.5/9 )}, xscale=-1,yscale=-1 ]
			\draw [thick,black] plot  coordinates {(8/12,-8/8) (-1/6,1/4) };
			\fill [black] (0,0) circle (0.06);
			\begin{scope}[shift={(2/9,-3/9)}]
				\draw [thick,black] plot  coordinates {(3/32,2/32) (-3/32,-2/32) };
			\end{scope}
			\begin{scope}[shift={(4/9,-6/9)}]
				\draw [thick,black] plot  coordinates {(3/32,2/32) (-3/32,-2/32) };
			\end{scope}
		\end{scope}
		\draw [thick,orange] plot  coordinates {(-1/4,1/6) (6/4,-6/6) };
		\draw [thick,orange] plot  coordinates {(5/2+1/4,1/6) (5/2-6/4,-6/6) };
		\fill [black] (0,0) circle (0.06);
		\fill [black] (5/4,-5/6) circle (0.06);
		\begin{scope}[shift={(3/5,-2/5)}]
			\draw [thick,blue] plot coordinates {(1/12,1/3) (-7/24,-7/6) };
			\fill [black] (0,0) circle (0.06);
			\draw (-6.5/24,-6.5/6) node [below] {\hbox{\footnotesize$\geqone$}};
		\end{scope}
		\begin{scope}[shift={(5/2-3/5,-2/5)}]
			\draw [thick,blue] plot coordinates {(-1/12,1/3) (7/24,-7/6) };
			\fill [black] (0,0) circle (0.06);
			\draw (6.5/24,-6.5/6) node [below] {\hbox{\footnotesize$\geqone$}};
		\end{scope}
		\begin{scope}[shift={(5/2,0)}]
			\draw [thick,black] plot  coordinates {(8/12,-8/8) (-1/6,1/4) };
			\fill [black] (0,0) circle (0.06);
			\begin{scope}[shift={(2/9,-3/9)}]
				\draw [thick,black] plot  coordinates {(3/32,2/32) (-3/32,-2/32) };
			\end{scope}
			\begin{scope}[shift={(4/9,-6/9)}]
				\draw [thick,black] plot  coordinates {(3/32,2/32) (-3/32,-2/32) };
			\end{scope}
		\end{scope}
\end{tikzpicture}}

\def\FigTwoCaseDxv{\begin{tikzpicture}[baseline=13pt]
		\draw [thick,black] plot  coordinates {(-14/12,-14/8) (1/6,1/4) };
		\begin{scope} [shift={(-4/9,-6/9 )}, xscale=-1, yscale=-1]
			\draw [thick,black] plot  coordinates {(8/12,-8/8) (-1/6,1/4) };
			\fill [black] (0,0) circle (0.06);
			\begin{scope}[shift={(2/9,-3/9)}]
				\draw [thick,black] plot  coordinates {(3/32,2/32) (-3/32,-2/32) };
			\end{scope}
			\begin{scope}[shift={(4/9,-6/9)}]
				\draw [thick,black] plot  coordinates {(3/32,2/32) (-3/32,-2/32) };
			\end{scope}
		\end{scope}
		\begin{scope} [shift={(-7/9,-10.5/9 )}, xscale=-1,yscale=-1 ]
			\draw [thick,black] plot  coordinates {(8/12,-8/8) (-1/6,1/4) };
			\fill [black] (0,0) circle (0.06);
			\begin{scope}[shift={(2/9,-3/9)}]
				\draw [thick,black] plot  coordinates {(3/32,2/32) (-3/32,-2/32) };
			\end{scope}
			\begin{scope}[shift={(4/9,-6/9)}]
				\draw [thick,black] plot  coordinates {(3/32,2/32) (-3/32,-2/32) };
			\end{scope}
		\end{scope}
		\draw [thick,orange] plot  coordinates {(-1/4,1/6) (6/4,-6/6) };
		\draw [thick,orange] plot  coordinates {(5/2+1/4,1/6) (5/2-6/4,-6/6) };
		\fill [black] (0,0) circle (0.06);
		\fill [black] (5/4,-5/6) circle (0.06);
		\begin{scope}[shift={(3/5,-2/5)}]
			\draw [thick,blue] plot coordinates {(1/12,1/3) (-7/24,-7/6) };
			\fill [black] (0,0) circle (0.06);
			\draw (-6.5/24,-6.5/6) node [right] {\hbox{\footnotesize$\geqone$}};
		\end{scope}
		\draw (1+4/4+1/12,-1+4/6) node [below] {\hbox{\footnotesize$\geqone$}};
		\begin{scope}[shift={(5/2,0)}]
			\draw [thick,black] plot  coordinates {(8/12,-8/8) (-1/6,1/4) };
			\fill [black] (0,0) circle (0.06);
			\begin{scope}[shift={(2/9,-3/9)}]
				\draw [thick,black] plot  coordinates {(3/32,2/32) (-3/32,-2/32) };
			\end{scope}
			\begin{scope}[shift={(4/9,-6/9)}]
				\draw [thick,black] plot  coordinates {(3/32,2/32) (-3/32,-2/32) };
			\end{scope}
		\end{scope}
\end{tikzpicture}}

\def\FigTwoCaseDxvi{\begin{tikzpicture}[baseline=13pt]
		\draw [thick,black] plot  coordinates {(-1-14/12,-14/8) (-1+1/6,1/4) };
		\begin{scope} [shift={(-4/9-1,-6/9 )}, xscale=-1, yscale=-1]
			\draw [thick,black] plot  coordinates {(8/12,-8/8) (-1/6,1/4) };
			\fill [black] (0,0) circle (0.06);
			\begin{scope}[shift={(2/9,-3/9)}]
				\draw [thick,black] plot  coordinates {(3/32,2/32) (-3/32,-2/32) };
			\end{scope}
			\begin{scope}[shift={(4/9,-6/9)}]
				\draw [thick,black] plot  coordinates {(3/32,2/32) (-3/32,-2/32) };
			\end{scope}
		\end{scope}
		\begin{scope} [shift={(-7/9-1,-10.5/9 )}, xscale=-1, yscale=-1]
			\draw [thick,black] plot  coordinates {(8/12,-8/8) (-1/6,1/4) };
			\fill [black] (0,0) circle (0.06);
			\begin{scope}[shift={(2/9,-3/9)}]
				\draw [thick,black] plot  coordinates {(3/32,2/32) (-3/32,-2/32) };
			\end{scope}
			\begin{scope}[shift={(4/9,-6/9)}]
				\draw [thick,black] plot  coordinates {(3/32,2/32) (-3/32,-2/32) };
			\end{scope}
		\end{scope}
		\draw [thick,orange] plot [smooth, tension=1] coordinates {(4/12,1/4) (-6/12,-1/4) (-16/12,1/4) };
		\begin{scope}[shift={(3/4,0)}]
			\draw [thick,orange] plot [smooth, tension=1] coordinates {(-3/3,2/12) (1/4,-6/12) (-1/3,-14/12) };
			\draw [thick,orange] plot [smooth, tension=1] coordinates {(3/3,2/12) (-1/4,-6/12) (1/3,-14/12) };
			\fill [black] (0,-0.26) circle (0.06);
			\fill [black] (0,-0.99) circle (0.06);
			\fill [black] (3/4-0.14/3,0.14/2) circle (0.06);
		\end{scope}
		\fill [black] (-5/6-0.6/3,1/4-0.6/2) circle (0.06);
		\fill [black] (0.14,0.036) circle (0.06);
		\begin{scope}[shift={(3/2,0)}]
			\draw [thick,black] plot  coordinates {(8/12,-8/8) (-1/6,1/4) };
			\begin{scope}[shift={(2/9,-3/9)}]
				\draw [thick,black] plot  coordinates {(3/32,2/32) (-3/32,-2/32) };
			\end{scope}
			\begin{scope}[shift={(4/9,-6/9)}]
				\draw [thick,black] plot  coordinates {(3/32,2/32) (-3/32,-2/32) };
			\end{scope}
		\end{scope}
		\draw [thick,blue] plot  coordinates {(-6/12,-15/9+7/4) (-6/12,-15/9+1/4) };
		\fill [black] (-6/12,-.25) circle (0.06);
		\draw (-6/6,-10/6) node [right] {\hbox{\footnotesize$\geqone$}};
\end{tikzpicture}}

\def\FigTwoCaseDxvii{\begin{tikzpicture}[baseline=13pt]
		\draw [thick,black] plot  coordinates {(-14/12,-14/8) (1/6,1/4) };
		\begin{scope} [shift={(-4/9,-6/9 )}, xscale=-1, yscale=-1]
			\draw [thick,black] plot  coordinates {(8/12,-8/8) (-1/6,1/4) };
			\fill [black] (0,0) circle (0.06);
			\begin{scope}[shift={(2/9,-3/9)}]
				\draw [thick,black] plot  coordinates {(3/32,2/32) (-3/32,-2/32) };
			\end{scope}
			\begin{scope}[shift={(4/9,-6/9)}]
				\draw [thick,black] plot  coordinates {(3/32,2/32) (-3/32,-2/32) };
			\end{scope}
		\end{scope}
		\begin{scope} [shift={(-7/9,-10.5/9 )},xscale=-1, yscale=-1]
			\draw [thick,black] plot  coordinates {(8/12,-8/8) (-1/6,1/4) };
			\fill [black] (0,0) circle (0.06);
			\begin{scope}[shift={(2/9,-3/9)}]
				\draw [thick,black] plot  coordinates {(3/32,2/32) (-3/32,-2/32) };
			\end{scope}
			\begin{scope}[shift={(4/9,-6/9)}]
				\draw [thick,black] plot  coordinates {(3/32,2/32) (-3/32,-2/32) };
			\end{scope}
		\end{scope}
		\begin{scope}[shift={(3/4,0)}]
			\draw [thick,orange] plot [smooth, tension=1] coordinates {(-3/3,2/12) (1,-8/12) (-1/2,-20/12) };
			\draw [thick,orange] plot [smooth, tension=.8] coordinates {(1/2+3/3,2/12) (5/12,-8/12) (11/12,-16/12) };
			\fill [black] (0.69,-0.36) circle (0.06);
			\fill [black] (0.66,-1.16) circle (0.06);
			\draw (-1/2+1/4,-20/12+1/8) node [below] {\hbox{\footnotesize$\geqone$}};
		\end{scope}
		\fill [black] (0.195/3,0.195/2) circle (0.06);
		\begin{scope}[shift={(2,0)}]
			\draw [thick,black] plot  coordinates {(8/12,-8/8) (-1/6,1/4) };
			\fill [black] (-0.035/3,0.035/2) circle (0.06);
			\begin{scope}[shift={(2/9,-3/9)}]
				\draw [thick,black] plot  coordinates {(3/32,2/32) (-3/32,-2/32) };
			\end{scope}
			\begin{scope}[shift={(4/9,-6/9)}]
				\draw [thick,black] plot  coordinates {(3/32,2/32) (-3/32,-2/32) };
			\end{scope}
		\end{scope}
\end{tikzpicture}}

\def\FigTwoCaseDxviii{\begin{tikzpicture}[baseline=15pt]
		\draw [thick,black] plot  coordinates {(-14/12,-14/8) (1/6,1/4) };
		\begin{scope} [shift={(-4/9,-6/9 )}, xscale=-1, yscale=-1]
			\draw [thick,black] plot  coordinates {(8/12,-8/8) (-1/6,1/4) };
			\fill [black] (0,0) circle (0.06);
			\begin{scope}[shift={(2/9,-3/9)}]
				\draw [thick,black] plot  coordinates {(3/32,2/32) (-3/32,-2/32) };
			\end{scope}
			\begin{scope}[shift={(4/9,-6/9)}]
				\draw [thick,black] plot  coordinates {(3/32,2/32) (-3/32,-2/32) };
			\end{scope}
		\end{scope}
		\begin{scope} [shift={(-7/9,-10.5/9 )},xscale=-1, yscale=-1]
			\draw [thick,black] plot  coordinates {(8/12,-8/8) (-1/6,1/4) };
			\fill [black] (0,0) circle (0.06);
			\begin{scope}[shift={(2/9,-3/9)}]
				\draw [thick,black] plot  coordinates {(3/32,2/32) (-3/32,-2/32) };
			\end{scope}
			\begin{scope}[shift={(4/9,-6/9)}]
				\draw [thick,black] plot  coordinates {(3/32,2/32) (-3/32,-2/32) };
			\end{scope}
		\end{scope}
		\draw [thick,orange] plot  coordinates {(-1/2,0) (5/2,0) };
		\fill [black] (0,0) circle (0.06);
		\begin{scope}[shift={(4/2,0)}]
			\draw [thick,black] plot  coordinates {(8/12,-8/8) (-1/6,1/4) };
			\fill [black] (0,0) circle (0.06);
			\begin{scope}[shift={(2/9,-3/9)}]
				\draw [thick,black] plot  coordinates {(3/32,2/32) (-3/32,-2/32) };
			\end{scope}
			\begin{scope}[shift={(4/9,-6/9)}]
				\draw [thick,black] plot  coordinates {(3/32,2/32) (-3/32,-2/32) };
			\end{scope}
		\end{scope}
		\begin{scope}[shift={(2/3,0)}]
			\draw [thick,blue] plot  coordinates {(0,1/3) (0,-6/4) };
			\fill [black] (0,0) circle (0.06);
			\draw (0,-6/4) node [left] {\hbox{\footnotesize$\geqone$\kern-1pt}};
		\end{scope}
		\begin{scope}[shift={(4/3,0)}]
			\draw [thick,blue] plot  coordinates {(0,1/3) (0,-6/4) };
			\fill [black] (0,0) circle (0.06);
			\draw (0,-6/4) node [right] {\hbox{\footnotesize$\geqone$}};
		\end{scope}
\end{tikzpicture}}

\def\FigTwoCaseDxix{\begin{tikzpicture}[baseline=15pt]
		\draw [thick,black] plot  coordinates {(-14/12,-14/8) (1/6,1/4) };
		\begin{scope} [shift={(-4/9,-6/9 )}, xscale=-1, yscale=-1]
			\draw [thick,black] plot  coordinates {(8/12,-8/8) (-1/6,1/4) };
			\fill [black] (0,0) circle (0.06);
			\begin{scope}[shift={(2/9,-3/9)}]
				\draw [thick,black] plot  coordinates {(3/32,2/32) (-3/32,-2/32) };
			\end{scope}
			\begin{scope}[shift={(4/9,-6/9)}]
				\draw [thick,black] plot  coordinates {(3/32,2/32) (-3/32,-2/32) };
			\end{scope}
		\end{scope}
		\begin{scope} [shift={(-7/9,-10.5/9 )}, xscale=-1, yscale=-1]
			\draw [thick,black] plot  coordinates {(8/12,-8/8) (-1/6,1/4) };
			\fill [black] (0,0) circle (0.06);
			\begin{scope}[shift={(2/9,-3/9)}]
				\draw [thick,black] plot  coordinates {(3/32,2/32) (-3/32,-2/32) };
			\end{scope}
			\begin{scope}[shift={(4/9,-6/9)}]
				\draw [thick,black] plot  coordinates {(3/32,2/32) (-3/32,-2/32) };
			\end{scope}
		\end{scope}
		\draw [thick,orange] plot  coordinates {(-1/2,0) (5/2,0) };
		\fill [black] (0,0) circle (0.06);
		\begin{scope}[shift={(4/2,0)}]
			\draw [thick,black] plot  coordinates {(8/12,-8/8) (-1/6,1/4) };
			\fill [black] (0,0) circle (0.06);
			\begin{scope}[shift={(2/9,-3/9)}]
				\draw [thick,black] plot  coordinates {(3/32,2/32) (-3/32,-2/32) };
			\end{scope}
			\begin{scope}[shift={(4/9,-6/9)}]
				\draw [thick,black] plot  coordinates {(3/32,2/32) (-3/32,-2/32) };
			\end{scope}
		\end{scope}
		\begin{scope}[shift={(1/2,0)}]
			\draw [thick,green] plot  coordinates {(0,1/3) (0,-8/4) };
			\fill [black] (0,0) circle (0.06);
		\end{scope}
		\begin{scope}[shift={(1/2,-5/6)}]
			\draw [thick,blue] plot  coordinates {(-1/3,0) (4/3,0)};
			\fill [black] (0,0) circle (0.06);
			\draw (2/2,0) node [below] {\hbox{\footnotesize$\geqone$\kern-1pt}};
		\end{scope}
		\begin{scope}[shift={(1/2,-5/3)}]
			\draw [thick,blue] plot  coordinates {(-1/3,0) (4/3,0)};
			\fill [black] (0,0) circle (0.06);
			\draw (2/2,0) node [below] {\hbox{\footnotesize$\geqone$\kern-1pt}};
		\end{scope}
\end{tikzpicture}}

\def\FigTwoCaseDxx{\begin{tikzpicture}[baseline=15pt]
		\draw [thick,black] plot  coordinates {(-14/12,-14/8) (1/6,1/4) };
		\begin{scope} [shift={(-4/9,-6/9 )}, xscale=-1, yscale=-1]
			\draw [thick,black] plot  coordinates {(8/12,-8/8) (-1/6,1/4) };
			\fill [black] (0,0) circle (0.06);
			\begin{scope}[shift={(2/9,-3/9)}]
				\draw [thick,black] plot  coordinates {(3/32,2/32) (-3/32,-2/32) };
			\end{scope}
			\begin{scope}[shift={(4/9,-6/9)}]
				\draw [thick,black] plot  coordinates {(3/32,2/32) (-3/32,-2/32) };
			\end{scope}
		\end{scope}
		\begin{scope} [shift={(-7/9,-10.5/9 )}, xscale=-1, yscale=-1]
			\draw [thick,black] plot  coordinates {(8/12,-8/8) (-1/6,1/4) };
			\fill [black] (0,0) circle (0.06);
			\begin{scope}[shift={(2/9,-3/9)}]
				\draw [thick,black] plot  coordinates {(3/32,2/32) (-3/32,-2/32) };
			\end{scope}
			\begin{scope}[shift={(4/9,-6/9)}]
				\draw [thick,black] plot  coordinates {(3/32,2/32) (-3/32,-2/32) };
			\end{scope}
		\end{scope}
		\draw [thick,orange] plot  coordinates {(-1/2,0) (5/2,0) };
		\fill [black] (0,0) circle (0.06);
		\begin{scope}[shift={(4/2,0)}]
			\draw [thick,black] plot  coordinates {(8/12,-8/8) (-1/6,1/4) };
			\fill [black] (0,0) circle (0.06);
			\begin{scope}[shift={(2/9,-3/9)}]
				\draw [thick,black] plot  coordinates {(3/32,2/32) (-3/32,-2/32) };
			\end{scope}
			\begin{scope}[shift={(4/9,-6/9)}]
				\draw [thick,black] plot  coordinates {(3/32,2/32) (-3/32,-2/32) };
			\end{scope}
		\end{scope}
		\begin{scope}[shift={(1,0)}]
			\draw [thick,blue] plot  coordinates {(0,1/3) (0,-6/4) };
			\fill [black] (0,0) circle (0.06);
			\draw (0,-6/4) node [right] {\hbox{\footnotesize$\geqtwo$}};
		\end{scope}
\end{tikzpicture}}

\def\FigTwoCaseDxxi{\begin{tikzpicture}[baseline=15pt]
		\draw [thick,black] plot  coordinates {(-14/12,-14/8) (1/6,1/4) };
	\begin{scope} [shift={(-4/9,-6/9 )}, xscale=-1, yscale=-1]
		\draw [thick,black] plot  coordinates {(8/12,-8/8) (-1/6,1/4) };
		\fill [black] (0,0) circle (0.06);
		\begin{scope}[shift={(2/9,-3/9)}]
			\draw [thick,black] plot  coordinates {(3/32,2/32) (-3/32,-2/32) };
		\end{scope}
		\begin{scope}[shift={(4/9,-6/9)}]
			\draw [thick,black] plot  coordinates {(3/32,2/32) (-3/32,-2/32) };
		\end{scope}
	\end{scope}
	\begin{scope} [shift={(-7/9,-10.5/9 )}, xscale=-1, yscale=-1]
		\draw [thick,black] plot  coordinates {(8/12,-8/8) (-1/6,1/4) };
		\fill [black] (0,0) circle (0.06);
		\begin{scope}[shift={(2/9,-3/9)}]
			\draw [thick,black] plot  coordinates {(3/32,2/32) (-3/32,-2/32) };
		\end{scope}
		\begin{scope}[shift={(4/9,-6/9)}]
			\draw [thick,black] plot  coordinates {(3/32,2/32) (-3/32,-2/32) };
		\end{scope}
	\end{scope}
		\draw [thick,orange] plot  coordinates {(-1/2,0) (5/2,0) };
		\fill [black] (0,0) circle (0.06);
		\begin{scope}[shift={(4/2,0)}]
			\draw [thick,black] plot  coordinates {(8/12,-8/8) (-1/6,1/4) };
			\fill [black] (0,0) circle (0.06);
			\begin{scope}[shift={(2/9,-3/9)}]
				\draw [thick,black] plot  coordinates {(3/32,2/32) (-3/32,-2/32) };
			\end{scope}
			\begin{scope}[shift={(4/9,-6/9)}]
				\draw [thick,black] plot  coordinates {(3/32,2/32) (-3/32,-2/32) };
			\end{scope}
		\end{scope}
		\draw (1,0) node [below] {\hbox{\footnotesize$\geqtwo$}};
\end{tikzpicture}}

\def\FigTwoCaseDxxii{\begin{tikzpicture}[baseline=15pt]
		\begin{scope} [shift={(-3/4,0 )}, xscale=-1, yscale=1]
			\draw [thick,black] plot  coordinates {(8/12,-8/8) (-1/6,1/4) };
			\fill [black] (0,0) circle (0.06);
			\begin{scope}[shift={(2/9,-3/9)}]
				\draw [thick,black] plot  coordinates {(3/32,2/32) (-3/32,-2/32) };
			\end{scope}
			\begin{scope}[shift={(4/9,-6/9)}]
				\draw [thick,black] plot  coordinates {(3/32,2/32) (-3/32,-2/32) };
			\end{scope}
		\end{scope}
		\begin{scope} [shift={(0,0 )}, xscale=-1, yscale=1]
		\draw [thick,black] plot  coordinates {(8/12,-8/8) (-1/6,1/4) };
		\fill [black] (0,0) circle (0.06);
		\begin{scope}[shift={(2/9,-3/9)}]
			\draw [thick,black] plot  coordinates {(3/32,2/32) (-3/32,-2/32) };
		\end{scope}
		\begin{scope}[shift={(4/9,-6/9)}]
			\draw [thick,black] plot  coordinates {(3/32,2/32) (-3/32,-2/32) };
		\end{scope}
	\end{scope}
		\draw [thick,black] plot  coordinates {(-3/2+1/4,0) (5/2,0) };
		\begin{scope}[shift={(4/2,0)}]
			\draw [thick,black] plot  coordinates {(8/12,-8/8) (-1/6,1/4) };
			\fill [black] (0,0) circle (0.06);
			\begin{scope}[shift={(2/9,-3/9)}]
				\draw [thick,black] plot  coordinates {(3/32,2/32) (-3/32,-2/32) };
			\end{scope}
			\begin{scope}[shift={(4/9,-6/9)}]
				\draw [thick,black] plot  coordinates {(3/32,2/32) (-3/32,-2/32) };
			\end{scope}
		\end{scope}
		\begin{scope}[shift={(2/3,0)}]
			\draw [thick,blue] plot  coordinates {(0,1/3) (0,-6/4) };
			\fill [black] (0,0) circle (0.06);
			\draw (0,-6/4) node [left] {\hbox{\footnotesize$\geqone$\kern-1pt}};
		\end{scope}
		\begin{scope}[shift={(4/3,0)}]
			\draw [thick,blue] plot  coordinates {(0,1/3) (0,-6/4) };
			\fill [black] (0,0) circle (0.06);
			\draw (0,-6/4) node [right] {\hbox{\footnotesize$\geqone$}};
		\end{scope}
\end{tikzpicture}}

\def\FigTwoCaseDxxiii{\begin{tikzpicture}[baseline=15pt]
		\begin{scope} [shift={(-3/4,0 )}, xscale=-1, yscale=1]
			\draw [thick,black] plot  coordinates {(8/12,-8/8) (-1/6,1/4) };
			\fill [black] (0,0) circle (0.06);
			\begin{scope}[shift={(2/9,-3/9)}]
				\draw [thick,black] plot  coordinates {(3/32,2/32) (-3/32,-2/32) };
			\end{scope}
			\begin{scope}[shift={(4/9,-6/9)}]
				\draw [thick,black] plot  coordinates {(3/32,2/32) (-3/32,-2/32) };
			\end{scope}
		\end{scope}
		\begin{scope} [shift={(0,0 )}, xscale=-1, yscale=1]
			\draw [thick,black] plot  coordinates {(8/12,-8/8) (-1/6,1/4) };
			\fill [black] (0,0) circle (0.06);
			\begin{scope}[shift={(2/9,-3/9)}]
				\draw [thick,black] plot  coordinates {(3/32,2/32) (-3/32,-2/32) };
			\end{scope}
			\begin{scope}[shift={(4/9,-6/9)}]
				\draw [thick,black] plot  coordinates {(3/32,2/32) (-3/32,-2/32) };
			\end{scope}
		\end{scope}
		\draw [thick,black] plot  coordinates {(-3/2+1/4,0) (5/2,0) };
		\begin{scope}[shift={(4/2,0)}]
			\draw [thick,black] plot  coordinates {(8/12,-8/8) (-1/6,1/4) };
			\fill [black] (0,0) circle (0.06);
			\begin{scope}[shift={(2/9,-3/9)}]
				\draw [thick,black] plot  coordinates {(3/32,2/32) (-3/32,-2/32) };
			\end{scope}
			\begin{scope}[shift={(4/9,-6/9)}]
				\draw [thick,black] plot  coordinates {(3/32,2/32) (-3/32,-2/32) };
			\end{scope}
		\end{scope}
\begin{scope}[shift={(1/2,0)}]
	\draw [thick,green] plot  coordinates {(0,1/3) (0,-8/4) };
	\fill [black] (0,0) circle (0.06);
\end{scope}
\begin{scope}[shift={(1/2,-5/6)}]
	\draw [thick,blue] plot  coordinates {(-1/3,0) (4/3,0)};
	\fill [black] (0,0) circle (0.06);
	\draw (2/2,0) node [below] {\hbox{\footnotesize$\geqone$\kern-1pt}};
\end{scope}
\begin{scope}[shift={(1/2,-5/3)}]
	\draw [thick,blue] plot  coordinates {(-1/3,0) (4/3,0)};
	\fill [black] (0,0) circle (0.06);
	\draw (2/2,0) node [below] {\hbox{\footnotesize$\geqone$\kern-1pt}};
\end{scope}
\end{tikzpicture}}

\def\FigTwoCaseDxxiv{\begin{tikzpicture}[baseline=15pt]
		\begin{scope} [shift={(-3/4,0 )}, xscale=-1, yscale=1]
			\draw [thick,black] plot  coordinates {(8/12,-8/8) (-1/6,1/4) };
			\fill [black] (0,0) circle (0.06);
			\begin{scope}[shift={(2/9,-3/9)}]
				\draw [thick,black] plot  coordinates {(3/32,2/32) (-3/32,-2/32) };
			\end{scope}
			\begin{scope}[shift={(4/9,-6/9)}]
				\draw [thick,black] plot  coordinates {(3/32,2/32) (-3/32,-2/32) };
			\end{scope}
		\end{scope}
		\begin{scope} [shift={(0,0 )}, xscale=-1, yscale=1]
			\draw [thick,black] plot  coordinates {(8/12,-8/8) (-1/6,1/4) };
			\fill [black] (0,0) circle (0.06);
			\begin{scope}[shift={(2/9,-3/9)}]
				\draw [thick,black] plot  coordinates {(3/32,2/32) (-3/32,-2/32) };
			\end{scope}
			\begin{scope}[shift={(4/9,-6/9)}]
				\draw [thick,black] plot  coordinates {(3/32,2/32) (-3/32,-2/32) };
			\end{scope}
		\end{scope}
		\draw [thick,black] plot  coordinates {(-3/2+1/4,0) (5/2,0) };
		\begin{scope}[shift={(4/2,0)}]
			\draw [thick,black] plot  coordinates {(8/12,-8/8) (-1/6,1/4) };
			\fill [black] (0,0) circle (0.06);
			\begin{scope}[shift={(2/9,-3/9)}]
				\draw [thick,black] plot  coordinates {(3/32,2/32) (-3/32,-2/32) };
			\end{scope}
			\begin{scope}[shift={(4/9,-6/9)}]
				\draw [thick,black] plot  coordinates {(3/32,2/32) (-3/32,-2/32) };
			\end{scope}
		\end{scope}
		\begin{scope}[shift={(1,0)}]
			\draw [thick,blue] plot  coordinates {(0,1/3) (0,-6/4) };
			\fill [black] (0,0) circle (0.06);
			\draw (0,-6/4) node [right] {\hbox{\footnotesize$\geqtwo$}};
		\end{scope}
\end{tikzpicture}}

\def\includefigD#1{
\makebox(130pt,75pt){#1}}

\def\FigTwoCaseDone{
\begin{tabular}{|c|c|c|}
\hline
{\large\strut} (D1) & (D2) & (D3) \\
\hline
${\Large\strut}
\alpha=\beta=\gamma=2$ 
& $\alpha>\beta=\gamma=2$ 
& $\beta>\gamma=2$ 
\\
\hline
\includefigD{\FigTwoCaseDi} & 
\includefigD{\FigTwoCaseDii} & 
\includefigD{\FigTwoCaseDiii}
 \\
\hline\hline
{\large\strut} (D4) & (D5) & (D6) \\
\hline
${\Large\strut}
\gamma >2$ 
& $\beta>\gamma \ \land \ \alpha+\gamma<4$ 
& $\beta>\gamma \ \land \ \beta+\gamma<4=\alpha+\gamma$ 
\\
\hline
\includefigD{\FigTwoCaseDiv} &
\includefigD{\FigTwoCaseDv}&
\includefigD{\FigTwoCaseDvi}
\\
\hline\hline
{\large\strut} (D7) & (D8) & (D9) \\
\hline
${\Large\strut}
\beta>\gamma \ \land \ \beta+\gamma=4=\alpha+\gamma$
& $\beta>\gamma \ \land \ \beta+\gamma<4<\alpha+\gamma$
& $2>\gamma \ \land \ \beta+\gamma=4<\alpha+\gamma$ 
\\
\hline
\includefigD{\FigTwoCaseDvii} &
\includefigD{\FigTwoCaseDviii}&
\includefigD{\FigTwoCaseDix}
\\
\hline\hline
{\large\strut} (D10) & (D11) & (D12) \\
\hline
${\Large\strut}
2>\gamma<\beta \ \land \ \beta+\gamma>4$ 
& $\beta=\gamma=\delta<\alpha \ \land \ \alpha+\gamma<4$ 
& $\beta=\gamma=\delta<\alpha \ \land \ \alpha+\gamma=4$ 
\\
\hline
\includefigD{\FigTwoCaseDx} &
\includefigD{\FigTwoCaseDxi} &
\includefigD{\FigTwoCaseDxii}
\\
\hline 

\end{tabular}}

\def\FigTwoCaseDtwo{
\begin{tabular}{|c|c|c|}
	\hline
	{\large\strut} (D13) & (D14) & (D15) \\
	\hline
	${\Large\strut}
	\beta=\gamma=\delta<2 \ \land \ \alpha+\gamma>4$ 
	& $ \gamma <\delta<\frac{\alpha+\gamma}{2} \ \land \ \alpha +\delta<4
	 $ 
	& $ \gamma <\delta<\frac{\alpha+\gamma}{2} \ \land \ \alpha +\delta=4
	$ 
	\\
	\hline
	\includefigD{\FigTwoCaseDxiii} & 
	\includefigD{\FigTwoCaseDxiv} & 
	\includefigD{\FigTwoCaseDxv}
	\\
	\hline\hline
	{\large\strut} (D16) & (D17) & (D18) \\
	\hline
	${\Large\strut}
	\gamma<\delta<\frac{4+\gamma}{3} \ \land \ \alpha+\delta>4
	$ 
	& $\delta\geq \frac{4+\gamma}{3} \ \land \ \alpha>\frac{8-\gamma}{3}>2  
	$ 
	& $  3\alpha+\gamma < 8 \ \land \  \delta = \frac{\alpha+\gamma}{2}>\gamma$
	\\
	\hline
	\includefigD{\FigTwoCaseDxvi} & 
	\includefigD{\FigTwoCaseDxvii} & 
	\includefigD{\FigTwoCaseDxviii}
	\\
	\hline\hline
		{\large\strut} (D19) & (D20) & (D21) \\
	\hline
	{\Large\strut}
	$\gamma<\alpha \ \land \   \frac{\alpha+\gamma}{2}<\delta<\frac{\alpha+2\gamma+4}{5} $ 
	& $\gamma< \frac{3\alpha+\gamma}{4} < 2 \ \land \  \delta\geq \frac{\alpha+2\gamma+4}{5} $
	& $  3\alpha+\gamma = 8 \ \land \   \delta\geq\frac{\alpha+\gamma}{2}>\gamma $
	\\
	\hline
	\includefigD{\FigTwoCaseDxix} & 
	\includefigD{\FigTwoCaseDxx} & 
	\includefigD{\FigTwoCaseDxxi}
	\\
	\hline\hline
		{\large\strut} (D22) & (D23) & (D24) \\
	\hline
	{\Large\strut}
	$\alpha=\gamma=\delta<2 $ 
	& $2>\alpha=\gamma<\delta<\frac{\alpha+2\gamma +4}{5}$
	&$2>\alpha=\gamma \ \land \ \delta\geq\frac{\alpha+2\gamma +4}{5}$
	\\
	\hline
	\includefigD{\FigTwoCaseDxxii} & 
	\includefigD{\FigTwoCaseDxxiii} & 
	\includefigD{\FigTwoCaseDxxiv}
	\\
	\hline
\end{tabular}}

\def\FigTwoCaseEFGi{\begin{tikzpicture}[baseline=12pt]
\draw (7/4,0/4) node {};
\draw [thick,black] plot  coordinates {(13/8,-13/8) (-1/4,1/4) };
\begin{scope}[shift={(4/3,-4/3)}]
	\draw [thick,blue] plot  coordinates {(8/12,8/8) (-1/6,-1/4) };
	\fill [black] (0,0) circle (0.06);
	\draw (2/4,4/4) node [above] {\hbox{\footnotesize$\geqone$}};
\end{scope}
\begin{scope}[xscale=-1]
\draw [thick,black] plot  coordinates {(13/8,-13/8) (-1/4,1/4) };
\begin{scope}[shift={(1/3,-1/3)}]
	\draw [thick,black] plot  coordinates {(-1/12,-1/12) (1/12,1/12) };
\end{scope}
\begin{scope}[shift={(2/3,-2/3)}]
	\draw [thick,black] plot  coordinates {(-1/12,-1/12) (1/12,1/12) };
\end{scope}
\begin{scope}[shift={(3/3,-3/3)}]
	\draw [thick,black] plot  coordinates {(-1/12,-1/12) (1/12,1/12) };
\end{scope}

\begin{scope}[shift={(4/3,-4/3)}]
	\draw [thick,blue] plot  coordinates {(8/12,8/8) (-1/6,-1/4) };
	\fill [black] (0,0) circle (0.06);
	\draw (2/4,4/4) node [above] {\hbox{\footnotesize$\geqone$}};
\end{scope}
\end{scope}

\begin{scope}[shift={(0,0)}]
\fill [white] (0,0) circle (0.07);
\draw [thick] (0,0) circle (0.07);
\end{scope}

\begin{scope}[shift={(1/3,-1/3)}]
\draw [thick,black] plot  coordinates {(-1/12,-1/12) (1/12,1/12) };
\end{scope}
\begin{scope}[shift={(2/3,-2/3)}]
\draw [thick,black] plot  coordinates {(-1/12,-1/12) (1/12,1/12) };
\end{scope}
\begin{scope}[shift={(3/3,-3/3)}]
\draw [thick,black] plot  coordinates {(-1/12,-1/12) (1/12,1/12) };
\end{scope}
\end{tikzpicture}}

\def\FigTwoCaseEFGii{\begin{tikzpicture}[baseline=12pt]
		\draw (7/4,0/4) node {};
\draw [thick,black] plot  coordinates {(3/3,-3/3) (-1/4,1/4) };
\begin{scope}[shift={(-1/3,+1/3)}]
\draw [thick,orange] plot  coordinates {(2/3,-1) (8/3,-1) };
\fill [black] (+3/3,-1) circle (0.06);
\begin{scope}[shift={(-1/3,0)}]
\draw [thick,blue] plot  coordinates {(2,-4/3) (2,1/8) };
\fill [black] (2,-1) circle (0.06);
\draw (2,0) node [right] {\hbox{\footnotesize$\geqone$}};
\end{scope}
\end{scope}

\begin{scope} [shift={(2,-2/3 )}, xscale=-1, yscale=1]
	\draw [thick,black] plot  coordinates {(0,-8/8) (0,1/4) };
	\fill [black] (0,0) circle (0.06);
	\begin{scope}[shift={(0,-3/9)}]
		\draw [thick,black] plot  coordinates {(3.6/32,0) (-3.6/32,0) };
	\end{scope}
	\begin{scope}[shift={(0,-6/9)}]
		\draw [thick,black] plot  coordinates {(3.6/32,0) (-3.6/32,0) };
	\end{scope}
\end{scope}

		\begin{scope}[xscale=-1]
			\draw [thick,black] plot  coordinates {(13/8,-13/8) (-1/4,1/4) };
			\begin{scope}[shift={(1/3,-1/3)}]
				\draw [thick,black] plot  coordinates {(-1/12,-1/12) (1/12,1/12) };
			\end{scope}
			\begin{scope}[shift={(2/3,-2/3)}]
				\draw [thick,black] plot  coordinates {(-1/12,-1/12) (1/12,1/12) };
			\end{scope}
			\begin{scope}[shift={(3/3,-3/3)}]
				\draw [thick,black] plot  coordinates {(-1/12,-1/12) (1/12,1/12) };
			\end{scope}
			
			\begin{scope}[shift={(4/3,-4/3)}]
				\draw [thick,blue] plot  coordinates {(8/12,8/8) (-1/6,-1/4) };
				\fill [black] (0,0) circle (0.06);
				\draw (2/4,4/4) node [above] {\hbox{\footnotesize$\geqone$}};
			\end{scope}
		\end{scope}
		
		\begin{scope}[shift={(0,0)}]
			\fill [white] (0,0) circle (0.07);
			\draw [thick] (0,0) circle (0.07);
		\end{scope}
		
		\begin{scope}[shift={(1/3,-1/3)}]
			\draw [thick,black] plot  coordinates {(-1/12,-1/12) (1/12,1/12) };
		\end{scope}
\end{tikzpicture}}

\def\FigTwoCaseEFGiii{\begin{tikzpicture}[baseline=12pt]
		\draw (7/4,0/4) node {};
		\draw [thick,black] plot  coordinates {(3/3,-3/3) (-1/4,1/4) };
		\begin{scope}[shift={(-1/3,+1/3)}]
			\draw [thick,orange] plot  coordinates {(2/3,-1) (8/3,-1) };
			\fill [black] (+3/3,-1) circle (0.06);
			\begin{scope}[shift={(-1/3,0)}]
				\draw [thick,blue] plot  coordinates {(2,-4/3) (2,1/8) };
				\fill [black] (2,-1) circle (0.06);
				\draw (2,0) node [right] {\hbox{\footnotesize$\geqone$}};
			\end{scope}
		\end{scope}
		
	\begin{scope} [shift={(2,-2/3 )}, xscale=-1, yscale=1]
		\draw [thick,black] plot  coordinates {(0,-8/8) (0,1/4) };
		\fill [black] (0,0) circle (0.06);
		\begin{scope}[shift={(0,-3/9)}]
			\draw [thick,black] plot  coordinates {(3.6/32,0) (-3.6/32,0) };
		\end{scope}
		\begin{scope}[shift={(0,-6/9)}]
			\draw [thick,black] plot  coordinates {(3.6/32,0) (-3.6/32,0) };
		\end{scope}
	\end{scope}
		
		\begin{scope}[xscale=-1]
					\draw [thick,black] plot  coordinates {(3/3,-3/3) (-1/4,1/4) };
			\begin{scope}[shift={(-1/3,+1/3)}]
				\draw [thick,orange] plot  coordinates {(2/3,-1) (8/3,-1) };
				\fill [black] (+3/3,-1) circle (0.06);
				\begin{scope}[shift={(-1/3,0)}]
					\draw [thick,blue] plot  coordinates {(2,-4/3) (2,1/8) };
					\fill [black] (2,-1) circle (0.06);
					\draw (2,0) node [left] {\hbox{\footnotesize$\geqone$}};
				\end{scope}
			\end{scope}
			
			\begin{scope} [shift={(2,-2/3 )}, xscale=-1, yscale=1]
				\draw [thick,black] plot  coordinates {(0,-8/8) (0,1/4) };
				\fill [black] (0,0) circle (0.06);
				\begin{scope}[shift={(0,-3/9)}]
					\draw [thick,black] plot  coordinates {(3.6/32,0) (-3.6/32,0) };
				\end{scope}
				\begin{scope}[shift={(0,-6/9)}]
					\draw [thick,black] plot  coordinates {(3.6/32,0) (-3.6/32,0) };
				\end{scope}
			\end{scope}
					\begin{scope}[shift={(1/3,-1/3)}]
				\draw [thick,black] plot  coordinates {(-1/12,-1/12) (1/12,1/12) };
			\end{scope}
		\end{scope}
		
		\begin{scope}[shift={(0,0)}]
			\fill [white] (0,0) circle (0.07);
			\draw [thick] (0,0) circle (0.07);
		\end{scope}
		
		\begin{scope}[shift={(1/3,-1/3)}]
			\draw [thick,black] plot  coordinates {(-1/12,-1/12) (1/12,1/12) };
		\end{scope}
\end{tikzpicture}}

\def\FigTwoCaseEFGiv{\begin{tikzpicture}[baseline=12pt]
		\draw (7/4,0/4) node {};
		\draw [thick,black] plot  coordinates {(3/3,-3/3) (-1/4,1/4) };
		\begin{scope}[shift={(-1/3,+1/3)}]
			\draw [thick,orange] plot  coordinates {(2/3,-1) (8/3,-1) };
			\fill [black] (+3/3,-1) circle (0.06);
			\begin{scope}[shift={(-1/3,0)}]
				\draw (2,-1) node [above] {\hbox{\footnotesize$\geqone$}};
			\end{scope}
		\end{scope}
		
\begin{scope} [shift={(2,-2/3 )}, xscale=-1, yscale=1]
	\draw [thick,black] plot  coordinates {(0,-8/8) (0,1/4) };
	\fill [black] (0,0) circle (0.06);
	\begin{scope}[shift={(0,-3/9)}]
		\draw [thick,black] plot  coordinates {(3.6/32,0) (-3.6/32,0) };
	\end{scope}
	\begin{scope}[shift={(0,-6/9)}]
		\draw [thick,black] plot  coordinates {(3.6/32,0) (-3.6/32,0) };
	\end{scope}
\end{scope}
		
		\begin{scope}[xscale=-1]
			\draw [thick,black] plot  coordinates {(13/8,-13/8) (-1/4,1/4) };
			\begin{scope}[shift={(1/3,-1/3)}]
				\draw [thick,black] plot  coordinates {(-1/12,-1/12) (1/12,1/12) };
			\end{scope}
			\begin{scope}[shift={(2/3,-2/3)}]
				\draw [thick,black] plot  coordinates {(-1/12,-1/12) (1/12,1/12) };
			\end{scope}
			\begin{scope}[shift={(3/3,-3/3)}]
				\draw [thick,black] plot  coordinates {(-1/12,-1/12) (1/12,1/12) };
			\end{scope}
			
			\begin{scope}[shift={(4/3,-4/3)}]
				\draw [thick,blue] plot  coordinates {(8/12,8/8) (-1/6,-1/4) };
				\fill [black] (0,0) circle (0.06);
				\draw (2/4,4/4) node [above] {\hbox{\footnotesize$\geqone$}};
			\end{scope}
		\end{scope}
		
		\begin{scope}[shift={(0,0)}]
			\fill [white] (0,0) circle (0.07);
			\draw [thick] (0,0) circle (0.07);
		\end{scope}
		
		\begin{scope}[shift={(1/3,-1/3)}]
			\draw [thick,black] plot  coordinates {(-1/12,-1/12) (1/12,1/12) };
		\end{scope}
\end{tikzpicture}}

\def\FigTwoCaseEFGv{\begin{tikzpicture}[baseline=12pt]
		\draw (7/4,0/4) node {};
		\draw [thick,black] plot  coordinates {(3/3,-3/3) (-1/4,1/4) };
		\begin{scope}[shift={(-1/3,+1/3)}]
			\draw [thick,orange] plot  coordinates {(2/3,-1) (8/3,-1) };
			\fill [black] (+3/3,-1) circle (0.06);
			\begin{scope}[shift={(-1/3,0)}]
				\draw (2,-1) node [above] {\hbox{\footnotesize$\geqone$}};
			\end{scope}
		\end{scope}
		\begin{scope}[shift={(1/3,-1/3)}]
			\draw [thick,black] plot  coordinates {(-1/12,-1/12) (1/12,1/12) };
		\end{scope}
\begin{scope} [shift={(2,-2/3 )}, xscale=-1, yscale=1]
	\draw [thick,black] plot  coordinates {(0,-8/8) (0,1/4) };
	\fill [black] (0,0) circle (0.06);
	\begin{scope}[shift={(0,-3/9)}]
		\draw [thick,black] plot  coordinates {(3.6/32,0) (-3.6/32,0) };
	\end{scope}
	\begin{scope}[shift={(0,-6/9)}]
		\draw [thick,black] plot  coordinates {(3.6/32,0) (-3.6/32,0) };
	\end{scope}
\end{scope}
		
		\begin{scope}[xscale=-1]
			\draw [thick,black] plot  coordinates {(3/3,-3/3) (-1/4,1/4) };
			\begin{scope}[shift={(-1/3,+1/3)}]
				\draw [thick,orange] plot  coordinates {(2/3,-1) (8/3,-1) };
				\fill [black] (+3/3,-1) circle (0.06);
				\begin{scope}[shift={(-1/3,0)}]
					\draw [thick,blue] plot  coordinates {(2,-4/3) (2,1/8) };
					\fill [black] (2,-1) circle (0.06);
					\draw (2,0) node [left] {\hbox{\footnotesize$\geqone$}};
				\end{scope}
			\end{scope}
			\begin{scope}[shift={(1/3,-1/3)}]
				\draw [thick,black] plot  coordinates {(-1/12,-1/12) (1/12,1/12) };
			\end{scope}
			\begin{scope} [shift={(2,-2/3 )}, xscale=-1, yscale=1]
				\draw [thick,black] plot  coordinates {(0,-8/8) (0,1/4) };
				\fill [black] (0,0) circle (0.06);
				\begin{scope}[shift={(0,-3/9)}]
					\draw [thick,black] plot  coordinates {(3.6/32,0) (-3.6/32,0) };
				\end{scope}
				\begin{scope}[shift={(0,-6/9)}]
					\draw [thick,black] plot  coordinates {(3.6/32,0) (-3.6/32,0) };
				\end{scope}
			\end{scope}
		\end{scope}
		
		\begin{scope}[shift={(0,0)}]
			\fill [white] (0,0) circle (0.07);
			\draw [thick] (0,0) circle (0.07);
		\end{scope}

\end{tikzpicture}}

\def\FigTwoCaseEFGvi{\begin{tikzpicture}[baseline=12pt]
		\draw (7/4,0/4) node {};
		\draw [thick,black] plot  coordinates {(3/3,-3/3) (-1/4,1/4) };
		\begin{scope}[shift={(-1/3,+1/3)}]
			\draw [thick,orange] plot  coordinates {(2/3,-1) (8/3,-1) };
			\fill [black] (+3/3,-1) circle (0.06);
			\begin{scope}[shift={(-1/3,0)}]
				\draw (2,-1) node [above] {\hbox{\footnotesize$\geqone$}};
			\end{scope}
		\end{scope}
		
		\begin{scope} [shift={(2,-2/3 )}, xscale=-1, yscale=1]
			\draw [thick,black] plot  coordinates {(0,-8/8) (0,1/4) };
			\fill [black] (0,0) circle (0.06);
			\begin{scope}[shift={(0,-3/9)}]
				\draw [thick,black] plot  coordinates {(3.6/32,0) (-3.6/32,0) };
			\end{scope}
			\begin{scope}[shift={(0,-6/9)}]
				\draw [thick,black] plot  coordinates {(3.6/32,0) (-3.6/32,0) };
			\end{scope}
		\end{scope}
		
		\begin{scope}[xscale=-1]
			\draw [thick,black] plot  coordinates {(3/3,-3/3) (-1/4,1/4) };
			\begin{scope}[shift={(-1/3,+1/3)}]
				\draw [thick,orange] plot  coordinates {(2/3,-1) (8/3,-1) };
				\fill [black] (+3/3,-1) circle (0.06);
				\begin{scope}[shift={(-1/3,0)}]
					\draw (2,-1) node [above] {\hbox{\footnotesize$\geqone$}};
				\end{scope}
			\end{scope}
			\begin{scope}[shift={(1/3,-1/3)}]
				\draw [thick,black] plot  coordinates {(-1/12,-1/12) (1/12,1/12) };
			\end{scope}
			\begin{scope} [shift={(2,-2/3 )}, xscale=-1, yscale=1]
				\draw [thick,black] plot  coordinates {(0,-8/8) (0,1/4) };
				\fill [black] (0,0) circle (0.06);
				\begin{scope}[shift={(0,-3/9)}]
					\draw [thick,black] plot  coordinates {(3.6/32,0) (-3.6/32,0) };
				\end{scope}
				\begin{scope}[shift={(0,-6/9)}]
					\draw [thick,black] plot  coordinates {(3.6/32,0) (-3.6/32,0) };
				\end{scope}
			\end{scope}
		\end{scope}
		
		\begin{scope}[shift={(0,0)}]
			\fill [white] (0,0) circle (0.07);
			\draw [thick] (0,0) circle (0.07);
		\end{scope}
		
		\begin{scope}[shift={(1/3,-1/3)}]
			\draw [thick,black] plot  coordinates {(-1/12,-1/12) (1/12,1/12) };
		\end{scope}
\end{tikzpicture}}

\def\FigTwoCaseEFGvii{\begin{tikzpicture}[baseline=12pt]
		\draw (7/4,0/4) node {};
		\draw [thick,black] plot  coordinates {(3/3,-3/3) (-1/4,1/4) };
		\begin{scope}[shift={(-1/3,+1/3)}]
			\draw [thick,orange] plot [smooth, tension=1] coordinates {(2/3,-9/8) (2,-4/8) (4/3,1/8) };
		\draw [thick,orange] plot [smooth, tension=1] coordinates {(7/2-2/3,-9/8) (7/2-2,-4/8) (7/2-4/3,1/8) };
		\fill [black] (7/4,-0.075) circle (0.06);
		\fill [black] (7/4,-0.72) circle (0.06);
		\fill [black] (1.025,-1.025) circle (0.06);
		\end{scope}
	
				\begin{scope} [shift={(2+1/32,-2/3 )}, xscale=-1, yscale=1]
			\draw [thick,black] plot  coordinates {(0,-8/8) (0,1/4) };
			\fill [black] (0,0) circle (0.06);
			\begin{scope}[shift={(0,-3/9)}]
				\draw [thick,black] plot  coordinates {(3.6/32,0) (-3.6/32,0) };
			\end{scope}
			\begin{scope}[shift={(0,-6/9)}]
				\draw [thick,black] plot  coordinates {(3.6/32,0) (-3.6/32,0) };
			\end{scope}
		\end{scope}
		
		\begin{scope}[xscale=-1]
			\draw [thick,black] plot  coordinates {(13/8,-13/8) (-1/4,1/4) };
			\begin{scope}[shift={(1/3,-1/3)}]
				\draw [thick,black] plot  coordinates {(-1/12,-1/12) (1/12,1/12) };
			\end{scope}
			\begin{scope}[shift={(2/3,-2/3)}]
				\draw [thick,black] plot  coordinates {(-1/12,-1/12) (1/12,1/12) };
			\end{scope}
			\begin{scope}[shift={(3/3,-3/3)}]
				\draw [thick,black] plot  coordinates {(-1/12,-1/12) (1/12,1/12) };
			\end{scope}
			
			\begin{scope}[shift={(4/3,-4/3)}]
				\draw [thick,blue] plot  coordinates {(8/12,8/8) (-1/6,-1/4) };
				\fill [black] (0,0) circle (0.06);
				\draw (2/4,4/4) node [above] {\hbox{\footnotesize$\geqone$}};
			\end{scope}
		\end{scope}
		
		\begin{scope}[shift={(0,0)}]
			\fill [white] (0,0) circle (0.07);
			\draw [thick] (0,0) circle (0.07);
		\end{scope}
		
		\begin{scope}[shift={(1/3,-1/3)}]
			\draw [thick,black] plot  coordinates {(-1/12,-1/12) (1/12,1/12) };
		\end{scope}
\end{tikzpicture}}

\def\FigTwoCaseEFGviii{\begin{tikzpicture}[baseline=12pt]
		\draw (7/4,0/4) node {};
		\draw [thick,black] plot  coordinates {(3/3,-3/3) (-1/4,1/4) };
		\begin{scope}[shift={(-1/3,+1/3)}]
			\draw [thick,orange] plot [smooth, tension=1] coordinates {(2/3,-9/8) (2,-4/8) (4/3,1/8) };
		\draw [thick,orange] plot [smooth, tension=1] coordinates {(7/2-2/3,-9/8) (7/2-2,-4/8) (7/2-4/3,1/8) };
		\fill [black] (7/4,-0.075) circle (0.06);
		\fill [black] (7/4,-0.72) circle (0.06);
		\fill [black] (1.025,-1.025) circle (0.06);
		\end{scope}
		\begin{scope}[shift={(1/3,-1/3)}]
			\draw [thick,black] plot  coordinates {(-1/12,-1/12) (1/12,1/12) };
		\end{scope}
		\begin{scope} [shift={(2+1/32,-2/3 )}, xscale=-1, yscale=1]
			\draw [thick,black] plot  coordinates {(0,-8/8) (0,1/4) };
			\fill [black] (0,0) circle (0.06);
			\begin{scope}[shift={(0,-3/9)}]
				\draw [thick,black] plot  coordinates {(3.6/32,0) (-3.6/32,0) };
			\end{scope}
			\begin{scope}[shift={(0,-6/9)}]
				\draw [thick,black] plot  coordinates {(3.6/32,0) (-3.6/32,0) };
			\end{scope}
		\end{scope}
		
		\begin{scope}[xscale=-1]
			\draw [thick,black] plot  coordinates {(3/3,-3/3) (-1/4,1/4) };
			\begin{scope}[shift={(-1/3,+1/3)}]
				\draw [thick,orange] plot  coordinates {(2/3,-1) (8/3,-1) };
				\fill [black] (+3/3,-1) circle (0.06);
				\begin{scope}[shift={(-1/3,0)}]
					\draw [thick,blue] plot  coordinates {(2,-4/3) (2,1/8) };
					\fill [black] (2,-1) circle (0.06);
					\draw (2,0) node [left] {\hbox{\footnotesize$\geqone$}};
				\end{scope}
			\end{scope}
			\begin{scope}[shift={(1/3,-1/3)}]
				\draw [thick,black] plot  coordinates {(-1/12,-1/12) (1/12,1/12) };
			\end{scope}
			\begin{scope} [shift={(2,-2/3 )}, xscale=-1, yscale=1]
				\draw [thick,black] plot  coordinates {(0,-8/8) (0,1/4) };
				\fill [black] (0,0) circle (0.06);
				\begin{scope}[shift={(0,-3/9)}]
					\draw [thick,black] plot  coordinates {(3.6/32,0) (-3.6/32,0) };
				\end{scope}
				\begin{scope}[shift={(0,-6/9)}]
					\draw [thick,black] plot  coordinates {(3.6/32,0) (-3.6/32,0) };
				\end{scope}
			\end{scope}
		\end{scope}
		
		\begin{scope}[shift={(0,0)}]
			\fill [white] (0,0) circle (0.07);
			\draw [thick] (0,0) circle (0.07);
		\end{scope}
\end{tikzpicture}}

\def\FigTwoCaseEFGix{\begin{tikzpicture}[baseline=12pt]
		\draw (7/4,0/4) node {};
		\draw [thick,black] plot  coordinates {(3/3,-3/3) (-1/4,1/4) };
		\begin{scope}[shift={(-1/3,+1/3)}]
			\draw [thick,orange] plot [smooth, tension=1] coordinates {(2/3,-9/8) (2,-4/8) (4/3,1/8) };
			\draw [thick,orange] plot [smooth, tension=1] coordinates {(7/2-2/3,-9/8) (7/2-2,-4/8) (7/2-4/3,1/8) };
			\fill [black] (7/4,-0.075) circle (0.06);
			\fill [black] (7/4,-0.72) circle (0.06);
			\fill [black] (1.025,-1.025) circle (0.06);
		\end{scope}
		\begin{scope}[shift={(1/3,-1/3)}]
			\draw [thick,black] plot  coordinates {(-1/12,-1/12) (1/12,1/12) };
		\end{scope}
		\begin{scope} [shift={(2+1/32,-2/3 )}, xscale=-1, yscale=1]
			\draw [thick,black] plot  coordinates {(0,-8/8) (0,1/4) };
			\fill [black] (0,0) circle (0.06);
			\begin{scope}[shift={(0,-3/9)}]
				\draw [thick,black] plot  coordinates {(3.6/32,0) (-3.6/32,0) };
			\end{scope}
			\begin{scope}[shift={(0,-6/9)}]
				\draw [thick,black] plot  coordinates {(3.6/32,0) (-3.6/32,0) };
			\end{scope}
		\end{scope}
		
		\begin{scope}[xscale=-1]
			\draw [thick,black] plot  coordinates {(3/3,-3/3) (-1/4,1/4) };
			\begin{scope}[shift={(-1/3,+1/3)}]
				\draw [thick,orange] plot  coordinates {(2/3,-1) (8/3,-1) };
				\fill [black] (+3/3,-1) circle (0.06);
				\begin{scope}[shift={(-1/3,0)}]
					\draw (2,-1) node [above] {\hbox{\footnotesize$\geqone$}};
				\end{scope}
			\end{scope}
			\begin{scope}[shift={(1/3,-1/3)}]
				\draw [thick,black] plot  coordinates {(-1/12,-1/12) (1/12,1/12) };
			\end{scope}
			\begin{scope} [shift={(2,-2/3 )}, xscale=-1, yscale=1]
				\draw [thick,black] plot  coordinates {(0,-8/8) (0,1/4) };
				\fill [black] (0,0) circle (0.06);
				\begin{scope}[shift={(0,-3/9)}]
					\draw [thick,black] plot  coordinates {(3.6/32,0) (-3.6/32,0) };
				\end{scope}
				\begin{scope}[shift={(0,-6/9)}]
					\draw [thick,black] plot  coordinates {(3.6/32,0) (-3.6/32,0) };
				\end{scope}
			\end{scope}
		\end{scope}
		
		\begin{scope}[shift={(0,0)}]
			\fill [white] (0,0) circle (0.07);
			\draw [thick] (0,0) circle (0.07);
		\end{scope}
\end{tikzpicture}}

\def\FigTwoCaseEFGx{\begin{tikzpicture}[baseline=12pt]
		\draw (7/4,0/4) node {};
		\draw [thick,black] plot  coordinates {(3/3,-3/3) (-1/4,1/4) };
		\begin{scope}[shift={(-1/3,+1/3)}]
			\draw [thick,orange] plot [smooth, tension=1] coordinates {(2/3,-9/8) (2,-4/8) (4/3,1/8) };
			\draw [thick,orange] plot [smooth, tension=1] coordinates {(7/2-2/3,-9/8) (7/2-2,-4/8) (7/2-4/3,1/8) };
			\fill [black] (7/4,-0.075) circle (0.06);
			\fill [black] (7/4,-0.72) circle (0.06);
			\fill [black] (1.025,-1.025) circle (0.06);
		\end{scope}
		\begin{scope}[shift={(1/3,-1/3)}]
			\draw [thick,black] plot  coordinates {(-1/12,-1/12) (1/12,1/12) };
		\end{scope}
			\begin{scope} [shift={(2+1/32,-2/3 )}, xscale=-1, yscale=1]
			\draw [thick,black] plot  coordinates {(0,-8/8) (0,1/4) };
			\fill [black] (0,0) circle (0.06);
			\begin{scope}[shift={(0,-3/9)}]
				\draw [thick,black] plot  coordinates {(3.6/32,0) (-3.6/32,0) };
			\end{scope}
			\begin{scope}[shift={(0,-6/9)}]
				\draw [thick,black] plot  coordinates {(3.6/32,0) (-3.6/32,0) };
			\end{scope}
		\end{scope}
		
		\begin{scope}[xscale=-1]
			\draw [thick,black] plot  coordinates {(3/3,-3/3) (-1/4,1/4) };
			\begin{scope}[shift={(-1/3,+1/3)}]
				\draw [thick,orange] plot [smooth, tension=1] coordinates {(2/3,-9/8) (2,-4/8) (4/3,1/8) };
				\draw [thick,orange] plot [smooth, tension=1] coordinates {(7/2-2/3,-9/8) (7/2-2,-4/8) (7/2-4/3,1/8) };
				\fill [black] (7/4,-0.075) circle (0.06);
				\fill [black] (7/4,-0.72) circle (0.06);
				\fill [black] (1.025,-1.025) circle (0.06);
			\end{scope}
			\begin{scope}[shift={(1/3,-1/3)}]
				\draw [thick,black] plot  coordinates {(-1/12,-1/12) (1/12,1/12) };
			\end{scope}
				\begin{scope} [shift={(2+1/32,-2/3 )}, xscale=-1, yscale=1]
				\draw [thick,black] plot  coordinates {(0,-8/8) (0,1/4) };
				\fill [black] (0,0) circle (0.06);
				\begin{scope}[shift={(0,-3/9)}]
					\draw [thick,black] plot  coordinates {(3.6/32,0) (-3.6/32,0) };
				\end{scope}
				\begin{scope}[shift={(0,-6/9)}]
					\draw [thick,black] plot  coordinates {(3.6/32,0) (-3.6/32,0) };
				\end{scope}
			\end{scope}
		\end{scope}
		
		\begin{scope}[shift={(0,0)}]
			\fill [white] (0,0) circle (0.07);
			\draw [thick] (0,0) circle (0.07);
		\end{scope}
\end{tikzpicture}}

\def\includefigEFG#1{
\makebox(140pt,70pt){#1}}

\def\FigTwoCaseE{
\begin{tabular}{|c|}
\hline{\large\strut} (E1) \\
\hline
\includefigEFG{\FigTwoCaseEFGi} \\
\hline\end{tabular}}

\def\FigTwoCaseF{
\begin{tabular}{|c|c|c|}
\hline{\large\strut} (F1) & (F2) & (F3) \\
\hline{\Large\strut}
$\alpha<4$ & $\alpha=4$ & $\alpha>4$\\
\hline
\includefigEFG{\FigTwoCaseEFGii} &
\includefigEFG{\FigTwoCaseEFGiv} & 
\includefigEFG{\FigTwoCaseEFGvii}\\
\hline\end{tabular}}

\def\FigTwoCaseG{
\begin{tabular}{|c|c|c|}
\hline{\large\strut} (G1) & (G2) & (G3) \\
\hline{\Large\strut}
$4>\alpha\ge\beta$ & $4=\alpha>\beta$ & $\alpha=4=\beta$\\
\hline
\includefigEFG{\FigTwoCaseEFGiii} &
\includefigEFG{\FigTwoCaseEFGv} & 
\includefigEFG{\FigTwoCaseEFGvi}\\
\hline\hline
{\large\strut} (G4) & (G5) & (G6) \\ 
\hline
{\Large\strut} $\alpha>4>\beta$ & $\alpha>4=\beta$ & $\alpha\ge\beta>4$ \\
\hline
\includefigEFG{\FigTwoCaseEFGviii} & 
\includefigEFG{\FigTwoCaseEFGix} &
\includefigEFG{\FigTwoCaseEFGx} \\
\hline 
\end{tabular}}

\def\FigTwoUnmarkedi{\begin{tikzpicture}[baseline=13pt]
		\draw [thick,black] plot  coordinates {(-1/2,0) (4/2,0) };
		\draw (-2/4,0) node [left] {\hbox{\footnotesize$\geqtwo$}};
\end{tikzpicture}}

\def\FigTwoUnmarkedii{\begin{tikzpicture}[baseline=13pt]
				\draw [thick, black] (0,0) arc (-45:225:0.4);
				\draw [thick, black] plot coordinates {(0,0) (-1,-1)};
				\begin{scope}[shift={(-0.707106781*0.8,0)},xscale=-1]
									\draw [thick, black] plot coordinates {(0,0) (-1,-1)};
				\end{scope}
		\fill [black] (-0.707106781*0.4,-0.4*0.707106781) circle (0.06);
		\draw (-3/4,-1/2) node [left] {\hbox{\footnotesize$\geqone$}};
\end{tikzpicture}}

\def\FigTwoUnmarkediii{\begin{tikzpicture}[baseline=13pt]
		\draw [thick, black] (0,0) arc (-45:225:0.4);
		\draw [thick, black] plot coordinates {(0,0) (-1,-1)};
		\begin{scope}[shift={(-0.707106781*0.4+1/2,-0.4*0.707106781-1/4)},xscale=-1]
					\draw[scale=0.5, domain=-1:1, smooth, variable=\x, thick, black] plot ({\x}, {\x*\x/2});
		\end{scope}
		
		\begin{scope}[shift={(-0.707106781*0.8,0)},xscale=-1]
			\draw [thick, black] plot coordinates {(0,0) (-0.707106781*0.4,-0.707106781*0.4)};
		\end{scope}
		\fill [black] (-0.707106781*0.4,-0.4*0.707106781) circle (0.06);
		\begin{scope}[shift={(1-0.8*0.707106781,0)},xscale=-1]
	\draw [thick, black] (0,0) arc (-45:225:0.4);
\draw [thick, black] plot coordinates {(0,0) (-1,-1)};
	\begin{scope}[shift={(-0.707106781*0.8,0)},xscale=-1]
	\draw [thick, black] plot coordinates {(0,0) (-0.707106781*0.4,-0.707106781*0.4)};
\end{scope}
		\fill [black] (-0.707106781*0.4,-0.4*0.707106781) circle (0.06);

		\end{scope}

\end{tikzpicture}}

\def\FigTwoUnmarkediv{\begin{tikzpicture}[baseline=13pt]
		
							\draw[scale=0.4, domain=0:3*pi, smooth, variable=\x, thick, black] plot ({\x}, {cos(\x r)});
							\draw[scale=0.4, domain=0:3*pi, smooth, variable=\x, thick, black] plot ({\x}, {-cos(\x r)});
							\fill [black] (0.4*1/2*pi,0) circle (0.06);
							\fill [black] (0.4*3/2*pi,0) circle (0.06);
							\fill [black] (0.4*5/2*pi,0) circle (0.06);
\end{tikzpicture}}

\def\FigTwoUnmarkedv{\begin{tikzpicture}[baseline=13pt]
				\draw [thick,black] plot  coordinates {(1/4,1/4) (-1,-1) };
				\draw [thick,black] plot  coordinates {(-1/4,1/4) (1,-1) };
						\fill [black] (0,0) circle (0.06);
				\draw (-1,-1) node [left] {\hbox{\footnotesize$\geqone$}};
					\draw (1,-1) node [right] {\hbox{\footnotesize$\geqone$}};
\end{tikzpicture}}

\def\FigTwoUnmarkedvi{\begin{tikzpicture}[baseline=13pt]
		\draw [thick,black] plot  coordinates {(1/4,1/4) (-1,-1) };
			\draw [thick,black] plot  coordinates {(-1/4,1/4) (0,0) };
\begin{scope}[yscale=-1, shift={(1,1)}]
	\draw [thick, black] (0,0) arc (-45:225:0.4);
	\draw [thick, black] plot coordinates {(0,0) (-1,-1)};
	\begin{scope}[shift={(-0.707106781*0.8,0)},xscale=-1]
		\draw [thick, black] plot coordinates {(0,0) (-1/2,-1/2)};
	\end{scope}
	\fill [black] (-0.707106781*0.4,-0.4*0.707106781) circle (0.06);
\end{scope}
			
		\fill [black] (0,0) circle (0.06);
		
		\draw (-1,-1) node [left] {\hbox{\footnotesize$\geqone$}};
\end{tikzpicture}}

\def\FigTwoUnmarkedvii{\begin{tikzpicture}[baseline=13pt]
		\draw [thick,black] plot  coordinates {(-1/4,1/4) (0,0) };
\begin{scope}[yscale=-1, shift={(1,1)}]
	\draw [thick, black] (0,0) arc (-45:225:0.4);
	\draw [thick, black] plot coordinates {(0,0) (-1,-1)};
	\begin{scope}[shift={(-0.707106781*0.8,0)},xscale=-1]
		\draw [thick, black] plot coordinates {(0,0) (-1/2,-1/2)};
	\end{scope}
	\fill [black] (-0.707106781*0.4,-0.4*0.707106781) circle (0.06);
\end{scope}
		\begin{scope}[xscale=-1]
				\draw [thick,black] plot  coordinates {(-1/4,1/4) (0,0) };
\begin{scope}[yscale=-1, shift={(1,1)}]
	\draw [thick, black] (0,0) arc (-45:225:0.4);
	\draw [thick, black] plot coordinates {(0,0) (-1,-1)};
	\begin{scope}[shift={(-0.707106781*0.8,0)},xscale=-1]
		\draw [thick, black] plot coordinates {(0,0) (-1/2,-1/2)};
	\end{scope}
	\fill [black] (-0.707106781*0.4,-0.4*0.707106781) circle (0.06);
\end{scope} 
		\end{scope}
		
		\fill [black] (0,0) circle (0.06);
\end{tikzpicture}}

\def\includefigUnmarked#1{
	\makebox(130pt,75pt){#1}}
\def\FigTwoUnmarked{
\begin{tabular}{|c|c|c|}
	\hline
	{\large\strut} (I) & (II) & (III) \\
	\hline
	{\Large\strut}
$3\alpha +\gamma\leq 8 \ \land  \  \delta \geq \frac{4+\alpha+2\gamma}{5}$ 
	&  \vphantom{\huge\strut} \makecell{ $3 \alpha+\gamma >8 \ \ \land  \ $ \\ $3\delta-\gamma>4\geq \beta+\gamma$ } 
	& $\gamma=2<\beta$ \\
	\hline
	\includefigUnmarked{\FigTwoUnmarkedi} & 
	\includefigUnmarked{\FigTwoUnmarkedii} & 
	\includefigUnmarked{\FigTwoUnmarkediii} \\
	\hline\hline
	{\large\strut} (IV) & (V) & (VI) \\
	\hline
	{\Large\strut}
	 $\gamma>2$
	& $\alpha+\delta\leq 4 \ \land \  \delta<\frac{4+\alpha+2\gamma}{5}$ 
	& \vphantom{\huge\strut} \makecell{$3\delta-\gamma<4\geq \beta+\gamma$ \\   $ \land\ \ \delta+\alpha>4$} \\
	\hline
	\includefigUnmarked{\FigTwoUnmarkediv} &
	\includefigUnmarked{\FigTwoUnmarkedv} & 
	\includefigUnmarked{\FigTwoUnmarkedvi}   \\
	\hline 
	\noalign{\global\arrayrulewidth=2.0pt}
	\arrayrulecolor{white}\hline
	\noalign{\global\arrayrulewidth=0.4pt}
	\arrayrulecolor{black} \cline{1-1}
	{\large\strut} (VII) & \multicolumn{1}{c}{} & \multicolumn{1}{c}{}  \\
	\cline{1-1}
	{\Large\strut}

	$\gamma<2 \ \land \  \beta+\gamma>4$
	& \multicolumn{1}{c}{}
	& \multicolumn{1}{c}{}  \\
	\cline{1-1}
	\includefigUnmarked{\FigTwoUnmarkedvii}   &
	\multicolumn{1}{c}{} & 
	\multicolumn{1}{c}{} \\
	\cline{1-1}

	\end{tabular}}


%% file: GP23.bbl
\begin{thebibliography}{99}
\setlength{\parskip}{0pt}
\addcontentsline{toc}{section}{References}
\footnotesize

\bibitem{ArzdorfWewers}
Arzdorf, K., Wewers, S.:
{ \it Another proof of the Semistable Reduction Theorem}
\url{arXiv:1211.4624} November 2012, 31p. 

\bibitem{Bosch1980}
Bosch, S.:
Formelle Standardmodelle hyperelliptischer Kurven. 
{\it Mathematische Annalen} {\bf 251}, no.~1 (1980), 19--42.

\bibitem{BoschLuetkebohmertRaynaud1990}
Bosch, S.; L\"{u}tkebohmert, W.; Raynaud, M.:
{\it N\'{e}ron models.}
Ergebnisse der Mathematik und ihrer Grenzgebiete (3), vol.~21. 
Berlin: Springer-Verlag, 1990.




\bibitem{Cuzub2018}
Cuzub, A.:
Properties of Models of Algebraic Curves
{\it Mediterr. J. Math.} {\bf 15.75} (2018) 22p.

\bibitem{deJongIHES1996}
de Jong, A. J.: 
Smoothness, semi-stability and alterations. 
{\it Inst. Hautes \'Etudes Sci. Publ. Math.} No. 83 (1996), 51--93.

\bibitem{DDMM}
Dokchitser, T.; Dokchitser, V.; Maistret, C.; Morgan, A.:
{\it Arithmetic of hyperelliptic curves over local fields.}
\url{arXiv:1808.02936v2} October 2018, 93pp.







\bibitem{EGA4}
Grothendieck, A., 
{\it \'Etudes locale des sch\'emas et des morphismes de sch\'emas},
\'Elements de G\'eom\'etrie Alg\'ebrique IV, EGA4,
{\it Publ. Math. IHES} 
{\bf 20} (1964),
{\bf 24} (1965), 
{\bf 28} (1966), 
{\bf 32} (1967).

\bibitem{Fiore2018}
Fiore, L.:
{\it Semistable models of hyperelliptic curves over residue characteristic 2},
Tesi di laurea, University of Milan (2018) 210p.

\bibitem{FioreYelton}
Fiore, L., Yelton, J.:
{\it Clusters and semistable models of hyperelliptic curves in the wild case},
\url{arXiv:2207.12490v4} August 2023, 85p.  

\bibitem{Gehrunger2020}
Gehrunger, T.:
{\it Reduction of Hyperelliptic Curves.}
Master Thesis, ETH Z\"urich, October 2020, 58pp.

\bibitem{GehrungerPink2021}
Gehrunger, T., Pink, R,:
{\it Reduction of Hyperelliptic Curves in Characteristic $\neq2$.}
 \url{arXiv:2112.05550} December 2021, 19p.  

\bibitem{GPComp}
Gehrunger, T., Pink, R,:
Worksheets for this paper.
\url{https://doi.org/10.3929/ethz-b-000669537}


\bibitem{Gehrunger2024}
Gehrunger, T.:
{\it Effectively computing the Stable Reduction of Hyperelliptic Curves in the wild case.}
in preparation.




\bibitem{Igusa1960}
Igusa, J.-I.:
Arithmetic Variety of Moduli for Genus Two.
{\it Annals of Mathematics} {\bf 72}, no. 3 (1960), 612--649.


\bibitem{Knudsen1983}
Knudsen, F. F.:
The projectivity of the moduli space of stable curves, II: 
The stacks $M_{g,n}$. 
{\it Mathematica Scandinavica} {\bf 52} (1983), 161--199.

\bibitem{Lehr2001ReductionOP}
Lehr, C.:
Reduction of p-cyclic covers of the projective line
{\it Manuscripta Mathematica} {\bf 106} (2001) 151--175.

\bibitem{LehrMatignon2006}
Lehr, C.; Matignon, M.:
Wild monodromy and automorphisms of curves.
{\it Duke Math. J.} {\bf 135} no.~3 (2006) 569--586.

\bibitem{Liu1993}
Liu, Q.:
Courbes stables de genre 2 et leur sch\'ema de modules. 
{\it Mathematische Annalen} {\bf 295}, no.~1 (1993), 201--222.

\bibitem{LiuAlgGeo2002}
Liu, Q.:
{\it Algebraic geometry and arithmetic curves.}
Oxford Graduate Texts in Mathematics, vol.~6.
Oxford: Oxford University Press, 2002.

\bibitem{Liu2006}
Liu, Q.:
Stable reduction of finite covers of curves.
{\it Compos. Math.} {\bf 142} no.~1 (2006), 101--118. 

\bibitem{Arithmetica}
J.Sesiano.:
Books IV to VII of Diophantus’ Arithmetica.
{\it Sources in the History of Mathematics and Physical Science}.
Springer, 1982.


\bibitem{MuellerPink2022}
M\"uller, N., Pink, R.:
Hyperelliptic curves with many automorphisms.
{\it International Journal of Number Theory} {\bf 18}, no.04 (2022) 913--930.

\bibitem{Raynaud1970}
Raynaud, M.:
Sp\'ecialisation du foncteur de Picard.
Inst. Hautes Études Sci. Publ. Math. No. 38 (1970), 27--76. 

\bibitem{Raynaud1990}
Raynaud, M.:
$p$-groupes et r\'{e}duction semi-stable des courbes,
in: {\it The Grothendieck Festschrift}, vol. III,
Progr. Math. 88, Boston: Birkh\"{a}user (1990) 179--197.

\bibitem{Raynaud1999}
Raynaud, M.:
Sp\'{e}cialisation des rev\^{e}tements en caract\'{e}ristique $p>0$,
{\it Ann. Sci. \'{E}cole Norm. Sup.} (4),
{\bf 32} no. 1 (1999) 87--126.

\bibitem{SGA7}
Grothendieck, A., et al.
{\it Groupes de monodromie en g\'eom\'etrie alg\'ebrique},
S\'eminaire de G\'eom\'etrie Alg\'ebrique du Bois Marie,
SGA7, Lect. Notes Math. {\bf 288}, {\bf 340},
Berlin etc.: Springer (1972, 1973).



\bibitem{Temkin2010}
Temkin, M.:
Stable Modification of Relative Curves,
{\it J. Algebraic Geometry} {\bf 19} (2010) 603--677.

\bibitem{Yelton2021}
Yelton, J.:
Semistable models of elliptic curves over residue characteristic 2,
{\it Canadian Mathematical Bulletin} {61} no. 1 (2021), 154–162.


\end{thebibliography}
